\font\Bigrm = cmr10 scaled \magstep 2
\font\bigrm = cmr10 scaled \magstep 1
\font\Bigmath = cmmi10 scaled \magstep 2
 1
\font\Bigbb = msbm10 scaled \magstep 2
 1
\input amssym

\centerline{\Bigrm WEIGHTED SPECTRAL LARGE SIEVE INEQUALITIES}
\vskip 4pt 
\centerline{{\Bigrm FOR HECKE CONGRUENCE SUBGROUPS OF
{\Bigmath SL}({\Bigrm 2},$\,${\Bigbb Z}{\bigrm [{\Bigmath i}]})} }

\bigskip\centerline{by}\bigskip 

\centerline{\bigrm NIGEL WATT}

\vskip 11mm  

\footnote{{}}{\hskip -2 mm
{\bf 2010 Mathematics Subject Classification:}\quad
11F$30^{*}$, 11F37, 11F70, 11F72, 11L05, 11L07, 11L40, 11M41, 11N13, 
11N35, 11R42, 11R44, 22E30, 33C10, 44A15.}

\footline={\hss{\vbox to 2cm{\vfil\hbox{\rm\folio}}}\hss}

\noindent{\bf Abstract:}\ 
We prove new bounds for weighted mean values of sums involving 
Fourier coefficients of cusp forms that 
are automorphic with respect to a Hecke congruence subgroup 
$\Gamma\leq SL(2,{\Bbb Z}[i])$, 
and correspond to exceptional 
eigenvalues of the Laplace operator on the space   
$L^2(\Gamma\backslash SL(2,{\Bbb C})/SU(2))$. These results are, 
for certain applications, an effective substitute for the generalised 
Selberg eigenvalue conjecture. We give a proof of one such application, which is 
an upper bound for a sum of generalised Kloosterman sums (of significance  
in the study of 
certain mean values of Hecke zeta-functions with groessencharakters). 

Our proofs make extensive use of Lokvenec-Guleska's generalisation of the 
Bruggeman-Motohashi summation formulae for 
$PSL(2,{\Bbb Z}[i])\backslash PSL(2,{\Bbb C})$. 
We also employ a bound of Kim and Shahidi for 
the first eigenvalues 
of the relevant Laplace operators, and an `unweighted'  
spectral large sieve inequality (our proof of which is to appear separately). 

\medskip 


\vskip 11mm

\noindent{\bf Contents}

\vskip 6mm  

{\settabs\+\quad &5..\ &1.5.\ 
&Schwartz Spaces, Fourier Integrals, Poisson Summation and the 
Analytic Large Sieve\qquad &\cr  
\+&&&&Page\cr 
\+&Outline of Results and Methods&&&2\cr 
\smallskip 
\+&Acknowledgements&&&5\cr 
\smallskip
\+&1.&Definitions and Statements of the Results&&6\cr
\+&&1.1&The space $L^2(\Gamma\backslash G)$, 
Kloosterman Sums and Fourier Coefficients of Cusp Forms&6\cr  
\+&&1.2&A Kloosterman to Spectral Sum Formula and Other Key Ingredients&12\cr  
\+&&1.3&New Results on Sums over Exceptional Eigenvalues&17\cr
\+&&1.4&A Sum of Kloosterman Sums and an Application&20\cr  
\+&&1.5&Notation and Conventions&23\cr  
\smallskip 
\+&2.&Upper and Lower Bounds for the ${\bf K}$-Transform&&26\cr 
\smallskip 
\+&3.&A Bound in respect of a Single Level: the Proof of Theorem~4&&36\cr 
\smallskip 
\+&4.&Averaging over the Level&&39\cr 
\smallskip 
\+&5.&Schwartz Spaces, Fourier Integrals, Poisson Summation and the 
Analytic Large Sieve&&47\cr 
\smallskip 
\+&6.&An Elementary Bound for a Sum of Kloosterman Sums&&58\cr 
\smallskip 
\+&7.&Switching to Levels of Greater Modulus&&74\cr 
\smallskip 
\+&8.&The Proof of Theorem~9&&84\cr 
\smallskip 
\+&9.&The Proofs of Theorems~10 and~11&&98\cr  
\smallskip 
\+&References&&&115\cr}

\bigskip 
\bigskip 

\goodbreak\centerline{\bf Outline of Results and Methods}

\bigskip 

In [4] Bruggeman and Motohashi have obtained summation 
formulae for $PSL(2,{\Bbb Z}[i])\backslash PSL(2,{\Bbb C})$, 
analogous to the summation formulae for $PSL(2,{\Bbb Z})\backslash PSL(2,{\Bbb R})$ 
of Bruggeman [2] and Kuznetsov [15,16]. 
Their first formula shows that a certain wide class of sums involving 
Fourier coefficients of modular forms
may be expressed in terms of sums of `generalised' Kloosterman sums; 
their second formula does the reverse (and is a partial inverse of their first formula). 
In order to distinguish these two types of summation formula, 
we shall use the terminology `spectral to Kloosterman summation formula', 
and `Kloosterman to spectral summation formula'.   
Subsequently, in [19], Lokvenec-Guleska succeeded in generalising 
this work of Bruggeman and Motohashi, so as to obtain, for  
each imaginary quadratic number field $F={\Bbb Q}(\sqrt{D})$, and 
for each Hecke congruence subgroup $\Gamma$ of the special linear group 
$SL(2,{\frak O}_F)$   
(where ${\frak O}_F$ denotes the ring of integers of $F$), the corresponding  
summation formulae for $\Gamma\backslash SL(2,{\Bbb C})$. Our aim in this paper is to describe several applications 
of the case $F={\Bbb Q}(i)$ of these generalised summation formulae. 
Hence we assume, in what follows, that $\Gamma$ is one of the 
Hecke congruence  subgroups of the group 
$SL(2,{\Bbb Z}[i])$; these subgroups, and the `levels' 
with which they are associated, are 
defined at the beginning of the next section. 
Our applications of the summation formulae include new bounds for mean values 
of sums involving Fourier coefficients of modular forms, and new bounds for 
sums of generalised Kloosterman sums. Our results 
depend, in part, upon the best available lower bound for the absolute 
value $\lambda_1=\lambda_1(\Gamma)$ of the first non-zero  
eigenvalue  of the Laplace operator ${\bf \Delta}$ on the space 
$L^2(\Gamma\backslash SL(2,{\Bbb C})/SU(2))$. 
This is currently the bound $\lambda_1>77/81$ of Kim and Shahidi [13,14]. 

In the spectral theory of $L^2(\Gamma\backslash SL(2,{\Bbb C}))$ utilised 
in [19] the fundamental building blocks are not so much individual modular forms 
as whole subspaces $V$ that are irreducible and invariant 
with respect to the right-actions of 
all elements of $SL(2,{\Bbb C})$, and that are also `cuspidal' (in that the
$\Gamma$-automorphic functions they contain 
have, at each cusp of $\Gamma$,  a Fourier expansion  
in which the `constant' term is zero). 
Associated with each of these subspaces $V$ there is a pair 
of `spectral parameters' $(\nu_V , p_V)\in{\Bbb C}\times{\Bbb Z}$ 
satisfying either $\nu_V\in i[0,\infty)$, 
or else $\nu_V\in (0,2/9)$ and $p_V=0$ (these are the parameters appearing 
in the equations~(1.1.4) below). Moreover, for each $V$, and each cusp  
${\frak c}\in{\Bbb Q}(i)\cup\{\infty\}$, there are associated 
Fourier coefficients $c_V^{\frak c}(m)$ ($0\neq m\in{\Bbb Z}[i]$), 
which occur in the Fourier expansion at ${\frak c}$ of every one of a certain 
system of generators of the space $V$ (for details see (1.1.5)-(1.1.9) below). 
To state our results it is convenient to define the `modified Fourier coefficient' 
$c_V^{\frak c}(m;\nu_V,p_V)$ to be equal to $(\pi |m|)^{\nu_V}(m/|m|)^{-p_V}c_V^{\frak c}(m)$. 

A countably infinite set of pairwise orthogonal cuspidal subspaces $V$ 
arise in the spectral decomposition of $L^2(\Gamma\backslash SL(2,{\Bbb C}))$, 
but do not generate the whole of that space: indeed, the complete form of this  
spectral decomposition involves, 
in addition to the subspaces $V$, certain subspaces generated by 
continuously weighted mean values of Eisenstein series  
(for more details see (1.1.2), (1.1.3), (1.1.12) and (1.1.19) below). On one side 
of the summation formulae (of Bruggeman and Motohashi, or Lokvenec-Guleska) 
stand sums and integrals defined in terms of the spectral data and Fourier coefficients 
discussed above; 
on the other side are sums of `generalised Kloosterman sums' 
$S_{{\frak a},{\frak b}}(m,n;c)$, where the cusps ${\frak a}$, ${\frak b}$ 
and $m,n\in{\Bbb Z}[i]-\{ 0\}$ are fixed, 
and where the summation is over all $c$ lying in a certain 
countably infinite set ${}^{\frak a}{\cal C}^{\frak b}\subset {\Bbb C}^{*}$ with 
no point of accumulation in ${\Bbb C}$ (see (1.1.13)-(1.1.15) and (1.1.1) 
for the relevant definitions).  

In the paper [22] (to appear) the case $F={\Bbb Q}(i)$ of  the 
spectral to Kloosterman formula for $\Gamma_0(q)$ obtained in [19] is 
slightly generalised, so as to apply for arbitrary pairs of cusps 
${\frak a},{\frak b}$ (rather than  just for ${\frak a}={\frak b}=\infty$); 
by means of that generalised formula, and bounds for the relevant generalised 
Kloosterman sums, we obtain, in [22, Theorem~1], the spectral large  sieve inequality 
which is reproduced as `Theorem~2' in this paper. 
In the present paper it is instead the Kloosterman to spectral summation formula,  
Theorem~1 below, that has the more prominent part to play 
(though we use the spectral to Kloosterman 
summation formula in proving Theorem~11).  

To understand what motivates much of our work one has only to consider 
the application of this Kloosterman to spectral summation formula 
with the `test-function' $f=F_X$ given by 
$$F_X(z)=\Phi\!\left( X |z|^2\right)\qquad\qquad\hbox{($z\in{\Bbb C}^{*}$),}$$ 
where $X\geq 2$, and the function $\Phi$ is infinitely differentiable 
on $(0,\infty)$, with support $[1/2,2]$ and range $[0,1]$, say. 
By (1.2.1), we have 
$$\sum_{c\in {}^{\frak a}{\cal C}^{\frak b}}^{(\Gamma)}   
\,{S_{{\frak a},{\frak b}}\left( m , n ; c\right)\over |c|^2}\,
F_X\!\left( {2\pi\sqrt{mn}\over c}\right) 
=\pi\sum_V^{(\Gamma)}\,\overline{c_V^{\frak a}\left( m;\nu_V,p_V\right)}\,
c_V^{\frak b}\left( n;\nu_V,p_V\right) {\bf K}F_X\!\left(\nu_V , p_V\right)  
+ \cdots\;,\eqno(0.1)$$
where the transform ${\bf K}f(\nu,p)$ is given by (1.2.2)-(1.2.4), 
while the ellipsis `$\cdots$' signifies a sum involving modified 
Fourier coefficients of Eisenstein series, and the 
suffix `($\Gamma$)', placed above both summation signs, serves only 
as a reminder that the set 
${}^{\frak a}{\cal C}^{\frak b}$, the Kloosterman sum $S_{{\frak a},{\frak b}}(m,n;c)\,$  
and the relevant set of cuspidal subspaces $V$  are dependent on the 
Hecke congruence subgroup $\Gamma\subset SL(2,{\Bbb Z}[i])$.
Note that, since $\Phi$ has support $[1/2,2]$, the first summation in (0.1)  
is effectively a sum over the finitely many 
$c\in{}^{\frak a}{\cal C}^{\frak b}$ which satisfy $X'/2<|c|^2<2X'$, where 
$X'=4\pi^2 |mn|X$. By Lemma~2.2 of this paper, one moreover has 
$${\bf K}F_{X}(\nu ,p)\ll_{\Phi} (\log X) X^{-|p|}(1+|\nu|)^{-4}\qquad\quad  
\hbox{for $\quad (\nu,p)\in i{\Bbb R}\times{\Bbb Z}$.}\eqno(0.2)$$
Since the relevant spectral parameters $(\nu_V,p_V)$ are either contained in
$i{\Bbb R}\times{\Bbb Z}$, or else have $p_V=0$ and $0<\nu\leq 2/9$, 
it follows by (0.1), (0.2) and the spectral large sieve inequality 
(Theorem~2 below) that    
$$\sum_{\scriptstyle c\in {}^{\frak a}{\cal C}^{\frak b}\atop\scriptstyle 
{\textstyle{X'\over 2}}<|c|^2<2X'}^{(\Gamma)}   
\!\!\!\!\!\!\!\!\!{S_{{\frak a},{\frak b}}\!\left( m , n ; c\right)\over\pi |c|^2} 
\,F_X\!\!\left( {2\pi\sqrt{mn}\over c}\right) 
=O_{\Phi ,m,n}(\log X)+\sum_{\scriptstyle V\atop\scriptstyle\nu_V >0}^{(\Gamma)}
\!\overline{c_V^{\frak a}\!\left( m;\nu_V,0\right)}\,
c_V^{\frak b}\!\left( n;\nu_V,0\right) {\bf K}F_X\!\left(\nu_V , 0\right) ,\eqno(0.3)$$
where (see the discussion around (1.1.11) below) the last summation is 
certainly finite, and each relevant $V$ (if there be any) corresponds 
to an eigenvalue $\lambda_V$ of 
the (symmetric and positive) operator $-{\bf \Delta}$ on $L^2(\Gamma\backslash SL(2,{\Bbb C})/SU(2))$ satisfying 
$\lambda_V=1-\nu_V^2<1$; this is, moreover, a bijective correspondence.
If the generalised Selberg eigenvalue conjecture 
(for which see [7, Chapter~7, Part~6]) is correct, then 
the sum over $V$ in (0.3) is empty. 
In the absence of any proof of this conjecture it remains relevant to note 
that, by Lemmas~2.2 and~2.3 of this paper,  
$$\nu^{-1}X^{\nu}\ll_{\Phi} {\bf K}F_{X}(\nu ,0)\ll_{\Phi} \nu^{-1}X^{\nu}\qquad\quad  
\hbox{for $\quad 0<\nu\leq 1/2\ $ and $\ X\rightarrow +\infty$.}\eqno(0.4)$$
Hence, if it is (for example) the case that the eigenvalue $\lambda_1(\Gamma)$  
is both `exceptional' (i.e. less than $1$) and  of multiplicity $1$, then 
by (0.3) and (0.4) one will have 
$$\Biggl|\sum_{\scriptstyle c\in {}^{\frak a}{\cal C}^{\frak b}\atop\scriptstyle 
{\textstyle{X'\over 2}}<|c|^2<2X'}^{(\Gamma)}   
\!\!\!\!\!\!\!\!\!{S_{{\frak a},{\frak b}}\!\left( m , n ; c\right)\over\pi |c|^2} 
\,F_X\!\!\left( {2\pi\sqrt{mn}\over c}\right)\Biggr|  
\asymp_{\Phi}\,{\left| c_{V_1}^{\frak a}\!\left( m;\nu_1,0\right)  
c_{V_1}^{\frak b}\!\left( n;\nu_1,0\right)\right|\over\nu_1}\ X^{\nu_1}\qquad\quad  
\hbox{as $\quad X\rightarrow +\infty$,}\eqno(0.5)$$
where $\nu_1=\sqrt{1-\lambda_1(\Gamma)}$, and where $V_1$ is that 
cuspidal subspace $V$ which occurs in the spectral decomposition of 
$L^2(\Gamma\backslash SL(2,{\Bbb C}))$ and has $(\nu_V,p_V)=(\nu_1,0)$. 
It is therefore reasonable to expect that, in the event that the sum 
over $V$ in (0.3) being non-empty, that sum will be the crucial 
determinant of the asymptotic behaviour (as $X\rightarrow +\infty$) 
of the sum of Kloosterman sums appearing in (0.3). By the 
bound $\lambda_1>77/81$ of Kim and Shahidi, one has 
$0<\nu_1<2/9$ in (0.5), and $0<\nu_V<2/9$ in the sum over $V$ in (0.3). 

We follow the pattern set by Deshouillers and Iwaniec [5], 
in considering weighted 
mean values (over $m$ and $n$)  
of the sum of Kloosterman sums 
appearing on the left-hand sides of (0.5). Their results on 
sums of generalised Kloosterman sums associated  with Hecke congruence 
subgroups of $SL(2,{\Bbb Z})$ had (see [5, Section~1.5]) 
numerous applications to problems concerning the 
multiplicative number theory of the rational integers:  
so one motivation for this paper is to obtain results 
that may help to similarly 
advance the multiplicative number theory of the Gaussian integers.  
Bearing in mind the equation (0.3), and the bounds in (0.4), 
we are led (via the Cauchy-Schwarz inequality) to investigate  
what upper bounds may be obtained for the sums 
$$\sigma_q^{\frak a}({\bf b},N;X)
=\sum_{\scriptstyle V\atop\scriptstyle\nu_V>0}^{\left(\Gamma_0(q)\right)} 
X^{\nu_V}
\Biggl|\sum_{{\textstyle{N\over C}}<|n|^2\leq N} 
b_n c_V^{\frak a}\left( n;\nu_V,0\right)\Biggr|^2\qquad\qquad  
\hbox{($0\neq q\in{\Bbb Z}[i]$, ${\frak a}\in{\Bbb Q}(i)\cup\{\infty\}$),}\eqno(0.6)$$
where $N,X\geq 1$, $C>1$ and $\Gamma_0(q)$ is the Hecke congruence subgroup 
of $SL(2,{\Bbb Z}[i])$ of level $q$, while  
the coefficients $b_n$ ($0\neq n\in{\Bbb Z}[i]$) are arbitrary complex 
numbers (collectively represented by the symbol `${\bf b}$'). 
Usually we have either $C=2$, or $C=4$.  Note, moreover, that 
`$n$', in the above, is a Gaussian integer variable of summation 
(i.e. it ranges over all values in ${\Bbb Z}[i]$ permitted 
by the conditions attached to the summation sign). 
Indeed, since most of the summations that appear in this paper 
are summations over ${\Bbb Z}[i]$, it has suited us to 
make it our convention that, where there is nothing to 
indicate the contrary, variables of summation are  
understood to be Gaussian integer variables. 

Our principal results 
presuppose a uniform bound of the form 
$${\rm Re}\!\left(\nu_V\right)\leq\vartheta\leq{2\over 9}\;,\eqno(0.7)$$
such as follows (with $\vartheta =2/9$) from the lower bound 
$\lambda_1>77/81$ of Kim and Shahidi. 
By combining the bound (0.7) with a 
carefully targeted application of the Kloosterman to spectral summation formula, 
we obtain, in Theorem~4 below, an upper bound for 
$\sigma_q^{\frak a}({\bf b},N;X)$. The proof of this result, and those of 
Theorems~5,~6 and~7 are modelled on the proofs of the 
analogous results in [5]. 

In Theorems~5-9 we specialise to the case ${\frak a}=\infty$, 
and consider the mean value, over levels $q\in{\Bbb Z}[i]$ satisfying a 
condition of the form $Q/2<|q|^2\leq Q$, of the sum 
$\sigma_q^{\infty}({\bf a},N;X)$: note that Theorems~8 and~9 apply only when   
the relevant coefficients $a_n$ ($0\neq n\in{\Bbb Z}[i]$) are of a special type. 
Theorems~6 and~7 are a key tool in our proofs (by induction) 
of Theorems~5 and~9: they enable one to relate the mean value 
$$S(Q,X,N)=\sum_{{\textstyle{Q\over 2}}<|q|^2\leq Q}\sigma_q^{\infty}({\bf a},N;X)$$ 
to other mean values 
of the same form, but with $X$ and $Q$ replaced by other numbers (and, in the case of 
Theorem~7, with each coefficient $a_n$ replaced by the corresponding 
product $a_n |n|^{2it}$,  
where $t$ is some real number independent of $n$).
In cases where both $Q/N$ and $X/(Q/N)$ are sufficiently large, 
the result (1.3.7) of Theorem~5 
is a sharper upper bound than that which follows directly from Theorem~4 and (0.7). 

We consider Theorem~9, which is an analogue of [21, Theorem~2], 
to be the foremost achievement of this paper. 
Indeed, both Theorems~8 and~9 are play a crucial part in a significant 
application that we will come to shortly (after some discussion of Theorem~8, and 
of the proof of Theorem~9). 

Theorem~8 is (as shown towards the end of Subsection~1.3) an easy corollary of 
Theorem~9. There is, of course, a direct proof of Theorem~8 
(one considerably shorter than 
that of Theorem~9), but we have not included it in this paper. 
A measure of the strength of Theorem~8 is that, in respect of 
cases where $1\leq X\leq Q^2/N$ and the coefficients $a_n$ satisfy the 
required hypotheses (i.e. with $H=N$), it yields the same bound for 
the mean value $S(Q,X,N)$
as would follow (by Theorems~3 and~4), were it known to be the case that, 
for every $\varepsilon >0$, one has 
$\lambda_1(\Gamma)\geq 1-\varepsilon$ for all but 
finitely many of the Hecke congruence subgroups 
$\Gamma =\Gamma_0(q)\leq SL(2,{\Bbb Z}[i])$. 

Our proof of Theorem~9 resembles that of the analogous  
result [21, Theorem~2]  in having three distinct phases. 
In the first phase (to which Sections~5 and~6 are devoted) 
we obtain, ultimately through Lemmas~6.1-6.3, a bound for a sum of the form  
$${\cal R}
=\sum_{p\neq 0}\theta_p|p|^{-2}\sum_{q\neq 0}|q|^{-2}\sum_h\phi_h\sum_k\sum_{\ell}
\,S(hk,\ell;pq)\,\varphi(h,k,\ell,p,q)\,\Upsilon_{\ell}\;,$$
where $S(u,v;w)$ denotes the `simple Kloosterman sum' defined in (1.3.6),  
and where it is supposed that $\theta_m,\phi_m=O(1)$ for $0\neq m\in{\Bbb Z}[i]$;  
that $\Upsilon_{\ell}\in{\Bbb C}$ for ${\ell}\in{\Bbb Z}[i]$; 
and that the 
function $\varphi : {\Bbb C}^5\rightarrow{\Bbb C}$ 
is `sufficiently smooth' (in the sense made clear at the start of Section~6) and, 
for some $Z_1,\ldots Z_5\geq 1$, 
has its support 
${\rm Supp}(\varphi)$ contained in the set $\{ {\bf z}\in{\Bbb C}^5 : 
Z_j/2\leq |z_j|^2\leq Z_j\ {\rm for}\ j=1,\ldots ,5\}$. 

The steps in the proof of Lemma~6.1 are similar to steps 
in the initial part of the proof of [21, Proposition~2.1]. The only   
(rather minor) novelty there is the use of Poisson summation 
over ${\Bbb Z}[i]$, instead of Poisson summation over ${\Bbb Z}$. 
Although it would require some additional hypotheses concerning the 
coefficients $\phi_h$ and $\Upsilon_{\ell}$, 
the entire proof 
of [21, Proposition~2.1] could be adapted for the current context: 
this would not yield identically the same bound for ${\cal R}$  
as that obtained in Lemma~6.3, but would nevertheless produce  
a result from which Theorem~9 could be deduced. Rather than do this,   
we instead take the opportunity 
to try out some new ideas, in the hope of achieving a 
proof in which the key features are less obscured by 
lengthy computations than is the case in respect of the proof 
of [21, Proposition~2.1].

Our primary innovation (in the estimation of ${\cal R}$) 
is to be found in the proof of Lemma~6.2. There we apply, in conjunction with 
the Cauchy-Schwarz inequality,  a `special analytic 
large sieve for ${\Bbb Z}[i]$', which is obtained in Lemmas~5.8 and~5.9 
(as a corollary of Huxley's more general large sieve estimates in 
[9, Theorem~1]). This ultimately results in the bound for ${\cal R}$ 
that we obtain in Lemma~6.3. That bound, however, 
is not quite adequate for our purposes, for  
it is only obtained subject to quite stringent conditions 
(these being the same conditions as appear in the hypotheses of Lemma~6.1). 
The hypothesis that $Q\gg\max\{ HK , L\}$, is the most irksome of these conditions;   
it causes Lemma~8.4 to be conditional upon having $R\geq N$; 
and if we had no means of setting aside this last constraint, then 
we would be unable to deduce the case $1\leq Q\leq N$ of the result 
in Lemma~8.5, which would (at best) make the deduction 
of the result of Theorem~9 more difficult. These considerations prompt 
our work in Section~7, which is an application of the Parseval identity 
[19, Theorem~8.1] pertaining to a 
certain subspace of $L^2(\Gamma_0(q)\backslash SL(2,{\Bbb C}))$. 
In Lemma~7.3 we find that 
$$S(R,X,N)\ll (\log R_1)\left( S\left( R_1,X,N\right) + S\left( 2R_1,X,N\right)\right)\qquad\quad  
\hbox{for $\quad 1\leq R\leq {2\over 5}\,R_1$.}\eqno(0.8)$$
Hence, at the (acceptable) cost of increasing our final bounds on $S(Q,X,N)$ 
by a factor $O(\log N)$, we are effectively able to nullify the condition  
$R\geq N$ of Lemma~8.4, and so compensate for the above mentioned 
inadequacy of the bound for ${\cal R}$ obtained in Lemma~6.3.

Lemma~8.5 marks the end of the second phase of our proof of Theorem~9. 
At the end of Section~8 comes the third and final phase, in which 
it is shown 
(with the help of Lemma~4.2, a corollary of Lemmas~6 and~7) 
that Theorem~9 follows by induction from Lemma~8.5. 

In the paper [23] (to appear), an analysis of the contribution 
of `off-diagonal terms' to a certain  
mean value of groessencharakter zeta-functions 
(a smoothly weighted majorant of the mean value $J(D,N)$ 
defined in (1.4.22) below)
leads to a sum of generalised Kloosterman sums $S_{{\frak a},{\frak b}}(m,n;c)$, 
in which ${\frak a}=1/s$, with $0\neq s\in{\Bbb Z}[i]$, and  
${\frak b}=\infty$, while the relevant discrete subgroups of 
$SL(2,{\Bbb C})$ are Hecke congruence subgroups $\Gamma_0(rs)\leq SL(2,{\Bbb Z}[i])$,   
with $0\neq r\in{\Bbb Z}[i]$ and $r$ coprime to $s$. 
This is analogous to the situation which obtained in respect of 
the proofs of both [6, Theorems~1~and~2] and the later result  
[21, Theorem~1]; and it provides the motivation for Theorem~10 
of the present paper, in which we obtain a bound for the sum of generalised Kloosterman 
sums in question that is     
(if one allows for the stronger 
lower bounds for $\lambda_1(\Gamma)$ now available) 
analogous to the bound obtained in 
[21, Proposition~4.1] . 

See the end of Section~1.4 for a brief description    
(with some history) of the main result obtained in [23]. 
Note, in particular, that the proof, in [23], of the bound (1.4.23) for $J(D,N)$  
depends critically on the result that we obtain in Theorem~10.  
Since Theorem~10 is essentially a corollary of 
Theorems~3,~4,~8 and~9, the part it plays in [23] therefore constitutes 
a significant application of those results. 

\bigskip 
\bigskip 

\goodbreak\centerline{\bf Acknowledgements}

\bigskip 

This work was begun while the author 
was employed at Royal Holloway, University of London, 
through a fellowship associated with the EPSRC funded project, 
`The Development and Application of Mean Value Results in Multiplicative
Number Theory' (GR/T20236/01), led by Glyn Harman. 
The author thanks the EPSRC and Royal Holloway for their support, and 
good working conditions. He is grateful to Glyn for his advice and 
encouragement during 2004-6, and in the years since. 

The author's research work in connection with the above EPSRC 
project has benefitted greatly from correspondence 
with Professor Roelof Bruggeman (at Utrecht University) 
and Professor Yoichi Motohashi (at Nihon University); 
the author wishes to record his gratitude for their 
comprehensive answers to many questions, and for their 
correspondence on various interesting related subjects.  

A significant part of the work described in this paper 
was completed in the period 2006-9, 
during which the author held a Lectureship 
at Cardiff University. The author thanks Professor Martin Huxley 
(at Cardiff University) for useful discussions related to this paper. 

The author wishes to thank his parents for their 
encouragement, and for providing very helpful support 
while this paper was being written. 

\bigskip
\bigskip

\goodbreak\centerline{\bf \S 1. Definitions and Statements of the Results}

\bigskip 
\smallskip 

\goodbreak\centerline{\bf \S 1.1. The Space $L^2(\Gamma\backslash G)$, 
Kloosterman Sums and Fourier Coefficients of Cusp Forms}

\bigskip 

Let ${\frak O}={\Bbb Z}[i]$ (the ring of Gaussian integers). 
Then, for each non-zero $q\in{\frak O}$, the Hecke congruence subgroup of 
$SL(2,{\frak O})$ of level $q$ is the group 
$$\Gamma_0(q)
=\left\{ \pmatrix{a&b \cr c&d}\in SL(2,{\frak O}) : c\in q{\frak O}\right\}\;,$$
endowed with the associative binary operation of matrix multiplication;  
and all Hecke congruence subgroups of $SL(2,{\frak O})$ 
are contained in the set  
$\{ \Gamma_0(q) : 0\neq q\in{\frak O}\}=\{ \Gamma_0(q) : q\in{\frak O},\ 
{\rm Re}(q)> 0\ {\rm and}\ {\rm Im}(q)\geq 0\}$. 

Let $\Gamma$ be a Hecke congruence subgroup of $SL(2,{\frak O})$.  
Then $\Gamma$ is a discrete and cofinite (but not cocompact) subgroup of 
the Lie group $G=SL(2,{\Bbb C})$. 
A function $f : G \rightarrow {\Bbb C}$ is said to be $\Gamma$-automorphic 
if and only if it is such that  
$$f(\gamma g)=f(g)\qquad\quad\hbox{for $\quad\gamma\in\Gamma\ $ and $\ g\in G\,$.}$$
In preparation for further discussion of $\Gamma$-automorphic functions we 
next define a coordinate system and measure for $G$. 

The maximal compact subgroup of $G$ is 
$K=SU(2)=\left\{ k[\alpha,\beta] : 
\alpha,\beta\in{\Bbb C}\ {\rm and}\ |\alpha|^2+|\beta|^2=1\right\}$ 
where 
$$k[\alpha,\beta]=\pmatrix{\alpha&\beta\cr -\overline{\beta}&\overline{\alpha}}\qquad\qquad 
\hbox{($\alpha,\beta\in{\Bbb C}$).}$$
One has $G=NAK$ (the Iwasawa decomposition)
where 
$N=\{ n[z] : z\in{\Bbb C}\}$ and 
$A=\{ a[r] : r>0\}$, with 
$$n[z]=\pmatrix{1&z\cr 0&1}\qquad\hbox{($z\in{\Bbb C}$)}\qquad\quad\hbox{and}\qquad\quad   
a[r]=\pmatrix{\sqrt{r}&0\cr 0&1/\sqrt{r}}\qquad\hbox{($r>0$).}$$
Moreover, each $k=k[\alpha,\beta]\in K$ has a factorisation of the form 
$k=h\bigl[ e^{i\varphi /2}\bigr] v[i\theta] h\bigl[ e^{i\psi /2}\bigr]$, 
where $\varphi,\theta,\psi\in{\Bbb R}$, 
$$h[u]=\pmatrix{u&0\cr 0&u^{-1}}\qquad\quad\hbox{and}\qquad\quad  
v[i\theta]=\pmatrix{\cos(\theta /2)&i\sin(\theta /2)\cr i\sin(\theta /2)&\cos(\theta /2)}\;.$$
Each $g\in G$ has unique Iwasawa coordinates 
$(z,r,\theta,\varphi,\psi)\in{\Bbb C}\times(0,\infty)\times{\Bbb R}^3$ such that
$\theta\in[0,\pi)$, $\varphi\pm\psi\in[0,4\pi)$ and 
$n[z]a[r]h\bigl[ e^{i\varphi /2}\bigr] v[i\theta] h\bigl[ e^{i\psi /2}\bigr]=g$. 
In terms of these coordinates 
(and with $x={\rm Re}(z)$, $y={\rm Im}(z)$) 
the subgroups $N$, $A$ and $K$ have left and right 
Haar measures 
${\rm d}n={\rm d}_{+}z={\rm d}x\,{\rm d}y$,  
${\rm d}a=r^{-1}\,{\rm d}r$ and  
${\rm d}k=2^{-3}\pi^{-2}\sin(\theta)\,{\rm d}\varphi\,{\rm d}\theta\,{\rm d}\psi$,
respectively. 
Note that ${\rm d}k$ here is normalised so that 
$\int_{K}{\rm d}k=2$. 
Similarly 
$${\rm d}g=r^{-2}\,{\rm d}n\,{\rm d}a\,{\rm d}k=r^{-3}\,{\rm d}x\,{\rm d}y\,{\rm d}r\,{\rm d}k
\qquad\qquad\hbox{($g=n[x+iy]a[r]k$, $x,y\in{\Bbb R}$, $r>0$ and $k\in K$)}$$
is a left and right Haar measure for $G$. 

By [7, Chapter~7, Proposition~3.9], 
a fundamental domain for the action of $SL(2,{\frak O})$ upon 
$G$ is the set 
$${\cal F}_{{\Bbb Q}(i)}^G={\cal F}_{{\Bbb Q}(i)}^{NA} K^{+}\;,$$ 
where 
$${\cal F}_{{\Bbb Q}(i)}^{NA}
=\left\{ n[z]a[r] : z\in{\Bbb C}, r>0, |z|^2+r^2\geq 1\ {\rm and}\ 
|{\rm Re}(z)|,{\rm Im}(z)\in[0,1/2]\right\}$$
and 
$$K^{+}
=\left\{ 
h\left[ e^{i\varphi /2}\right] v\left[ i\theta\right] h\left[ e^{i\psi /2}\right] : 
0\leq\theta <\pi, 0\leq\varphi-\psi<4\pi\ {\rm and}\ 0\leq\varphi+\psi<2\pi\right\}$$
(the latter set being a fundamental domain for 
the action of the group $\{ h[1] , h[-1]\}$ on $K$\/). 
Since the group $\Gamma =\Gamma_0(q)$ is of finite index in 
$SL(2,{\frak O})=\Gamma_0(1)$, there exist representatives  
$\gamma_1,\ldots ,\gamma_{[SL(2,{\frak O}) : \Gamma]}$ of the 
right cosets of $\Gamma$ in $SL(2,{\frak O})$
such that the set
$$\bigcup_{k=1}^{[SL(2,{\frak O}) : \Gamma]}\gamma_k {\cal F}_{{\Bbb Q}(i)}^G
=\bigcup_{k=1}^{[SL(2,{\frak O}) : \Gamma]}\gamma_k {\cal F}_{{\Bbb Q}(i)}^{NA} 
K^{+}={\cal F}_{\Gamma\backslash G}\qquad\ \hbox{(say)}$$
is a fundamental domain for the action of $\Gamma$ upon $G$.

Since $\Gamma\ni h[-1]$, the $\Gamma$-automorphic functions $f$ are 
even  (i.e. they satisfy $f(h[-1]g)=f(g h[-1])=f(g)$, for $g\in G$\/). 
Given any measurable $\Gamma$-automorphic function $f : G\rightarrow{\Bbb C}$, 
one defines 
$$\int_{\Gamma\backslash G} f(g) {\rm d}g 
=\int_{{\cal F}_{\Gamma\backslash G}} f(g) {\rm d}g$$
if the latter integral exists (note that this integral is independent of 
our particular choice of fundamental domain ${\cal F}_{\Gamma\backslash G}$). 
Such a function  $f$ is said to be `square integrable' if and only if 
$$\int_{\Gamma\backslash G} |f(g)|^2 {\rm d}g<\infty\;.$$
We define $L^2(\Gamma\backslash G)$ to be the 
set of all square integrable $\Gamma$-automorphic functions 
$f : G\rightarrow{\Bbb C}$. This set $L^2(\Gamma\backslash G)$ 
is a Hilbert space with inner product 
$$\langle f , h\rangle_{\Gamma\backslash G}
=\int_{\Gamma\backslash G} 
f(g)\overline{h(g)} {\rm d}g\qquad\qquad\hbox{($f,h\in L^2(\Gamma\backslash G)$\/).}$$

We now define what is meant by 
`Fourier expansion at a cusp': this concept will prove useful in discussing the 
decomposition of the space $L^2(\Gamma\backslash G)$. 
When ${\frak z}=\bigl[ z_1,z_2\bigr]\in{\Bbb P}^{1}({\Bbb C})={\Bbb C}\cup\{\infty\}$ 
(the Riemann sphere) and 
$$\ g=\pmatrix{a&b\cr c&d}\in G\;,$$
one may define
$g{\frak z}=\bigl[ a z_1+b z_2,c z_1+d z_2\bigr]\in{\Bbb P}^{1}({\Bbb C})$
(so that $g\infty =\infty=[1,0]$ if and only if $c=0$\/). 
This (since $\Gamma<G$) determines an action of $\Gamma$ on the Riemann sphere. 
The `cusps' of $\Gamma$ are the points ${\frak c}\in P^{1}({\Bbb Q}(i))={\Bbb Q}(i)\cup\{\infty\}$. 
For a cusp ${\frak c}$ of $\Gamma$, the corresponding stabiliser and 
`parabolic stabiliser' subgroups, 
$$\Gamma_{\frak c}=\{ \gamma\in\Gamma : \gamma{\frak c} ={\frak c}\}\qquad\quad
\hbox{and}\qquad\quad  
\Gamma_{\frak c}'=\{ \gamma\in\Gamma_{\frak c} : {\rm Tr}(\gamma)=2\}\;,$$
(where ${\rm Tr}(\gamma)$ denotes the matrix trace of $\gamma$\/) are both infinite, 
with $\left[\Gamma_{\frak c}:\Gamma_{\frak c}'\right]\in\{ 2,4\}$; 
and it is possible to choose a `scaling matrix' $g_{\frak c}\in G$ such that  
$g_{\frak c}\infty ={\frak c}$ and   
$$g_{{\frak c}}^{-1}\Gamma_{{\frak c}}' g_{{\frak c}}
=\left\{ n[\alpha] : \alpha\in{\frak O}\right\}=B^{+}\qquad\ \hbox{(say).}\eqno(1.1.1)$$
We assume henceforth that each cusp ${\frak c}$ of $\Gamma$ has  
assigned to it just such a scaling matrix $g_{\frak c}$. 
When $f : G\rightarrow{\Bbb C}$ is $\Gamma$-automorphic one has   
$f\left( g_{\frak c} n[\alpha] g\right)
=f\left( g_{\frak c} n[\alpha]g_{\frak c}^{-1}g_{\frak c}g\right)
=f\left( g_{\frak c}g\right)$, for $g\in G$ and 
$\alpha\in{\frak O}$.  Hence if $f$ is (for example) a $\Gamma$-automorphic 
function that is continuous on $G$, then 
one has a Fourier expansion at the cusp ${\frak c}$:
$$f\left( g_{\frak c}g\right) 
=\sum_{\omega\in{\frak O}} \left( F_{\omega}^{\frak c} f\right) (g)\qquad\qquad 
\hbox{($g\in G$),}$$
where, for $\omega\in{\frak O}$, the function $F_{\omega}^{\frak c}f : G\rightarrow{\Bbb C}$ 
is continuous on $G$ and satisfies  
$$\left( F_{\omega}^{\frak c}f\right)\left( n[z] g\right) 
={\rm e}\left( {\rm Re}(\omega z)\right) \left( F_{\omega}^{\frak c}f\right)(g)\qquad\qquad 
\hbox{($g\in G$, $z\in{\Bbb C}$),}$$
with `${\rm e}(x)$' being a convenient notation for $\exp(2\pi ix)$. 

Let ${}^0L^2(\Gamma\backslash G)$ denote the 
closure in $L^2(\Gamma\backslash G)$ of the subspace spanned 
by cusp forms (we define the term `cusp form' below (1.1.10)). 
Then, by the discussion in [19, Chapter~8], the Hilbert space $L^2(\Gamma\backslash G)$ 
has a decomposition into mutually orthogonal subspaces, 
$$L^2(\Gamma\backslash G)
={\Bbb C}\oplus{}^{0}L^2(\Gamma\backslash G)\oplus{}^{\rm e}L^2(\Gamma\backslash G)\;,\eqno(1.1.2)$$
where ${\Bbb C}$ denotes the $1$-dimensional space of constant functions, while 
$${}^{0}L^2(\Gamma\backslash G)
=\overline{\bigoplus V}\;,\eqno(1.1.3)$$ 
with $V$ running over a countably infinite set of mutually orthogonal `cuspidal' proper subspaces (in 
the terminology of representation theory each $V$ here is both invariant and irreducible with respect to the right-actions 
of the elements of $G$\/). To classify the spaces $V$ we need the 
two Casimir operators associated with $G$, which are $\Omega_{+}$ and 
$\Omega_{-}=\overline{\Omega_{+}}$, where 
in terms of the Iwasawa coordinates (and with 
$\partial/\partial z=(1/2)(\partial/\partial x-i\partial/\partial y)$ and 
$\partial/\partial\overline{z}=\overline{\partial/\partial z}
=(1/2)(\partial/\partial x+i\partial/\partial y)$\/) one has:  
$$\eqalignno{\Omega_{+}
 &={1\over 2}\,r^2{\partial\over\partial z}{\partial\over\partial\overline{z}}
+{1\over 2}\,r e^{i\varphi}\cot(\theta){\partial\over\partial z}{\partial\over\partial\varphi}
-{1\over 2}\,ir e^{i\varphi}{\partial\over\partial z}{\partial\over\partial\theta}
-{1\over 2}\,re^{i\varphi}\csc(\theta){\partial\over\partial z}{\partial\over\partial\psi} + {}\qquad\cr
 &\qquad +{1\over 8}\,r^2{\partial^2\over\partial r^2}  
-{1\over 4}\,ir {\partial\over\partial r}\,{\partial\over\partial\varphi} 
-{1\over 8}\,{\partial^2\over\partial\varphi^2} -{1\over 8}\,r\,{\partial\over\partial r} 
+{1\over 4}\,i{\partial\over\partial\varphi}\;. 
}$$
By the discussion in [19, Subsection~3.2.2], 
each $V$ occurring in the decomposition (1.1.3) has associated with it 
a unique pair of `spectral parameters' 
$\left(\nu_V,p_V\right)\in \left( i[0,\infty)\times{\Bbb Z}\right)\cup\left( (0,1)\times\{ 0\}\right)$ 
such that 
$$\Omega_{\pm} f
={1\over 8}\,\left(\left(\nu_V\mp p_V\right)^2 -1\right) f\qquad\quad  
\hbox{for $\quad f\in V$;}\eqno(1.1.4)$$
and each has, itself, a decomposition into mutually 
orthogonal proper subspaces:  
$$V=\overline{\bigoplus_{\ell =\left| p_V\right|}^{\infty}\bigoplus_{q=-\ell}^{\ell} 
V_{K,\ell,q}}\eqno(1.1.5)$$
with $V_{K,\ell,q}\subseteq\left\{ 
f\in V : \Omega_K f=\textstyle{1\over 2}(\ell +1)\ell f\ {\rm and}\ 
(\partial/\partial\psi)f=-iqf\right\}$ for 
$q,\ell\in{\Bbb Z}$, $\ell\geq |p_V|$ and $|q|\leq\ell$, 
where $\Omega_{K}$ (the Casimir operator associated with $K=SU(2)$\/) 
is given by  
$$\Omega_{K}
={1\over 2}\,\csc^2(\theta)\left({\partial^2\over\partial\varphi^2}
+{\partial^2\over\partial\psi^2}\right) 
-\csc(\theta)\cot(\theta){\partial^2\over\partial\varphi\partial\psi}
+{1\over 2}\,{\partial^2\over\partial\theta^2} 
+{1\over 2}\,\cot(\theta){\partial\over\partial\theta}\;.$$
We follow [4] and [19] in our use of the symbols `$\ell$', `$p_V$' and `$q$'  
in the above: it may therefore be worth clarifying that `$q$', 
in the context of the spaces $V_{K,\ell,q}$ in (1.1.5), 
denotes a rational integer valued 
variable that is independent of the (as yet unspecified) 
level of the group $\Gamma$.  But, from Subsection~1.2 onwards  
(where there is 
little need to discuss the spaces $V_{K,\ell,q}$ or related matters,   
unless it be in respect of the case $\ell =p_V =q=0$)  
we generally have $\Gamma =\Gamma_0(q)$,  
so that $q$ then denotes some non-zero Gaussian integer `level'.  
As explained below [19, Equation (8.3)], each factor $V_{K,\ell,q}$ in the decomposition (1.1.5)  
is a $1$-dimensional space over ${\Bbb C}$, and so 
contains some generator $f_{\ell,q}^V\neq 0$ such that  
$$V_{K,\ell,q}={\Bbb C} f_{\ell,q}^V\eqno(1.1.6)$$
(we take this generator $f_{\ell,q}^V$ to equal the 
`$T_V\varphi_{\ell,q}(\nu_V,p_V)$' of [19, Chapter~8]). 

The Fourier expansions at cusps of the generators of the above spaces $V_{K,\ell,q}$ 
are a central concern of this paper. In order to best describe (and compare) 
these Fourier expansions we now define certain `Jacquet integrals'. 
For $k=k[\alpha ,\beta]\in K$ and $p,q,\ell\in{\Bbb Z}$ with $\ell\geq\max\{ |p| , |q|\}$, 
let $\Phi_{p,q}^{\ell}\left( k[\alpha ,\beta]\right)$ be the coefficient of 
$X^{\ell -p}$ in the polynomial 
$(\alpha X-\overline{\beta})^{\ell -q} (\beta X+\overline{\alpha})^{\ell +q}$. 
Then the system 
$\bigl\{\Phi_{p,q}^{\ell} : 
\hbox{$p,q,\ell\in{\Bbb Z}$ and $\ell\geq \max\{ |p|,|q|\}$}\bigr\}$ 
is a complete orthogonal basis of the Hilbert space $L^2_{\rm even}(K)$ 
endowed with the inner product  
$\langle f , h\rangle_K =\int_K f(k)\overline{h(k)} {\rm d}k$ 
(i.e. the space of even functions $f : K\rightarrow{\Bbb C}$ 
such that $\int_K |f|^2 {\rm d}k <\infty$\/). One has 
$$\left\|\Phi_{p,q}^{\ell}\right\|_K^2
=\int_K \left|\Phi_{p,q}^{\ell}(k)\right|^2 {\rm d}k 
={1\over (\ell +\textstyle{1\over 2})}\,
\pmatrix{2\ell\cr\ell -p}\pmatrix{2\ell\cr\ell -q}^{-1}\qquad   
\hbox{for $\ \,p,q,\ell\in{\Bbb Z}\ $ and $\ \ell\geq \max\{ |p|,|q|\}$.}$$
For $\omega,\nu\in{\Bbb C}$, with ${\rm Re}(\nu)>0$, and $\ell,p,q\in{\Bbb Z}$ with 
$\ell\geq\max\{ |p| , |q|\}$, one defines $\varphi_{\ell ,q}(\nu ,p) : G\rightarrow{\Bbb C}$ 
and the corresponding Jacquet integral ${\bf J}_{\omega}\varphi_{\ell ,q}(\nu ,p) : G\rightarrow{\Bbb C}$ by: 
$$\varphi_{\ell ,q}(\nu ,p)\left( n a[r] k\right) =r^{1+\nu}\Phi_{p,q}^{\ell}(k)\qquad\qquad 
\hbox{($n\in N$, $r>0$, $k\in K$\/);}\eqno(1.1.7)$$
$$\left( {\bf J}_{\omega}\varphi_{\ell ,q}(\nu ,p)\right) (g)
=\int_{\Bbb C} \varphi_{\ell ,q}(\nu ,p)\left( k[0,-1] n[z] g\right)
{\rm e}(-{\rm Re}(\omega z)) {\rm d}n[z]\qquad\quad
\hbox{for $\quad g\in G$.}$$
The last integral converges absolutely when ${\rm Re}(\nu)>0$: though 
it fails to do so when ${\rm Re}(\nu)\leq 0$, it 
is shown by [4, Lemma~5.1] that if $\ell,p,q,g$ and $\omega\neq 0$ are given  
then the function $\nu\mapsto\left( {\bf J}_{\omega}\varphi_{\ell ,q}(\nu ,p)\right) (g)$  
has an entire analytic continuation. Through this one defines  
the function ${\bf J}_{\omega}\varphi_{\ell,q}\left(\nu_V ,p_V\right) : G\rightarrow{\Bbb C}$  
when $\left(\nu_V,p_V\right)$ are the spectral parameters of an arbitrary irreducible 
subspace $V\subset{}^{0}L^2(\Gamma\backslash G)$.
As noted in [22, Subsection~1.7] (see, in particular, [22, Relations (1.7.10) and (1.5.17)]\/), each term $F_{\omega}^{\frak c} f$ 
in the Fourier expansion at any cusp ${\frak c}$ of any function $f\in V_{K,\ell,q}$ 
is a constant multiple of the corresponding Jacquet integral, 
${\bf J}_{\omega}\varphi_{\ell,q}\left(\nu_V ,p_V\right) : G\rightarrow{\Bbb C}$. 
Indeed, it is even possible to choose, for 
the subspace factors $V_{K,\ell,q}$ in (1.1.5), a system of generators, 
$${\cal B}^V_K=\left\{ f^V_{\ell,q}\in V_{K,\ell,q}-\{ 0\} :
\ell,q\in{\Bbb Z}, \ell\geq |p_V|\ {\rm and}\  |q|\leq\ell\right\}\qquad\ \hbox{(say),}
$$ 
such that at 
each cusp ${\frak c}$ of $\Gamma$ one has Fourier expansions  
$$f^V_{\ell,q}\left( g_{\frak c} g\right) 
=\sum_{0\neq\omega\in{\frak O}} c_V^{\frak c}(\omega) 
\left( {\bf J}_{\omega}\varphi_{\ell,q}\left(\nu_V,p_V\right)\right) (g)\qquad\qquad  
\hbox{($g\in G$, $\ell\geq |p_V|$ and $|q|\leq\ell$\/),}\eqno(1.1.8)$$
with coefficients 
$c_V(\omega)$ that, in addition to being independent of $g$, are also  
independent of $\ell$ and $q$.
The system ${\cal B}^V_K$ may be normalised so that,  
for $\ell,q\in{\Bbb Z}$ with $\ell\geq|p_V|$ and $|q|\leq\ell$, one has  
$$\left\| f^V_{\ell,q}\right\|_{\Gamma\backslash G}^2
=\cases{\left\|\Phi_{p_V,q}^{\ell}\right\|_K^2 
&if $\left(\nu_V,p_V\right)\in i[0,\infty)\times{\Bbb Z}\,$,\cr 
 &\hbox{\quad} \cr
{\Gamma(1+\ell -\nu_V)\over\Gamma(1+\ell+\nu_V)}
\left\|\Phi_{0,q}^{\ell}\right\|_K^2 &if $0<\nu_V<1$ and $p_V=0\,$.}\eqno(1.1.9)$$
Subject to this normalisation,  the function $c_V : {\frak O}\rightarrow{\Bbb C}$  
is determined, up to a constant multiplier of absolute value $1$, 
by $V$ and $g_{\frak c}$ alone (the same being true of the system ${\cal B}^V_K$\/). 

It is implicit in the equation (1.1.8) that at all 
cusps ${\frak c}$ of $\Gamma$ one has $F_0^{\frak c} f_{\ell,q}^V(g)=0$ 
for $g\in G$.  Moreover, it follows by [19, Lemma~5.2.1]
that each $f_{\ell,q}^V\in{\cal B}_K^V$ satisfies, 
at every cusp ${\frak c}$ of $\Gamma$,  a growth condition 
$$f_{\ell,q}^V\left( g_{\frak c}n a[r] k\right)\ll_{V,\ell,q,{\frak c}}\ 
r^{\ell +1/2}e^{-\pi r}\qquad\qquad  
\hbox{($\,n\in N,\,k\in K\ $ and 
$\ r\geq R( f_{\ell ,q}^V , {\frak c})\,$),}\eqno(1.1.10)$$
where $R(f,{\frak c})>0$ depends only upon $f$ and ${\frak c}$. 
Any such $\Gamma$-automorphic eigenfunction of 
both Casimir operators $\Omega_{\pm}$ is commonly termed a `cusp form'
(hence the designation of $V$ as a `cuspidal' subspace).

The spectral parameters $\left(\nu_V,p_V\right)$ associated with the decomposition 
(1.1.3) merit some further consideration. Let $V$ be one of the 
relevant cuspidal subspaces. Then, as indicated prior to (1.1.4), 
$V$ is either of the 
`unitary principal series' (i.e. has $\nu_V\in i[0,\infty)$ and $p_V\in{\Bbb Z}$\/), 
or else is of the `complementary series'
(having $0<\nu_V<1$ and $p_V=0$). If $p_V=0$, then the generator 
$f^V_{0,0}$ of $V_{K,0,0}$ satisfies 
$-{\bf \Delta} f^V_{0,0}=\lambda_V f^V_{0,0}$, where $\lambda_V=1-\nu_V^2$ and 
${\bf \Delta}$ is the hyperbolic Laplacian 
operator: 
$${\bf \Delta} =4\left(\Omega_{+}+\Omega_{-}\right)\bigl|_{C^{\infty}(G/K)}\ 
=r^{2}\left( {\partial^2\over\partial x^2}+{\partial^2\over\partial y^2}
+{\partial^2\over\partial r^2}\right) 
-r\,{\partial\over\partial r}\;,$$
with $C^{\infty}(G/K)$ signifying the space of infinitely differentiable functions 
$f : G\rightarrow{\Bbb C}$ which, for 
$k\in K$, $r>0$ and $z=x+iy\in{\Bbb C}$, satisfy $f(n[z]a[r]k)=f(n[z]a[r])$.
By [7, Theorem~1.7] the operator $-{\bf \Delta}$ is symmetric and positive 
on the space $\bigl\{ f\in L^2(\Gamma\backslash G)\cap C^{\infty}(G/K) : 
{\bf \Delta} f\in L^2(\Gamma\backslash G)\bigr\}\supset V_{K,0,0}$, which of itself 
implies $1-\nu_V^2=\lambda_V>0$ (partially explaining why we have $\nu_V\in i[0,\infty)\cup(0,1)$ when 
$p_V=0$\/). Recent work of Kim and Shahidi [14], 
[13, Theorem~4.10] has shown that $\lambda_V >77/81$, so that one has  
$$\left(\nu_V , p_V\right)\in(0 , 2/9)\times\{ 0\}\qquad\quad\hbox{if $\quad V$ is 
of the complementary series.}\eqno(1.1.11)$$
Eigenvalues $\lambda_V<1$ (and these only) are termed `exceptional'. 
Since the group $\Gamma$ here is 
(in the terminology of [7, Chapter~2, Definition~2.3]) 
cofinite but non-cocompact, at most finitely many of the factors $V$ in 
the decomposition (1.1.3) 
correspond to such exceptional eigenvalues of $-{\bf \Delta}$. 
Indeed, by [7, Chapter~4, Corollary~5.3]  
one has 
$\sum_V^* \lambda_V^{-2}<\infty$ 
(the asterisk indicating summation over those of the 
cuspidal subspaces $V$ occurring in (1.1.3) that  
have $p_V=0$).

The generalised Selberg eigenvalue conjecture, if true, would 
(in the present context) entail the complete  
absence of any exceptional eigenvalues:  
so that all $V$ occurring in the decomposition (1.1.3) 
would necessarily be of the unitary principal series.  
Whilst the generalised Selberg eigenvalue conjecture has neither been proved, nor disproved, 
it is known that certain discrete groups, such as $\Gamma_0(1)<SL(2,{\Bbb C})$, 
are not associated (in the manner described above) with any exceptional eigenvalues: 
see [7, Chapter~7, Proposition~6.2] for other examples. 
Yet the current state of knowledge does not, for example, 
rule out the possibility that there may exist an infinite sequence 
of distinct Gaussian primes, $\varpi_1,\varpi_2,\ldots\ $, 
such that each group in the sequence 
$\Gamma_0(\varpi_1),\Gamma_0(\varpi_2),\ldots\ $  
is associated with at least one exceptional eigenvalue of $-{\bf \Delta}$. 
It is fair to say that any proof of the generalised Selberg eigenvalue conjecture 
(or of just those cases of it that are relevant) 
would render much of this paper obsolete.   

The subspace ${}^{\rm e}L^2(\Gamma\backslash G)$ in (1.1.2) 
is a special case 
of the space 
referred to in [19, Chapter~8] as `$L^{2,{\rm cont}}(\Gamma\backslash G ,\chi)$',   
and (as noted there) is generated by integrals of certain Eisenstein series. 
In determining a suitable set of such generators it helps to note  
that, by  
the relation of $\Gamma$-equivalence of cusps 
(whereby ${\frak a}$ is deemed
$\Gamma$-equivalent to ${\frak b}$ if and only if $\gamma{\frak a}={\frak b}$ for some 
$\gamma\in\Gamma$), 
the set ${\Bbb P}^1({\Bbb Q}(i))$ of cusps of $\Gamma$  
is partitioned into finitely many $\Gamma$-equivalence classes, each of form 
$\{\gamma{\frak c} : \gamma\in\Gamma\}$ for some ${\frak c}\in {\Bbb P}^1({\Bbb Q}(i))$. 
We shall use the notation 
${\frak a}\sim^{\!\!\!\!\Gamma}{\frak b}$ to signify that 
${\frak a}$ is $\Gamma$-equivalent to ${\frak b}$. 
Let ${\frak C}$ be a complete set of representatives 
of the $\Gamma$-equivalence classes of cusps in ${\Bbb P}^1({\Bbb Q}(i))$. 
Then, for ${\frak c}\in{\frak C}$, $\ell,p,q\in{\Bbb Z}$ with $\ell\geq |p|,|q|$    
and $\nu\in{\Bbb C}$ with ${\rm Re}(\nu)>1$, the Eisenstein series
$E_{\ell,q}^{\frak c}(\nu,p) : G\rightarrow{\Bbb C}$ is given by: 
$$E_{\ell,q}^{\frak c}(\nu,p)(g)
={1\over\left[\Gamma_{\frak c} : \Gamma_{\frak c}'\right]} 
\sum_{\Gamma_{\frak c}'\gamma\in\Gamma_{\frak c}'\backslash\Gamma}
\varphi_{\ell,q}(\nu,p)\left( g_{\frak c}^{-1}\gamma g\right)\qquad\qquad 
\hbox{($g\in G$),}\eqno(1.1.12)$$
where $\varphi_{\ell,q}(\nu,p)$ is as defined in (1.1.7). 
By virtue of (1.1.1), the sum in (1.1.12) is well-defined.  
Moreover, a property of the function $\varphi_{\ell,q}(\nu,p) : G\rightarrow{\Bbb C}$ 
ensures that if 
$\Gamma_{\frak c}\neq\Gamma_{\frak c}'$, and if $p$ is odd, 
then the terms of that sum cancell one another out; since 
$\left[\Gamma_{\frak c}:\Gamma_{\frak c}'\right]\in\{ 2,4\}$ for 
${\frak c}\in{\Bbb P}^1({\Bbb Q}(i))$, one therefore has
$$E_{\ell,q}^{\frak c}(\nu,p) \neq 0\qquad\ \hbox{only if}\qquad\  
p\in\textstyle{1\over 2}\left[\Gamma_{\frak c} : \Gamma_{\frak c}'\right]{\Bbb Z}\;.$$
The condition ${\rm Re}(\nu)>1$ ensures absolute convergence of 
the sum in (1.1.12): this, and more delicate issues of convergence, are discussed in 
[22, Subsection~1.8] (but see also [7, Chapter~3], [4, Section~5] or [19, Section~3.3]). 
Here it suffices to record that the Eisenstein series given by (1.1.12) are 
infinitely differentiable $\Gamma$-automorphic functions on $G$, and 
inherit from $\varphi_{\ell,q}(\nu,p)$  
the property of being 
eigenfunctions of both Casimir operators $\Omega_{\pm}$ with 
corresponding eigenvalues $\,{\textstyle{1\over 8}}((\nu\mp p)^2-1)\,$.

In parallel with the Fourier expansions (1.1.8) we shall need also 
the Fourier expansions of the Eisenstein series. Preparatory to this we now 
define certain `generalised Kloosterman sums'. Given any pair of 
cusps ${\frak a},{\frak b}\in{\Bbb P}^1({\Bbb Q}(i))$, let 
$${}^{\frak a}\Gamma^{{\frak b}}(c)
=\left\{\gamma\in\Gamma : 
g_{\frak a}^{-1}\gamma g_{{\frak b}}=\pmatrix{*&*\cr c&*}\right\}\qquad\quad 
\hbox{for $\quad c\in{\Bbb C}\,$,}\eqno(1.1.13)$$
and put 
$${}^{\frak a}{\cal C}^{{\frak b}}
=\left\{ c\in{\Bbb C}-\{ 0\} : {}^{\frak a}\Gamma^{{\frak b}}(c)\neq\emptyset\right\}\;.\eqno(1.1.14)$$
Then, for $c\in {}^{\frak a}{\cal C}^{{\frak b}}$ and $\omega,\omega'\in{\frak O}$, 
the generalised Kloosterman sum 
$S_{{\frak a} , {\frak b}}\left(\omega , \omega' ; c\right)$
is given by: 
$$S_{{\frak a} , {\frak b}}\left(\omega , \omega' ; c\right) 
=\sum_{\scriptstyle \Gamma_{\frak a}'\gamma\Gamma_{{\frak b}}'\in\Gamma_{\frak a}'\backslash 
{}^{\frak a}\Gamma^{{\frak b}}(c)/\Gamma_{{\frak b}}'\atop\scriptstyle 
g_{\frak a}^{-1} \gamma g_{{\frak b}}
=\left( {\scriptstyle\!\!\!\!\!s(\gamma)\quad *\atop\scriptstyle\  c\quad\ d(\gamma)}\right) }
{\rm e}\!\left({\rm Re}\!\left(\omega\,{s(\gamma)\over c}+\omega'\,{d(\gamma)\over c}\right)\right) ,\eqno(1.1.15)$$
where ${\rm e}(x)=\exp(2\pi ix)$. 
If ${\frak a},{\frak b}\in{\frak C}$, and if  $\ell,p,q\in{\Bbb Z}$ and $\nu\in{\Bbb C}$
are such that $\ell\geq |p|,|q|$ and  ${\rm Re}(\nu)>1$, then 
at the cusp ${\frak b}$ the Eisenstein series $E_{\ell,q}^{\frak a}(\nu,p)$
has the Fourier expansion 
$$\eqalign{
E_{\ell,q}^{\frak a}(\nu,p)\left( g_{\frak b} g\right) 
 &=\delta^{*}_{{\frak a},{\frak b}}\varphi_{\ell,q}(\nu,p)( g) 
+{1\over \left[\Gamma_{\frak a} : \Gamma_{\frak a}'\right]}\,D_{\frak a}^{\frak b}(0;\nu,p)
\,{\pi\Gamma(|p|+\nu)\over\Gamma(\ell+1+\nu)}\,{\Gamma(\ell+1-\nu)\over\Gamma(|p|+1-\nu)}\,
\varphi_{\ell,q}(-\nu,-p)(g) + {}\cr 
 &\qquad +{1\over \left[\Gamma_{\frak a} : \Gamma_{\frak a}'\right]}
\sum_{0\neq\psi\in{\frak O}}
D_{\frak a}^{\frak b}\left(\psi;\nu,p\right) 
{\bf J}_{\psi}\varphi_{\ell,q}(\nu,p)(g)\;,}\eqno(1.1.16)$$
where
$$\delta^{*}_{{\frak a},{\frak b}}=\cases{1 &if ${\frak a}\sim^{\!\!\!\!\Gamma}{\frak b}$,\cr 
0 &otherwise,}\eqno(1.1.17)$$
and
$$D_{\frak a}^{\frak b}(\psi;\nu,p)
=\sum_{c\in\,{}^{\frak a}{\cal C}^{\frak b}}^{(\Gamma)}
S_{{\frak a} , {\frak b}}\!\left( 0 , \psi ; c\right) |c|^{-2(1+\nu)}\left( c/|c|\right)^{2p}\eqno(1.1.18)$$ 
(all the sums here being absolutely convergent). 
Using an evaluation
of  $S_{{\frak a} , {\frak b}}\!\left( 0 , \psi ; c\right)\,$ 
(analogous to the classical evaluation [8, Theorem~271] of Ramanujan's sum)
it can be shown that 
when $0\neq\psi\in{\frak O}$ the right-hand side of (1.1.18) converges absolutely
for ${\rm Re}(\nu)>0$. It may, on the other hand, be deduced from (1.1.12)
that when ${\frak a}$, $\ell$, $p$, $q$ and $g$ are given, the function 
$\nu\mapsto E_{\ell,q}^{\frak a}(\nu,p)(g)$ is holomorphic for ${\rm Re}(\nu)>1$
(see [7, Chapter~3, Proposition~2.5] for the case $p=q=\ell =0$\/); 
and it is known that this function of $\nu$ has a meromorphic 
continuation to all of ${\Bbb C}$, with a simple pole at $\nu =1$ 
if and only if $p=q=\ell =0$, and with no other poles in 
the closed half plane $\{\nu\in{\Bbb C} : {\rm Re}(\nu)\geq 0\}$. 
This may be shown by application of 
Langlands' general theory [18], or by expressing the coefficients  
$D_{\frak a}^{\infty}(\psi;\nu,p)$ ($\psi\in{\frak O}$\/) 
in terms of Hecke zeta-functions: [4, Lemma~5.2] being a prototypic example of the latter approach. 
Applying this meromorphic continuation one obtains, when 
${\rm Re}(\nu)\geq 0$ and $(\nu ,p)\neq (1,0)$, 
an infinitely differentiable $\Gamma$-automorphic function 
$E_{\ell,q}^{\frak a}(\nu,p) : G\rightarrow{\Bbb C}$ satisfying 
$\Omega_{\pm}E_{\ell,q}^{\frak a}(\nu,p)={1\over 8}\bigl( (\nu\mp p)^2-1\bigr)
E_{\ell,q}^{\frak a}(\nu,p)$.

Because of the behaviour 
(as the Iwasawa coordinate $r$ tends to $\infty$) 
of first two terms on the right-hand side of (1.1.16), one has 
$E_{\ell,q}^{\frak a}(\nu,p)\not\in L^2(\Gamma\backslash G)$.  
Nevertheless, by averaging $E_{\ell,q}^{\frak a}(\nu,p)$ over a range of values 
of $\nu\in i{\Bbb R}$ one can obtain a suitable generator in the space 
${}^{e}L^{2}(\Gamma\backslash G)$. Indeed, by an extension of 
[7, Chapter~6, Theorem~3.2], one has:
$${}^{e}L^{2}(\Gamma\backslash G)
=\overline{\bigoplus_{{\frak c}\in{\frak C}}^{(\Gamma)} 
\bigoplus_{\scriptstyle\ell,q\in{\Bbb Z}\atop\scriptstyle 
\ell\geq |q|}
\bigoplus_{\scriptstyle 
p\in{1\over 2}\left[\Gamma_{\frak c} : \Gamma_{\frak c}'\right]{\Bbb Z}\atop\scriptstyle 
|p|\leq\ell} 
\Biggl\{\ \int\limits_0^{\infty} E_{\ell,q}^{\frak c}(it,p)(g) H(t) {\rm d}t : 
\,H\in L^2(0,\infty)\Biggr\}}\;.\eqno(1.1.19)$$
Equations~(1.1.2), (1.1.3), (1.1.5), (1.1.6) and~(1.1.19)  
describe the spectral decomposition of $L^2(\Gamma\backslash G)$:  
for the subspace 
$$L^2(\Gamma\backslash G;\ell,q)=
\left\{ f\in L^2(\Gamma\backslash G)\, : 
\,\Omega_K f=-{(\ell +1)\ell\over 2}\,f\ \ {\rm and}\ \   
{\partial\over\partial\psi}\,f=-iqf\right\}\eqno(1.1.20)$$ 
(where we assume that $\ell,q\in{\Bbb Z}$ and $\ell\geq |q|$), one has  
a corresponding Parseval identity [22, Theorem~A],  
which is a special case of [19, Theorem~8.1]. 

Apart from the coefficient $\delta^{*}_{{\frak a},{\frak b}}$ in (1.1.16) possibly 
being replaced by $\epsilon^p \delta^{*}_{{\frak a},{\frak b}}$, for some $\epsilon\in{\frak O}^{*}$, 
the Fourier expansion (1.1.16)-(1.1.18) is valid for arbitrary cusps 
${\frak a},{\frak b}\,$ of $\Gamma$ (i.e. not only for 
${\frak a},{\frak b}\in{\frak C}$). Indeed, when 
${\frak c}\sim^{\!\!\!\!\Gamma}{\frak d}$, there will exist  
a unit $\epsilon\in{\frak O}^{*}$ such that 
$E_{\ell,q}^{\frak c}(\nu,p)=\epsilon^p E_{\ell,q}^{\frak d}(\nu,p)$ 
for $\ell,p,q\in{\Bbb Z}$ with $\ell\geq\max\{ |p|,|q|\}$ and all $\nu\in{\Bbb C}$ 
that are not poles of $E_{\ell,q}^{\frak c}(\nu,p)$ (the set of 
poles of the function $\nu\mapsto E_{\ell,q}^{\frak c}(\nu,p)(g)$ 
being independent of the variable $g$\/). 

The meromorphic continuation of the function $\nu\mapsto E_{\ell,q}^{\frak a}(\nu,p)(g)$ 
implies a corresponding meromorphic continuation of each term 
$( F_{\psi}^{\frak b}E_{\ell,q}^{\frak a}(\nu,p))(g)\,$  
occurring in the Fourier expansion 
of the Eisenstein series $E_{\ell,q}^{\frak a}(\nu,p)$ at the (arbitrary) cusp ${\frak b}$, and 
hence the meromophic continuation over ${\Bbb C}$ of the function $\nu\mapsto D_{\frak a}^{\frak b}(\psi;\nu,p)$, 
given (for ${\rm Re}(\nu)>1$) by (1.1.18). 
Let this meromorphic continuation define $D_{\frak a}^{\frak b}(\psi;\nu,p)$ 
when ${\rm Re}(\nu)\leq 1$ and $\nu$ is not a pole. Then (from the above discussion) 
the function 
$\nu\mapsto D_{\frak a}^{\frak b}(\psi;\nu,p)$ either has no poles   
in the closed half plane $\{ \nu\in{\Bbb C} : {\rm Re}(\nu)\geq 0\}$, or 
has there just the one simple pole, at $\nu =1$: it is, in particular, 
holomorphic at all points $\nu\in i{\Bbb R}$. 

In the next three subsections we present our main results. 
These may be more concisely expressed in terms of modified 
Fourier coefficients, $c_V^{\frak c}\left(\omega;\nu_V,p_V\right)$ and 
$B_{\frak a}^{\frak b}(\psi;\nu,p)$, which, 
for ${\frak a},{\frak b},{\frak c}\in{\Bbb P}^1({\Bbb Q}(i))$, 
$0\neq\omega\in{\frak O}$, any cuspidal subspace $V$ occurring as 
a factor in (1.1.3), any $p\in{\Bbb Z}$ and 
any $\nu\in{\Bbb C}$ that is not a pole of $D_{\frak a}^{\frak b}(\omega;\nu,p)$, 
are given by:
$$c_V^{\frak c}\!\left(\omega;\nu_V,p_V\right) 
=\left(\pi|\omega|\right)^{\nu_V}\left({\omega/|\omega|}\right)^{-p_V}
c_V^{\frak c}(\omega)\quad\,\hbox{and}\quad\,  
B_{\frak a}^{\frak b}(\omega;\nu,p)
=\left(\pi|\omega|\right)^{\nu}\left(\omega/|\omega|\right)^{-p}
D_{\frak a}^{\frak b}(\omega;\nu,p)\;.\eqno(1.1.21)$$ 
Note that these modified coefficients, and the generalised Kloosterman sums 
defined by (1.1.13)-(1.1.15), are to a large extent determined 
by the $\Gamma$-equivalence classes of the relevant cusps 
(rather than by the cusps themselves, or by the choice of scaling matrices). 
This follows from 
the fact that ${\frak a}\sim^{\!\!\!\!\Gamma}{\frak b}$ if and only if
the subset ${}^{\frak a}\Gamma^{\frak b}(0)\subset\Gamma$ given by (1.1.13) is non-empty. 
For, given the requirement that (1.1.1) holds for all cusps~${\frak c}$,  
one can (by a calculation) show that if ${}^{\frak a}\Gamma^{\frak u}(0)\ni\gamma_1$ 
and ${}^{\frak b}\Gamma^{\frak v}(0)\ni\gamma_2$ then, 
for some $\beta_1,\beta_2\in{\Bbb C}$, some $\epsilon_1,\epsilon_2\in{\frak O}^{*}$ and 
some $\eta_1,\eta_2\in{\Bbb C}^{*}$  with $\eta_j^2 =\epsilon_j$ ($j=1,2$\/), 
one will have   
$$g_{{\frak a}}^{-1}\gamma_1 g_{\frak u}=h\left[\eta_1\right] n\left[\beta_1\right]\;,\qquad\quad   
g_{{\frak b}}^{-1}\gamma_2 g_{\frak v}=h\left[\eta_2\right] n\left[\beta_2\right]\qquad\quad    
{\rm and}\qquad\    
{}^{\frak u}{\cal C}^{\frak v}=\eta_1\eta_2 {}^{{\frak a}}{\cal C}^{{\frak b}}$$
(where the sets ${}^{\frak a}{\cal C}^{\frak b} , {}^{\frak u}{\cal C}^{\frak v}$ are defined by (1.1.14))   
and, for $\omega_1,\omega_2,\in{\frak O}$ and $0\neq\omega\in{\frak O}$, the identities:    
$$S_{{\frak u},{\frak v}}\left(\omega_1,\omega_2;c_1\right) 
={\rm e}\left( {\rm Re}\left(\beta_2\omega_2 -\beta_1\omega_1\right)\right) 
S_{{\frak a},{\frak b}}\left(\overline{\epsilon_1}\,\omega_1 , \overline{\epsilon_2}\,\omega_2 , 
\overline{\eta_1\eta_2}\,c_1\right)\qquad\qquad   
\hbox{($c_1\in {}^{\frak u}{\cal C}^{\frak v}$\/),}$$ 
$$c_V^{\frak u}(\omega;\nu_V,p_V)
={\rm e}\left({\rm Re}(\beta_1\omega)\right) 
c_V^{{\frak a}}\left( \overline{\epsilon_1}\,\omega ;\nu_V,p_V\right)\qquad\  
{\rm and}\qquad\  
B_{\frak u}^{\frak v}(\omega;\nu,p)
=\epsilon_1^p\,{\rm e}\left( {\rm Re}\left(\beta_2\omega\right)\right) 
B_{{\frak a}}^{{\frak b}}\left(\overline{\epsilon_2}\,\omega ;\nu,p\right)$$
(the latter pair being valid for any cuspidal subspace $V$ occurring 
in (1.1.3), and any $(\nu,p)\in{\Bbb C}\times{\Bbb Z}$ 
such that $\nu$ is not a pole of $B_{\frak u}^{\frak v}(\omega;\nu,p)$\/).

\bigskip 
\bigskip 

\goodbreak\centerline{\bf \S 1.2. A Kloosterman to Spectral Sum Formula and Other Key Ingredients}

\bigskip 

An essential underlying 
component of the proofs of the principal new results of this paper 
is the following 
`Kloosterman to spectral' summation formula for $SL(2,{\Bbb C})$, which is 
analogous to (though in some ways simpler than) the summation formula 
for $SL(2,{\Bbb R})$ of Kuznetsov [15,16].  
Before stating this formula     
it is worth clarifying, firstly,    
that when $D$ is an open subset of $\,{\Bbb C}$, a function $f : D\rightarrow{\Bbb C}$ 
may be termed `smooth' if and only if  
each of the functions $u(x,y)={\rm Re}(f(x+iy))$ and 
$v(x,y)={\rm Im}(f(x+iy))$ (both having the set $D'=\{ (x,y)\in{\Bbb R}^2 : x+iy\in D\}$
as their domain) is such that, for all $n\in{\Bbb N}$, 
every one of its $2^n$ partial derivatives of order $n$ is a continuous 
real-valued function on $D'$.  
Secondly,  
for $\sigma\in{\Bbb R}$, we 
use the subscript `$(\sigma)$' to denote integration 
from $\sigma -i\infty$ to $\sigma +i\infty$ along the contour   
${\rm Re}(z)=\sigma$, so that  if 
$f$ is a complex function such that the function $t\mapsto f(\sigma +it)$ 
is Lebesgue integrable on $(-\infty ,\infty)$ then 
$$\int\limits_{(\sigma)} f(z) {\rm d}z
=i\int\limits_{-\infty}^{\infty} f(\sigma +it) {\rm d}t\;.$$

\goodbreak\proclaim Theorem~1 (A Kloosterman to spectral sum formula). 
Let $f : {\Bbb C}^{*}\rightarrow{\Bbb C}$ be an even smooth function 
compactly supported in ${\Bbb C}^{*}$. 
Suppose moreover that $q,\omega_1,\omega_2\in {\frak O}={\Bbb Z}[i]$, 
with $q\omega_1\omega_2\neq 0$; and that 
${\frak C}$ is a complete set of representatives for 
the $\Gamma$-equivalence classes of cusps for the Hecke congruence subgroup 
$\Gamma=\Gamma_0(q)\leq SL(2,{\frak O})$. Then, for all pairs of cusps ${\frak a},{\frak b}$ of  
$\Gamma$ (and all associated pairs of 
scaling matrices $g_{\frak a}, g_{\frak b}\in SL(2,{\Bbb C})$ such that (1.1.1)
holds for ${\frak c}={\frak a}$ and for ${\frak c}={\frak b}$\/), one has   
$$\eqalignno{ 
&\sum_{c\in {}^{\frak a}{\cal C}^{\frak b}}^{(\Gamma)}  
\,{S_{{\frak a},{\frak b}}\left(\omega_1 , \omega_2 ; c\right)\over |c|^2}\,
f\left( {2\pi\sqrt{\omega_1\omega_2}\over c}\right) = {}\qquad 
\qquad\qquad\qquad\qquad\qquad\qquad\qquad\qquad\qquad\qquad\qquad\qquad\qquad\quad &(1.2.1)\cr
 &\qquad\qquad\qquad\qquad\qquad        
=\pi\sum_V^{(\Gamma)} \,\overline{c_V^{\frak a}\left(\omega_1;\nu_V,p_V\right)}\,
c_V^{\frak b}\left(\omega_2;\nu_V,p_V\right) {\bf K}f\left(\nu_V , p_V\right) + {}\cr 
 &\qquad\qquad\qquad\qquad\qquad\quad\ {} +(-i/4)\sum_{{\frak c}\in{\frak C}}^{(\Gamma)} 
{1\over\left[\Gamma_{\frak c} : \Gamma_{\frak c}'\right]} 
\sum_{p\in{1\over 2}\left[\Gamma_{\frak c} : \Gamma_{\frak c}'\right]{\Bbb Z}}
\ \int\limits_{(0)} \overline{B_{\frak c}^{\frak a}\left(\omega_1;\nu,p\right)}\,
B_{\frak c}^{\frak b}\left(\omega_2;\nu,p\right) {\bf K}f(\nu,p)\,{\rm d}\nu\;,}$$
where ${}^{\frak a}{\cal C}^{\frak b}$ and the generalised 
Kloosterman sums $S_{{\frak a},{\frak b}}\left(\omega_1 , \omega_2 ; c\right)$ 
are as defined in (1.1.13)-(1.1.15); 
where the system of irreducible cuspidal subspaces $V\subset L^2(\Gamma\backslash G)$,  
spectral parameters $(\nu_V,p_V)\in{\Bbb C}\times{\Bbb Z}$ and  
modified Fourier coefficients $c_V^{\frak d}\left(\omega;\nu_V,p_V\right)$ and  
$B_{\frak c}^{\frak d}\left(\omega;\nu,p\right)$ 
are as described in (1.1.2)-(1.1.11) and (1.1.16)-(1.1.21), while the subgroups 
$\Gamma_{\frak c}'\leq\Gamma_{\frak c}<\Gamma$ 
are as defined above (1.1.1); and where (as in [4, Theorem~10.1]) 
one defines the ${\bf K}$-transform by:   
$${\bf K}f(\nu,p)
=\int_{{\Bbb C}^{*}} {\cal K}_{\nu,p}(z) f(z)\,{\rm d}_{\times}z\;,\eqno(1.2.2)$$
with $\,{\rm d}_{\times}z=|z|^{-2}{\rm d}_{+}z=|z|^{-2}{\rm d}x{\rm d}y\,$ 
(for $0\neq z=x+iy$ and $x,y\in{\Bbb R}$\/), 
$${\cal K}_{\nu,p}(z)={{\cal J}_{-\nu,-p}(z)-{\cal J}_{\nu,p}(z)\over\sin(\pi\nu)}\;,\qquad\quad  
{\cal J}_{\mu,k}(z)
=\left| {z\over 2}\right|^{2\mu}\left( {z\over |z|}\right)^{\!\!-2k}J_{\mu -k}^{*}(z)
J_{\mu +k}^{*}\left( \overline{z}\right)\;,\eqno(1.2.3)$$
and
$$J_{\xi}^{*}(w)
=\sum_{m=0}^{\infty} {(-1)^m\,(w/2)^{2m}\over\Gamma(m+1)\Gamma(\xi +m+1)}\;.\eqno(1.2.4)$$
Nothing more than (1.1.1) need be assumed in respect of the  
scaling matrices $g_{\frak c}\in SL(2,{\Bbb C})$ chosen for ${\frak c}\in{\frak C}$, 
even when  ${\frak C}\cap\{ {\frak a},{\frak b}\}\neq\emptyset$. Similarly, 
$g_{\frak a}$ may differ from $g_{\frak b}$, even when ${\frak a}={\frak b}$.

\medskip 

\goodbreak
{\bf Proof.}\quad This theorem is a minor extension of 
Lokvenec-Guleska's result [19, Theorem~12.3.2], which applies only to 
the case ${\frak a}={\frak b}=\infty$ (though being, in other important respects 
considerably more general than our theorem). 
The proof is a straightforward application of [22, Theorem~B] 
(a spectral to Kloosterman summation formula, generalising [4, Theorem~10.1] of 
Bruggeman and Motohashi, and 
extending [19, Theorem~11.3.3] of Lokvenec-Guleska), in combination with 
Bruggeman and Motohashi's one-sided ${\bf B}$-transform inversion formula [4, Theorem~11.1] and 
`annihilation lemma' [4, Lemma~11.1].
The ${\bf B}$-transform in question maps any suitable complex-valued two-variable 
function $h(\nu,p)$ to 
the function ${\bf B}h : {\Bbb C}^{*}\rightarrow{\Bbb C}$ given by:  
$${\bf B} h(z)
=\sum_{p\in{\Bbb Z}} {1\over 4\pi i}\int\limits_{(0)} 
{\cal K}_{\nu,p}(z) h(\nu,p)\left( p^2 -\nu^2\right) {\rm d}\nu\;.$$
Subject to our hypotheses concerning $f$, we 
have (see [4, Theorem~11.1]) the one-sided inversion formula:  
$$\pi {\bf BK}f = f\;.\eqno(1.2.5)$$
In addition, 
[4, Lemma~11.1] shows that 
$$\sum_{p\in{\Bbb Z}}\ \int\limits_{(0)} 
{\bf K}f (\nu,p)\left( p^2 -\nu^2\right) {\rm d}\nu =0\;,\eqno(1.2.6)$$

To prove our theorem we need only verify that, for some 
$\sigma >1/2$, the function $h={\bf K}f$ satisfies the 
hypotheses (i)-(iii) of [22, Theorem~B]: for then the 
equation (1.2.1) follows by the direct use of (1.2.6) and (1.2.5) 
to effect appropriate substitutions 
in the case $h={\bf K}f$ of [22, Theorem~B, Equation~(1.9.1)].
Those hypotheses are satisfied by $h={\bf K}f\,$ if, 
when $S_{\sigma}=\{\nu\in{\Bbb C} : |{\rm Re}(\nu)|\leq\sigma\}$, 
one has all of the following:\smallskip 
\hskip 15mm (i)\quad ${\bf K}f(\nu,p)={\bf K}f(-\nu ,-p)$\quad 
for $(\nu,p)\in S_{\sigma}\times{\Bbb Z}\;$;\hfill\smallskip 
\hskip 14.5mm (ii)\quad $\nu\mapsto {\bf K}f(\nu,p)$ is holomorphic on a neighbourhood of the strip 
$\,S_{\sigma}\;$;\smallskip 
\hskip 14mm (iii)\quad ${\bf K}f(\nu,p)\ll_{f,\sigma} 
(1+|{\rm Im}(\nu)|)^{-4} (1+|p|)^{-4}\;$\quad (say)\quad  
for $(\nu,p)\in S_{\sigma}\times{\Bbb Z}\;$. \medskip 

Note firstly that by (1.2.3) and (1.2.4) the functions 
$\mu\mapsto{\cal J}_{\mu,k}(z)$  are entire. 
Moreover, by using the relations 
$J_{\xi}=(-1)^{\xi} J_{-\xi}\,$ ($\xi\in{\Bbb Z}$\/) 
satisfied by the $J$-Bessel function 
$J_{\xi}(z)=(z/2)^{\xi} J_{\xi}^{*}(z)$, 
one may show that ${\cal J}_{-\nu,-p}={\cal J}_{\nu,p}$ when both $p$ and $\nu$ 
are integers. By this and the first equation in (1.2.3) it follows 
that the functions $\nu\mapsto{\cal K}_{\nu,p}(z)$ are entire 
(the singularities at $\nu\in{\Bbb Z}$ being removable). 
In addition, since it is also the case that the functions 
$z\mapsto J_{\xi}^{*}(z)$ are entire, each function 
$z\mapsto{\cal J}_{\mu,k}(z)$ is continuous on ${\Bbb C}^{*}$; and so the 
same is true of the 
functions $z\mapsto{\cal K}_{\nu,p}(z)$. Therefore (given 
that $f$ is compactly supported in ${\Bbb C}^{*}$\/)  
it follows by the definition (1.2.2) and the holomorphicity 
of the functions $\nu\mapsto{\cal K}_{\nu,p}(z)$ that, 
for each $p\in{\Bbb Z}$ the function $\nu\mapsto{\bf K}f(\nu,p)$ is entire. 
This has verified that the condition~(ii) above is satisfied. 

By [19, Lemma~12.1.1, Estimate (12.24)], the condition~(iii) is 
satisfied for all $\sigma >0$. 
Finally, since the condition~(ii) has already been verified, 
the condition~(i) is a trivial consequence of (1.2.2) and 
the relation ${\cal K}_{\nu,p}={\cal K}_{-\nu,-p}$ 
implicit in the first equation of (1.2.3). The proof is now complete, 
for it has been shown that the conditions~(i), (ii) and~(iii) 
hold for all $\sigma >0$ (and so certainly for some $\sigma >1/2$\/) 
\quad$\blacksquare$   

\bigskip

\goodbreak 
\noindent{\bf Remark~1.}\quad 
The above inversion of the summation formula
[22, Theorem~B] is one-sided (i.e. non-surjective):  
for it contains no `diagonal term' (i.e. no counterpart of the term 
in [22, Equation~(1.9.1)] with 
coefficient $\delta_{\omega_1,\omega_2}^{{\frak a},{\frak b}}$\/),   
whereas, as is pointed out in [4, Section~11], there exist test functions $h$ 
satisfying the conditions~(i)-(iii) of [22, Theorem~B] that do produce a non-zero 
diagonal term on the right-hand side of [22, Equation~(1.9.1)]. 

\bigskip

We next state the principal new result of [22], followed by a very useful corollary. 

\bigskip

\goodbreak\proclaim Theorem~2. 
Let $\varepsilon >0$, $0\neq q\in{\frak O}={\Bbb Z}[i]$, 
$\Gamma =\Gamma_0(q)\leq SL(2,{\frak O})$ and 
$K,P,N\geq 1$.  
Suppose further that $b_{n}\in{\Bbb C}$ for $n\in{\frak O}-\{ 0\}$, and that 
$u,w\in{\frak O}$ satisfy $w\neq 0$ and $(u,w)\sim 1$ (i.e. that 
$u$ and $w$ are coprime). 
Then, when
${\frak a}$ is a cusp of $\Gamma$ 
with ${\frak a}\sim^{\!\!\!\!\Gamma}u/w$, and when 
$E_{0}^{\frak a}\left( q,P,K;N,{\bf b}\right)$ and $E_{1}^{\frak a}\left( q,P,K;N,{\bf b}\right)$
are given by  
$$E_{0}^{\frak a}\left( q,P,K;N,{\bf b}\right) 
=\sum_{\scriptstyle V\atop\scriptstyle 
\left| p_V\right|\leq P ,\ \left|\nu_V\right|\leq K}^{(\Gamma)} 
\Bigl|\sum_{\scriptstyle n\in{\frak O}\atop\scriptstyle N/2<|n|^2\leq N}
b_{n} c_V^{\frak a}\left( n;\nu_V,p_V\right)\Bigr|^2\;,\eqno(1.2.7)$$
$$E_{1}^{\frak a}\left( q,P,K;N,{\bf b}\right)=\sum_{{\frak c}\in{\frak C}}^{(\Gamma)} 
{1\over 4\pi\left[\Gamma_{\frak c} : \Gamma_{\frak c}'\right]}
\sum_{\scriptstyle p\in{1\over 2}\left[\Gamma_{\frak c} : \Gamma_{\frak c}'\right]{\Bbb Z}\atop\scriptstyle 
|p|\leq P}
\int_{-K}^{K}\Bigl|\sum_{\scriptstyle n\in{\frak O}\atop\scriptstyle N/2<|n|^2\leq N}
b_{n} B_{\frak c}^{\frak a}(n;it,p)\Bigr|^2 {\rm d}t\eqno(1.2.8)$$
(where the terminology used has the same meaning as in Theorem~1), one has the upper bounds:    
$$E_{j}^{\frak a}\left( q,P,K;N,{\bf b}\right)
\ll \left( P^2+K^2\right)\left(PK+
O_{\varepsilon}\!\!\left( {N^{1+\varepsilon}\over (PK)^{1/2}}\,|\mu({\frak a})|^2\right)\right)
\left\|{\bf b}_N\right\|_2^2
\qquad\quad\ \hbox{($j=0,1$)}\,,\eqno(1.2.9)$$
where $\,\mu({\frak a})\in\{ 1/\alpha : 0\neq \alpha\in{\frak O}\}$,  
$${1\over\mu({\frak a})}\sim{(w,q)q\over \left( w^2 , q\right)}
\sim{q\over\bigl( (w,q) , q/(w,q)\bigr)}\eqno(1.2.10)$$
and 
$$\left\|{\bf b}_N\right\|_2
=\Biggl( \sum_{\scriptstyle n\in{\frak O}\atop\scriptstyle 
0<|n|^2\leq N}\!\!\!\!\!\!|b_{n}|^2\Biggr)^{1/2}\,.\eqno(1.2.11)$$

\bigskip

\goodbreak 
\noindent{\bf Proof.}\quad This is [22, Theorem~1]: the modification, in (1.2.11), of 
the notation defined in [22, (1.9.17)] is of no significance here, but does 
help in stating other results below (Theorems~5, 7, 10 and~11 for example)\quad$\blacksquare$

\bigskip

\goodbreak 
\noindent{\bf Remark~2.}\quad 
Since $\left[\Gamma_{\frak c} : \Gamma_{\frak c}'\right]\in\{ 2 , 4\}$ for 
all cusps ${\frak c}$ of $\Gamma$, one may omit the factor 
$(4\pi\left[\Gamma_{\frak c} : \Gamma_{\frak c}'\right])^{-1}$ in (1.2.8).

\medskip

\goodbreak
\noindent{\bf Remark~3.}\quad 
The factor $|\mu({\frak a})|^2$ in the bound (1.2.9) has its origin 
in [22, Proposition~1, (1.9.18)], where it is established that, for each cusp 
${\frak a}$ of $\Gamma$, the set ${}^{\frak a}{\cal C}^{\frak a}$ 
defined by (1.1.14) satisfies 
$${}^{\frak a}{\cal C}^{\frak a}
\subset {1\over\mu({\frak a})}\,{\frak O} -\{ 0\}\;.\eqno(1.2.12)$$
For $\mu({\frak a})$ as in (1.2.10), the ideal $(1/\mu({\frak a})){\frak O}$  
and absolute value $|\mu({\frak a})|$ are   
determined by the $\Gamma$-equivalence class of the cusp ${\frak a}$. 
Since $\infty\sim^{\!\!\!\!\Gamma}1/q\ $ (for $\Gamma =\Gamma_0(q)$), 
one has in particular $1/\mu(\infty)\sim 1/\mu(1/q)\sim q$. 

\bigskip

\goodbreak\proclaim Corollary to Theorems~1 and~2. Let all the hypotheses of the 
case $\omega_1=\omega_2=1$ of Theorem~1 hold. 
Suppose, moreover, that $A>1$, $\varepsilon >0$, $M,N\geq 1$ and $X\geq 2$; 
and suppose that one has  
$$f(z)=\varphi( |z| )\qquad\qquad\hbox{($z\in{\Bbb C}^{*}$\/),}\eqno(1.2.13)$$
where the function $\varphi : (0,\infty)\rightarrow{\Bbb C}$ is infinitely differentiable,  
and has its support contained in the interval 
$\bigl[ A^{-1}X^{-1/2} , A X^{-1/2}\bigr]$. 
Put
$$Y=X^{-3/2}\max_{r>0}\left|\varphi^{(3)}(r)\right|\;.\eqno(1.2.14)$$
Then, for all pairs of cusps ${\frak a},{\frak b}$ of  
$\Gamma$, 
for all choices of scaling matrices $g_{\frak a}, g_{\frak b}\in SL(2,{\Bbb C})$ 
such that (1.1.1)
holds for ${\frak c}={\frak a}$ and for ${\frak c}={\frak b}$, 
and for arbitrary complex coefficients $\overline{a_m} , b_n$ 
($0\neq m,n\in{\frak O}$\/), one has   
$$\eqalignno{ 
&\sum_{M/2<|m|^2\leq M}\overline{a_m}\sum_{N/2<|n|^2\leq N}
b_n\sum_{c\in {}^{\frak a}{\cal C}^{\frak b}}^{(\Gamma)} 
\,{S_{{\frak a},{\frak b}}\left( m , n ; c\right)\over |c|^2}\,
f\left( {2\pi\sqrt{mn}\over c}\right) = {}\qquad 
\qquad\qquad\qquad\qquad\qquad\qquad\qquad &(1.2.15)\cr
 &\qquad\qquad           
=\pi\!\sum_{\scriptstyle V\atop\scriptstyle \nu_V>0}^{(\Gamma)}   
\!{\bf K}f\left(\nu_V , 0\right)
\overline{\sum_{M/2<|m|^2\leq M} a_m
c_V^{\frak a}\left( m;\nu_V,0\right)}  
\sum_{N/2<|n|^2\leq N} b_n c_V^{\frak b}\left( n;\nu_V,0\right) + {}\cr 
 &\qquad\qquad\quad\     
+O_{A}\left( (\log X)\,Y
\left( 1+O_{\varepsilon}\left( |\mu({\frak a})| M^{\varepsilon+1/2}\right)\right) 
\!\left( 1+O_{\varepsilon}\left( |\mu({\frak b})| N^{\varepsilon+1/2}\right)\right) 
\left\| {\bf a}_M\right\|_2 
\left\| {\bf b}_N\right\|_2
\right) ,}$$
where, in the first sum on the right-hand side, 
one sums over just those factors $V$ of the orthogonal decomposition 
(1.1.3) that  lie in the complementary series (i.e. have spectral parameters 
$\nu_V\in (0,1)$ and $p_V=0$\/); and where all other 
terminology either has the same meaning as in Theorem~1, or else 
is defined by the relations (1.2.10) and (1.2.11) of Theorem~2. 

\bigskip

\goodbreak
\noindent{\bf Proof.}\quad 
For any non-zero Gaussian integers $m,n$, an application
of Theorem~1 yields  
the case $\omega_1=m$, $\omega_2=n$ of the summation formula~(1.2.1). 
Upon multiplying both sides of this summation formula by 
$\overline{a_m}\,b_n$, and then summing over all pairs $m,n\in{\frak O}$ 
such that $M/2<|m|^2\leq M$ and $N/2<|n|^2\leq N$, one arrives at an 
expression for the left-hand side of (1.2.15) in terms of a sum  
involving Fourier coefficients $c_V^{\frak a}(m;\nu_V,p_V)$, 
$c_V^{\frak b}(n;\nu_V,p_V)$, $B_{\frak c}^{\frak a}(m;\nu,p)$ and 
$B_{\frak c}^{\frak b}(n;\nu,p)$, and transforms 
${\bf K}f(\nu_V,p_V)$ and ${\bf K}f(\nu,p)$. 
The result (1.2.15) is deduced from this expression by 
applying the upper bound
$${\bf K}f(\nu,p)\ll (\log X)\,Y\left( A^2/X\right)^{|p|} (|p|!)^{-2} (1+|\nu|)^{-4}\qquad\quad\   
\hbox{($\nu\in i{\Bbb R}$, $p\in{\Bbb Z}$)}\eqno(1.2.16)$$
in combination with bounds for the sums $S_j(H,r)$ which, 
for $r\in{\Bbb N}\cup\{ 0\}$, $H\in\bigl\{ 2^k : k\in{\Bbb N}\bigr\}$ and $j=0,\pm i$,   
are given by: 
$$S_j(H,r)=\cases{\sum\limits_{\scriptstyle V\atop{\scriptstyle \nu_V\in i{\Bbb R},\ 
|p_V|=r\atop\scriptstyle H/2\leq |\nu_V|+1<H}} 
\Biggl| \sum\limits_{{\textstyle{M\over 2}}<|m|^2\leq M} a_m
c_V^{\frak a}\left( m;\nu_V,p_V\right)\Biggr|
\,\Biggl|\sum\limits_{{\textstyle{N\over 2}}<|n|^2\leq N} b_n 
c_V^{\frak b}\left( n;\nu_V,p_V\right)\Biggr|\;, &if $j=0$; \cr 
{}^{\quad}_{\quad} &\hbox{\quad} \cr 
\sum\limits_{\scriptstyle {\frak c}\in{\frak C}\atop\scriptstyle 
\left[\Gamma_{\frak c} : \Gamma_{\frak c}'\right] {\textstyle |}\,2r}^{(\Gamma)} 
\sum\limits_{p=\pm r}
\int\limits_{{\textstyle{H\over 2}}-1}^{H-1} 
\Biggl| \sum\limits_{{\textstyle{M\over 2}}<|m|^2\leq M} a_m
B_{\frak c}^{\frak a}\left( m;jt,p\right)\Biggr| 
\,\Biggl|\sum\limits_{{\textstyle{N\over 2}}<|n|^2\leq N} b_n 
B_{\frak c}^{\frak b}\left( n;jt,p\right)\Biggr|\,{\rm d}t\;, &otherwise.\cr}$$
The bound (1.2.16) is proved in Section~2 (see the remark following Lemma~2.2 there). 
As for the relevant bounds on the above sums $S_j(H,r)$: it 
follows by the Cauchy-Schwarz inequality, Remark~2 (above) and the case 
$P=r+1$, $K=H-1$ of the bounds (1.2.9) of Theorem~2 that
$$S_j(H,r)
\ll\left( r^2 +H^2\right)\left( r+1\right) H 
\left( 1+O_{\varepsilon}\left( M^{1+\varepsilon}|\mu({\frak a})|^2\right)\right)^{1/2} 
\left( 1+O_{\varepsilon}\left( N^{1+\varepsilon}|\mu({\frak b})|^2\right)\right)^{1/2}
\left\| {\bf a}_M\right\|_2 
\left\| {\bf b}_N\right\|_2\;,$$
for  $r\in{\Bbb N}\cup\{ 0\}$, $H=2^1,2^2,2^3,\ldots\ $ and $j =0,\pm i$.
On the other hand, for $(\nu,p)\in i{\Bbb R}\times{\Bbb Z}$ with 
$|\nu|+1\geq H/2\geq 1$ and $|p|=r$, the bound (1.2.16) implies 
${\bf K}f(\nu,p)\ll (\log X)Y (A^2/2)^{r}(r!)^{-1} (r+1)^{-3} H^{-4}$ (given that $X\geq 2$). 
Verification of the $O$-term in (1.2.15) may therefore 
be completed by noting that
$$\sum_{H=2^k\,:\,k\in{\Bbb N}} {r^2+H^2\over H^3 (r+1)^2}\ll 1\qquad   
\hbox{(for $r\geq 0$)}\qquad\quad\hbox{and}\qquad\quad 
\sum_{r=0}^{\infty}{\left( A^2/2\right)^{r}\over r!}=\exp\left( A^2/2\right)\quad\blacksquare$$ 

\bigskip

\goodbreak 
\noindent{\bf Remark~4.}\quad 
Let $\varphi : (0,\infty)\rightarrow{\Bbb C}$ 
be a function which is 
infinitely differentiable on $(0,\infty)$, and 
has compact support 
(i.e. support which is a compact subset of $(0,\infty)$). Suppose moreover that 
the function $f : {\Bbb C}^*\rightarrow{\Bbb C}$ satisfies 
$f(z)=\varphi(|z|)$, for all $z\in{\Bbb C}^*$. 
Then, as an almost immediate corollary  
of Lemma~9.4 (below),  it follows that 
the function $f$ is smooth and has compact support: to verify this,  
one has only to check, firstly, that 
the function $\Omega(u)=\varphi(\sqrt u)$ is infinitely 
differentiable on $(0,\infty)$,  and compactly supported, 
before then applying that lemma with $X=1$, $t=0$ and any $B>1$ such that 
$[B^{-1},B]\supseteq{\rm Supp}(\Omega)$. 

\bigskip 

The bound (1.2.16) for ${\bf K}f(\nu,p)$ does not apply in the 
`complementary series' case (i.e. when $p=0$ and $0<\nu<1$), so 
it is of no help in estimating the factors ${\bf K}f(\nu_V,0)$ 
which occur in the sum on the right-hand side of (1.2.15). 
In Section~2, Lemma~2.3, we show that if $A>1$ and $X\geq 2$, and if 
$f(z)=\varphi(|z|)$ ($z\in{\Bbb C}^*$), where the function 
$\varphi : (0,\infty)\rightarrow{\Bbb C}$ is continuous and supported 
in $\bigl[ A^{-1}X^{-1/2}, A X^{-1/2}\bigr]$, then 
$${\bf K}f(\nu,0)\ll_A\ \int_0^{\infty} |\varphi(r)|\,{{\rm d}r\over r}
\,\min\left\{ \log X\,,\,\nu^{-1}\right\} X^{\nu}\qquad\quad  
\hbox{for $\quad 0<\nu\leq 1/2$}\eqno(1.2.17)$$
(note that, by (1.1.11), we do not require information about 
${\bf K}f(\nu,0)$ for $\nu >1/2$). 
This upper bound is near to being best-possible: for if $A$, $X$, $f$ and $\varphi$ 
remain as just described, if 
$A X^{-1/2}\leq 1/c$, with $c>2 e^{\gamma}$ 
(where $\gamma$ is Euler's constant), 
and if the range of $\varphi$ is a subset of $[0,\infty)$, then,  
by the remark following Lemma~2.3 in Section~2,  one will have  
$${\bf K}f(\nu,0)\gg_{A,c}\ \int_0^{\infty} \varphi(r)\,{{\rm d}r\over r}
\,\min\left\{ \log X\,,\,\nu^{-1}\right\} X^{\nu}\qquad\quad  
\hbox{for $\quad 0<\nu\leq 1/2$.}\eqno(1.2.18)$$ 

Upon combining (1.2.17) with Theorem~2 and the Corollary to Theorems~1 and~2, 
one obtains 
(under the same conditions as those under which (1.2.15) is obtained) the upper bound  
$$\eqalignno{ 
&\sum_{M/2<|m|^2\leq M}\overline{a_m}\sum_{N/2<|n|^2\leq N}
b_n\sum_{c\in {}^{\frak a}{\cal C}^{\frak b}}^{(\Gamma)}  
\,{S_{{\frak a},{\frak b}}\left( m , n ; c\right)\over |c|^2}
\,f\!\left( {2\pi\sqrt{mn}\over c}\right) \ll_{A,\varepsilon} {} &(1.2.19)\cr
 &\qquad\qquad\qquad\quad   
\ll_{A,\varepsilon} Y X^{\Theta(q)} (\log X) 
\!\left( 1+|\mu({\frak a})| M^{\varepsilon+1/2}\right) 
\!\left( 1+|\mu({\frak b})| N^{\varepsilon+1/2}\right) 
\!\left\| {\bf a}_M\right\|_2 
\left\| {\bf b}_N\right\|_2\;,}$$
with the exponent $\Theta(q)$ being defined by: 
$$\Theta(q)=\max\left\{ {\rm Re}\left(\nu_V\right) : 
\hbox{$\,V$ occurs in the case $\Gamma =\Gamma_0(q)$ of (1.1.3)}\right\}\eqno(1.2.20)$$
(so that, by (1.1.11) and the points noted in the same paragraph,   
$\Theta(q)=\sqrt{\max\{ 0 , 1-\lambda_1(\Gamma_0(q))\}}$, where 
$\lambda_1(\Gamma)=\min\{ \lambda_V : \hbox{$\,V$ occurs in (1.1.3) and $p_V=0$}\}$).
Many of our new results depend on the constant 
$$\vartheta =\sup_{0\neq q\in{\frak O}}\Theta(q)\;\eqno(1.2.21)$$
If the generalised Selberg eigenvalue conjecture is correct then $\vartheta=0$. 
In this paper we seek unconditional results, and so make do with the following theorem, 
which is an immediate corollary of the result in (1.1.11).

\bigskip

\goodbreak\proclaim Theorem~3 (the Kim-Shahidi bound). For $0\neq q\in{\frak O}$, one has   
$$0\leq\Theta(q)\leq\vartheta\leq{2\over 9}\;.\eqno(1.2.22)$$

\bigskip 
\bigskip 

\goodbreak\centerline{\bf \S 1.3. New Results on Sums over Exceptional Eigenvalues}

\bigskip 

In this section we state our principal new results   
concerning estimates 
for mean values involving Fourier coefficients of cusp forms 
(relevant results from [22] and [13] having already been 
covered in the previous subsection).
In stating these results we have chosen not to include very much 
in the way of ad hoc comments on the definitions and terminology already  
introduced in Subsections~1.1 and~1.2: such definitions and terminology 
are taken as understood. 

We start with a theorem on a sum over exceptional eigenvalues 
pertaining to a single level, $q$. 

\bigskip

\goodbreak\proclaim Theorem~4. 
Let $\varepsilon >0$, $0\neq q\in{\frak O}={\Bbb Z}[i]$ and 
$N,X\geq 1$. Then, for each cusp ${\frak a}$ of the group 
$\Gamma=\Gamma_0(q)\leq SL(2,{\frak O})$, 
and arbitrary complex coefficients $b_n$ ($n\in{\frak O}-\{ 0\}$), one has 
$$\eqalignno{\!\!\!\!\sum_{\scriptstyle V\atop\scriptstyle\nu_V>0}^{(\Gamma)} 
X^{\nu_V}
\Biggl|\sum_{N/2<|n|^2\leq N}& b_n c_V^{\frak a}\left( n;\nu_V,0\right)\Biggr|^2 \ll {} 
 &(1.3.1)\cr 
 &\ll\left( 1+X M_{\frak a}N\right)^{\Theta(q)} 
\left( 1+O_{\varepsilon}\left( M_{\frak a}N^{1+\varepsilon}\right)\right)^{1-\Theta(q)}  
\left\|{\bf b}_N\right\|_2^2\log\left( 2+{1\over M_{\frak a} N}\right) ,}$$
where $M_{\frak a}=|\mu({\frak a})|^2$, while $\left\|{\bf b}_N\right\|_2$, $\mu({\frak a})$ 
and $\Theta(q)$  
are as indicated in (1.2.10), (1.2.11) and (1.2.20).

\bigskip

\goodbreak 
This theorem is analogous to a result [5, Theorem~5] of Deshouillers and Iwaniec; 
and is proved using the same basic idea (a choice of function $f$, in applying 
the Corollary to Theorems~1 and~2, by which the sum of Kloosterman sums in (1.2.15) 
is effectively made void). Full details of the proof appear in Section~3.

In Theorems~5, 8 and~9 below 
we obtain some improvement, on 
average over the level $q$, of the bound (1.3.1). These theorems 
are each concerned with 
estimating, for given complex coefficients $a_n$ ($n\in{\frak O}$), 
$t\in{\Bbb R}$, $Q,N>0$ and $X\geq 1$, the  
sum 
$$S_t(Q,X,N)=
\sum_{Q/2<|q|^2\leq Q}
\sum_{\scriptstyle V\atop\scriptstyle\nu_V>0}^{(\Gamma_0(q))} 
X^{\nu_V}
\left|\sum_{N/4<|n|^2\leq N}a_n |n|^{2it} c_V^{\infty}
\left( n;\nu_V,0\right)\right|^2\;,
\eqno(1.3.2)$$
where we now assume the fixed choice of scaling matrix  
$$g_{\infty}=\pmatrix{1&0\cr 0&1} .\eqno(1.3.3)$$
Note that, when  
$0\neq q\in{\frak O}$ and $\Gamma=\Gamma_0(q)$, 
the above choice of $g_{\infty}$ guarantees that (1.1.1) will hold for ${\frak c}=\infty$;  
and, with regard to the Kloosterman sum  
$S_{\infty,\infty}\left(\omega ,\omega';c\right)$ given by (1.1.13)-(1.1.15), 
the choice (1.3.3) also ensures that 
$${}^{\infty}{\cal C}^{\infty}=q{\frak O}-\{ 0\}\eqno(1.3.4)$$
and 
$$S_{\infty,\infty}\left(\omega ,\omega';\ell q\right) 
=S(\omega ,\omega';\ell q)\qquad\qquad\hbox{($0\neq\ell\in{\frak O}$),}\eqno(1.3.5)$$
where, for $u,v\in{\frak O}$ and $0\neq w\in{\frak O}$, we 
define the `simple Kloosterman sum' $S(u,v;w)$ by: 
$$S(u,v;w)=\sum_{\scriptstyle d\bmod w{\frak O}\atop\scriptstyle (d,w)\sim 1} 
{\rm e}\left( {\rm Re}\left( {u d^{*}+vd\over w}\right)\right) ,\eqno(1.3.6)$$
with $d^{*}$ denoting an arbitrary Gaussian integer solution of the 
congruence $d d^{*}\equiv 1\bmod w{\frak O}$ 
(and with ${\rm e}(x)=\exp(2\pi ix)$, as in the equation (1.1.15)). 
 
In Theorems~5, 6 and~7 
the relevant complex coefficients $a_n$ ($n\in{\frak O}$) may be arbitrary, but  
Theorems~8 and~9 require that these coefficients satisfy   
additional special hypotheses.

\bigskip

\goodbreak\proclaim Theorem~5. Let $\varepsilon >0$. Then, for $X\geq 1$, $Q,N>0$ and 
$t\in{\Bbb R}$, one has 
$$S_t(Q,X,N)\ll_{\varepsilon} 
(QN)^{\varepsilon}\left( Q+Q^{1-2\vartheta}N^{\vartheta}X^{\vartheta}
+N X^{\vartheta}\right)
\left\|{\bf a}_N\right\|_2^2\;,\eqno(1.3.7)$$
where $\vartheta$ is the absolute constant defined by (1.2.20) and (1.2.21).

\bigskip

\goodbreak Theorem~5 is analogous to [5, Theorem~6]. For its proof 
(at the end of Section~4) we need the 
next two results. 

\bigskip

\goodbreak\proclaim Theorem~6 (Change of Weight). Let $\vartheta$ be  
given by (1.2.10) and (1.2.22). Then, for $X,Y\geq 1$, $Q,N>0$ and $t\in{\Bbb R}$, 
one has 
$$S_t(Q,X,N)
\leq\max\left\{ 1\,,\,(X/Y)^{\vartheta}\right\} S_t(Q,Y,N)\;.\eqno(1.3.8)$$

\bigskip

\goodbreak 
\proclaim Theorem~7 (Swapping of Levels). Let $X\geq 1$, $Q,N>0$, $t\in{\Bbb R}$, 
$\varepsilon >0$ and $j\geq 2$; and let 
$$Q^{*}=64\pi^2 XN/Q\;.\eqno(1.3.9)$$ 
Then there exists an 
$L\in\left\{ Q^{*}\,,\,Q^{*}/2\,,\,Q^{*}/2^2\,,\,\ldots\,,\,Q^{*}/2^9\right\}$ such that 
$$S_t(Q,X,N)
\leq O_{\varepsilon, j}\!\left(\,\int\limits_{-\infty}^{\infty} 
S_{t+u}(L,X,N)\,{{\rm d}u\over(1+|u|)^j}\right) 
+O_{\varepsilon}\left( 
X^{\varepsilon}\left( Q+Q^{*}
+N^{1+\varepsilon}\right)\left\|{\bf a}_N\right\|_2^2\right) .\eqno(1.3.10)$$

\bigskip

Given (1.2.20), (1.2.21) and the definition of $S_t(Q,X,N)$ in (1.3.2), 
Theorem~6 is a trivial corollary of Theorem~3 (the Kim-Shahidi bound).   
For proof of this it suffices to note that  
when $0<\nu\leq\vartheta$ one has:    
$$X^{\nu}
=Y^{\nu} (X/Y)^{\nu}
\leq\cases{Y^{\nu}(X/Y)^{\vartheta} &if $0<Y<X$;\cr Y^{\nu} &if $0<X\leq Y$.}$$
Theorem~6, and the use subsequently made of it, 
are modelled on steps in some of the proofs in [5]. 
Theorem~7 is modelled on [5, Lemma~8.1] and (like that lemma) may be proved 
by exploiting the symmetrical nature of a relevant sum of Kloosterman sums. 
We give the proof of Theorem~7 in Section~4. 

Theorems~6 and~7 are useful for more than just the proof of Theorem~5:  
they also help us to deduce, from 
certain elementary estimates for sums of Kloosterman sums, 
the following two theorems (in which $\vartheta\in [0,2/9]$ is 
given by (1.2.20) and (1.2.21)). 

\bigskip 

\goodbreak\proclaim Theorem~8. Let $H\geq 1\geq\delta>0$. Suppose that 
$a_n=\alpha(n)$ for $0\neq n\in{\frak O}$, where the function 
$\alpha : {\Bbb C}\rightarrow{\Bbb C}$ is smooth, 
has its support contained within the annulus 
$\bigl\{ z\in{\Bbb C} : H/2\leq |z|^2\leq H\bigr\}$, and  
satisfies 
$$(\delta |x+iy|)^{j+k}
\,{\partial^{j+k}\over\partial x^j\partial y^k}\,\alpha(x+iy)=O_{j,k}(1)\;, 
\qquad\quad   
\hbox{for $\quad j,k\in\{ 0,1,2,\ldots \}\ $ and $\ x,y\in{\Bbb R}$.}\eqno(1.3.11)$$
Then, for all $Q,X\geq 1$, all $t\in{\Bbb R}$ and all $\varepsilon >0$,  
one has: 
$$S_t(Q,X,H) 
\ll_{\varepsilon} 
\left(\delta^{-1}+|t|\right)^{11} 
\left( 1+{X\over Z(Q,H)}\right)^{\!\vartheta} \!(Q+H)^{1+\varepsilon} H\;,\eqno(1.3.12)$$
where
$$Z(Q,H)={Q^2\over H}+H\;.\eqno(1.3.13)$$

\bigskip

\goodbreak\proclaim Theorem~9. 
Let $H,K\geq 1\geq\delta>0$. Suppose that $N=HK$ and that, 
for $0\neq n\in{\frak O}$, one has 
$$a_n=\sum_{\scriptstyle h\mid n}
\alpha(h)\beta\left( {n\over h}\right)\;,\eqno(1.3.14)$$
where the functions $\alpha,\beta : {\Bbb C}\rightarrow{\Bbb C}$ are smooth, 
have both 
${\rm Supp}(\alpha)\subseteq\bigl\{ z\in{\Bbb C} : H/2\leq |z|^2\leq H\bigr\}$ and 
${\rm Supp}(\beta)\subseteq\bigl\{ z\in{\Bbb C} : K/2\leq |z|^2\leq K\bigr\}$ 
(where $\ {\rm Supp}(f)$ denotes the support of $f$), and, 
at all points $(x,y)\in{\Bbb R}^2$, satisfy   
$$(\delta |x+iy|)^{j+k}
\max\left\{\left| {\partial^{j+k}\over\partial x^j\partial y^k}\,\alpha(x+iy)\right|\,, 
\,\left| {\partial^{j+k}\over\partial x^j\partial y^k}\,\beta(x+iy)\right|\right\} 
=O_{j,k}(1)\qquad\   
\hbox{($j,k\in{\Bbb N}\cup\{ 0\}$).}\eqno(1.3.15)$$
Then, for all $Q,X\geq 1$, all $t\in{\Bbb R}$ and all $\varepsilon>0$, one has: 
$$S_t(Q,X,N) 
\ll_{\varepsilon} 
\left(\delta^{-1}+|t|\right)^{11} 
\left( \left( 1+{X\over Q^2 N^{-1}}\right)^{\!\vartheta}\!Q 
+\left( 1+{X\over H+K}\right)^{\!\vartheta}\!N\right) Q^{\varepsilon} N^{1+\varepsilon}\;.
\eqno(1.3.16)$$

\bigskip

\goodbreak 
Theorem~8 is an analogue of [21, Theorem~3]. 
We think it worth noting that [21, Theorem~3] is a corollary of 
a deeper result [5, Theorem~7] obtained by Deshouillers and Iwaniec. 
By analogy with [5, Theorem~7],  
one might expect the results (1.3.12)-(1.3.13)  also to 
be valid in cases where, for some $N\geq 1$, one has: 
$$a_n=\cases{1 &if $N/2<|n|^2\leq N$;\cr 0 &otherwise.}\eqno(1.3.17)$$
Although good bounds for $S_t(Q,X,N)$ in these cases would be interesting, 
we do not require (or prove, or claim) any such bounds in this paper: 
what work we have done in this area falls well short of 
giving (1.3.12)-(1.3.13), subject to (1.3.17), and is in too much of an 
unfinished state to be worth recording here. 

Though it might have been instructive to have included  
an independent proof of Theorem~8, we prefer 
just to point out that Theorem~8 is trivially implied by Theorem~9. 
To see that this is indeed so, observe firstly that if $K=\sqrt{2}$ (say),  
and if $0<\delta\leq 1$, then there exists a function $\beta : {\Bbb C}\rightarrow{\Bbb C}$ 
which, while satisfying the hypotheses of Theorem~9, is also 
such that ${\frak O}\cap{\rm Supp}(\beta)=\{ 1\}$. 
For such $K$ and $\beta$ the definition (1.3.14) simplifies 
to give just $a_n=\alpha(n)$, for $n\in{\frak O}-\{ 0\}$. 
This shows that the hypotheses concerning $S_t(Q,X,N)$ in 
Theorem~8 justify the application of Theorem~9 
(i.e. with $K=\sqrt{2}$, $N=\sqrt{2} H$ and $\beta$ as just described).  
Hence (after simplifying the relevant case of (1.3.16)) one obtains:  
$$S_t(Q,X,H)=S_t(Q,X,\sqrt{2}H)
\ll_{\varepsilon}\left(\delta^{-1}+|t|\right)^{11}  
\left( Q+H+\sqrt{XH}\right)^{2\vartheta} (Q+H)^{1+2\varepsilon -2\vartheta} H\;.$$
Upon substituting $\varepsilon /2$ for $\varepsilon$ here (as one may), 
the result (1.3.12)-(1.3.13) of Theorem~8 follows immediately, since  
$(1+\sqrt{XH}/(Q+H))^2\ll 1+XH/(Q+H)^2<1+XH/(Q^2+H^2)$.  

We prove Theorem~9 at the end of Section~8, following 
extensive preparation undertaken in the first part of that section, and (before that) in  
Sections~5, 6 and~7. 
This (given the observations of the preceding paragraph) 
makes it unnecessary to include a separate proof of Theorem~8. 
It is nevertheless worth mentioning that we could prove Theorem~8 
independently of Theorem~9, by taking as a starting point the case 
$\beta =\overline{\alpha}$ 
of the elementary upper bound
$$\sum_{H/2<|h|^2\leq H}\,\sum_{K/2<|k|^2\leq K} 
\alpha(h)\beta(k)|hk|^{iu}S(h,k;c)
\ll_{\varepsilon} 
|c|^{\varepsilon} HK +|c|^{\varepsilon +2}(\log H)(\log K)\left(\delta^{-1}+|u|\right)^8\eqno(1.3.18)$$
(valid when $u\in{\Bbb R}$, $\varepsilon >0$ and $\alpha$, $\beta$, $H$, $K$ and $\delta$ are 
as in Theorem~9). Our proof of (1.3.18) is 
omitted from this paper, since there is nothing very novel about it,
and since the other results in this paper are obtained 
independently of (1.3.18).
The corresponding starting point for the proof of Theorem~9 is the 
estimate for a sum of Kloosterman sums obtained in Lemma~6.3. 
Our work in Section~7 enables us to compensate for 
the inconveniently restrictive conditions under which 
the result of Lemma~6.3 is obtained.  

Although Theorem~9 is analogous to our result in [21, Theorem~2], 
the proof we give of Theorem~9 is not obtained 
by adapting, in its entirety, the corresponding proof in [21];   
the relevant innovations are 
discussed in our `Outline of Results and Methods', above (see 
the paragraph containing the bound (0.8), and the two paragraphs preceding it). 

\bigskip 
\bigskip 

\goodbreak\centerline{\bf \S 1.4. A Sum of Kloosterman Sums and an Application}

\bigskip 

As just mentioned, estimates for sums of Kloosterman sums play 
a part in proving Theorems~8 and~9.  Conversely, the Corollary to Theorems~1 and~2
makes it possible to deduce, from our new results on 
sums over exceptional eigenvalues, some results on 
sums of generalised Kloosterman sums that are genuinely new 
(in that they do not follow directly from (1.2.19) and Theorem 3). 
In this paper we obtain  
just one such result, which is Theorem~10 below;   
for its proof we require also an auxilliary result, Theorem~11. 
The complete proofs of Theorems~10 and~11
appear in Section~9 (here those proofs are only outlined briefly). 

We work, as before, with Hecke congruence subgroups 
$\Gamma=\Gamma_0(q)<SL(2,{\Bbb C})$, where $0\neq q\in{\frak O}={\Bbb Z}[i]$. 
Given $q$ (and hence the group $\Gamma$), 
the associated generalised Kloosterman sums 
that Theorem~10 relates to are those of the form 
$S_{1/s,\infty}(\omega,\omega';c)$, 
where $s$ divides $q$ and is coprime to the Gaussian integer $r=q/s$. 
To completely determine the values of these sums 
one must specify scaling matrices for the cusps ${1/s}$ and $\infty$. 
We choose the scaling matrix $g_{\infty}$ as in (1.3.3);   
and for each pair $r,s$ of non-zero coprime Gaussian integers such that $rs=q$,  
we choose the scaling matrix for the cusp $1/s$ of $\Gamma_0(q)$ 
to be 
$$g_{1/s}=\pmatrix{\sqrt{r}&t/\sqrt{r}\cr s\sqrt{r}&u\sqrt{r}} ,\eqno(1.4.1)$$
where the square root is chosen arbitrarily, while  
$t=t(r,s)$ and $u=u(r,s)$ may be any pair of Gaussian integers with  
$$ru-st=1\eqno(1.4.2)$$
(so that $g_{1/s}\in SL(2,{\Bbb C})$\/). 
A suitable choice of  
$u$ and $t$ may be determined by means of 
the Euclidean algorithm for ${\Bbb Z}[i]$. 
Then (see Lemmas~9.1 and~9.2) 
the condition~(1.1.1) holds both for ${\frak c}=\infty$ and 
for ${\frak c}=1/s$; while  the 
set ${}^{1/s}{\cal C}^{\infty}$ and generalised Kloosterman sum 
$S_{{1/s},\infty}\bigl(\omega,\omega';c\bigr)$ 
defined in (1.1.13)-(1.1.15) satisfy
$${}^{1/s}{\cal C}^{\infty}
=\left\{ ps\sqrt{r} : 0\neq p\in{\frak O}\ {\rm and}\ (p,r)\sim 1\right\}\eqno(1.4.3)$$
and 
$$S_{{1/s},\infty}\bigl(\omega,\omega';ps\sqrt{r}\bigr)
=S\left(r^{*}\omega,\omega';ps\right)\qquad\qquad    
\hbox{($\omega,\omega'\in{\frak O}\ $ and $\ 0\neq p\in{\frak O}\,$ with $\,(p,r)\sim 1$),}\eqno(1.4.4)$$
with $r^{*}=r^{*}(r,ps)$ being an 
arbitrary Gaussian integer such that $r r^{*}\equiv 1\bmod ps{\frak O}$, 
and with $S(u,v;w)$ being the Kloosterman sum defined in (1.3.6). 

We are now ready to state Theorem~10: a
new result concerning 
the Kloosterman sums in (1.4.4).  

\bigskip

\goodbreak\proclaim Theorem~10. Let $\vartheta$ be the real absolute constant 
given by (1.2.20) and (1.2.21); let $\varepsilon >0$; and let 
$N,L\geq 1\geq\delta >0$. Let $a_n\in{\Bbb C}$ for $n\in{\frak O}-\{ 0\}$; 
and let $A : {\Bbb C}\rightarrow{\Bbb C}$ be a smooth function 
which satisfies 
$$(\delta |x+iy|)^{j+k}{\partial^{j+k}\over\partial x^j\partial y^k}\,A(x+iy) 
=O_{j,k}(1)\qquad\qquad 
\hbox{($j,k\in{\Bbb N}\cup\{ 0\}\ $ and $\ x,y\in{\Bbb R}$),}\eqno(1.4.5)$$
and which has 
$\,{\rm Supp}(A)\subseteq\left\{ z\in{\Bbb C} : L/2\leq |z|^2\leq L\right\}$. 
Let $P,Q,R,S\geq 1$ and $X>0$ satisfy 
$$Q=RS\geq\max\left\{ \sqrt{N}\,,\,\sqrt{L}\right\}\eqno(1.4.6)$$
and 
$$X={PS\sqrt{R}\over 4\pi^2\sqrt{L N}}\geq 2\;.\eqno(1.4.7)$$ 
Let also $b$ be a complex-valued function with domain 
$${\cal B}(R,S)
=\left\{ (r,s)\in{\frak O}\times{\frak O} : 
\,R/2<|r|^2\leq R\,,\ S/2<|s|^2\leq S\ {\rm and}\ (r,s)\sim 1\right\} ;\eqno(1.4.8)$$
and, for each pair 
$(r,s)\in{\cal B}(R,S)$, 
let $g_{r,s} : (0,\infty)\rightarrow{\Bbb C}$ 
be an infinitely differentiable function which satisfies 
$$g_{r,s}^{(j)}(x)\ll_j x^{-j}\qquad\qquad\hbox{($\,(r,s)\in{\cal B}(R,S)$,  
$\,j\in{\Bbb N}\cup\{ 0\}\,$ and $\,x>0$),}\eqno(1.4.9)$$
and which has 
${\rm Supp}\bigl( g_{r,s}\bigr)\subseteq [P/2,P]$. Put 
$${\Lambda}=\sum_{(r,s)\in{\cal B}(R,S)} b(r,s) 
\sum_{N/4<|n|^2\leq N} a_n \sum_{L/2<|\ell|^2\leq L} A(\ell)  
K_{r,s}(n,\ell)\;,\eqno(1.4.10)$$
where 
$$K_{r,s}(n,\ell) 
=\sum_{\scriptstyle 0\neq p\in{\frak O}\atop\scriptstyle (p,r)\sim 1} 
g_{r,s}\!\left( |p|^2\right) 
S\!\left( r^{*}n , \ell ; ps\right)\eqno(1.4.11)$$
(with $r^{*}=r^{*}(r,ps)$ and $S(u,v;w)$ as described below (1.4.4)). 
Then 
$${\Lambda}^2\ll_{\varepsilon} 
Q^{\varepsilon}\left\|b\right\|_2^2\,\left\| {\bf a}_N\right\|_2^2 L P^2 S\log^2(X)\,(L+Q)(N+Q) 
\!\left( 1+
{X^2\over \left( 1+Q N^{-1}\right)\!\left( 1+Q L^{-1}\right)^{2}\!L}\right)^{\!\!\vartheta}
\delta^{-11}\;,\eqno(1.4.12)$$
where 
$$\left\| b\right\|_2
=\biggl( \sum_{(r,s)\in{\cal B}(R,S)} |b(r,s)|^2\biggr)^{1/2}\;,\eqno(1.4.13)$$
and where 
the terminology `$\|{\bf x}_M\|_2$' is that introduced in (1.2.11). 
If it is moreover the case that the hypotheses of Theorem~9 
concerning $a_n$ ($n\in{\frak O}-\{ 0\}$), $N,H,K\in[1,\infty)$ and 
$\alpha,\beta : {\Bbb C}\rightarrow{\Bbb C}$ are satisfied, then 
one has also 
$$\eqalignno{{\Lambda}^2
 &\ll_{\varepsilon} 
Q^{1+\varepsilon}\left\| b\right\|_{\infty}^2 
N L P^2 S\log^2(X) \times {}  &(1.4.14)\cr 
 &\qquad\times (L+Q)  
\!\left(\!\Biggl( 1+{X^2\over (H+K)\left( 1+Q L^{-1}\right)^{2} L}\Biggr)^{\!\!\vartheta} N 
+\Biggl( 1+{X^2\over Q^2 N^{-1}\left( 1+Q L^{-1}\right)^{2} L}\Biggr)^{\!\!\vartheta} 
Q\!\right)\!\delta^{-22}\,,\qquad\quad}$$
where 
$$\left\| b\right\|_{\infty}
=\max_{(r,s)\in{\cal B}(R,S)} |b(r,s)|\;.\eqno(1.4.15)$$

\medskip

\goodbreak 
\noindent{\bf Remark.}\quad 
The implicit constants in (1.4.12) and (1.4.14) may of course depend on 
those in the conditions~(1.4.5) and~(1.4.9); and the one in (1.4.14) 
may also depend on the implicit constants in (1.3.15).

\bigskip

\goodbreak 
To prove Theorem~10 we first use (1.4.4) and the Corollary 
to Theorems~1 and~2 to bound the sum $\Lambda$, given by (1.4.10) and (1.4.11),  
in terms of an acceptable $O$-term plus a sum over exceptional eigenvalues 
$\lambda_V=1-\nu_V^2$.  
It is then almost (but not quite) straightforward to deduce the results in (1.4.12)-(1.4.15) 
from Theorems~3, 4, 8 and~9, via the 
Cauchy-Schwarz inequality:  
the one (minor) problem that we encounter, in carrying this out, 
is the necessity of 
dealing with sums involving Fourier coefficients 
at cusps $1/s$, whereas Theorems~8 and~9 apply only to 
sums involving Fourier coefficients at the cusp $\infty\,$  
(which is $\Gamma_0(q)$-equivalent to $1/s$ only if $q\mid s$). 
In addressing the problem just mentioned we rely on an idea of Iwaniec,   
which  applies here through the observation that, if 
$\Gamma$, $q$, $r$, $s$, $g_{1/s}$ and $g_{\infty}$ are as assumed in (1.4.1)-(1.4.2) and (1.3.3), 
then, in addition to (1.4.3) and (1.4.4), one will have
$${}^{1/s}{\cal C}^{1/s}={}^{\infty}{\cal C}^{\infty}=q{\frak O}-\{ 0\}\eqno(1.4.16)$$
and, for $\omega,\omega'\in{\frak O}$ and $0\neq\ell\in{\frak O}$, 
$$S_{1/s,1/s}\left(\omega,\omega';\ell q\right) 
=S\left(\omega,\omega';\ell q\right) 
=S_{\infty,\infty}\left(\omega,\omega';\ell q\right)\eqno(1.4.17)$$
(see Lemmas~9.1 and~9.2 for the proofs). 
By combining these facts with the spectral to Kloosterman summation formula and 
spectral large sieve inequalities of [22, Theorems~B~and~1]  
we obtain the following theorem, which is analogous to a result of Iwaniec in [11], 
and through which we solve the above mentioned 
problem of dealing with sums involving Fourier coefficients at cusps  
that are not $\Gamma$-equivalent to $\infty$. 

\bigskip

\goodbreak 
\proclaim Theorem~11. For all $n\in{\frak O}-\{ 0\}$, let $b_n\in{\Bbb C}$. 
Let $X>0$ and $N\geq 1$; 
let $q,r,s\in{\frak O}-\{ 0\}$ satisfy 
$$q=rs\qquad\quad {\rm and}\qquad\quad (r,s)\sim 1\;.\eqno(1.4.18)$$
Let $\Gamma =\Gamma_0(q)\leq SL(2,{\frak O})$; 
let $g_{\infty},g_{1/s}\in SL(2,{\Bbb C})$ be as stated in 
(1.3.3) and (1.4.1)-(1.4.2); and, 
for ${\frak a}\in\{ \infty , 1/s\}$, let   
$$\rho^{\frak a}=\rho_q^{\frak a}({\bf b},N;X)
=\sum_{\scriptstyle V\atop\scriptstyle\nu_V>0}^{(\Gamma)} 
\!\!\left( X^{\nu_V}+X^{-\nu_V}\right)\exp\left(\nu_V^2\right) 
\!\Biggl|\sum_{{\textstyle{N\over 4}}<|n|^2\leq N} \!\!b_n 
c_V^{\frak a}\!\left( n;\nu_V,0\right)\Biggr|^{\,2}\;.\eqno(1.4.19)$$
Then, when $\varepsilon >0$, one has 
$$\rho^{1/s}-\rho^{\infty}
\ll\left( 1+O_{\varepsilon}\!\left( |q|^{-2}N^{1+\varepsilon}\right)\right)
\left\|{\bf b}_N\right\|_2^2\;,\eqno(1.4.20)$$ 
where $\left\|{\bf b}_N\right\|_2$ is as defined in (1.2.11). 

\bigskip 

\goodbreak 
In a forthcoming paper [23] we show that Theorem~10 has a significant application 
in respect of mean values involving certain Hecke zeta-functions. 
This echoes the way in which the results [5, Theorems~10~and~11], 
which bound sums involving the 
classical analogue of the Kloosterman sum $S(r^*\omega,\omega';ks)$ in 
(1.4.4), were used by Deshouillers and Iwaniec to obtain, in [6], 
new upper bounds for the mean value 
$$I(T,M)={1\over T}\int_0^T 
\Biggl|\sum_{0<m\leq M} a_m m^{-it}\Biggr|^2 
\left|\zeta\left(\textstyle{1\over 2}+it\right)\right|^4 {\rm d}t\;,$$
where $\zeta(s)$ is Riemann's zeta-function, and $(a_m)$ an arbitrary 
complex sequence. 

Ideas proposed  by Iwaniec in [11] inspired the 
work [21], where (by building upon the approach of Deshouillers and Iwaniec) 
it was shown that when $\varepsilon >0$ one has 
$$I(T,M)\ll_{\varepsilon} T^{\varepsilon} M \max_{m\leq M} |a_m|^2\qquad\quad  
\hbox{for $\quad M^4\ll T$.}\eqno(1.4.21)$$
Our Theorems~8, 9 and~10 are analogues of the results [21, Theorem~3, Theorem~2 and 
Proposition~4.1], respectively. In [23] we prove, with the help of Theorem~10, 
a new upper bound for the mean value 
$$J(D,N)
={1\over D^2}\sum_{-D\leq k\leq D}\int\limits_{-D}^D 
\Biggl|\sum_{0<|n|^2\leq N} \alpha_n \lambda^k(n) |n|^{-2it}\Biggr|^2 
\left|\zeta\left(\textstyle{1\over 2}+it,\lambda^k\right)\right|^4 {\rm d}t\;,\eqno(1.4.22)$$
where the coefficients $\alpha_n\,$ ($0\neq n\in{\Bbb Z}[i]$) are arbitrary complex numbers 
and $\lambda^k$ denotes the `gr\"{o}ssencharakter' 
given by $\lambda^k(n)=(n/|n|)^{4k}\ $ ($0\neq n\in{\frak O}$), 
while $\zeta(s , \lambda^k)$ 
denotes that Hecke zeta function which satisfies   
$$\zeta\bigl( s , \lambda^k\bigr)
={1\over 4}\sum_{n\neq 0}
\lambda^k(n) |n|^{-2s}\qquad\qquad   
\hbox{(${\rm Re}(s)>1$),}$$
and is thereby uniquely defined, 
through analytic continuation, 
for all $s\in{\Bbb C}-\{ 1\}\,$ 
(the point $s=1$ is excluded here solely due to it being 
a simple pole of the Dedekind zeta function 
$\,\zeta_{{\Bbb Q}(i)}(s)=\zeta(s,\lambda^0)$\/).  
The principal result of [23] contains the 
natural analogue of the result in (1.4.21), which is the bound 
$$J(D,N)\ll_{\varepsilon} D^{\varepsilon} N \max_{0<|n|^2\leq N} |\alpha_n|^2\qquad\quad  
\hbox{for $\quad N^2\leq D$.}\eqno(1.4.23)$$

\bigskip 
\bigskip 

\goodbreak\centerline{\bf \S 1.5. Notation and Conventions}

\bigskip 

Notation in common use requires no comment, so 
our Index (below) cover only the more unusual 
(or potentially ambiguous) notation and  
conventions; though we have not made it comprehensive, even in this respect.  
Some of the notation not listed in the index is discussed in   
supplementary paragraphs. 

\medskip 

\goodbreak 
\noindent{\bf Index of Notation:} 

\medskip 

{\settabs\+$E_j^{\frak a}(q,P,K;N,{\bf b})\,$\  
&a complete set of representatives of the $\Gamma$-equivalence classes of cusps\quad\  
&\cr 
\+{\bf Symbol}&{\bf Description}&{\bf Where Defined}\cr 

\smallskip 

\+${\bf v}\cdot{\bf w}$
&{\it equal to $v_1 w_1+\cdots+v_n w_n$, 
the inner product of vectors ${\bf v},{\bf w}\in{\Bbb C}^n$ 
}&--\cr
\+${\frak a}\sim^{\!\!\!\!\!\Gamma}{\frak b}$
&{\it the relation of $\Gamma$-equivalence (for cusps ${\frak a}, {\frak b}$)
}&above (1.1.12)\cr 
\+$m\mid n$
&{\it (when $m,n\in{\frak O}$): the relation `$n$ is divisible by $m$'  
}&--\cr
\+$m\sim n$
&{\it (when $m,n\in{\frak O}$): the relation `$n$ is an associate of $m$' 
}&--\cr 
\+$(m,n)$
&{\it a highest common factor (of $m,n\in{\frak O}$)
}&in Conventions\cr 
\+$[x]$
&{\it the greatest rational integer less than or equal to $x$
}&--\cr 
\+$\|\beta\|$
&{\it the distance from $\beta\in{\Bbb C}$ to the nearest 
Gaussian integer
}&above (5.7)\cr 
\+$(\alpha)_m$
&{\it equal to $\alpha(\alpha +1)\cdots(\alpha+m-1)$
}&in (2.1)\cr 
\+$\|{\bf b}_N\|_2$
&{\it the Euclidean norm of a vector involving coefficients $b_n$ ($0\neq n\in{\frak O}$) 
}&in (1.2.11)\cr 
\+$\| b\|_2$, $\| b\|_{\infty}$
&{\it the `Euclidean' and `Sup' norms of the function $b : {\cal B}(R,S)\rightarrow{\Bbb C}$  
}&(1.4.13), (1.4.15)\cr 

\smallskip 

\+$\int_{\Gamma\backslash G} f(g) {\rm d}g$
&{\it a right-invariant integral of a $\Gamma$-automorphic measurable function $f$
}&beginning of \S 1.1\cr 

\smallskip 

\+$\int\limits_{(\sigma)}f(z) {\rm d}z$
&{\it integral along a `vertical' contour line in ${\Bbb C}$
}&beginning of \S 1.2\cr 
\smallskip 

\+$\hat{F}({\bf y})$, $\hat{f}({\bf w})$
&{\it Fourier transforms 
for $F\in{\cal S}\bigl( {\Bbb R}^n\bigr)$, 
$f\in{\cal S}\bigl( {\Bbb C}^n\bigr)$
}&in (5.1)-(5.3)\cr 

\smallskip 

\+$\bigoplus\limits_{{\frak c}\in{\frak C}}^{(\Gamma)}$
&{\it direct sum over representatives of the $\Gamma$-equivalence classes of cusps
}&as in (1.2.1)\cr

\smallskip 

\+$\Gamma$, $\Gamma_0(q)$
&{\it Hecke congruence subgroups of $SL(2,{\frak O})$  
($\Gamma_0(q)$ is of `level' $q$) 
}&beginning of \S 1.1\cr
\+$\Gamma_{\frak c}$, $\Gamma_{\frak c}'$
&{\it `stabiliser' and 
`parabolic stabiliser' subgroups (for the cusp ${\frak c}$) 
}&above (1.1.1)\cr 
\+${}^{\frak a}\Gamma^{\frak b}(c)$
&{\it a `Bruhat cell'
}&in (1.1.13)\cr  
\+$\Gamma(z)$
&{\it Euler's Gamma function, defined for $z\in {\Bbb C}-\{0,-1,-2,\ldots\ \}$ 
}&--\cr 
\+$\gamma$
&{\it usually an element of $\Gamma$; sometimes Euler's constant, $0.5772157\ldots\ $  
}&--\cr
\+${\bf \Delta}$ 
&{\it the hyperbolic Laplacian operator on $L^2(\Gamma\backslash G/K)$
}&above (1.1.11)\cr 
\+$\partial/\partial z_m$, 
$\partial/\partial\overline{z_m}$
&{\it Complex partial differentiation operators
}&in (5.19)\cr 
\+$\delta^{*}_{{\frak a},{\frak b}}$ 
&{\it the `delta symbol' for $\Gamma$-equivalence of the cusps ${\frak a}$ and ${\frak b}$ 
}&in (1.1.17)\cr 
\+$\delta_{\omega,\omega'}^{{\frak a},{\frak b}}$
&{\it the `delta symbol' of the `spectral to Kloosterman' sum formula
}&in (9.11)\cr 
\+$\delta_{w,z}$
&{\it the `delta-symbol' for equality of the complex numbers $w$ and $z$
}&below (9.11)\cr 
\+$\Theta(q)$
&{\it equal to $\sqrt{\max\{ 0 , 1-\lambda_1(\Gamma_0(q))\}}$ 
(conjecturally zero for $q\in{\frak O}-\{ 0\}$) 
}&in (1.2.10)\cr 
\+$\vartheta$
&{\it the least upper bound for the set $\{\Theta(q) : 0\neq q\in{\frak O}\}$ 
}&in (1.2.21)\cr 
\+$\lambda_V$
&{\it equal to $1-\nu_V^2$ (an eigenvalue of the operator $-{\bf \Delta}$ on $L^2(\Gamma\backslash G/K)\,$) 
}&above (1.1.11)\cr
\+$\lambda_1$ (or $\lambda_1(\Gamma)$)
&{\it the first eigenvalue of the operator $-{\bf \Delta}$ on $L^2(\Gamma\backslash G/K)$
}&below (1.2.20)\cr 
\+$\mu({\frak a})$
&{\it $1/|\mu({\frak a})|$ is a useful lower bound for 
the set $\{ |c| : c\in{}^{\frak a}{\cal C}^{\frak a}\}$}&in (1.2.10)\cr 
\+$\mu_{\frak O}(n)$
&{\it the `Gaussian' M\"{o}bius function
}&in (5.4)\cr  
\+$(\nu_V,p_V)$
&{\it the spectral parameters of the cuspidal space $V$
}&in (1.1.4)\cr 

\smallskip 

\+$\sum
\limits_{c\in {}^{\frak a}{\cal C}^{\frak b}}^{(\Gamma)}
f(c)S(c)$
&{\it (with $S(c)=S_{{\frak a},{\frak b}}(m ,n ;c)$): here ${}^{\frak a}{\cal C}^{\frak b}$ 
and $S(c)$ are dependent on $\Gamma$
}&as in (1.2.1)\cr 

\smallskip 

\+$\sum\limits_V^{(\Gamma)}$
&{\it sum over irreducible cuspidal subspaces $V\subset L^2(\Gamma\backslash G)$
}&as in (1.2.1)\cr 

\smallskip 

\+$\sum\limits_{{\frak c}\in{\frak C}}^{(\Gamma)}$
&{\it sum over representatives of the $\Gamma$-equivalence classes of cusps
}&as in (1.2.1)\cr

\smallskip 

\+$\sigma_q^{\frak a}({\bf b},N;X)$
&{\it a weighted mean value for the group $\Gamma_0(q)$
}&(0.6), \S 3, (4.16)\cr 
\+$\Omega_{\pm}$
&{\it the Casimir operators associated with $G$
}&above (1.1.4)\cr 
\+${\bf \Omega}, {\bf A}, {\bf a}, {\bf B},\ \ldots\ $
&{\it vectors in ${\Bbb R}^n$ or ${\Bbb C}^n$; sets of coefficients 
(see `$b_n$' and `${\bf b}$' in (1.2.7))
}&--\cr
\+${\frak a},{\frak b},{\frak c},\ldots\ $
&{\it cusps of $\,\Gamma$, or (more generally) points 
in ${\Bbb P}^1({\Bbb C})={\Bbb C}\cup\{\infty\}$
}&above (1.1.1)\cr 
\+$A$, $a[r]$
&{\it $A=\{ a[r] : r>0\}<G$
}&beginning of \S 1.1\cr
\+${\cal A}^{\Bbb R}_m I$ 
&{\it a `generalised annulus' in ${\Bbb R}^{2m}$,  
determined by the set $I\subset [0,\infty)$  
}&in (8.1)\cr  
\+$B^{+}$
&{\it the group $\{ n[\alpha] : \alpha\in{\frak O}\}<N<G$
}& in (1.1.1)\cr 
\+$B_{\frak a}^{\frak b}\left(\omega;\nu,p\right)$
&{\it a modified Fourier coefficient of an Eisenstein series
}&in (1.1.21)\cr 
\+${\cal B}(R,S)$
&{\it a bounded subset of ${\frak O}\times{\frak O}-\{ (0,0)\}$
}&in (1.4.8)\cr 
\+${\frak C}$&{\it a complete set of representatives of 
the $\Gamma$-equivalence classes of cusps
}&above (1.1.12)\cr 
\+${}^{\frak a}{\cal C}^{\frak b}$
&{\it the set of arguments of a generalised Kloosterman sum
}&in (1.1.14)\cr
\+$c_V^{\frak c}\left(\omega\right)$ 
&{\it a Fourier coefficient of a cuspidal subspace
}&in (1.1.8)-(1.1.9)\cr 
\+$c_V^{\frak c}\left(\omega;\nu_V,p_V\right)$
&{\it a modified Fourier coefficient of a cuspidal subspace
}&in (1.1.21)\cr 
\+$c_q(b,h;k)$
&{\it a generalisation of the Ramanujan sum
}&in (5.31)\cr 
\+$D_{\frak a}^{\frak b}\left(\omega;\nu,p\right)$
&{\it a Fourier coefficient of an Eisenstein series
}&in (1.1.18)\cr
\+${\rm d}_{+}z$, ${\rm d}_{\times}z$
&{\it the standard Lebesgue measure on ${\Bbb C}$, 
and a Haar measure for ${\Bbb C}^*$
}&below (1.2.2)\cr 
\+${\rm d}g$
&{\it a normalised left and right Haar measure on $G$
}&beginning of \S 1.1\cr 
\+$E_j^{\frak a}(q,P,K;N,{\bf b})$ 
&{\it a spectral mean, for cusp forms ($j=0$), or 
Eisenstein series ($j=1$)
}&in (1.2.7)-(1.2.8)\cr 
\+${\rm e}(x)$
&{\it equal to $\exp(2\pi i x)$, a character for the additive group ${\Bbb R}/{\Bbb Z}$
}&below (1.1.15)\cr 
\+$(F_m^{\frak c} f)(g)$
&{\it the `$m$-th order' term in the Fourier expansion of 
$f$ at a cusp ${\frak c}$
}&below (1.1.1)\cr  
\+$g_{\frak c}$ 
&{\it a scaling matrix for the cusp ${\frak c}$
}&in Conventions\cr
\+$G$
&{\it the special linear group, $SL(2,{\Bbb C})$
}&beginning of \S 1.1\cr 
\+$h[u]$
&{\it $h[u]\in G$ for $u\in{\Bbb C}^*$; and $h[u]=k[u,0]$ when $|u|=1$
}&beginning of \S 1.1\cr 
\+$J_{\nu}^{*}(z)$
&{\it equal to $(z/2)^{-\nu}J_{\nu}(z)$ when $z>0$   
($J_{\nu}(z)$ being Bessel's $J$-function) 
}&in (1.2.4)\cr 
\+${\cal J}_{\mu,k}(z)$, ${\cal K}_{\nu,p}(z)$
&{\it functions related to Bessel functions of representations of 
$PSL(2,{\Bbb C})$ 
}&in (1.2.3)\cr 
\+$K$, $k[\alpha ,\beta]$
&{\it the special unitary group, $SU(2)<G$, and one of its elements  
}&beginning of \S 1.1\cr 
\+${\bf K}f(\nu,p)$
&{\it the ${\bf K}$-transform of $f$
}&in (1.2.2)-(1.2.4)\cr 
\+${\bf K_m}\varphi(\nu ,k)$
&{\it a component of the ${\bf K}$-transform
}&in (2.6)\cr 
\+$\log(x)$
&{\it equal to $\log_{e}(x)$, the natural logarithm
}&--\cr 
\+${\cal L}_m$
&{\it second order differential operators on ${\cal S}\bigl( {\Bbb C}^n\bigr)$
}&in (5.12)\cr 
\+$M_{\frak a}$
&{\it equal to $|\mu({\frak a})|^2$ 
}&below (1.3.1)\cr 
\+${\bf M}\varphi(s)$
&{\it the ${\bf M}$-transform of $\varphi$ (a variant of the Mellin transform)
}&in (2.2)\cr 
\+$N$, $n[z]$ 
&{\it $N=\{ n[z] : z\in{\Bbb C}\}<G$ 
}&beginning of \S 1.1\cr 
\+${\frak O}$
&{\it equal to ${\Bbb Z}[i]$, the ring of integers of the Gaussian number field 
${\Bbb Q}(i)$
}&beginning of \S 1.1\cr
\+${\Bbb P}^1({\Bbb C})$
&{\it a projective line, identified with the Riemann sphere, ${\Bbb C}\cup\{\infty\}$
}&above (1.1.1)\cr
\+${\Bbb P}^1\bigl( {\Bbb Q}(i)\bigr)$
&{\it a projective line, identified with ${\Bbb Q}(i)\cup\{\infty\}$, the set of all cusps
}&above (1.1.1)\cr  
\+${\cal R}$
&{\it a certain sum of simple Kloosterman sums
}&in (6.1)-(6.5)\cr 
\+$S(u,v;w)$
&{\it the simple Kloosterman sum
}&in (1.3.6)\cr
\+smooth $f$
&{\it a complex function, all partial derivatives of which are continuous 
}&beginning of \S 1.2\cr 
\+${\rm Supp}(f)$
&{\it the support of $f$, with respect to the topology of the Euclidean metric 
}&--\cr 
\+$S_{{\frak a},{\frak b}}\left(\omega_1 , \omega_2 ; c\right)$
&{\it a generalised Kloosterman sum
}&in (1.1.15)\cr  
\+$S_t(Q,X,N)$
&{\it a weighted mean value, with averaging over the level of the group $\Gamma$
}&(1.3.2)-(1.3.3)\cr 
\+${\cal S}\bigl( {\Bbb R}^n\bigr)$, ${\cal S}\bigl( {\Bbb C}^n\bigr)$ 
&{\it the Schwartz spaces on ${\Bbb R}^n$, and the Schwartz space on ${\Bbb C}^n$
}&beginning of \S 5\cr 
\+$V$
&{\it an irreducible cuspidal subspace of ${}^0L^2(\Gamma\backslash G)$
}&below (1.1.3)\cr 
\+$V_{K,\ell,q}$ 
&{\it a one dimensional subspace of $V$
}&below (1.1.5)\cr 
}

\medskip 

\goodbreak 
\proclaim The $L^2$-spaces. We define $L^2(\Gamma\backslash G)$ to be the Hilbert space of 
all square-integrable $\Gamma$-automorphic functions $f : G\rightarrow{\Bbb C}$. 
See the first four paragraphs of Subsection~1.1 
for the definitions of the terms `$\Gamma$-automorphic' 
and `square-integrable', and for the definition of the Hilbert-space inner product 
$\langle f , h\rangle_{\Gamma\backslash G}$ (for which the corresponding norm is 
$\| f\|_{\Gamma\backslash G}=\sqrt{\langle f , f\rangle_{\Gamma\backslash G}}\ $). 
\hfill\break $\hbox{\qquad}$ 
We define the term `cusp form' immediately below (1.1.10). 
The space ${}^0L^2(\Gamma\backslash G)$ is 
the closure of the subspace of $L^2(\Gamma\backslash G)$ 
spanned by cusp forms. The space ${}^{\rm e}L^2(\Gamma\backslash G)$ 
is the orthogonal complement in $L^2(\Gamma\backslash G)$ 
of the subspace ${\Bbb C}\oplus {}^0L^2(\Gamma\backslash G)$  
(with `${\Bbb C}$' here signifying the $1$-dimensional space of constant functions).  
More general spaces, $L^{2,{\rm cusp}}(\Gamma\backslash G ,\chi)$ and 
$L^{2,{\rm cont}}(\Gamma\backslash G ,\chi)$, are discussed    
in [19, Chapter~8]. When $\Gamma=\Gamma_0(q)$ and $\chi_0$ is the trivial character 
on $({\frak O}/q{\frak O})^*$, 
one has $L^{2,{\rm cusp}}(\Gamma\backslash G,\chi_0)={}^0L^2(\Gamma\backslash G)$ 
and 
$L^{2,{\rm cont}}(\Gamma\backslash G,\chi_0)={}^{\rm e}L^2(\Gamma\backslash G)$.

\goodbreak 
\proclaim The Scaling Matrices for Cusps. The notation `$g_{\frak c}$' denotes 
an element of $G$ satisfying both the equation $g_{\frak c}\infty ={\frak c}$ 
and the condition (1.1.1); we call any such element of $G$ 
a `scaling matrix for the cusp ${\frak c}$'. 

\goodbreak 
\proclaim Set-Theoretic Notation.   
We denote the cardinality of any set ${\cal A}$ by $|{\cal A}|$, so that 
$|\{ x\in{\Bbb R} : x^2=1\}|=2$ (for example). 
Given suitable functions $f$ and $g$, 
we define $g\circ f$ to be the function obtained by composing $f$ with $g$, 
so that $(g\circ f)(x)=g(f(x))$ whenever $g(f(x))$ is defined. 

\goodbreak 
\proclaim Algebraic Notation.  
If $R$ is a ring with identity, then $R^{*}$ denotes the group of units of $R$.  
When $U$, $V$ and $W$ are groups, the notation 
$U\leq W$ (resp. $U<W$) is used to indicate 
that $U$ is a subgroup (resp. proper subgroup) of $W$. 
If $U$ and $V$ are subgroups of the group $W$, then 
$W/V$, $U\backslash W$ and $U\backslash W/V$ denote 
the relevant sets of left cosets, right cosets and double cosets (respectively); 
and $\bigl[ W : U\bigr]$ denotes the index of $U$ in $W$, 
so that $\bigl[ W : U\bigr]=|W/U|$.
This notation for `quotients', such as $U\backslash W$ and $U\backslash W/V$, 
may apply in more general contexts. For example, if $U$ is a subgroup of $W$, and if 
$S$ is a subset of the elements of the group $W$ such that $uS\subseteq S$ 
for all $u\in U$, then $S$ can be expressed as a disjoint union of 
certain of the right cosets of $U$ in $W$, and so the notation $U\backslash S$ 
makes sense (as shorthand for the set of right cosets occurring in 
that disjoint union). Similar considerations apply in the 
case of quotients $S/V$ and $U\backslash S/V$, provided that the 
set $S$  is suitably invariant 
(either under left-multiplication by elements of $U\leq W$, or under 
right-multiplication by elements of the group $V\leq W$). 

\goodbreak 
\proclaim Notation for Upper and Lower Bounds.  
The greatest element of a set ${\cal X}\subset{\Bbb R}$ 
(where there is such an element) will be denoted by 
$\max{\cal X}$; similarly (and with a similar caveat) 
we use $\min{\cal X}$  to denote 
the least element of ${\cal X}$. Any notation of the form 
$\max_{A(x)} f(x)$, in which $A(x)$ is some statement about $x$, 
has the same meaning as $\max\{ f(x) : A(x)\ {\rm is\ true}\}$.  
Similarly, $\min_{A(x)} f(x)=\min\{ f(x) : A(x)\ {\rm is\ true}\}$.
\hfill\break $\hbox{\qquad}$ 
Where $B\geq 0$, we use the notation $O_{\alpha_1 ,\ldots,\alpha_n}(B)$  
to denote    
a complex-valued variable $\beta$  
satisfying a condition of the form $|\beta|\leq C(\alpha_1,\ldots,\alpha_n) B$, in which  
the `implicit constant' 
$C(\alpha_1,\ldots,\alpha_n)$ is positive and 
depends only on previously declared constants and 
$\alpha_1,\ldots,\alpha_n$.
As alternatives to an expression of the form 
`$\xi =O_{\alpha_1 ,\ldots,\alpha_n}(B)$', we may prefer
to follows Vinogradov in using either `$\xi\ll_{\alpha_1 ,\ldots,\alpha_n} B$', 
or `$B\gg_{\alpha_1 ,\ldots,\alpha_n} \xi$'. 
Where $A\geq 0$ and $B\geq 0$, the notation 
$A\asymp_{\alpha_1 ,\ldots,\alpha_n}B$ may be used to 
signify that one has both $A\ll_{\alpha_1 ,\ldots,\alpha_n}B$ 
and $B\ll_{\alpha_1 ,\ldots,\alpha_n}A$: we may also sometimes 
write this as `$A\ll_{\alpha_1 ,\ldots,\alpha_n}B\ll_{\alpha_1 ,\ldots,\alpha_n}A$'. 

\goodbreak 
\proclaim Epsilon.   
The part played by `$\varepsilon$' in our results is effectively  
that of an `arbitrarily small positive constant'. Indeed, although $\varepsilon$ 
is technically a variable, any practical application of our main 
results would involve a case in 
which $\varepsilon$ is assigned a value equal to some small absolute constant,  
such as the constant $10^{-10}$, for example (this is because 
implicit constants associated with the bounds in those results are 
dependent on $\varepsilon$).  
The value of $\varepsilon$ may vary from result to result, but
will generally remain fixed within each individual proof. 

\goodbreak 
\proclaim Complex Numbers.   
When $z\in{\Bbb C}$, the real and
imaginary parts of $z$, its absolute value and its complex conjugate 
are denoted, respectively, by  $\hbox{\rm Re}(z)$, $\hbox{\rm Im}(z)$, 
$|z|$ and $\overline{z}$ (so that $\overline{z}=\hbox{\rm Re}(z)-i\hbox{\rm Im}(z)$ 
and $|z|^2 =z\,\overline{z}\,$).

\goodbreak 
\proclaim Number-Theoretic Notation.  When $\alpha,\beta\in{\frak O}$ are not both zero, 
we may use the notation $(\alpha , \beta)$ 
to denote a highest common factor of $\alpha$ and $\beta$. This creates 
some ambiguity, for if $d$ is a highest common factor of 
$\alpha$ and $\beta$, then so too are the three other associates of $d$ 
(i.e. $id$, $-d$ and $-id$). This ambiguity does not, however, 
lead to any serious difficulties, since 
relations of the form $(\alpha ,\beta)\sim d$, or 
$|(\alpha ,\beta)|^2=n$, remain valid if 
the number $(\alpha ,\beta)$ is replaced by any one of its associates. 
If $\alpha$ and $\beta$ happen to be rational integer valued variables, 
then it is natural (and not inconsistent with the statements above) 
that we unambiguously put  
$(\alpha , \beta)=\max\{ d\in{\Bbb N} : d|\alpha\ {\rm and}\ d\mid\beta\}$.  
\hfill\break $\hbox{\qquad}$ 
When $b\in{\frak O}$, 
the ideal $\{ bm : m\in{\frak O}\}<{\frak O}$ 
is denoted by $b{\frak O}$; and, for $a,b\in{\frak O}$, we denote 
the coset 
$\{ a+n : n\in b{\frak O}\}\in {\frak O}/(b{\frak O})$ 
by $a+b{\frak O}$. Given $c\in{\frak O}$, we define  
$a,b\in{\frak O}$ to be `equivalent (to one another) modulo 
$c{\frak O}$' if and only if it is the case that $a+c{\frak O}=b+c{\frak O}$ 
(i.e. if and only if $c\mid (a-b)$). We write 
$a\equiv b\bmod c{\frak O}$ 
to signify that 
$a$ is equivalent to $b$ modulo $c{\frak O}$. 
\hfill\break $\hbox{\qquad}$ 
In relations such as $h m^*\equiv\ell\bmod c{\frak O}$, 
or in expressions such as the highest common factor $(h m^* , c)$, 
the rational expression $h m^*/c$, or 
(see (1.3.6)) the `simple Kloosterman sum' $S(h m^*,\ell;c)$,  
it is to be understood that $m^*$ denotes an
arbitrary element of ${\frak O}$ satisfying 
$m m^*\equiv 1\bmod c{\frak O}$.  
It is therefore implicit in such expressions that one has both $(m,c)\sim1$ and 
$(m^*,c)\sim 1$.

\goodbreak 
\proclaim Summation Related Conventions.   
Where there is no indication to the contrary, variables of 
summation range over all values in ${\frak O}$ 
consistent with all the conditions attached to the summation.
\hfill\break $\hbox{\qquad}$ 
When a `condition' of the form `$m\bmod c{\frak O}$' 
appears below the summation sign, 
it is to be understood that the variable of summation $m$ ranges 
(to the extent permitted 
by any other conditions of summation) 
over some fixed set of coset representatives 
$\{ m_1,\ldots ,m_{[{\frak O} : c{\frak O}]}\}$ of $c{\frak O}$ in ${\frak O}$. 
\hfill\break $\hbox{\qquad}$ 
If the very first condition of a summation is expressed in terms 
a certain set ${\cal X}$ (defined in terms of some variable parameters $x,y,\ldots\,$) 
having a certain fixed property, 
then it is that set itself (and not the variable parameters) which must 
be regarded as the variable of summation: such a summation 
may therefore only be considered well-defined when the value of its summand 
is uniquely determined, within the sum, by the set ${\cal X}$. 
This convention applies, for example, in (9.11), (9.16) and (9.18);  
and in (9.11) it results in the summation there effectively 
being such that the variable 
$\gamma$ runs (once) over the elements in a 
set $\{\rho_1\tau,\ldots ,\rho_{n}\tau\}$, 
where $\{\rho_1,\ldots ,\rho_{n}\}$  
is some complete set of representatives 
for the right-cosets of $\Gamma_{\frak a}'$ in $\Gamma_{\frak a}$, and 
$\tau$ is some element of $\Gamma$ such that $\tau{\frak b}={\frak a}$. 
We adopt a similar convention in respect of products. Hence in the equation~(5.35), 
for example, the final product may be expressed 
as $\prod_{p\in{\cal P}(q,k)}(1-|p|^{-2})$, 
where  ${\cal P}(q,k)$ denotes the set of those Gaussian primes $\varpi$ 
with ${\rm Re}(\varpi)>0$, ${\rm Im}(\varpi)\geq 0$,   
$\varpi\mid (q,k)$ and $\varpi\not{\mid}\,(q/(q,k))$. 
\hfill\break $\hbox{\qquad}$ 
Our notation for generalised Kloosterman sums and Fourier 
coefficients of Eisenstein series is ambiguous, in that it gives no 
indication of the dependence of  
those sums and coefficients on the group $\Gamma$. In order to compensate for this 
ambiguity, we
adopt the following conventions regarding summation: 
in sums involving Kloosterman sums 
$S_{{\frak a},{\frak b}}(m,n;c)$,  the symbol clarifying the 
relevant group $\Gamma$ 
is shown (within brackets) above the sign for summation over 
$c\in{}^{\frak a}{\cal C}^{\frak b}$;  
while, in  sums involving the Fourier coefficients  
$B_{\frak c}^{\frak a}(\omega;\nu ,p)$, the symbol for the relevant group 
$\Gamma$ appears (within brackets) above the sign for summation over the 
cusps ${\frak c}\in{\frak C}$. 

\bigskip 
\bigskip 

\goodbreak\centerline{\bf \S 2. Upper and Lower Bounds for the ${\bf K}$-Transform}

\bigskip 

This section concerns the ${\bf K}$-transform defined by
the equations~(1.2.2)-(1.2.4) of Theorem~1. In it we establish both upper and 
lower bounds for the ${\bf K}$-transform of suitable functions. 
Two new notational conventions are convenient for stating 
these results and their proofs. The first is the convention that 
$$(\alpha)_m=\alpha (\alpha+1)\cdots(\alpha +m-1)={\Gamma(\alpha +m)\over\Gamma(\alpha)}\;,\eqno(2.1)$$
for $\alpha\in{\Bbb C}$ and all non-negative integers $m$ 
(i.e. with $\Gamma(\alpha+m)/\Gamma(\alpha)$ being defined by 
analytic continuation at the removable singularities $\alpha =0,-1,-2,\ldots\ $\/). 
The second convention is that 
$${\bf M}\varphi(s) 
=\int_0^{\infty}\varphi(2\rho)\rho^{s-1}{\rm d}\rho\eqno(2.2)$$
(a variant of the Mellin transform) when the function $\varphi$ and $s\in{\Bbb C}$ are such 
that $\varphi(2x)x^{s-1}\in L^1(0,\infty)$.

\bigskip

\goodbreak\proclaim Lemma~2.1. Let 
$$f(z)=\varphi(|z|)\qquad\qquad\hbox{($z\in{\Bbb C}^{*}$),}\eqno(2.3)$$
where $\varphi : (0,\infty)\rightarrow{\Bbb C}$ is continuous, and 
compactly supported.  
Then, for $\nu\in{\Bbb C}$ and $p\in{\Bbb Z}$, one has  
$${\bf K}f(\nu,p)
={\bf K}f(\nu,|p|)
={\bf K}f(-\nu,|p|)\eqno(2.4)$$
and 
$${\bf K}f(\nu,p) 
=2\pi\sum_{m=0}^{\infty}{1\over m! (m+|p|)!}\,{\bf K}_m\varphi(\nu,|p|)\;,\eqno(2.5)$$
where, for non-negative integers $k$ and $m$,  
$${\bf K}_m\varphi (\nu, k)
={(-1)^k\over\sin(\pi\nu)}\left( 
{{\bf M}\varphi(-2\nu+4m+2k)\over\Gamma(-\nu +m+1)\Gamma(-\nu +m+1+k)}
-{{\bf M}\varphi(2\nu+4m+2k)\over\Gamma(\nu +m+1)\Gamma(\nu +m+1+k)}\right) .\eqno(2.6)$$
If, moreover, the function $\varphi$ has a continuous 
derivative $\varphi^{(j)} : (0,\infty)\rightarrow{\Bbb C}$, of order $j\in{\Bbb N}$, 
then 
$$(s)_j\,{\bf M}\varphi(s)
=(-2)^j\,{\bf M}\bigl(\varphi^{(j)}\bigr)(s+j)\qquad\qquad 
\hbox{($s\in{\Bbb C}$).}\eqno(2.7)$$

\medskip

\goodbreak
\noindent{\bf Proof.}\quad 
By (1.2.2)-(1.2.4),  
$${\bf K}f(\nu,p)
={1\over\sin(\pi\nu)}
\left({\bf J}f(-\nu,-p)-{\bf J}(\nu,p)\right) ,\eqno(2.8)$$
where 
$$\eqalignno{
{\bf J}f(\mu,k) &=\int_{{\Bbb C}^{*}}{\cal J}_{\mu,k}(z) f(z) {\rm d}_{\times} z = {}\cr
 &=\int\limits_0^{\infty}\int\limits_0^{2\pi} 
{\cal J}_{\mu,k}\left( r e^{i\theta}\right) f\left( r e^{i\theta}\right) {\rm d}\theta 
{{\rm d}r\over r} = {}\cr 
 &=\int\limits_0^{\infty}\int\limits_0^{2\pi} 
\left( {r\over 2}\right)^{2\mu} e^{-2ki\theta} 
\sum_{m=0}^{\infty} 
{(-1)^m (r/2)^{2m} e^{2mi\theta}\over m! \Gamma(\mu -k+m+1)}
\sum_{n=0}^{\infty} 
{(-1)^n (r/2)^{2n} e^{-2ni\theta}\over n! \Gamma(\mu +k+n+1)} 
\,f\left( r e^{i\theta}\right) {\rm d}\theta 
{{\rm d}r\over r}\;.\qquad &(2.9)}$$
By hypothesis, there exists an $R>0$ such that $f(z)=0$ for all $z\in{\Bbb C}$ 
with $|z|\geq R$. Therefore, and by virtue of the fact that 
the sums over $m$ and $n$ in (2.9) are uniformly absolutely convergent 
for $(r,\theta)\in (0,R)\times (0,2\pi)$, one may integrate term-by-term 
on the right-hand side of (2.9), and so obtain: 
$${\bf J}f(\mu,k) 
=2\pi\sum\sum_{\!\!\!\!\!\!\!\!\!\!\!\!m,n\geq 0} 
{(-1/4)^{m+n} 2^{-2\mu} {\bf I}f(\mu+m+n\,,\,k-m+n)\over 
m! n! \Gamma(\mu-k+m+1)\Gamma(\mu+k+n+1)}\;,\eqno(2.10)$$
where 
$${\bf I}f(\lambda ,\ell)
=\int\limits_0^{\infty} f_{\ell}(r) r^{2\lambda -1} {\rm d}r\qquad\quad  
{\rm with}\qquad\quad f_{\ell}(r)
={1\over 2\pi}\int\limits_0^{2\pi} f\left( r e^{i\theta}\right) e^{-2\ell i\theta} 
{\rm d}\theta\;.\eqno(2.11)$$ 
The expansion (2.10) is a result of Bruggeman and Motohashi 
[4, Equations~(11.8)~and~(11.9)].

By (2.3) and (2.11), 
$$f_{\ell}(r)={\varphi(r)\over 2\pi}
\int\limits_0^{2\pi}  e^{-2\ell i\theta} 
{\rm d}\theta
=\cases{\varphi(r) &if $\ell =0$;\cr 0 &if $0\neq 2\ell\in{\Bbb Z}$.}\eqno(2.12)$$ 
The summation in (2.10) is therefore effectively restricted to 
the pairs $m,n\in{\Bbb N}\cup\{ 0\}$ with $m-n=k$. 
If $k\geq 0$ then the equation $m-n=k$ implies 
that $m=n+k=n+|k|$; whereas if $k<0$, then it implies that 
$n=m-k=m+|k|$. Hence, by (2.10)-(2.12), one obtains 
(both for $k\geq 0$, and for $k<0$):
$${\bf J}f(\mu,k) 
=2\pi\sum_{m\geq 0} 
{(-1)^k 4^{-(\mu +2m+|k|)}\over 
m! (m+|k|)! \Gamma(\mu+m+1)\Gamma(\mu+m+|k|+1)}
\,\int\limits_0^{\infty}\varphi(r) r^{2(\mu+2m+|k|)-1} {\rm d}r\;.\eqno(2.13)$$

The equation~(2.13) shows that, subject to our hypotheses concerning $f$,  
we will have ${\bf J}f(\mu,k)={\bf J}f(\mu,|k|)$ for $\mu\in{\Bbb C}$, $k\in{\Bbb Z}$.
By this observation and (2.8) we obtain an identity,   
$${\bf K}f(\nu,p)
={1\over\sin(\pi\nu)}
\left({\bf J}f(-\nu,|p|)-{\bf J}f(\nu,|p|)\right) ,\eqno(2.14)$$ 
from which both of the equations in (2.4) follow, as immediate corollaries. 

We next have to deduce the result (2.5)-(2.6). 
By the substitution $r=2\rho$, one has 
$$\int\limits_0^{\infty}\varphi(r) r^{2(\pm\nu+2m+|p|)-1} {\rm d}r
=2^{\pm 2\nu+4m+2|p|}{\bf M}\varphi(\pm 2\nu+4m+2|p|)\;,$$
with ${\bf M}\varphi(s)$ defined by (2.2).
Hence, and by (2.13), 
$${\bf J}f(\pm\nu,|p|) 
=2\pi\sum_{m\geq 0} 
{(-1)^{|p|} {\bf M}\varphi(\pm 2\nu+4m+2|p|)\over 
m! (m+|p|)!\,\Gamma(\pm\nu+m+1)\Gamma(\pm\nu+m+|p|+1)}
\;.$$
The result (2.5)-(2.6) of the lemma now follows: one has only to 
substitute, for each term ${\bf J}(\pm\nu,|p|)$ in the equation~(2.14), 
the corresponding expansion (above); and then 
note that $\sum a_m -\sum b_m=\sum (a_m-b_m)$, whenever the 
first two series are convergent. 

For proof of (2.7), note that if $j\in{\Bbb N}$, and if the $j$-th order derivative 
$\varphi^{(j)} : (0,\infty)\rightarrow{\Bbb C}$ is continuous, then, 
given that the support of $\varphi$ is contained in some 
closed bounded interval $[a,b]\subset(0,\infty)$, it will 
follow by (2.2) and integration by parts that, for all $s\in{\Bbb C}$,  
$$\eqalign{ 
{\bf M}\bigl(\varphi^{(j)}\bigr) (s+j) 
 &=\int\limits_{a/2}^{b/2}\varphi^{(j)}(2\rho)\rho^{s+j-1} {\rm d}\rho = {}\cr 
 &=\left( {1\over 2}\varphi^{(j-1)}(b)(b/2)^{s+j-1}
-{1\over 2}\varphi^{(j-1)}(a)(a/2)^{s+j-1}\!\right) 
-{(s+j-1)\over 2}\!\int\limits_{a/2}^{b/2}\!\!\varphi^{(j-1)}(2\rho)\rho^{s+j-2} {\rm d}\rho = {}\cr 
 &=\left( 0-0\right) +{(s+j-1)\over (-2)}{\bf M}\bigl(\varphi^{(j-1)}\bigr) (s+j-1)\;.
}$$
By means of this 
identity one obtains a proof by induction 
of (2.7) (which is trivially true for $j=0$) 
\quad$\blacksquare$

\bigskip

\goodbreak 
\proclaim Lemma~2.2. Let $A>1$, $X>0$, $j\in{\Bbb N}\cup\{ 0\}$ and $0<\delta <1$; 
and let $f : {\Bbb C}^{*}\rightarrow{\Bbb C}$ be given by the equation~(2.3) 
of Lemma~2.1, 
where the function $\varphi : (0,\infty)\rightarrow{\Bbb C}$ 
has a continuous $j$-th order derivative, $\varphi^{(j)} : (0,\infty)\rightarrow{\Bbb C}$, 
and has its support contained within the 
closed interval $\bigl[ A^{-1}X^{-1/2} , A X^{-1/2}\bigr]$. 
Suppose, moreover, that $\nu\in{\Bbb C}$ and $p\in{\Bbb N}\cup\{ 0\}$; that either 
$\nu\in i{\Bbb R}$ or $0<\nu\leq 1-\delta$ and $p=0$; and 
that $j=0$ if $\nu =p=0$. 
Then 
$${\bf K}f(\nu,p)
\ll_{A,\delta ,j} 
{ Y_j \left(\left| X^{\nu}\right| +\left| X^{-\nu}\right|\right)
\min\left\{ 1+|\log X|\,,\,|\nu|^{-1}\right\} 
\over  |\nu +p|^j
\,(p!)\left| (\nu+1)_{p}\right|} 
\,\left( {A^2\over 4X}\right)^{p}
\exp\left( {A^4\over 16 X^2}\right) ,\eqno(2.15)$$
where 
$$Y_k=X^{-k/2}\max_{x\in{\Bbb R}}\left|\varphi^{(k)}(x)\right|\qquad\qquad  
\hbox{($k=0,1,\ldots ,j$).}\eqno(2.16)$$

\medskip

\goodbreak 
\noindent{\bf Proof.}\quad 
We showed in the proof of Theorem~1 that, for $p\in{\Bbb Z}$, the 
function $\nu\mapsto{\bf K}f(\nu,p)$ is entire. Moroever, 
when $p^j\neq 0$, the upper bound in 
(2.15) is a continuous function of $\nu$ for $\nu\in{\Bbb C}-\{ -1,-2,-3,\ldots\ \}$. 
It will therefore suffice 
to prove Lemma~2.2 for cases where $f$, $\nu$ and $p$ satisfy both the 
stated hypotheses of the lemma and the additional hypothesis that $\nu\neq 0$ 
(the cases of the lemma in which $\nu =0$ will then follow by 
taking the limit as $t\rightarrow 0+$ of cases with $\nu =it$). 

We shall also assume (henceforth) that 
$$j=0\qquad\quad\hbox{if}\qquad |\nu +p|<1\;.\eqno(2.17)$$ 
This is permissible, given that when $|\nu +p|<1$ the bound 
(2.15) will be at its strongest for $j=0$. Indeed, if $j\geq 1$ and  
$k\in\{ 0,1,\ldots ,j-1\}$ then, for each $x>0$ there exists some $x_1\in(0,x)$ such that  
$$\varphi^{(k)}(x)=\int\limits_0^x\varphi^{(k+1)}(y)\,{\rm d}y=(x-0)\varphi^{(k+1)}(x_1)\;;$$ 
so that, by (2.16) and the hypotheses concerning the support of 
$\varphi$, one has $Y_k X^{k/2}\leq A X^{-1/2} Y_{k+1}X^{(k+1)/2}$ and, 
hence, $Y_k\leq A Y_{k+1}$.  Consequently  
$0\leq Y_0 |\nu +p|^{-0}=Y_0\leq A^j Y_j\leq A^j Y_j |\nu +p|^{-j}$ 
when $|\nu +p|<1$;  
which establishes that, for such $\nu$ and $p$, the bound  
(2.15) for $j=0$ implies the bound (2.15) for all $j\in{\Bbb N}\cup\{ 0\}$.
 
Suppose now that the hypotheses of the lemma are satisfied, with 
$\nu\neq 0$ and $j$ satisfying (2.17). We may complete the proof 
of the lemma by showing that if $m$ is a non-negative integer then one has  
$${\bf K}_m\varphi(\nu,p)
\ll_{A,\delta ,j}
{ Y_j \left(\left| X^{\nu}\right| +\left| X^{-\nu}\right|\right)
\min\left\{ 1+|\log X|\,,\,|\nu|^{-1}\right\} 
\over  |\nu +p|^j
\left| (\nu+1)_{p}\right|} 
\,\left( {A^2\over 4X}\right)^{p+2m} ,\eqno(2.18)$$
with $Y_j$ given by (2.16), and with 
${\bf K}_m\varphi(\nu,k)$ defined by the equation~(2.6) (as in Lemma~2.1). 
For then, since 
$$\sum_{m=0}^{\infty} {1\over m! (m+p)!}
\,O_{A,\delta ,j}\left(\left( {A^2\over 4X}\right)^{2m}\right)  
\ll_{A,\delta,j} \sum_{m=0}^{\infty} {1\over m! p!} 
\,\left( {A^4\over 16 X^2}\right)^{m}  
= {1\over p!}\exp\left( {A^4\over 16 X^2}\right) ,$$
the bound (2.15) will follow, by (2.18), from the result (2.5) of Lemma~2.1.

Let $m\in{\Bbb N}\cup\{ 0\}$. By (2.6) and (2.7) of Lemma~2.1, and the definition (2.2), 
$$\eqalignno{
{\bf K}_m\varphi(\nu, p)
 &={(-1)^p (-2)^j\over\sin(\pi\nu)}
\Biggl( {{\bf M}\left(\varphi^{(j)}\right)(-2\nu+4m+2p+j)\over 
(-2\nu +4m+2p)_j\Gamma(-\nu +m+1)\Gamma(-\nu +m+1+p)} + {}\cr 
 &\qquad\qquad\qquad\qquad\qquad 
{} + (-1){{\bf M}\left(\varphi^{(j)}\right)(2\nu+4m+2p+j)\over 
(2\nu +4m+2p)_j\Gamma(\nu +m+1)\Gamma(\nu +m+1+p)}\Biggr) = {}\qquad\qquad\cr 
 &\quad =
(-1)^{p+j} 2^j\int\limits_0^{\infty}\varphi^{(j)}(2\rho) 
\rho^{4m+2p+j-1}\,{\Delta_{m,j}(\nu,p;\rho)\over\sin(\pi\nu)}\,{\rm d}\rho\;, &(2.19)}$$ 
where 
$$\Delta_{m,j}(\nu,p;\rho)
=\sum_{\epsilon =\pm 1}{-\epsilon\rho^{2\epsilon\nu}\over 
\Gamma(\epsilon\nu +m+1)\Gamma(\epsilon\nu +m+1+p)(2\epsilon\nu +4m+2p)_j}\;.\eqno(2.20)$$
Postponing consideration of the cases where $0<\nu\leq 1-\delta$ and $p=0$, 
let it temporarily be supposed that $0\neq\nu\in i{\Bbb R}$. 
Then (2.20) may be written as: 
$$\Delta_{m,j}(\nu,p;\rho)
=-2i\,{\rm Im}\left( {\rho^{2\nu}\over 
\Gamma(\nu +m+1)\Gamma(\nu +m+1+p)(2\nu +4m+2p)_j}\right)\qquad\qquad\hbox{($\rho>0$).}\eqno(2.21)$$ 
Given that $0\neq\nu\in i{\Bbb R}$, it follows trivially from (2.21) that 
$$\left| \Delta_{m,j}(\nu,p;\rho)\right| 
\leq{2\over\left| \Gamma(\nu +m+1)\Gamma(\nu +m+1+p)(2\nu +4m+2p)_j\right|}\qquad\qquad  
\hbox{($\rho >0$).}$$
Here it is helpful to note that, by the three functional equations $z\Gamma(z)=\Gamma(z+1)$, 
$\Gamma(\overline{z})=\overline{\Gamma(z)}$ and 
$\Gamma(z)\Gamma(1-z)=\pi/\sin(\pi z)$, one has: 
$$\Gamma(\nu +k+1)=(\nu+1)_k \Gamma(\nu+1)\qquad\qquad\hbox{($k\in{\Bbb N}\cup\{ 0\}$)}\eqno(2.22)$$
and 
$$\left|\Gamma^2(\nu +1)\right|
=|\nu\Gamma(\nu)\Gamma(1-\nu)|={|\pi\nu|\over |\sin(\pi\nu)|}\;.\eqno(2.23)$$
Hence, for $\rho>0$, 
$$\left|{\Delta_{m,j}(\nu,p;\rho)\over\sin(\pi\nu)}\right| 
\leq 
{2\over\left| 
\pi\nu (\nu+1)_m (\nu+1)_{m+p} (2\nu +4m+2p)_j\right|}
\leq 
{2\pi^{-1}\over  
\left|\nu (\nu+1)_{p}\right| |\nu +p|^j}\;.\eqno(2.24)$$
Since this bound is independent of $\rho$, and since 
the hypotheses concerning $\varphi$ imply that, for $\sigma\geq 0$,  
$$\eqalignno{
\int\limits_0^{\infty}\left|\varphi^{(j)}(2\rho)\right|\rho^{\sigma -1}{\rm d}\rho 
 &\leq\int\limits_{1/A}^A Y_j X^{j/2}\left({\alpha\over 2 X^{1/2}}\right)^{\sigma} 
{{\rm d}\alpha\over\alpha} \leq {}\cr
 &\leq 2^{-j} Y_j (4X)^{(j-\sigma)/2} A^{\sigma}\min\left\{ 2\log(A)\,,\,\sigma^{-1}\right\} \ll_{A,j} {}\cr
 &\ll_{A,j} Y_j\left({A^2\over 4X}\right)^{(\sigma -j)/2} (1+\sigma)^{-1} &(2.25)}$$
(where $Y_j$ is as in (2.16)), it therefore follows by (2.19) and (2.24) that 
when $|\nu|\geq 1$  we do obtain the desired bound (2.18):  
for, if $\nu\in i{\Bbb R}$ is such that $|\nu|\geq 1$, then
$\left| X^{\nu}\right|+\left| X^{-\nu}\right|=2$ and $\min\left\{ 
1+|\log X|\,,\,|\nu|^{-1}\right\}=|\nu|^{-1}$. This therefore 
completes our proof in respect of cases where 
$|\nu|\geq 1\,$ (given that the hypotheses of the lemma ensure that 
one has $\nu\in i{\Bbb R}$ in such cases).

It only remains to prove the cases of the lemma in which either 
$p\in{\Bbb Z}$ and $\nu\in i{\Bbb R}$ with $0<|\nu|<1$, or 
$p=0$ and $\nu\in(0,1-\delta]$. Taking the former case first 
(i.e. supposing now that $\nu\in i{\Bbb R}$ and $0<|\nu|<1$), 
we deduce from (2.21) and (2.22) that 
$$\eqalignno{
{i\over 2}\,\Delta_{m,j}(\nu,p;\rho) 
 &={\rm Im}\left(\rho^{2\nu}\right) 
{\rm Re}\left({1\over\Gamma^2(\nu+1)(\nu+1)_m(\nu+1)_{m+p}(2\nu+4m+2p)_j}\right) + {} 
 &(2.26)\cr
 &\qquad +{\rm Re}\left(\rho^{2\nu}\right) 
{\rm Im}\left({1\over\Gamma^2(\nu+1)(\nu+1)_m(\nu+1)_{m+p}(2\nu+4m+2p)_j}\right) 
.}$$
Given that $0\neq\nu\in i{\Bbb R}$, and assuming that $\varphi^{(j)}(2\rho)\neq 0$ 
(so that $A^{-1}X^{-1/2}\leq 2\rho\leq A X^{-1/2}$), we will have, in (2.26),  
$$\left| {\rm Re}\left(\rho^{2\nu}\right)\right|\leq\left|\rho^{2\nu}\right|=1$$
and 
$$\eqalign{
\left|{\rm Im}\left(\rho^{2\nu}\right)\right| 
=\left|\sin(2|\nu|\log\rho)\right| 
 &\leq\min\left\{ 1\,,\,2|\nu\log\rho|\right\} \leq {}\cr 
 &\leq\min\left\{ 1\,,\,|\nu|(2|\log A|+|\log(4X)|)\right\} \ll_A {}\cr
 &\ll_A |\nu|\min\left\{ |\nu|^{-1}\,,\,1+|\log X|\right\} .}$$
Moreover, since $1/\Gamma(z)$ is entire, it is trivially the case that  
$${\rm Re}\left({1\over\Gamma^2(\nu+1)(\nu+1)_m(\nu+1)_{m+p}(2\nu+4m+2p)_j}\right)
\ll{1\over \left|(\nu+1)_p (\nu+p)^j\right|}\quad 
\hbox{for $\nu\in i{\Bbb R}$ with $0<|\nu|<1$.}$$

We now lack only a bound for the latter of the two imaginary parts 
that appear in (2.26). In order to obtain a suitable bound for that imaginary 
part, we make use of the expansion 
$$\eqalignno{
 &{\rm Im}\!\left({1\over\Gamma^2(\nu+1)(\nu+1)_m(\nu+1)_{m+p}(2\nu+4m+2p)_j}\right) = {} 
 &(2.27)\cr 
 &\qquad\qquad\qquad\qquad\qquad\qquad\quad  
={\rm Im}\left({1\over\Gamma^2(\nu+1)}\right) 
{\rm Re}\left({1\over (\nu+1)_m(\nu+1)_{m+p}(2\nu+4m+2p)_j}\right) + {}\cr
 &\qquad\qquad\qquad\qquad\qquad\qquad\qquad\ +{\rm Re}\left({1\over\Gamma^2(\nu+1)}\right) 
{\rm Im}\left({1\over (\nu+1)_m(\nu+1)_{m+p}\,(2\nu+4m+2p)_j}\right) ,\qquad\quad}$$
in which (assuming that $\nu\in i{\Bbb R}$ and $|\nu|< 1$) one has: 
$${\rm Im}\left({1\over\Gamma^2(\nu+1)}\right)
={1\over 2i}\left({1\over\Gamma^2(1+\nu)}-{1\over\Gamma^2(1-\nu)}\right)
=\left({\Gamma^{-2}(1+\nu)-\Gamma^{-2}(1-\nu)\over 2i\nu}\right)\nu\ll|\nu|\;,$$
$${\rm Re}\left({1\over\Gamma^2(\nu+1)}\right)
\ll\left|{1\over\Gamma^2(\nu+1)}\right|\ll 1$$
and
$${\rm Re}\left({1\over (\nu+1)_m(\nu+1)_{m+p}\,(2\nu+4m+2p)_j}\right)
\ll{1\over\left|(\nu+1)_p (\nu+p)^j\right|}\;.$$
With regard to the final imaginary part in (2.27), we note firstly 
that, if  $-1\leq\theta\leq 1$, then (given that $0\neq\nu\in i{\Bbb R}$) 
one has:  
$$\eqalign{
 &{{\rm d}\over{\rm d}z}(z+1)_m(z+1)_{m+p}\,(2z+4m+2p)_j\Bigr|_{z=\theta\nu} = {}\cr 
 &\qquad\quad\ =(z+1)_m(z+1)_{m+p}\,(2z+4m+2p)_j
\left(\sum_{h=1}^m{1\over z+h}+\sum_{k=1}^{m+p}{1\over z+k}
+\sum_{\ell=0}^{j-1}{2\over 2z+4m+2p+\ell}\right)\Biggr|_{z=\theta\nu} \ll {}\cr 
 &\qquad\quad\ \ll\left|(\nu+1)_m(\nu+1)_{m+p}\,(2\nu+4m+2p)_j\right| 
\left(\sum_{h=1}^m {1\over |\nu+h|}+\sum_{k=1}^{m+p}{1\over |\nu+k|}
+\sum_{\ell=0}^{j-1}{2\over |2\nu+4m+2p+\ell|}\right) ,
}$$
$$\sum_{h=1}^m {1\over |\nu+h|}\leq\sum_{k=1}^{m+p}{1\over |\nu+k|}
\leq\sum_{k=1}^{m+p}{1\over k}\leq 1+\log(1+m+p)$$
and, by virtue of the additional hypothesis (2.17), 
$$\sum_{\ell=0}^{j-1}{2\over |2\nu+4m+2p+\ell|}
\leq\sum_{\ell=0}^{j-1}{1\over |\nu+p|}
\leq\sum_{\ell=0}^{j-1} 1=j\;.$$
Therefore, given that $0\neq\nu\in i{\Bbb R}$, it follows by the mean value 
theorem of differential calculus that 
$$\eqalign{\left|{\rm Im}\left({1\over 
(\nu+1)_m(\nu+1)_{m+p}\,(2\nu+4m+2p)_j}\right)\right| 
 &={\left|{\rm Im}\left( (\nu+1)_m(\nu+1)_{m+p}\,(2\nu+4m+2p)_j\right)\right|\over 
\left|(\nu+1)_m(\nu+1)_{m+p}\,(2\nu+4m+2p)_j\right|^2} = {}\cr
 &={\left|\sum\limits_{\epsilon=\pm 1}\epsilon(\epsilon\nu+1)_m(\epsilon\nu+1)_{m+p}
\,(2\epsilon\nu+4m+2p)_j\right|\over 
2\left|(\nu+1)_m(\nu+1)_{m+p}\,(2\nu+4m+2p)_j\right|^2} \leq {}\cr 
 &\leq 
{|\nu -(-\nu)|\,|2+2\log(m+p+1)+j|\over 
2\left|(\nu+1)_m(\nu+1)_{m+p}\,(2\nu+4m+2p)_j\right|} \ll {}\cr 
 &\ll {|\nu|\log(m+p+2)\over \left|(\nu+1)_p\,(\nu+p)^j\right|}\;. 
}$$
By combining the above bound with those obtained just below (2.27), 
we may deduce from (2.27) that 
$${\rm Im}\left({1\over\Gamma^2(\nu+1)(\nu+1)_m(\nu+1)_{m+p}(2\nu+4m+2p)_j}\right)
\ll {|\nu|\log(m+p+2)\over \left|(\nu+1)_p\,(\nu+p)^j\right|}\;. 
$$

Since  
$\sin(\pi\nu)/(\pi\nu)=\sinh(|\pi\nu|)/|\pi\nu|\geq 1$ when 
$0\neq\nu\in i{\Bbb R}$ (as we currently suppose), 
it follows by this last bound above, and by (2.26) 
and the bounds found between (2.26) and (2.27),  
that one has 
$${\Delta_{m,j}(\nu,p;\rho)\over\sin(\pi\nu)}
\ll_A {\min\left\{|\nu|^{-1}\,,\,1+|\log X|\right\}+\log(m+p+2)\over 
\left|(\nu+1)_p\,(\nu+p)^j\right|}\;.$$
By this bound, and (2.19) and (2.25), it follows that 
in the cases where $\nu\in i{\Bbb R}$ and $|\nu|<1$ we 
do obtain the desired bound (2.18): for these cases one has 
$$\eqalign{
{\min\left\{|\nu|^{-1}\,,\,1+|\log X|\right\}+\log(m+p+2)\over 1+4m+2p+j}
 &\leq \min\left\{|\nu|^{-1}\,,\,1+|\log X|\right\}+{\log(m+p+2)\over m+p+1} < {}\cr 
 &<\min\left\{|\nu|^{-1}\,,\,1+|\log X|\right\}+1 \ll {}\cr 
 &\ll \min\left\{|\nu|^{-1}\,,\,1+|\log X|\right\}\;.}$$

Since (2.18) (for all $m\in{\Bbb N}\cup\{ 0\}$) implies the result (2.15), 
we now have disposed of all those cases of the 
lemma in which $\nu\in i{\Bbb R}$. 
It only remains to consider the cases where (as we shall henceforth suppose) 
one has and $p=0$ and $0<\nu\leq 1-\delta$. In these cases 
$|\nu +p|=\nu<1$, so that by (2.17) we are also to assume now that $j=0$. 
By (2.20) and (2.22), we have 
$$\eqalignno{\Delta_{m,0}(\nu,0;\rho) 
 &={\rho^{-2\nu}\over\Gamma^2(-\nu+m+1)}-{\rho^{2\nu}\over\Gamma^2(\nu+m+1)} = {}
 &(2.28)\cr 
 &={\rho^{-2\nu}\over (-\nu+1)_m^2\Gamma^2(-\nu+1)}
-{\rho^{2\nu}\over (\nu+1)_m^2\Gamma^2(\nu+1)}=D_1+D_2\;, &(2.29)}$$
where 
$$D_1={\rho^{-2\nu}\over\Gamma^2(-\nu+1)}\left( 
{1\over (-\nu+1)_m^2}-{1\over (\nu+1)_m^2}\right)\qquad\hbox{and}\qquad 
D_2={1\over (\nu+1)_m^2}\left( 
{\rho^{-2\nu}\over\Gamma^2(-\nu+1)}
-{\rho^{2\nu}\over\Gamma^2(\nu+1)}\right) .$$
Since the function $x\mapsto 1/\Gamma(x)$ bounded on the interval $(0,\infty)$, 
it follows trivially from (2.28) that 
$$\Delta_{m,0}(\nu,0;\rho)\ll \rho^{2\nu}+\rho^{-2\nu}\qquad\qquad\hbox{($\rho>0$).}\eqno(2.30)$$ 

The bound (2.30) will not, by itself, suffice. The required alternative 
bound will be obtained by considering the terms $D_1$ and $D_2$ in (2.29).

Since $\nu\in{\Bbb R}$,  the mean value theorem of 
real differential calculus implies that, for some $\theta\in (-1,1)$,  
$$D_2
={1\over (\nu+1)_m^2}\,\left( (-\nu)-\nu\right) 
\left( {2\log\rho\over\Gamma^2(\theta\nu +1)}
+{{\rm d}\over{\rm d}x}\,{1\over\Gamma^2(x+1)}\biggr|_{x=\theta\nu}\right)
\rho^{2\theta\nu}\;.$$
Given that $0<\nu\leq 1-\delta <1$, so that $-1<\theta\nu<1$ here, 
one therefore has 
$$D_2\ll \nu\left( |\log\rho|+1\right)\left(\rho^{2\nu}+\rho^{-2\nu}\right)\qquad\qquad  
\hbox{($\rho>0$).}$$

To estimate the term $D_1$ we begin by observing that, 
since $\Gamma(2-\nu)\asymp 1$  for $0<\nu\leq 1-\delta<1$,  
$$D_1={\rho^{-2\nu}(1-\nu)^2\over\Gamma^2(2-\nu)(1-\nu)_m^2}\left( 
1- {(1-\nu)_m^2\over (1+\nu)_m^2}\right)
\asymp {\rho^{-2\nu}(1-\nu)^2\over (1-\nu)_m^2}\left( 
1- {(1-\nu)_m^2\over (1+\nu)_m^2}\right) .$$
Since $0<\nu\leq 1-\delta <1$, one has (in the above): 
$$1>{(1-\nu)_m\over (1+\nu)_m}
=\prod_{k=1}^m {k-\nu\over k+\nu}
=\prod_{k=1}^m {\left( 1-{\nu\over k}\right)\over 
\left( 1+{\nu\over k}\right)}
>\prod_{k=1}^m \left( 1-{\nu\over k}\right)^2=P_m^2(\nu)>0\;,$$
where 
$$P_m(\nu)=\prod_{k=1}^m \left( 1-{\nu\over k}\right)
=\prod_{k=1}^m {(k-\nu)\over k}={(1-\nu)_m\over m!}\;.
$$
Hence, 
$$D_1
\ll {(1-\nu)^2\rho^{-2\nu}\over\left( m! P_m(\nu)\right)^2}
\left( 1-P_m^4(\nu)\right) 
<{2\rho^{-2\nu}\over (m!)^2}\left({P_m^{-2}(\nu)-P_m^2(\nu)\over 2}\right) 
\ll (m+1)^{-2}\rho^{-2\nu}\sinh\left( 2 Q_m(\nu)\right) ,$$
where 
$$\eqalign{0\leq Q_m(\nu) &=\log\!\left({1\over P_m(\nu)}\right) = {}\cr  
 &=\sum_{k=1}^m \left(\log(1)-\log\left( 1-{\nu\over k}\right)\right) \leq {}\cr
 &\leq\sum_{k=1}^m {\left({\nu\over k}\right)\over\left( 1-{\nu\over k}\right)} = {}\cr 
 &=\sum_{k=1}^m\left( {\nu\over k}+{\left({\nu\over k}\right)^2\over\left( 1-{\nu\over k}\right)}\right) 
<\nu\sum_{k=1}^m {1\over k} +{\nu^2\over (1-\nu)}\sum_{k=1}^{\infty}{1\over k^2} 
\leq\nu\!\left(\log(1+m) +O\!\left({1\over 1-\nu}\right)\!\right) .
}$$
Therefore, and since $\sinh(y)\leq y\exp(y)$ for $y\geq 0$, we have 
$$\eqalign{
D_1 &\ll (m+1)^{-2}\rho^{-2\nu} \nu\!\left(\log(1+m) +O\!\left({1\over 1-\nu}\right)\!\right) 
(m+1)^{2\nu}\exp\left( O\left({1\over 1-\nu}\right)\right) \ll {}\cr
 &\ll \exp\left( O\left({1\over 1-\nu}\right)\right) 
{(\log(m+1) +1)\over (m+1)^{2-2\nu}}\,\nu\rho^{-2\nu}\ll_{\delta}\ \nu\rho^{-2\nu} 
}$$
(bearing in mind that $0<\nu\leq 1-\delta <1$, and that 
$\log x\ll_{\delta} x^{\delta}$ for $x\geq 1$ and $\delta >0$).

By the bounds just obtained for $D_2$ and $D_1$,
and by the equation~(2.29), it follows that 
$$\Delta_{m,0}(\nu,0;\rho)\ll_{\delta} 
\,\nu\left(\rho^{2\nu}+\rho^{-2\nu}\right)\left( 1+|\log\rho|\right)\qquad\qquad 
\hbox{($\rho>0$).}$$
In combination with the bound (2.30), this shows that when 
$0<\nu\leq 1-\delta <1$ one has 
$${\Delta_{m,0}(\nu,0;\rho)\over\sin(\pi\nu)}\ll_{\delta,A} 
\,\left( X^{-\nu}+X^{\nu}\right)
\min\left\{ \nu^{-1}\,,\,1+|\log X|\right\}\qquad\    
{\rm for}\quad {1\over 2A X^{1/2}}\leq\rho\leq {A\over 2X^{1/2}}\;.\eqno(2.31)$$
By hypothesis, $\varphi(x)=0$ unless $1/A\leq X^{1/2}x\leq A$; 
it therefore follows by the last bound above, the case $j=p=0$ of (2.19), 
and the case $j=0$ of (2.25), that when $0<\nu\leq 1-\delta<1$ and $j=p=0$
we do obtain the desired bound (2.18), 
for all $m\in{\Bbb N}\cup\{ 0\}$. 

We have now shown the bound (2.18) to hold 
in all relevant cases;    
given what was established in the paragraph containing (2.18), this
completes the proof of the lemma\quad$\blacksquare$ 

\bigskip

\goodbreak 
\noindent{\bf Remark.}\quad 
Lemma~2.2 implies that the upper bound (1.2.16) holds, subject to 
the hypotheses of the Corollary to Theorems~1 and~2. 
To see this note firstly that those hypotheses 
ensure that $X\geq 2$, that $\nu\in i{\Bbb R}$, and that,  
by Lemma~2.2, one obtains the result stated in (2.15)-(2.16) for 
$j=0,1,2,3$. 
The cases of (1.2.16) in which $|\nu|\geq 1$ therefore follow 
immediately from the case $j=3$ of (2.15)-(2.16), 
since one has there: $\bigl| X^{\nu}\bigr| +\bigl| X^{-\nu}\bigr| =2$, 
$\,\exp\!\left( A^4 /16 X^2\right)\ll_A 1$ 
and  (see (2.1)) $\,\bigl| (\nu+1)_p\bigr|\geq p!$. 
The remaining cases (where $|\nu|<1$) are implied by the 
case $j=0$ of (2.15)-(2.16): for, by the observation 
following (2.17),  the terms $Y_0$ and $Y_3$ defined by 
(2.16) (for $j=3$, say) satisfy $Y_0\ll_A Y_3\ll (1+|\nu|)^{-4}\,Y_3$, given that 
$|\nu|<1$. 

\bigskip

\goodbreak\proclaim Lemma~2.3. Let the hypotheses of the previous lemma 
concerning $A$, $X$, $j$, $\delta$, $f$, $\varphi$ and $\varphi^{(j)}$ 
be satisfied. Let $j=0$. Suppose that $0<\nu\leq 1-\delta$. 
Then 
$${\bf K}f(\nu,0)
\ll_{A,\delta} 
\ \int\limits_0^{\infty}\left|\varphi(r)\right| {{\rm d}r\over r} 
\,\min\left\{ 1+|\log X|\,,\,\nu^{-1}\right\}
\left( X^{\nu}+X^{-\nu}\right) 
\exp\left( {A^4\over 16 X^2}\right) .\eqno(2.32)$$
Suppose, moreover, that $\varphi : (0,\infty)\rightarrow [0,\infty)$; and that, 
for some $\varepsilon >0$, one has 
$$AX^{-1/2}\leq 2\exp\left( {\Gamma'(\delta)\over\Gamma(\delta)}-\varepsilon\right) . 
\eqno(2.33)$$ 
Then 
$${\bf K}f(\nu,0)
\gg_{A,\delta,\varepsilon} 
\ \int\limits_0^{\infty}\varphi(r) {{\rm d}r\over r}\, 
\min\left\{ 1+|\log X|\,,\,\nu^{-1}\right\} X^{\nu}\;.\eqno(2.34)$$
 
\medskip

\goodbreak 
\noindent{\bf Proof.}\quad 
Let $p=0$. Since all the hypotheses of the case $p=j=0$ of Lemma~2.2 are satisfied, 
and since $\nu\neq 0$, we obtain (as in the proof of Lemma~2.2) 
the case $p=j=0$ of the results in (2.19). Therefore  
$${\bf K}_m\varphi(\nu ,0) 
=\int\limits_0^{\infty}\varphi(2\rho)\rho^{4m}
\,{\Delta_{m,0}(\nu,0;\rho)\over\sin(\pi\nu)}\,{{\rm d}\rho\over\rho}\qquad\qquad  
\hbox{($m\in{\Bbb N}\cup\{ 0\}$),}\eqno(2.35)$$
with ${\bf K}_m\varphi(\nu ,k)$ as in Lemma~2.1; and 
with $\Delta_{m,0}(\nu,0;\rho)$ as in (2.28). 
Similarly, we obtain the bound (2.31), for all $m\in{\Bbb N}\cup\{ 0\}$. 
Given that  
${\rm Supp}(\varphi)\subseteq\bigl[ A^{-1}X^{-1/2}\,,\,A X^{-1/2}\bigr]$, 
it follows by (2.31) and (2.35) that 
$$\eqalign{
{\bf K}_m\varphi(\nu,0) &\ll_{A,\delta} 
\left( X^{\nu}+X^{-\nu}\right)
\min\left\{ 1+|\log X|\,,\,\nu^{-1}\right\} 
\int\limits_{A^{-1}X^{-1/2}}^{A X^{-1/2}}
\left|\varphi(r)\right|\left({r\over 2}\right)^{4m}{{\rm d}r\over r} \leq {}\cr 
 &\leq\left({A^2\over 4X}\right)^{2m}
\left( X^{\nu}+X^{-\nu}\right)
\min\left\{ 1+|\log X|\,,\,\nu^{-1}\right\} 
\int\limits_0^{\infty}
\left|\varphi(r)\right|{{\rm d}r\over r}\qquad\qquad\qquad\hbox{($m\in{\Bbb N}\cup\{ 0\}$).}
}$$
This upper bound for ${\bf K}_m\varphi(\nu,0)$ implies the result (2.32): for 
$\exp\left( A^4/16 X^2\right)=\sum_{m=1}^{\infty}\left( A^2 /4X\right)^{2m}/m!\,$ and, 
by the case $p=0$ of the result (2.5) of Lemma~2.1, one has 
$\sum_{m=0}^{\infty}\left|{\bf K}_m\varphi(\nu,0)\right|/ m!\gg {\bf K}f(\nu ,0)$.

It remains for us to prove the conditional lower bound (2.34).  
We therefore suppose now (in addition to what has been assumed) 
that $\varepsilon >0$; that (2.33) holds; and that 
$\varphi : (0,\infty)\rightarrow [0,\infty)$. 
By the result (2.5) of Lemma~2.1, the bound (2.34) will follow if 
it can be shown that one has both 
$${\bf K}_m\varphi(\nu ,0)\geq 0\qquad\qquad\hbox{($m\in{\Bbb N}$)}\eqno(2.36)$$
and 
$${\bf K}_0\varphi(\nu ,0)\gg_{A,\delta,\varepsilon} 
\ \int\limits_0^{\infty}
\varphi(r) {{\rm d}r\over r}\,X^{\nu} 
\min\left\{ 1+|\log X|\,,\,\nu^{-1}\right\}\;.\eqno(2.37)$$
Moreover, given our assumptions concerning the function $\varphi$ and its support, 
and given the lower bound 
$$\sin(\pi\nu)\gg_{\delta}\nu>0$$
(implied by our assumption that $1>1-\delta\geq\nu >0$), it follows 
by the equation~(2.35) that, in order to prove the lower bounds (2.36) and (2.37),  
it will be enough to show that, 
for $\rho\in\bigl[{1\over 2}A^{-1}X^{-1/2}, {1\over 2}A X^{-1/2}\bigr]$, 
one has both 
$$\Delta_{m,0}(\nu,0;\rho)\geq 0\qquad\qquad\hbox{($m\in{\Bbb N}$)}\eqno(2.38)$$
and 
$$\Delta_{0,0}(\nu,0;\rho)\gg_{A,\delta,\varepsilon}
\ X^{\nu} \min\left\{ 1\,,\,(1+|\log X|)\nu\right\}\eqno(2.39)$$
(with $\Delta_{m,0}(\nu,0;\rho)$ as in the equation~(2.28), in the proof of Lemma~2.2). 

We now complete the proof of the lemma by establishing first 
(2.38), and then (2.39). It is henceforth assumed that $m\in{\Bbb N}\cup\{ 0\}$ and 
$${1\over 2}\,A^{-1}X^{-1/2}\leq\rho\leq{1\over 2}\,A X^{-1/2}\;.\eqno(2.40)$$ 

By (2.28) and the fundamental theorem of calculus, 
$$\eqalignno{ 
\Delta_{m,0}(\nu,0;\rho) &=-\int_{-\nu}^{\nu} 
\left( {{\rm d}\over{\rm d}x}\,{\rho^{2x}\over\Gamma^2(x+m+1)}\right) {\rm d}x = {}\cr
 &=2\int_{-\nu}^{\nu}
{\rho^{2x}\over\Gamma^2(x+m+1)} 
\left(\log\left({1\over\rho}\right) +{\Gamma'(x+m+1)\over\Gamma(x+m+1)}\right)
{\rm d}x\;. &(2.41)}$$
By [24, Section~12.16],  
$${{\rm d}\over{\rm d}y}\,{\Gamma'(y+1)\over\Gamma(y+1)} 
={{\rm d}^2\over{\rm d}y^2}\,\log\Gamma(y+1) 
=\sum_{n=1}^{\infty}{1\over (y+n)^2} >0\qquad\quad\hbox{for $\quad y>-1$.}$$
Since $0<\nu\leq 1-\delta <1$ and $m\geq 0$, it follows 
that one has 
$${\Gamma'(x+m+1)\over\Gamma(x+m+1)}
\geq{\Gamma'(x+1)\over\Gamma(x+1)}
\geq{\Gamma'(1-\nu)\over\Gamma(1-\nu)}
\geq{\Gamma'(\delta)\over\Gamma(\delta)}\qquad\quad\hbox{for $\quad x\geq -\nu$.}\eqno(2.42)$$
The hypothesis (2.33), and the inequalities in (2.40), therefore ensure that   
$$\log\left({1\over\rho}\right) 
+{\Gamma'(x+m+1)\over\Gamma(x+m+1)}
\geq\varepsilon >0\qquad\quad\hbox{for $\quad x\geq -\nu$.}\eqno(2.43)$$
Since (2.41) and (2.43) combine to give the inequality  
$\Delta_{m,0}(\nu,0;\rho)\geq 0$, this concludes our proof of (2.38). 

We begin the proof of (2.39) by observing that the case $m=0$ of (2.43) shows that  
$$\log\left({1\over\rho}\right)\geq\varepsilon -{\Gamma'(1-\nu)\over\Gamma(1-\nu)}\;,$$
where, by the (2.42) for $x=0$, one has: 
$$-{\Gamma'(\delta)\over\Gamma(\delta)}\geq -{\Gamma'(1-\nu)\over\Gamma(1-\nu)}
\geq -{\Gamma'(1)\over\Gamma(1)}=\gamma >0$$
($\gamma$ here being Euler's constant, by [24, Section~12.16] for example). 
Therefore one has  
$$\log\left({1\over\rho}\right)\geq -{\Gamma'(1-\nu)\over\Gamma(1-\nu)} 
\left( 1+{\varepsilon\over G(\delta)}\right) >0\;,$$
where $G(\delta)=-\Gamma'(\delta)/\Gamma(\delta)\geq\gamma$. 
This implies the inequalities  
$${G(\delta)\over (G(\delta)+\varepsilon)}\,\log\left({1\over\rho}\right) 
\geq -{\Gamma'(1-\nu)\over\Gamma(1-\nu)} >0 \;,$$
from which it follows that 
$$\log\left({1\over\rho}\right) 
\geq -{\Gamma'(1-\nu)\over\Gamma(1-\nu)}
+{\varepsilon\over (G(\delta)+\varepsilon)}\,\log\left({1\over\rho}\right) >0\;.$$
By this and (2.42), it follows that one has 
$${\Gamma'(x+1)\over\Gamma(x+1)}+\log\left({1\over\rho}\right) 
\geq {\varepsilon\over (G(\delta)+\varepsilon)}\,\log\left({1\over\rho}\right)
\gg_{\delta,\varepsilon}\ \log\left({1\over\rho}\right)\qquad\quad  
\hbox{for $\quad x\geq -\nu$.}\eqno(2.44)$$

Postponing our application of (2.44), we 
take the opportunity to note here that, similarly to the above, 
it is implied by the hypothesis (2.33) that 
$$\log\left({2 X^{1/2}\over A}\right)\geq\varepsilon -{\Gamma'(\delta)\over\Gamma(\delta)}
\geq\varepsilon +\gamma>\gamma >0\;.$$
This, since $A>1$, already shows that $\log(4X)>2\gamma >1$. 
Moreover, the inequality $\log(2 X^{1/2}/A)>\gamma$ is (for $A>1$) 
equivalent to the inequality 
$$\log\left({2X^{1/2}\over A}\right)
>{\gamma\over(\gamma+\log A)}\,\log\left( 2X^{1/2}\right) ,$$
which (since we assume (2.40)) implies that we have 
$$\log\left({1\over\rho}\right) 
\gg_A\ \log\left( 2X^{1/2}\right)\gg 1+|\log X|\eqno(2.45)$$
(the final bound following since 
$1<\log(4X)$ and $|\log X|=|\log(4X)-\log 4|\leq \log(4X)+\log 4\ll\log(4X)$).

By (2.44) and the case $m=0$ of (2.41), we obtain the lower bound 
$$\Delta_{0,0}(\nu,0;\rho)
\gg_{\delta,\varepsilon}\ \log\left({1\over\rho}\right) 
\int_{-\nu}^{\nu} 
{\rho^{2x}\over\Gamma^2(x+1)}\,{\rm d}x\;,$$
where, since $0<\nu\leq 1-\delta<1$, and since $1/\Gamma(x)$ is continuous and 
positive-valued for $x>0$, one has: 
$$\int_{-\nu}^{\nu} 
{\rho^{2x}\over\Gamma^2(x+1)}\,{\rm d}x
\gg_{\delta}
\ \int_{-\nu}^{\nu} 
\rho^{2x}\,{\rm d}x
=\int_{-\nu}^{\nu}\left({1\over\rho}\right)^{2y}\,{\rm d}y 
={(1/\rho)^{2\nu}-(1/\rho)^{-2\nu}\over 2\log(1/\rho)}\;.$$
It therefore follows that 
$$\Delta_{0,0}(\nu,0;\rho)
\gg_{\delta,\varepsilon}
\ {(1/\rho)^{2\nu}-(1/\rho)^{-2\nu}\over 2}
=\sinh\left( 2\nu\log\left({1\over\rho}\right)\right) .$$
Now $\sinh(x)\gg\min\{ 1 , x\}\exp(x)$, for $x>0$, so 
that by (2.40), (2.45), the hypothesis $\nu >0$ and the lower bound for 
$\Delta_{0,0}(\nu,0;\rho)$ just obtained, we have:  
$$\Delta_{0,0}(\nu,0;\rho)
\gg_{\delta,\varepsilon}\ \min\left\{ 1\,,\,2\nu\log\left({1\over\rho}\right)\right\} 
\left({1\over\rho}\right)^{2\nu} 
\gg_A\ \min\left\{ 1\,,\,(1+|\log X|)\nu\right\} 
\left({4X\over A^2}\right)^{\nu}\;.$$
These bounds imply (2.39) (given that $0<\nu<1$), so that both (2.38) and (2.39) 
have now been shown to hold when $\rho$ satisfies (2.40). As already noted, 
the bounds (2.36) and (2.37) follow from these cases of (2.38) and (2.39);  
our proof of the lemma 
is therefore complete, 
for, by (2.36) and (2.37), the result (2.34) is obtained\ $\blacksquare$

\bigskip 

\goodbreak 
\noindent{\bf Remark.}\quad 
By logarithmic differentiation of the duplication formula 
[24, Section~12.15], one finds that 
$\Gamma'(1/2)/\Gamma(1/2)=\Gamma'(1)/\Gamma(1)-2\log 2=-\gamma -\log 4$. 
Hence the case $\delta =1/2$ of the result (2.33)-(2.34) of Lemma~2.3 
implies the conditional lower bound (1.2.18) (where it is implicitly assumed that 
$X\geq 4 e^{2\gamma}$, so that $|\log X|=\log X\gg 1$). 

\bigskip 
\bigskip 

\goodbreak\centerline{\bf \S 3. A Bound in respect of a Single Level: 
the Proof of Theorem~4}

\bigskip 

In this section we prove Theorem~4. Let $\varepsilon$, $q$, $N$ and $X$ 
satisfy the hypotheses of that theorem; let $\Gamma =\Gamma_0(q)\leq SL(2,{\frak O})$; 
let ${\frak a}$ be a cusp of $\Gamma$; and let $b_n\in{\Bbb C}$ for 
$n\in{\frak O}-\{ 0\}$. Then on the left-hand side of the relation (1.3.1) 
one has the sum 
$$\sum_{\scriptstyle V\atop\scriptstyle\nu_V>0}^{(\Gamma)} 
X^{\nu_V}
\Biggl|\sum_{N/2<|n|^2\leq N} b_n c_V^{\frak a}\left( n;\nu_V,0\right)\Biggr|^2
=\sigma_q^{\frak a}({\bf b},N;X)\qquad\ \hbox{(say).}$$
Let $\sigma_q^{\frak a}({\bf b},N;Y)$ be defined similarly for $Y>0$ 
(i.e. by substitution of $Y$ for $X$ in the above equation). 
Then, since $X^{\nu}=U^{\nu}Y^{\nu}$ when $U=X/Y$, and since 
the real function $\nu\mapsto U^{\nu}$ is increasing if $U\geq 1$, and 
decreasing if $1>U>0$, it follows from (1.2.20) that, for $Y>0$, 
$$\sigma_q^{\frak a}({\bf b},N;X)
\leq\left(\max\left\{ 1\,,\,{X\over Y}\right\}\right)^{\!\Theta(q)}
\sigma_q^{\frak a}({\bf b},N;Y)\;.\eqno(3.1)$$
By (3.1), the hypothesis $X\geq 1$, and the description of 
spectral parameters preceding (1.1.4), one has, in particular:  
$$\sigma_q^{\frak a}({\bf b},N;X)
\leq X^{\Theta(q)}\sigma_q^{\frak a}({\bf b},N;1)
\leq X^{\Theta(q)} E_0^{\frak a}(q,1,1;N,{\bf b})\;,$$
where $E_0^{\frak a}(q,P,K;N,{\bf b})$ is as defined by the equation~(1.2.7) of Theorem~2. 
Here the case 
$j=0$, $P=K=1$ of the result (1.2.9) of Theorem~2 may be applied 
with $\varepsilon /2$ substituted for $\varepsilon$: one thereby obtains the bound 
$$\sigma_q^{\frak a}({\bf b},N;X)
\ll X^{\Theta(q)}\left( 1
+O_{\varepsilon}\left( M_{\frak a} N^{1+(\varepsilon /2)}\right)\right)  
\left\|{\bf b}_N\right\|_2^2\;,$$ 
where $M_{\frak a}=|\mu({\frak a})|^2$; and where 
$\mu({\frak a})$ and $\left\|{\bf b}_N\right\|_2$ are as indicated 
by (1.2.10), (1.2.11) and Remark~3 (below Theorem~2). 
Moreover, since $N\geq 1$, and since the inequalities (1.2.22) of 
Theorem~3 imply that $\Theta(q)\geq 0$ and $1-\Theta(q)>1/2$, it 
therefore follows that, in cases where $M_{\frak a}N^{1+(\varepsilon/2)}\gg 1$, one has: 
$$\eqalign{
\sigma_q^{\frak a}({\bf b},N;X)
=O_{\varepsilon} \left( X^{\Theta(q)} M_{\frak a} N^{1+(\varepsilon /2)}
\left\|{\bf b}_N\right\|_2^2\right) 
 &\leq\left( X M_{\frak a} N\right)^{\Theta(q)} 
\left( O_{\varepsilon}\left( M_{\frak a} N^{1+\varepsilon}\right)\right)^{1-\Theta(q)} 
\left\|{\bf b}_N\right\|_2^2 < {}\cr 
 &<\left( 1+X M_{\frak a} N\right)^{\Theta(q)} 
\left( 1+O_{\varepsilon}\left( M_{\frak a} N^{1+\varepsilon}\right)\right)^{1-\Theta(q)} 
\left\|{\bf b}_N\right\|_2^2\;.}$$
This shows that, in cases where $M_{\frak a} N\gg 1$, the bound (1.3.1) holds:  
that completes proof of the theorem in those cases, so we 
may suppose, henceforth, that 
$$2^{10}\pi^2 M_{\frak a} N\leq 1\,.\eqno(3.2)$$

In order to complete this proof it will suffice to show that, subject to (3.2) 
(and our prior hypotheses), the bound (1.3.1) holds when  
$$X={1\over 16\pi^2 M_{\frak a} N}\;.\eqno(3.3)$$
For, if that particular case of the bound (1.3.1) holds, then 
$$\eqalign{
\sigma_q^{\frak a}\left( {\bf b},N;\left( 16\pi^2 M_{\frak a} N\right)^{-1}\right) 
 &\ll\left( 1+{1\over 16\pi^2}\right)^{\Theta(q)} 
\left( 1+O_{\varepsilon}\left( M_{\frak a}N^{1+\varepsilon}\right)\right)^{1-\Theta(q)}  
\left\|{\bf b}_N\right\|_2^2\log\left( 2+{1\over M_{\frak a} N}\right) \ll {}\cr 
 &\ll\left( 1+O_{\varepsilon}\left( M_{\frak a}N^{1+\varepsilon}\right)\right)^{1-\Theta(q)}  
\left\|{\bf b}_N\right\|_2^2\log\left( 2+{1\over M_{\frak a} N}\right) 
}$$
(by Theorem~3), and so, since the inequality (3.1) and Theorem~3 imply that one has 
$$\sigma_q^{\frak a}({\bf b},N;X)
\ll\left( 1+M_{\frak a} N X\right)^{\Theta(q)} 
\sigma_q^{\frak a}\left( {\bf b},N;\left( 16\pi^2 M_{\frak a} N\right)^{-1}\right) 
\qquad\qquad\hbox{($X\geq 1$),}$$
it follows that the bound (1.3.1) holds for all $X\geq 1$.  
Accordingly, we assume (3.3) for the remainder of this proof. 
Since we also assume (3.2), this implies that 
$$X\geq 64\;.\eqno(3.4)$$

Our aim is to show that, in the cases now being considered, the 
bound (1.3.1) results from the comparison of 
$\sigma_q^{\frak a}({\bf b},N;X)$ with the sum 
$$\tilde\sigma_q^{\frak a}( {\bf b},N;f) 
=\sum_{\scriptstyle V\atop\scriptstyle\nu_V>0}^{(\Gamma)} 
{\bf K}f\left(\nu_V ,0\right) 
\Biggl|\sum_{N/2<|n|^2\leq N} b_n c_V^{\frak a}\left( n;\nu_V,0\right)\Biggr|^2\;,\eqno(3.5)$$
where the ${\bf K}$-transform is as defined in Theorem~1, while 
the even `test-function' $f : {\Bbb C}^{*}\rightarrow{\Bbb C}$ is given by 
$$f(z)=\Phi\left( {\log\left( X^{1/2} |z|\right)\over\log 2}\right)\qquad\qquad\quad  
\hbox{($z\in{\Bbb C}^*$),}$$
with $\Phi : {\Bbb R}\rightarrow [0,\infty)$ satisfying:  
$$\Phi(t)=\cases{\exp\left( -{1\over 1-t^2}\right) &if $-1<t<1$; \cr 0 &otherwise.}$$

Given the fact that $\lim_{\delta\rightarrow 0+}\delta^{-k}\exp(-1/\delta)=0$, 
for all $k\in{\Bbb N}$, it may be seen by elementary calculus that the above 
function $\Phi$ is infinitely differentiable on ${\Bbb R}$. 
Hence the function $\varphi : (0,\infty)\rightarrow [0,\infty)$ given by 
$\varphi(r)=\Phi\bigl(\log\bigl( X^{1/2} r\bigr)/\log 2\bigr)\ $ ($r>0$), 
which has range $[0,1/e]$ and support $\bigl[ 2^{-1} X^{-1/2}\,,\,2 X^{-1/2}\bigr]$, 
is also infinitely differentiable, and (as may, for example, be deduced from the case 
$X=64$) satisfies:  
$$\varphi^{(j)}(r)\ll_j X^{j/2}\qquad\qquad
\hbox{($j\in{\Bbb N}\cup\{ 0\}\,$ and $\,r>0$).}\eqno(3.6)$$
It follows (see Remark~4, in Subsection~1.2)  
that $X$ and the function $f$, just defined, satisfy all 
of the relevant hypotheses of both Theorem~1 and the 
case $A=2$ of the Corollary to Theorems~1 and~2. 
Since (3.4) ensures that 
we have $A X^{-1/2}=2 X^{-1/2}\leq 1/c$ with $c=4>2 e^{\gamma}$, 
it is moreover the case that the conditions found sufficient to imply  
the lower bound (1.2.18) for ${\bf K}f(\nu,0)$ are satisfied. 

By (1.2.20) and Theorem~3, each term of the sum over subspaces $V$ 
in the equation~(3.5) certainly has $0<\nu_V<1/2$. 
We may therefore infer from the lower bound (1.2.18) that 
$$\tilde\sigma_q^{\frak a}({\bf b},N;f)
\gg\int_0^{\infty}\varphi(r)\,{{\rm d}r\over r}\,\min\{\log X\,,\,2\} 
\,\sigma_q^{\frak a}({\bf b},N;X)\;.$$
Here 
$$\int_0^{\infty}\varphi(r)\,{{\rm d}r\over r}
=(\log 2)\int_{-\infty}^{\infty}\Phi(t)\,{\rm d}t\gg 1$$
and, by (3.4), $\min\{\log X\,,\,2\}=2$; so it follows that we have  
$$\sigma_q^{\frak a}({\bf b},N;X) 
\ll \tilde\sigma_q^{\frak a}({\bf b},N;f)\;.\eqno(3.7)$$

To complete this proof we now 
deduce, from the Corollary to Theorems~1 and~2, an estimate for 
the sum $\tilde\sigma_q^{\frak a}({\bf b},N;f)$. 
By applying the result (1.2.15) of that corollary, with 
${\frak b}={\frak a}$, $M=N$, $a_m=b_m\ $ ($0\neq m\in{\frak O}$) and $A=2$, and 
with $\varepsilon /4$ substituted for $\varepsilon$, 
one obtains an equation with the sum of an $O$-term 
and the term $\pi\tilde\sigma_q^{\frak a}({\bf b},N;f)$ 
on one side, and a sum of Kloosterman sums on the other side. 
Each term of the latter sum of Kloosterman sums involves a factor  
$f\bigl( 2\pi\sqrt{mn}/c\bigr)=\varphi\bigl(\bigl| 2\pi\sqrt{mn}/c\bigr|\bigr)$, 
where, by (1.2.15) and the result (1.2.12) noted in remarks following Theorem~2, 
the constraints on the variables of summation $m$, $n$ and $c$ ensure that 
$$\left| {2\pi\sqrt{mn}\over c}\right|\leq {2\pi\sqrt{N}\over |1/\mu({\frak a})|}  
=2\pi\left( M_{\frak a} N\right)^{1/2}=2^{-1}X^{-1/2}$$
(the last equation following by (3.3)). 
Since $\varphi(r)=0$ for $x\leq 2^{-1} X^{-1/2}$, each term of the 
sum of Kloosterman sums equals zero; and so, in the cases being considered, 
the equation~(1.2.15) reduces to:  
$$0=\pi\tilde\sigma_q^{\frak a}({\bf b},N;f) 
+O_A\left( (\log X) Y\left( 1+O_{\varepsilon}\left( |\mu({\frak a})| 
N^{(1/2)+(\varepsilon /4)}\right)\right)^2 \left\| {\bf b}_N\right\|_2^2\right) ,$$
where $A=2$, and where, by (1.2.14) and (3.6), we have $Y\ll 1$. 
Therefore, and by (3.3) and (3.7), we obtain the bound 
$$\sigma_q^{\frak a}({\bf b},N;X) 
\ll\left( 1+O_{\varepsilon}\left( M_{\frak a} N^{1+(\varepsilon/2)}\right)\right) 
\left\|{\bf b}_N\right\|_2^2\log\left( {1\over M_{\frak a} N}\right) .\eqno(3.8)$$
Here, since Theorem~3 shows that $1\geq 1-\Theta(q)>1/2$, 
and since we have (by (3.4)) $X\geq 1$ and (by hypothesis) $\varepsilon >0$ and $N\geq 1$, it follows  
from (3.3) that one has  
$$O_{\varepsilon}\left( M_{\frak a} N^{1+(\varepsilon/2)}\right) 
<O_{\varepsilon}\left( {N^{\varepsilon/2}\over X}\right) 
\leq O_{\varepsilon}\left( {N^{(1-\Theta(q))\varepsilon}\over X^{1-\Theta(q)}}\right) 
\leq\left( O_{\varepsilon}\left({N^{\varepsilon}\over X}\right)\right)^{1-\Theta(q)} 
\ll\left( O_{\varepsilon}\left( M_{\frak a} N^{1+\varepsilon}\right)\right)^{1-\Theta(q)}\;.$$
Therefore, and since the real exponents $\Theta(q)$ and $1-\Theta(q)$ in the result 
(1.3.1) of Theorem~4 are non-negative, 
it follows that, subject to (3.2) and (3.3) holding, 
the bound in (3.8) implies that in (1.3.1). This, as was observed in the paragraph 
containing (3.3), is all that was needed to complete this proof\quad$\blacksquare$

\bigskip 
\bigskip 

\goodbreak\centerline{\bf \S 4. Averaging over the Level}

\bigskip 

In this section we are concerned with the sum $S_t(Q,X,N)$ defined in (1.3.2): 
the relevant arbitrary complex coefficients $a_n\ $ ($n\in{\frak O}$) 
are assumed to be fixed throughout. We first prove Theorem~7, and then apply it 
(together with Theorem~6) 
in our proof of Theorem~5. We omit the proof of Theorem~6, since  
the paragraph below (1.3.10) gives sufficient details of that (very simple) proof.

\bigskip 
\bigskip 

\goodbreak 
\centerline{\it The Proof of Theorem~7.} 

\bigskip 

Let $X$, $Q$, $N$, $t$, $\varepsilon$, $j$ and $Q^{*}$ satisfy the 
hypotheses of the theorem. If $1\leq X\ll 1$ then, by Theorem~6, Theorem~3, 
the definition (1.3.2) and 
the Cauchy-Schwarz inequality, one has 
$$S_t(Q,X,N)\leq X^{2/9} S_t(Q,1,N)\ll S_t(Q,1,N)
\leq\sum_{Q/2<|q|^2\leq Q}
2\left( E_0^{\infty}(q,1,1;N,{\bf b})+E_0^{\infty}(q,1,1;N/2,{\bf b})\right) ,$$
where $E_0^{\frak a}(q,P,K;N,{\bf b})$ is given by the equation~(1.2.7), in Theorem~2, 
with 
$$b_n=a_n |n|^{2it}\qquad\qquad\hbox{($0\neq n\in{\frak O}$).}\eqno(4.1)$$
By the results (1.2.9)-(1.2.11) 
of Theorem~2, and Remark~3 (below Theorem~2), it therefore follows that 
$$S_t(Q,X,N)\ll \sum_{Q/2<|q|^2\leq Q}
\left( 1+O_{\varepsilon}\left( {N^{1+\varepsilon}\over |q|^2}\right)\right)
\left\|{\bf a}_N\right\|_2^2 
\ll\left( Q+O_{\varepsilon}\left( N^{1+\varepsilon}\right)\right) 
\left\|{\bf a}_N\right\|_2^2\qquad\quad\hbox{if $\quad X\ll 1$.}\eqno(4.2)$$

Since (4.2) implies that (1.3.10) holds when $1\leq X<64$, 
we assume henceforth that $X\geq 64$ 
(as was the case, after (3.4), in our  proof of Theorem~4).  
We now define $g : {\Bbb C}^{*}\rightarrow [0,1/e]$ by setting 
$g(z)=\psi(|z|)$, for $0\neq z\in{\Bbb C}$, where, for $r>0$, 
one has $\psi(r)=\Phi\left(\log\bigl( Q^{-1/2}r\bigr)/\log 2\right)$, 
with $\Phi : {\Bbb R}\rightarrow [0,1/e]$ defined as in the proof of Theorem~4, 
below (3.5). The function $\psi : (0,\infty)\rightarrow [0,1/e]$ 
is infinitely differentiable; its support is the interval 
$\bigl[ Q^{1/2}/2\,,\,2Q^{1/2}\bigr]$; and it moreover satisfies 
$\psi(r)>\Phi(1/2)=e^{-4/3}$ for $(Q/2)^{1/2}<r<(2Q)^{1/2}$.  
It follows that, for $0\neq q\in{\frak O}$, one has 
$g(q)\gg 1$ if $Q/2<|q|^2<2Q$; and $g(q)\geq 0$ otherwise. 
Hence, and by the definition (1.3.2) of $S_t(Q,X,N)$, one has  
$$S_t(Q,X,N) 
\ll\sum_{0\neq q\in{\frak O}} g(q)
\sum_{\scriptstyle V\atop\scriptstyle\nu_V>0}^{\left(\Gamma_0(q)\right)} 
X^{\nu_V}
\Biggl|\sum_{N/4<|n|^2\leq N} b_n c_V^{\frak a}\left( n;\nu_V,0\right)\Biggr|^2\;,$$
where the coefficients 
$b_n$ ($0\neq n\in{\frak O}$) are given by (4.1). 
Since $g(q)\geq 0$ for $0\neq q\in{\frak O}$, and since $X\geq 64$, 
it follows from 
this last bound (similarly to how, in the proof of Theorem~4, we obtained (3.7))  
that 
$$S_t(Q,X,N)
\ll\sum_{0\neq q\in{\frak O}} g(q)
\sum_{\scriptstyle V\atop\scriptstyle\nu_V>0}^{\left(\Gamma_0(q)\right)} 
{\bf K}f\left(\nu_V ,0\right) 
\Biggl|\sum_{N/4<|n|^2\leq N} b_n c_V^{\frak a}\left( n;\nu_V,0\right)\Biggr|^2\;,\eqno(4.3)$$
with $f : {\Bbb C}^{*}\rightarrow [0,\infty)$ 
as defined, below (3.5), in the proof of Theorem~4; 
and with ${\bf K}f(\nu,p)$ as in Theorem~1. 

Our next step is to apply the Corollary to 
Theorems~1 and~2 to the inner sum (over spaces $V$) on the right-hand side of (4.3). 
That corollary does not apply directly to such sums; but 
one can obtain the required result from four distinct applications 
of the equation~(1.2.15), in which one substitutes, for
the pair $(M,N)$, the pairs $(N,N)$, $(N,N/2)$, $(N/2,N)$ and $(N/2,N/2)$, 
respectively. Then, by applying the result obtained, for ${\frak b}={\frak a}=\infty$ 
(so that $|\mu({\frak a})|=|\mu({\frak b})|=1/|q|$), 
and with $\varepsilon /2$ substituted for $\varepsilon$, one finds that, 
for $0\neq q\in{\frak O}$, 
$$\eqalignno{ 
\sum_{\scriptstyle V\atop\scriptstyle\nu_V>0}^{\left(\Gamma_0(q)\right)} 
{\bf K}f\left(\nu_V ,0\right) 
\Biggl|\sum_{{N\over 4}<|n|^2\leq N}&\!b_n c_V^{\frak a}\left( n;\nu_V,0\right)\Biggr|^2 = {} 
 &(4.4)\cr
 &=\sum_{{N\over 4}<|m|^2\leq N}\overline{b_m}\sum_{{N\over 4}<|n|^2\leq N}b_n 
\sum_{c\in{}^{\infty}{\cal C}^{\infty}}^{\left(\Gamma_0(q)\right)}  
\!{S_{\infty,\infty}(m,n;c)\over\pi |c|^2}\,f\!\left({2\pi\sqrt{mn}\over c}\right) + {}\cr
 &\quad\  +O_A\left( (\log X)Y\left( 1+O_{\varepsilon}\left( N^{1+\varepsilon}|q|^{-2}\right)\right) 
\left\|{\bf b}_N\right\|_2^2\right) ,}$$
where, since $f$ is the function defined in the proof of Theorem~4, 
we have $A=2$ and $Y\ll 1$ (just as is noted prior to (3.8)). 
Since the function $g$ has range $[0,1/e]$ and support 
$\big\{ z\in{\Bbb C} : Q/4\leq |z|^2\leq 4Q\bigr\}$, 
it follows by (4.4), combined with (1.3.4), (1.3.5), (4.1) and (4.3), 
that one has the upper bound 
$$\eqalign{
S_t(Q,X,N)
 &\ll\left|\sum_{0\neq q\in{\frak O}} g(q)
\qquad\sum\!\!\!\!\!\!\!\!\sum_{\!\!\!\!\!\!\!\!\!\!\!\!\!\!\!{N/4<|m|^2,|n|^2\leq N}}\overline{a_m}\,a_n \left|{m\over n}\right|^{-2it}
\sum_{0\neq\ell\in{\frak O}} 
\!{S(m,n;\ell q)\over\pi |\ell q|^2}\,f\!\left({2\pi\sqrt{mn}\over \ell q}\right)\right| + {}\cr 
 &\quad\ +(\log X)\left( Q+O_{\varepsilon}\left( N^{1+\varepsilon}\right)\right) 
\left\|{\bf a}_N\right\|_2^2\;,}$$
where $S(u,v;w)$ is the `simple Kloosterman sum' defined in (1.3.6). 
By our choice of $f$ and $g$, we have here that $f(w)g(z)\neq 0$ if and only if  
$1/4<X|w|^2,|z|^2/Q<4$. Therefore, in the above sum over $q$, $m$, $n$ and $\ell$ 
(where the constraints on $m$ and $n$ imply $N/4<|mn|\leq N$), the summand is 
zero whenever $|\ell|^2\not\in\bigl( 2^{-10}Q^{*}\,,\,Q^{*}\bigr)$, 
where, as in (1.3.9), $Q^{*}=64\pi^2XN/Q$. Consequently, one either has 
$$S_t(Q,X,N)\ll
(\log X)\left( Q+O_{\varepsilon}\left( N^{1+\varepsilon}\right)\right) 
\left\|{\bf a}_N\right\|_2^2\;,\eqno(4.5)$$ 
or else, for some $k\in\{ 0,1,2,\ldots\ ,9\}$, and $L=2^{-k}Q^{*}$, 
one has: 
$$S_t(Q,X,N)
\ll\left|\sum_{L/2<|\ell|^2\leq L} 
\,\qquad\sum\!\!\!\!\!\!\!\!\sum_{\!\!\!\!\!\!\!\!\!\!\!\!\!\!\!{N/4<|m|^2,|n|^2\leq N}}
\overline{a_m}\,a_n \left|{m\over n}\right|^{-2it}
\sum_{0\neq q\in{\frak O}} 
{S(m,n;\ell q)\over |\ell q|^2}\,f\!\left({2\pi\sqrt{mn}\over \ell q}\right) g(q)\right|\;.\eqno(4.6)$$ 

Since the bound (4.5) is satisfactory (i.e. it would imply the result (1.3.10)), let it 
now be supposed that $L\asymp Q^{*}$ is such that (4.6) holds. 
If it is shown that in this case one again obtains  (1.3.10),  
then by the conclusions reached in (4.5)-(4.6) (and after (4.2)) 
the proof of the theorem will be complete.

We seek to 
apply the Corollary to Theorems~1 and~2 to the sum over $m$, $n$ and $q$ in (4.6): 
for, by (1.3.4)-(1.3.5), the sums $S(m,n;\ell q)$ occurring in the inner summation 
in (4.6) are the 
generalised Kloosterman sums $S_{\infty,\infty}(m,n;c)$ associated with 
the Hecke congruence subgroup $\Gamma_0(\ell)\leq SL(2,{\frak O})$.
In order that (1.2.15) may be applied, we must first 
effect a replacement of the factors $f(2\pi\sqrt{mn}/\ell q)$ and $g(q)$, in (4.6), 
by single factor of the form $F(|2\pi\sqrt{mn}/(\ell q)|)$, where $F$ is a 
suitable complex-valued function.
We achieve this in the steps between (4.7) and (4.11) below. 

Recalling that $g(z)=\psi(|z|)$, for $0\neq z\in{\Bbb C}$, 
where the function $\psi : (0,\infty)\rightarrow [0,1/e]$ is both infinitely 
differentiable and of compact support, we have, by
Mellin's inversion formula [10, Appendix, Equation~(A.2)],     
$$g(z)=\psi(|z|)={1\over 2\pi i}\int\limits_{\sigma-i\infty}^{\sigma+i\infty} 
\Psi(s) |z|^{-s} {\rm d}s\qquad\hbox{($0\neq z\in{\Bbb C}$, $\sigma\in{\Bbb R}$),}\eqno(4.7)$$
where, for $s\in{\Bbb C}$,  
$$\Psi(s)=\int_0^{\infty} r^{s-1} \psi(r) {\rm d}r\;.\eqno(4.8)$$
Note that the integral in (4.7) is absolutely convergent in all relevant cases; indeed, by 
(4.8) and our particular choice of function $\psi$, one has, for $s\in{\Bbb C}$,  
$$\Psi(s)=\int\limits_0^{\infty} r^{s-1} \Phi\left({\log\left( Q^{-1/2}r\right)\over\log 2}\right) {\rm d}r
=(\log 2) Q^{s/2}\int\limits_{-\infty}^{\infty}2^{su}\Phi(u) {\rm d}u$$
where, by the definition of $\Phi$ (below (3.5)) and integration by parts, 
$$\int\limits_{-\infty}^{\infty}2^{su}\Phi(u) {\rm d}u
=\int\limits_{-1}^{1}2^{su}\Phi(u) {\rm d}u
=(-s\log 2)^{-j}\int_{-1}^1 2^{su}\Phi^{(j)}(u) {\rm d}u
\ll 2^{|{\rm Re}(s)|}\min\left\{ 1\,,\,O_j\left( |s|^{-j}\right)\right\}\qquad 
\hbox{($j\in{\Bbb N}$).}$$

We apply the case $\sigma =0$ of the inversion formula~(4.7), for $z=q$ (the variable of 
summation in (4.6)). 
Since the summation over $\ell$, $m$, $n$ and $q$ in (4.6) is 
(by virtue of the fact that 
$f(w)=0$ unless $X|w|^2>1/4$) effectively finite, our application of (4.7) 
allows us to deduce that  
$$S_t(Q,X,N)\ll\left|\,\int\limits_{-\infty}^{\infty}\!\Psi(i\tau) 
\!\sum_{{\textstyle{L\over 2}}<|\ell|^2\leq L} 
\,\qquad\sum\!\!\!\!\!\sum_{
\!\!\!\!\!\!\!\!\!\!\!\!\!\!\!{{\textstyle{N\over 4}}<|m|^2,|n|^2\leq N}}
\!\overline{a_m}\,a_n \left|{m\over n}\right|^{-2it}
\!\sum_{0\neq q\in{\frak O}} 
\!{S(m,n;\ell q)\over |\ell q|^2}\,f\!\left({2\pi\sqrt{mn}\over \ell q}\right) |q|^{-i\tau}
{\rm d}\tau\right|\;,$$
where, by the results of the preceding paragraph, 
$$\Psi(i\tau)\ll \bigl| Q^{i\tau/2}\bigr|\min\left\{ 1\,,\,O_j\left( |\tau|^{-j}\right)\right\}
\ll_j (1+|\tau|)^{-j}\qquad\qquad\hbox{($\tau\in{\Bbb R}$, $j=0,1,2,\ldots\ $).}\eqno(4.9)$$
It therefore follows by the substitution $\tau =4u$ that 
$$S_t(Q,X,N)\ll\int\limits_{-\infty}^{\infty} 
\left|\Psi(4iu)\!\!\sum_{{L\over 2}<|\ell|^2\leq L} |\ell|^{4iu}
\quad\sum\!\!\!\!\!\!\sum_{\!\!\!\!\!\!\!\!\!\!\!\!\!\!\!{{N\over 4}<|m|^2,|n|^2\leq N}}
\overline{a_m}\,|m|^{-2i(t+u)}
a_n |n|^{2i(t-u)} P_{m,n}\left(\ell; f_u\right)\right| {\rm d}u\;,\eqno(4.10)$$
for $j\geq 2$, where
$$f_u(z)=f(z)|z|^{4iu}=\varphi(|z|) |z|^{4iu}\eqno(4.11)$$
(with $\varphi : (0,\infty)\rightarrow[0,1/e]$ as defined 
just prior to (3.6)), 
and where, by (1.3.4) and (1.3.5),
$$P_{m,n}\left(\ell; f_u\right)
=\sum_{0\neq q\in{\frak O}} 
{S(m,n;\ell q)\over |\ell q|^2}\,f_u\left({2\pi\sqrt{mn}\over \ell q}\right)
=\sum_{c\in{}^{\infty}{\cal C}^{\infty}}^{\left(\Gamma_0(\ell)\right)} 
{S_{\infty,\infty}(m,n;c)\over |c|^2}\,f_u\left({2\pi\sqrt{mn}\over c}\right) .$$

Given arbitrary coefficients $a_n^{+},a_n^{-}\in{\Bbb C}$ ($0\neq n\in{\frak O}$), 
and with $f_u(z)$ and $P_{m,n}(\ell;f_u)$ as above, 
it follows by the Corollary to Theorems~1 and~2 
(similarly to how (4.4) was obtained) that, for $0\neq\ell\in{\frak O}$ and $u\in{\Bbb R}$, 
$$\eqalign{
\sum_{N/4<|m|^2,|n|^2\leq N}\overline{a_n^{+}}\,a_n^{-}
 &\,P_{m,n}\left(\ell; f_u\right) = {}\cr 
 &\quad =\pi\!\sum_{\scriptstyle V\atop\scriptstyle \nu_V>0}^{(\Gamma_0(\ell))}
\!{\bf K}f_u\left(\nu_V , 0\right)
\overline{\sum_{N/4<|m|^2\leq N}^{\hbox{\ }}  a_n^{+} c_V^{\infty}\left( m;\nu_V,0\right)}  
\sum_{N/4<|n|^2\leq N} a_n^{-} c_V^{\infty}\left( n;\nu_V,0\right) + {}\cr 
 &\qquad\  +O_{A}\left( (\log X)\,Y_u
\left( 1+O_{\varepsilon}\left( N^{1+\varepsilon}|\ell|^{-2}\right)\right) 
\left\| {\bf a}_N^{+}\right\|_2 
\left\| {\bf a}_N^{-}\right\|_2
\right) ,}$$
with $A=2$ as previously ($f_u$ having the same support as $f$), and, by 
(1.2.14), (4.11) and (3.6), with 
$$\eqalign{
Y_u &=X^{-3/2}\max_{r>0}\left|{{\rm d}^3\over{\rm d}r^3}\,\varphi(r)r^{4iu}\right| \ll {}\cr 
 &\ll X^{-3/2}\max_{r>(4X)^{-1/2}}\ \max_{k=0,1,2,3} X^{(3-k)/2}
\left|{{\rm d}^k\over{\rm d}r^k}\,r^{4iu}\right| 
\ll\max_{r>0}\max_{k=0,1,2,3}\left|r^k{{\rm d}^k\over{\rm d}r^k}\,r^{4iu}\right|
\asymp (1+|u|)^3
}$$
(since $\varphi(r)=0$ for $0<r\leq 2^{-1}X^{-1/2}$). 
Since one may take here
$a_n^{\pm}=a_n |n|^{2i(t\pm u)}$, for $0\neq n\in{\frak O}$ (when $u,t$ are any 
given real numbers), we are therefore able to deduce from (4.10) that 
$$S_t(Q,X,N)
\ll\int\limits_{-\infty}^{\infty}|\Psi(4iu)|  
\Biggl( (1+|u|)^3 E+\sum_{L/2<|\ell|^2\leq L} F_{\ell}(u)\Biggr) {\rm d}u\;,\eqno(4.12)$$
where
$$E\ll (\log X)\left( L+O_{\varepsilon}\left( N^{1+\varepsilon}\right)\right) 
\left\|{\bf a}_N\right\|_2^2
\ll_{\varepsilon} X^{\varepsilon}\left( L+ N^{1+\varepsilon}\right) 
\left\|{\bf a}_N\right\|_2^2\eqno(4.13)$$
and 
$$\eqalignno{
F_{\ell}(u) &=\sum_{\scriptstyle V\atop\scriptstyle \nu_V>0}^{(\Gamma_0(\ell))}
\left|{\bf K}f_u\left(\nu_V , 0\right)
\overline{\sum_{N/4<|m|^2\leq N}^{\hbox{\ }}  a_m |m|^{2i(t+u)} c_V^{\infty}\left( m;\nu_V,0\right)}  
\sum_{N/4<|n|^2\leq N} a_n |n|^{2i(t-u)} c_V^{\infty}\left( n;\nu_V,0\right)\right| \leq {}\cr 
 &\leq {1\over 2}\sum_{\sigma =\pm 1} 
\sum_{\scriptstyle V\atop\scriptstyle \nu_V>0}^{(\Gamma_0(\ell))}
\left|{\bf K}f_u\left(\nu_V , 0\right)\right|  
\,\left|\sum_{N/4<|n|^2\leq N} a_n |n|^{2i(t+\sigma u)} 
c_V^{\infty}\left( n;\nu_V,0\right)\right|^2&(4.14)}$$
(by the arithmetic-geometric mean inequality).

For $u\in{\Bbb R}$, let $\varphi_u : (0,\infty)\rightarrow{\Bbb C}$ 
be the function satisfying $\varphi_u(r) =\varphi(r) r^{4iu}$ for $r>0$, 
where $\varphi$ (as in (4.11)) is the function defined just above (3.6). 
The functions $\varphi_u$ just defined inherit from 
$\varphi$ both the property of having support 
$\bigl[ 2^{-1} X^{-1/2}\,,\,2 X^{-1/2}\bigr]$ 
and the property of being infinitely differentiable; they are, in particular, 
continuous functions on $(0,\infty)$. Given (4.11), and given that $X\geq 64\geq 2$, 
it therefore follows by (1.2.17) that, for $u\in{\Bbb R}$ and $0<\nu\leq 1/2$, 
$${\bf K}f_u(\nu,0)\ll \int\limits_0^{\infty}\left|\varphi_u(r)\right|{{\rm d}r\over r} 
\min\left\{\log X\,,\,\nu^{-1}\right\} X^{\nu}\;,$$
where 
$$\int\limits_0^{\infty}\left|\varphi_u(r)\right|{{\rm d}r\over r}
=\int\limits_0^{\infty}\varphi(r){{\rm d}r\over r}
=\int\limits_0^{\infty}
\Phi\left({\log\left( X^{1/2}r\right)\over\log 2}\right) {{\rm d}r\over r}
=(\log 2)\int\limits_{-1}^{1}\Phi(y){\rm d}y\ll 1$$
(by the definitions of $\varphi$ and $\Phi$ prior to (3.6)).
Since $X\geq 1$, and since one consequently has 
$\log X<(2/\varepsilon)X^{\varepsilon /2}$, one may 
deduce from the above results that, for $u\in{\Bbb R}$,  
$${\bf K} f_u(\nu ,0)\ll\cases{\varepsilon^{-1}X^{\nu}  
&if $\varepsilon /2<\nu\leq 1/2$; \cr
\varepsilon^{-1} X^{\varepsilon}  
&if $0<\nu\leq\min\{\varepsilon/2\,,\,1/2\}$.}$$ 
Therefore it follows by (4.14), (1.2.20) and Theorem~3, and the case $K=P=1$ 
of Theorem~2 that, in (4.12),
$$F_{\ell}(u)\ll_{\varepsilon}\sum_{\sigma =\pm 1} 
\sum_{\scriptstyle V\atop\scriptstyle \nu_V>\varepsilon /2}^{(\Gamma_0(\ell))}
X^{\nu}\,\left|\sum_{N/4<|n|^2\leq N} a_n |n|^{2i(t+\sigma u)} 
c_V^{\infty}\left( n;\nu_V,0\right)\right|^2
+X^{\varepsilon}\left( 1+{N^{1+\varepsilon}\over |\ell|^{2}}\right)\left\|{\bf a}_N\right\|_2^2\eqno(4.15)$$
(as the variable of summation $n$ 
is subject to stricter conditions in (1.2.7) than it is in (4.14),
the application here of Theorem~2 depends on prior use of the inequality 
$|S(N)+S(N/2)|^2\leq 2|S(N)|^2+2|S(N/2)|^2$, where $S(M)$
denotes the sum of those terms of the sum over $n$ in (4.14) 
for which $|n|^2$ lies in the interval $(M/2,M]$). 
By (4.12), (4.13) and (4.15), and the bounds on $\Psi(i\tau)$ in (4.9), 
we have 
$$\eqalign{
S_t(Q,X,N)
 &\ll_{\varepsilon}\ \int\limits_{-\infty}^{\infty} |\Psi(4iu)| 
\left( (1+|u|)^3 X^{\varepsilon}\left( L+N^{1+\varepsilon}\right)\left\|{\bf a}_N\right\|_2^2 
+\sum_{\sigma =\pm 1} S_{t+\sigma u}(L,X,N)\right) {\rm d}u = {}\cr 
 &=X^{\varepsilon}\left( L+N^{1+\varepsilon}\right)\left\|{\bf a}_N\right\|_2^2
\!\int\limits_{-\infty}^{\infty}\!\!O\!\left( (1+|u|)^{-2} \right) 
{\rm d}u 
+\int\limits_{-\infty}^{\infty}\!\!\left(|\Psi(4iu)|+|\Psi(-4iu)|\right) 
S_{t+u}(L,X,N)\,{\rm d}u = {}\cr 
 &=O\left( X^{\varepsilon}\left( L+N^{1+\varepsilon}\right)\left\|{\bf a}_N\right\|_2^2\right) 
+O_j\left(\,\int\limits_{-\infty}^{\infty} (1+|u|)^{-j}\,S_{t+u}(L,X,N)\,{\rm d}u\right)}$$
(both $S_t(Q,X,N)$ and $S_{t+\sigma u}(L,X,N)$ being 
examples of the sums defined in (1.3.2)). 
Since we have $L\asymp Q^{*}$, with $Q^{*}$ as in (1.3.9), 
this upper bound just obtained for 
$S_t(Q,X,N)$ therefore implies the result (1.3.10) of Theorem 7; 
as explained below (4.6), this completes the proof of that theorem 
\ $\blacksquare$

\bigskip 
\smallskip 

We use the remainder of this section to prepare for and 
present our proof of Theorem~5. 
Considering firstly the sum 
$$\sigma_q^{\infty}({\bf b},N;X)
=\sum_{\scriptstyle V\atop\scriptstyle\nu_V>0}^{(\Gamma_0(q))} 
X^{\nu_V}\left|\sum_{N/4<|n|^2\leq N} b_n c_V^{\infty}\left( n;\nu_V,0\right)\right|^2\;,\eqno(4.16)$$ 
it follows by Theorem~4 
(in conjunction with an inequality similar to that mentioned, in parenthesis, below (4.15)),  
and by the latter part of Remark~3 (below Theorem~2), 
that for $0\neq q\in{\frak O}$, $X,N\geq 1$, arbitrary complex 
coefficients $b_n\ $ ($0\neq n\in{\frak O}$) and any $\varepsilon >0$  
this sum satisfies 
$$\sigma_q^{\infty}({\bf b},N;X)
\ll \left( 1+{XN\over |q|^2}\right)^{\Theta(q)} 
\left( 1 +O_{\varepsilon}\left({N^{1+\varepsilon}\over |q|^2}\right)\right)^{1-\Theta(q)} 
\left\| {\bf b}_N\right\|_2^2\,\log\!\left( 2+{|q|^2\over N}\right) ,$$
where $\left\| {\bf b}_N\right\|_2$ and $\Theta(q)$ 
are as indicated by (1.2.11) and (1.2.20). 
By Theorem~3, and given the condition $\varepsilon >0$, and the conditions $X,N\geq 1$ 
(which ensure that 
$(1+XN|q|^{-2})/(1+N|q|^{-2})\geq 1$), 
it is implied by the above bound on $\sigma_q^{\infty}({\bf b},N;X)$ that one has 
$$\sigma_q^{\infty}({\bf b},N;X)
\ll_{\varepsilon} 
(QN)^{\varepsilon}\left( 1+{XN\over Q}\right)^{\vartheta} 
\left( 1+{N\over Q}\right)^{1-\vartheta}
\left\|{\bf b}_N\right\|_2^2\qquad\quad\hbox{if $\quad Q/2<|q|^2\leq Q$}$$
(with $\vartheta\in[0,2/9]$ given by (1.2.20) and (1.2.21)).
By summing this over the relevant 
$q\in{\frak O}$, one obtains:
$$\eqalignno{
\sum_{Q/2<|q|^2\leq Q}
\sigma_q^{\infty}({\bf b},N;X)
 &\ll_{\varepsilon} 
Q(QN)^{\varepsilon}\left( 1+{XN\over Q}\right)^{\vartheta} 
\left( 1+{N\over Q}\right)^{1-\vartheta}
\left\|{\bf b}_N\right\|_2^2 = {}\cr
 &=(QN)^{\varepsilon}\left( Q+XN\right)^{\vartheta} 
\left( Q+N\right)^{1-\vartheta}
\left\|{\bf b}_N\right\|_2^2\;, 
&(4.17)}$$
for $Q\geq 1$. Given $t\in{\Bbb R}$, and arbitrary coefficients 
$a_n\in{\Bbb C}$ ($0\neq n\in{\frak O}$), one 
may apply (4.17) with the particular coefficients 
$b_n=a_n |n|^{2it}$ ($0\neq n\in{\frak O}$). 
Then, by the definitions in (1.3.2) and (4.16), 
the sum bounded in (4.17) equals $S_t(Q,X,N)$; 
and so, given the equation $\left\|{\bf b}_N\right\|_2^2=\left\|{\bf a}_N\right\|_2^2$, 
and the inequalities 
$$\left( Q+XN\right)^{\vartheta} 
\left( Q+N\right)^{1-\vartheta}
\leq 2\left( Q^{\vartheta}+X^{\vartheta}N^{\vartheta}\right) 
\left( Q^{1-\vartheta}+N^{1-\vartheta}\right)\quad\ \hbox{and}\quad\   
Q^{\vartheta}N^{1-\vartheta}\leq\max\{ Q,N\}\leq Q+X^{\vartheta}N\;,$$
this application of (4.17) yields the bound 
$$S_t(Q,X,N)\ll_{\varepsilon} 
(QN)^{\varepsilon}\left( Q+X^{\vartheta}N^{\vartheta}Q^{1-\vartheta}+X^{\vartheta}N\right) 
\left\|{\bf a}_N\right\|_2^2\;,\eqno(4.18)$$
for $t\in{\Bbb R}$, $Q,X,N\geq 1$ and $\varepsilon >0$.  
Note that the term $X^{\vartheta}N^{\vartheta}Q^{1-\vartheta}$, in brackets, on 
the right-hand side of (4.18), is greater by a factor $Q^{\vartheta}$ than the 
corresponding term occurring in the result (1.3.7) of Theorem~5.  
From the bound (4.18) we deduce the first of the next two lemmas. The second of 
these lemmas is proved by a straightforward application of Theorems~6 and~7. 

\bigskip

\goodbreak\proclaim Lemma~4.1. Let $1/2\geq\varepsilon>0$; let $Q_0\geq 1$; 
and let $\vartheta$ be given by (1.2.20) and (1.2.21).   
Then there exists a number $C_1=C_1\!\left(\varepsilon , Q_0\right)\in[1,\infty)$ 
(depending only upon $\varepsilon$ and $Q_0$)  such that, 
for all $t\in{\Bbb R}$, and for  
all $Q,N,X\geq 1$ satisfying 
$Q^{1-\varepsilon}\leq N$, or $Q^{1+\varepsilon}\geq XN$, or $Q\leq Q_0$, 
one has 
$$S_t(Q,X,N)\leq C_1\!\left(\varepsilon,Q_0\right)(QN)^{\varepsilon}\left( Q+X^{\vartheta}N\right)
\left\|{\bf a}_N\right\|_2^2\;.\eqno(4.19)$$

\medskip

\goodbreak 
\noindent{\bf Proof.}\quad 
Assume the hypotheses of the lemma; and let $t\in{\Bbb R}$. 
By applying (4.18), with $\varepsilon /3$ substituted for~$\varepsilon$, 
one obtains the bound 
$$S_t(Q,X,N)\ll_{\varepsilon} (QN)^{\varepsilon}\left( Q
+X^{\vartheta}N\right)
\left\|{\bf a}_N\right\|_2^2\;,\eqno(4.20)$$
whenever $Q,N,X\geq 1$ are such that 
$X^{\vartheta}N^{\vartheta}Q^{1-\vartheta}\leq 
(QN)^{2\varepsilon /3}\max\left\{ Q\,,\,X^{\vartheta}N\right\}$. 
This condition holds if and only if one has either 
$(XN/Q)^{\vartheta}\leq (QN)^{2\varepsilon /3}$ 
or $(Q/N)^{1-\vartheta}\leq (QN)^{2\varepsilon /3}$. 
Here we may assume $Q,N,X\geq 1$; 
so by Theorem~3, the former of the two inequalities involving $(QN)^{2\varepsilon /3}$
will hold if $(XN/Q)^{2/9}\leq Q^{2\varepsilon /3}$; while, since $-\vartheta\leq 0$ 
and $\varepsilon\leq 1/2$, the 
latter inequality (which is equivalent to having  
$(Q/N)^{1+(2\varepsilon /3)-\vartheta}\leq Q^{4\varepsilon /3}$)  
will hold if $(Q/N)^{4/3}\leq Q^{4\varepsilon /3}$. 
Hence one obtains (4.20) for $Q,N,X\geq 1$ satisfying either $XN/Q\leq Q^{3\varepsilon}$ 
or $Q/N\leq Q^{\varepsilon}$ (so certainly in all the cases where 
$Q,X,N\geq 1$ and either $Q^{1+\varepsilon}\geq XN$ or $Q^{1-\varepsilon}\leq N$). 

When $1\leq Q\leq Q_0$, we may, instead of the above, simply apply 
(4.18) as it stands (i.e. with $\varepsilon$ there as it is here). 
Indeed, since $\vartheta<1$ and $1-\vartheta\leq 1$, it follows by 
(4.18) that for  $Q,X,N\geq 1$ one has 
$$S_t(Q,X,N)\ll_{\varepsilon} 
(QN)^{\varepsilon}\left( Q+X^{\vartheta}NQ+X^{\vartheta}N\right) 
\left\|{\bf a}_N\right\|_2^2
\ll Q (QN)^{\varepsilon}\left( Q+X^{\vartheta}N\right) 
\left\|{\bf a}_N\right\|_2^2\;,$$
so that   
$$S_t(Q,X,N)\ll_{\varepsilon , Q_0} (QN)^{\varepsilon}\left( Q+X^{\vartheta}N\right) 
\left\|{\bf a}_N\right\|_2^2\qquad\qquad\hbox{if $\quad Q\leq Q_0$.}\eqno(4.21)$$ 
Given both (4.21) and the conclusion of the preceding paragraph, it 
has now been shown that 
if, for example, 
$C_1(\varepsilon, Q_0)=\max\left\{ 1\,,\,A(\varepsilon)\,,\,B(\varepsilon,Q_0)\right\}$, 
where $A(\varepsilon)$ and $B(\varepsilon,Q_0)$ are any of the 
positive numbers that may serve as implicit constants 
in (4.20) and (4.21), respectively, 
then the bound (4.19) will hold in all the cases 
referred to by the lemma\quad$\blacksquare$ 

\bigskip

\goodbreak
\proclaim Lemma~4.2. Let $j\geq 2$; let $1/2\geq\varepsilon >0$; and let 
$\vartheta$ be given by (1.2.20) and (1.2.21). Then there exists 
a number $C_j=C_j(\varepsilon)\in[1,\infty)$ (depending only upon $\varepsilon$ and $j$) 
which is such that, for all $t\in{\Bbb R}$, and all $Q,X,N\geq 1$ and $Y\in{\Bbb R}$ 
satisfying both 
$$X\geq{Q\over N}\geq 1\;,\eqno(4.22)$$
and 
$$\min\left\{ X\,,\,{Q^{2-\varepsilon}\over N}\right\} =Y\;,\eqno(4.23)$$
there exist $L,v\in{\Bbb R}$ satisfying both 
$${YN\over Q}<L<{2^{10}YN\over Q}\;,\eqno(4.24)$$
and 
$$C_j(\varepsilon)\left( {X\over Y}\right)^{\!\!\vartheta} 
\!\left( (1+|v-t|)^{-(j-2)} S_{v}(L,Y,N) + Q^{1+(2-\varepsilon)\varepsilon /3}
\left\|{\bf a}_N\right\|_2^2\right)\geq S_t(Q,X,N)\;.\eqno(4.25)$$

\medskip

\goodbreak 
\noindent{\bf Proof.}\quad 
Let $j$, $\varepsilon$ and $\vartheta$ satisfy the stated hypotheses. 
Suppose, moreover, that $t\in{\Bbb R}$; that the parameters $Q,X,N\geq 1$ 
satisfy (4.22); and that $Y$ is given by the equation~(4.23). 
Then, since $Q\geq 1$ and $2-\varepsilon\geq 3/2>1$, 
it follows by (4.23) and (4.22) that $X\geq Y\geq Q/N\geq 1$. 
Hence (and since $\vartheta\geq 0$), 
the result (1.3.8) of Theorem~6 implies that 
$$S_t(Q,X,N)\leq\left({X\over Y}\right)^{\vartheta}S_t(Q,Y,N)\;.\eqno(4.26)$$

By substitution of $\varepsilon /3$ for $\varepsilon$ 
in the result (1.3.10) supplied by Theorem~7, one has, in (4.26),  
$$S_t(Q,Y,N)
\ll_{\varepsilon,j}\int\limits_{-\infty}^{\infty}S_{t+u}(L,Y,N) (1+|u|)^{-(j-2)} 
\,{{\rm d}u\over\left( 1+u^2\right)} 
+Y^{\varepsilon /3}\left( Q+{YN\over Q}+N^{1+\varepsilon /3}\right) 
\left\|{\bf a}_N\right\|_2^2\;,\eqno(4.27)$$
for some $L\in\left[ 2^{-3}\pi^2 YN/Q\,,\,2^6\pi^2 YN/Q\right]$. 
Choose such an $L$: this ensures (since $8<\pi^2<16$) that the inequalities in (4.24) 
are satisfied. Then, in considering the definition of the sum 
$S_{t+u}(L,Y,N)$, we may note that, on the right-hand side of the 
defining equation~(1.3.2), the inner 
sum (over spaces $V$) always has 
a finite number of terms: see 
the discussion of exceptional eigenvalues below (1.1.11). 
Consequently all the summations in (1.3.2) are finite; and so 
it is evident (from an inspection of (1.3.2)) that
the function $u\mapsto (1+|u|)^{-(j-2)} S_{t+u}(L,Y,N)$ 
has its range contained in $[0,\infty)$, 
and is both bounded and continuous on~${\Bbb R}$. 
Given the well-known evaluation of $\int_{\Bbb R} (1+u^2)^{-1}{\rm d}u$, 
it may therefore be deduced that there exists some $v\in{\Bbb R}$ 
such that the integral appearing in (4.27) is equal to 
$\pi (1+|v-t|)^{-(j-2)} S_{v}(L,Y,N)$. Hence, and by (4.26), it suffices for completion of 
this proof that  
we note (with regard to the rightmost terms in (4.27)) the three inequalities  
$$YN/Q\leq Q^{1-\varepsilon}\leq Q,\qquad 
\max\bigl\{ Q\,,\,N^{1+\varepsilon /3}\bigr\}\leq Q N^{\varepsilon /3}\qquad
\hbox{and}\qquad 
Y^{\varepsilon /3} N^{\varepsilon /3}=(YN)^{\varepsilon /3}
\leq Q^{(2-\varepsilon)\varepsilon /3}\;,$$
which (given that $\varepsilon >0$ and $Y,Q,N\geq 1$)
are implied by (4.23), (4.22) and (4.23), respectively\quad$\blacksquare$ 

\bigskip 
\bigskip  

\goodbreak 
\centerline{\it The Proof of Theorem~5.} 

\bigskip 

By (1.3.2) one has $S_t(Q,X,N)=0$ whenever either $Q$ or $N$ is less than $1$;  
so it will suffice to prove Theorem~5 in cases where $Q,N\geq 1$. 
Given the nature of the result (1.3.7), it may therefore also henceforth be 
supposed that $1/2\geq\varepsilon >0$ (the result for $\varepsilon =1/2$ 
implying the result for all $\varepsilon >1/2$, when $X,Q,N\geq 1$).
Taking now $C_2(\varepsilon)\in[1,\infty)$ to be one of those numbers 
shown to exist by the case $j=2$ of Lemma~4.2, we put 
$$Q_0=Q_0(\varepsilon)=\left( 2^{17} C_2(\varepsilon)\right)^{(3\varepsilon^{-2})}\qquad 
\hbox{and}\qquad 
C_0(\varepsilon)=C_1\!\left(\varepsilon\,,\,Q_0\right)\;,\eqno(4.28)$$
where $C_1\left(\varepsilon\,,\,Q_0\right)$ is any one of those numbers 
whose existence is established in Lemma~4.1 (note that we certainly have  
$Q_0\geq 1$ here). The values of $\varepsilon$, $C_2(\varepsilon)$, 
$Q_0(\varepsilon)$ and $C_0(\varepsilon)$ are to remain fixed 
throughout this proof. We assume also a fixed choice of 
the arbitrary complex coefficients $a_n\ $ ($0\neq n\in{\frak O}$) 
appearing in the definition (1.3.2) of the sums $S_t(Q,N,X)$.

For each $Q\in[1,\infty)$, let $A(Q)$ denote the proposition that, 
for all $t\in{\Bbb R}$ and all $X,N\geq 1$, one has 
$$S_t(Q,X,N)\leq C_0(\varepsilon)(QN)^{\varepsilon} 
\left( Q+X^{\vartheta}N+Q^{1-2\vartheta}(XN)^{\vartheta}\right)
\left\|{\bf a}_N\right\|_2^2\;,\eqno(4.29)$$
with $\vartheta$ being the constant defined in (1.2.20) and (1.2.21). 
In what follows `$A(Q)$' may be used as shorthand for either `Proposition~$A(Q)$' or
`the truth of Proposition~$A(Q)$' (which is the case should be clear from the context).
 
Since the bound (4.29) implies the result (1.3.7) of the theorem, 
we have only to show that $A(Q)$ is true for all $Q\in[1,\infty)\,$ (that will 
prove the theorem). Since $0\leq\vartheta\leq 2/9\leq 1/2$ (by Theorem~3), the 
right-hand side of (4.29) is an increasing function of $Q$. 
Therefore, and since the definition (1.3.2) implies that 
$$S_t(P,X,N)=S_t(Q,X,N)\qquad\ \hbox{if $\quad Q+1>P\geq Q\in{\Bbb N}$,}\eqno(4.30)$$
we have 
$$A(Q)\ \,{\rm implies}\ \,A(P)\qquad\ {\rm if}\quad Q+1>P\geq Q\in{\Bbb N}\;;\eqno(4.31)$$
and so may in fact complete this proof simply by showing that  
$A(Q)$ is true for all $Q\in{\Bbb N}$ (with it then following by (4.31) 
that $A(P)$ is true for all $P\in[1,\infty)$). 
The method of proof by contradiction is suited to this task. 

Suppose that 
$A(Q)$ is false for some $Q\in{\Bbb N}$. Then the set 
${\cal F}=\{ Q\in{\Bbb N} : A(Q)\ {\rm is\ false}\}$ 
is a non-empty subset of ${\Bbb N}$, and so contains a unique least element, $R=\min{\cal F}\in{\cal F}\subseteq{\Bbb N}$. 
Since $R\in{\cal F}$, we have: 
$$A(R)\ \,\hbox{is false.}\eqno(4.32)$$
On the other hand, for $Q=1,2,\ldots ,R-1$, the hypothesis $A(Q)$ must be true:  
else one would have $Q\in{\cal F}$, and so $Q\geq\min{\cal F}=R$, 
which is impossible when $Q<R$. Upon combining this with (4.31), we deduce:
$$A(Q)\ \,\hbox{is true}\qquad\ \hbox{for all}\quad Q\in[1,R)\;.\eqno(4.33)$$ 

Given the definition of $C_0(\varepsilon)$ 
in (4.28), it follows by Lemma~4.1 that $A(Q)$ is true for 
all $Q\in\bigl[1,Q_0(\varepsilon)\bigr]$ 
(the bound in (4.19) implying that in (4.29)). By this result and 
(4.32), it must be the case that 
$$R>Q_0(\varepsilon)\;.\eqno(4.34)$$
We aim to deduce from (4.28), (4.33) and (4.34) that 
$A(R)$ is true.  Such a deduction would directly contradict (4.32);  
would thereby establish the falsity of the premise that $A(Q)$ is false 
for some $Q\in{\Bbb N}$; and so would prove that 
$A(Q)$ is true for all $Q\in{\Bbb N}$. 
To achieve this we must show that, for $Q=R$, all $t\in{\Bbb R}$, and 
all $X,N\geq 1$, the inequality (4.29) holds. By Lemma~4.1, 
the inequality (4.29) does hold if it is the case that 
$Q=R$, $t\in{\Bbb R}$, $X,N\geq 1$ 
and either $R^{1-\varepsilon}\leq N$, or $R^{1+\varepsilon}\geq XN$. 
In all remaining cases (that are relevant to our purpose) one has  
$$R^{-\varepsilon}X>{R\over N}>R^{\varepsilon}\;.\eqno(4.35)$$
Therefore we may establish that $A(R)$ is true by showing that 
(4.29) holds if 
$Q=R$, $t\in{\Bbb R}$ and $X,N\geq 1$ are such that (4.35) holds. 
Accordingly, we assume henceforth  
that $t\in{\Bbb R}$ and $X,N\geq 1$; and that 
(4.35) holds. 

Since $\varepsilon >0$ and $R\geq 1$, the inequalities in (4.35) 
imply the case $Q=R$ of the condition~(4.22) of Lemma~4.2. 
Therefore, given our choice of $C_2(\varepsilon)$ (prior to (4.28)), 
it follows by the case $j=2$ of Lemma~4.2 that we have 
$$S_t(R,X,N)\leq C_2(\varepsilon)\left({X\over Y}\right)^{\vartheta} 
\left( S_v(L,Y,N)+R^{1+(2-\varepsilon)\varepsilon /3}\left\|{\bf a}_N\right\|_2^2\right) 
,\eqno(4.36)$$
with $Y=\min\left\{ X\,,\,R^{2-\varepsilon}/N\right\}$, 
for some $v\in{\Bbb R}$, and some $L$ satisfying $YN/R<L<2^{10}YN/R$. 
Here we have, by (4.35), 
$$Y\geq\min\{ X,R\}\geq 1$$ 
and (since $0<\varepsilon\leq 1/2$)
$$L>{YN\over R}=\min\left\{ {XN\over R}\,,\,R^{1-\varepsilon}\right\}
\geq\min\left\{ R^{\varepsilon}\,,\,R^{1-\varepsilon}\right\}\geq 1\;.$$ 
Moreover, since 
$$L<{2^{10}YN\over R}=2^{10}\min\left\{ {XN\over R}\,,\,R^{1-\varepsilon}\right\}\;,\eqno(4.37)$$
and since (given that $0<\varepsilon\leq 1/2<3$ and $C_2(\varepsilon)\geq 1$) 
it is implied by (4.34) and (4.28) that
$$R^{\varepsilon}>\left( 2^{17} C_2(\varepsilon)\right)^{3/\varepsilon}>2^{17}>2^{10}\;,\eqno(4.38)$$
it must consequently be the case that we have here: 
$$L<2^{10}R^{1-\varepsilon}<R\;.\eqno(4.39)$$

By the above, we have $1\leq L<R$; so that it follows by (4.33) that 
the proposition~$A(L)$ is true. 
Therefore, and since $v\in{\Bbb R}$, $N\geq 1$ and (as just shown) 
$Y\geq 1$, the inequality (4.29) holds when one replaces  
$t$, $Q$ and $X$ there by $v$, $L$ and $Y$, respectively. 
This means that we have: 
$$S_v(L,Y,N)
\leq C_0(\varepsilon)(LN)^{\varepsilon} 
\left( L+Y^{\vartheta}N+L^{1-2\vartheta}(YN)^{\vartheta}\right)
\left\|{\bf a}_N\right\|_2^2\;,$$
where, since $1\geq 1-\vartheta\geq 1-2\vartheta\geq 5/9>1/3>0$, 
it follows by (4.37) that 
$$2^{-10}L^{1-2\vartheta}(YN)^{\vartheta}
<\left( {YN\over R}\right)^{1-\vartheta} R^{\vartheta}
\leq R^{(1-\varepsilon)(1-\vartheta)+\vartheta}
=R^{1-(1-\vartheta)\varepsilon}
\leq R^{1-\varepsilon /3}\;.$$
This, with the  first part of (4.39) (and the hypothesis  
that $0<\varepsilon\leq 1/2$), allows us to conclude that 
$$S_v(L,Y,N)
\leq 2^{16}C_0(\varepsilon)\left( R^{1-\varepsilon}N\right)^{\varepsilon} 
\left( R^{1-\varepsilon /3}+Y^{\vartheta}N\right)\left\|{\bf a}_N\right\|_2^2\;.$$
Therefore, and since  
$1+(2-\varepsilon)\varepsilon /3=(1-\varepsilon /3)(\varepsilon+1)$,  
$\,\varepsilon >0$, and $R,N,C_0(\varepsilon)\geq 1$,  
it follows by (4.36) that 
$$\eqalign{
S_t(R,X,N)
 &\leq 2^{17}C_2(\varepsilon)C_0(\varepsilon)\left({X\over Y}\right)^{\vartheta} 
\left( R^{1-\varepsilon /3}N\right)^{\varepsilon}
\left( R^{1-\varepsilon /3}+Y^{\vartheta}N\right)\left\|{\bf a}_N\right\|_2^2 = {}\cr 
 &=2^{17}C_2(\varepsilon)C_0(\varepsilon) 
\left( R^{1-\varepsilon /3}N\right)^{\varepsilon}
\left( \left({X\over Y}\right)^{\vartheta}R^{1-\varepsilon /3}
+X^{\vartheta}N\right)\left\|{\bf a}_N\right\|_2^2\;.
}$$
Since $Y=\min\left\{ X\,,\,R^{2-\varepsilon}/N\right\}$, $\,\varepsilon >0$ and $R\geq 1$, 
we have, in the above, 
$$\left({X\over Y}\right)^{\vartheta}R^{1-\varepsilon /3}
\leq R^{1-\varepsilon /3}
+\left( {XN\over R^{2-\varepsilon}}\right)^{\vartheta} R^{1-\varepsilon /3}  
\leq R+(XN)^{\vartheta} R^{1-2\vartheta-(1/3-\vartheta)\varepsilon}$$
and (given that $1/3-\vartheta\geq 1/9>0$) therefore obtain: 
$$S_t(R,X,N)
\leq 2^{17}C_2(\varepsilon)C_0(\varepsilon) 
R^{-\varepsilon^2/3}(RN)^{\varepsilon}
\left( R+(XN)^{\vartheta} R^{1-2\vartheta} 
+X^{\vartheta}N\right)\left\|{\bf a}_N\right\|_2^2\;.$$
By (4.38), we have 
$$2^{17}C_2(\varepsilon)R^{-\varepsilon^2 /3}<1\;,$$
so that in obtaining the above bound for $S_t(R,X,N)$ we 
have achieved the objective of showing that 
(4.29) holds if $Q=R$, $t\in{\Bbb R}$ and $X,N\geq 1$ 
are such that (4.35) holds. This (as noted below (4.35)) is sufficient to
establish that Proposition~$A(R)$ is true, which contradicts what is 
stated in (4.32). Consequently, as explained below (4.34), we have 
proof by contradiction that $A(Q)$ is true for all $Q\in{\Bbb N}$. 
The theorem therefore follows (as discussed in the paragraph 
containing (4.31))\quad$\blacksquare$  

\bigskip 
\bigskip 

\goodbreak\centerline{\bf \S 5. Schwartz Spaces, Fourier Integrals, 
Poisson Summation and the 
Analytic Large Sieve}

\bigskip 

This section is where we begin our dedicated preparation for the proof of 
Theorem~9. In it we have collected together certain definitions, remarks and 
lemmas, for use in the sections which follow.
 
\medskip

\goodbreak\proclaim Definition (Schwartz Spaces for ${\Bbb R}^n$ and ${\Bbb C}^n$). For $n\in{\Bbb N}$, the `Schwartz space', 
${\cal S}\bigl( {\Bbb R}^n\bigr)$,  
is the space of all functions $F : {\Bbb R}^n\rightarrow{\Bbb C}$ such that, for 
each pair $(A,{\bf j})\in [0,\infty)\times\bigl( {\Bbb N}\cup\{ 0\}\bigr)^n$, 
there exists a continuous and bounded function  
$F_{A,{\bf j}} : {\Bbb R}^n\rightarrow{\Bbb C}$ 
such that one has  
$$(x_1^2+\cdots +x_n^2)^A 
\,{\partial^{j_1+\cdots +j_n}\over\partial x_1^{j_1}\cdots\partial x_n^{j_n}}
\,F({\bf x}) = F_{A,{\bf j}}({\bf x})\qquad\qquad   
\hbox{for all $\quad {\bf x}\in{\Bbb R}^n$.}$$  
For $n\in{\Bbb N}$, we define ${\cal S}\bigl( {\Bbb C}^n\bigr)$ to be 
the space of all functions $f : {\Bbb C}^n\rightarrow{\Bbb C}$ such that 
the space ${\cal S}\bigl({\Bbb R}^{2n}\bigr)$ contains  
the function $F : {\Bbb R}^{2n}\rightarrow{\Bbb C}$    
given by 
$F({\bf x})=f( x_1+ix_2, x_3+ix_4,\ldots , x_{2n-1}+ix_{2n})\ $ 
(${\bf x}\in{\Bbb R}^{2n}$).

\medskip 

\goodbreak\proclaim Definition (Fourier Transforms). Let $n\in{\Bbb N}$. 
For $F\in{\cal S}\bigl( {\Bbb R}^n\bigr)$ we define the corresponding `Fourier transform', 
$\hat F : {\Bbb R}^n\rightarrow{\Bbb C}$ by  
$$\hat F({\bf y})
=\int_{-\infty}^{\infty}\cdots\int_{-\infty}^{\infty}F({\bf x})
\,{\rm e}\!\left( -{\bf y}\cdot{\bf x}\right) {\rm d}x_1\cdots {\rm d}x_n\qquad\qquad  
\hbox{(${\bf y}\in{\Bbb R}^n$),}\eqno(5.1)$$
where ${\bf y}\cdot{\bf x}=y_1 x_1+\cdots +y_n x_n$. 
The existence of the integral in the equation~(5.1) is guaranteed by the 
continuity of $F$, and the boundedness of the function 
${\bf x}\mapsto (1+x_1^2)\cdots(1+x_n^2)F({\bf x})$ 
(which follow from the definition of the Schwartz space for ${\Bbb R}^n$, 
with the help of some 
elementary ineqalities). Indeed, although 
Fourier transforms are defined slightly differently in [17, Chapter~13, Section~4],   
it is effectively shown by [17, Chapter~13, Theorem~4.1] that (5.1) defines 
a linear mapping of ${\cal S}\bigl( {\Bbb R}^n\bigr)$ into itself. \hfill\break  
\indent For $f\in{\cal S}\bigl( {\Bbb C}^n\bigr)$, we define the Fourier transform 
$\hat f : {\Bbb C}^n\rightarrow{\Bbb C}$ by 
$$\hat f({\bf w})=\hat F\Bigl( {\rm Re}\bigl( w_1\bigr) , -{\rm Im}\bigl( w_1\bigr) , 
{\rm Re}\bigl( w_2\bigr) , -{\rm Im}\bigl( w_2\bigr) , 
\ldots\ldots , {\rm Re}\bigl( w_n\bigr) , -{\rm Im}\bigl( w_n\bigr)\Bigr) ,\eqno(5.2)$$
where $F$ is the element of ${\cal S}\bigl( {\Bbb R}^{2n}\bigr)$  
given by 
$F({\bf x})=f( x_1+ix_2, x_3+ix_4,\ldots , x_{2n-1}+ix_{2n})\ $ 
(${\bf x}\in{\Bbb R}^{2n}$); 
this means that  
$$\hat f({\bf w})=\int_{\Bbb C}\cdots\int_{\Bbb C}
f({\bf z}) {\rm e}\left( -{\rm Re}( {\bf w}\cdot{\bf z} )\right) 
{\rm d}_{+}z_1\cdots {\rm d}_{+}z_n\qquad\qquad\hbox{(${\bf z}\in{\Bbb C}^n$),}\eqno(5.3)$$
where ${\bf w}\cdot{\bf z}=w_1 z_1+\cdots w_n z_n$ and  
${\rm d}_{+}z={\rm d}x\,{\rm d}y$ for $z\in{\Bbb C}$ with ${\rm Re}(z)=x$,  
${\rm Im}(z)=y$. Given the final remark of the preceding paragraph, it follows 
by (5.2) that,  for each  
$f\in{\cal S}\bigl( {\Bbb C}^n\bigr)$, one has 
$\hat f\in{\cal S}\bigl( {\Bbb C}^n\bigr)$. 

\medskip 

\goodbreak\proclaim Definition (The M\"{o}bius function for ${\Bbb Z}[i]$). 
For $n\in{\frak O}-\{ 0\}$, we define 
$$\mu_{\frak O}(n)
=\cases{0 &if there exists a Gaussian prime $\varpi$ such that $\varpi^2\mid n$, \cr 
(-1)^{\omega(n)} &otherwise,}\eqno(5.4)$$
where $\omega(n)$ is the number of 
prime ideals of the ring ${\frak O}$ that contain $n$ (i.e. one quarter of the 
number of Gaussian primes that divide $n$). This is a `multiplicative' function 
on ${\frak O}-\{ 0\}$, in the sense that it satisfies 
$$\mu_{\frak O}(mn)=\mu_{\frak O}(m)\mu_{\frak O}(n)\qquad\qquad 
\hbox{for $\quad m,n\in{\frak O}-\{0 \}\ $ with $\ (m,n)\sim 1$.}\eqno(5.5)$$
A useful property of this function is the identity 
$${1\over 4}\sum_{\scriptstyle d\in{\frak O}\atop\scriptstyle d\mid n}\mu_{\frak O}(d)
=\cases{1 &if $n\sim 1$, \cr 0 &otherwise,}\eqno(5.6)$$ 
which is valid for all non-zero $n\in{\frak O}$ (and may be deduced, in a very few steps, 
directly from (5.4), given that the ring of Gaussian integers is a principal ideal 
domain with $4$ units). 

\medskip 

\goodbreak\proclaim Definition (The Distance to the Nearest Gaussian Integer). For $\beta\in{\Bbb C}$ 
the `distance from $\beta$ to the nearest Gaussian integer' is the number $\|\beta\|\in[0,\infty)$ 
given by 
$$\|\beta\|=\min\{ |\beta -m| : m\in{\frak O}\}\;.$$
Since ${\Bbb C}={\Bbb R}+i{\Bbb R}$, ${\frak O}={\Bbb Z}+i{\Bbb Z}$ 
and $|x+iy|^2=x^2+y^2\geq x^2$ ($x,y\in{\Bbb R}$), one has 
$$\|\beta\|^2=\|{\rm Re}(\beta)\|^2+\|{\rm Im}(\beta)\|^2\qquad\qquad 
\hbox{($\beta\in{\Bbb C}$).}\eqno(5.7)$$ 
Given that every real interval of form $(x-1/2,x+1/2]$ contains exactly one 
integer, it follows by (5.7) that, for all $\beta\in{\Bbb C}$, one has 
$\|\beta\|^2\leq (1/2)^2+(1/2)^2=1/2$. 

\medskip 

\goodbreak 
\noindent{\bf Remark.}$\,$   
In addition to the above definitions, we shall have cause to recall the  
orthogonality of the characters of the additive groups ${\frak O}/m{\frak O}$ 
($0\neq m\in{\frak O}$). Specifically, when $0\neq m\in{\frak O}$ and $a,b\in{\frak O}$, 
one has:  
$$\sum_{n\bmod m{\frak O}}{\rm e}\!\left( {\rm Re}\left( {an\over m}\right)\right)  
\,\overline{{\rm e}\!\left( {\rm Re}\left( {bn\over m}\right)\right) } 
=\sum_{n\bmod m{\frak O}}{\rm e}\!\left( {\rm Re}\left( {(a-b)n\over m}\right)\right)
=\cases{|m|^2 &if $a\equiv b\bmod m{\frak O}$, \cr 
0 &otherwise.}\eqno(5.8)$$
The first equality in (5.8) follows by Euler's formula for $e^{i\theta}$, and 
the identity $\exp(x)\exp(y)=\exp(x+y)$.  
In the case $a\equiv b\bmod m{\frak O}$ the second inequality of (5.8) 
follows since   
${\rm e}(k)=\exp(2\pi ik)=1$ for all $k\in{\Bbb Z}$, and since one has  
$[{\frak O}:m{\frak O}]^2=[{\frak O}:m{\frak O}][{\frak O}:\overline{m}{\frak O}]
=[{\frak O}:m{\frak O}][m{\frak O}:m\,\overline{m}{\frak O}]$, where  
${\frak O}={\Bbb Z}+{\Bbb Z}i$ and 
$m\,\overline{m}=|m|^2\in{\Bbb Z}$, which implies  that 
$[{\frak O}:m{\frak O}]^2=[{\frak O}:m\,\overline{m}{\frak O}]
=[{\frak O}:|m|^2{\frak O}]=\bigl|({\Bbb Z}/|m|^2{\Bbb Z})\times 
({\Bbb Z}/|m|^2{\Bbb Z})\bigr|=|m|^4$, 
and so shows that the additive group ${\frak O}/m{\frak O}$ has order 
$$\,|{\frak O}/m{\frak O}|=|m|^2\;.\eqno(5.9)$$  
In the case $a\not\equiv b\bmod m{\frak O}$, the second inequality of 
(5.8) follows by considering the effect of the substitutions $n=n'+1$ 
and $n=n''+i$ (applied to 
the variable of summation).

\bigskip

\goodbreak\proclaim Lemma~5.1 (Fourier's Inversion Formulae). Let $n\in{\Bbb N}$. 
Then, when $F\in{\cal S}\bigl( {\Bbb R}^n\bigr)$, one has 
$$\int_{-\infty}^{\infty}\cdots \int_{-\infty}^{\infty}
\hat F({\bf y})\,{\rm e}({\bf x}\cdot{\bf y})
\,{\rm d}y_1\,\cdots\,{\rm d}y_n
=F({\bf x})\qquad\qquad\hbox{(${\bf x}\in{\Bbb R}^n$).}\eqno(5.10)$$
For $f\in{\cal S}\bigl( {\Bbb C}^n\bigr)$, one has 
$$\int_{\Bbb C}\cdots\int_{\Bbb C} 
\hat f({\bf w})\,{\rm e}( {\rm Re}({\bf z}\cdot{\bf w}) )
\,{\rm d}_{+}w_1\,\cdots\,{\rm d}_{+}w_n
=f({\bf z})\qquad\qquad  
\hbox{(${\bf z}\in{\Bbb C}^n$).}\eqno(5.11)$$

\medskip

\goodbreak 
\noindent{\bf Proof.}\quad Given (5.1), the result (5.10) simply states 
that one has ${\hat G}(-{\bf x})=F({\bf x})$, where $G=\hat F$; as much 
may be deduced from  [17, Chapter~13, Theorem~5.1], by way of 
one linear substitution (applied to the relevant variables of integration).
The result (5.11) follows directly from (5.10): to see this one has only to 
apply the definition of ${\cal S}\bigl( {\Bbb C}^n\bigr)$, and 
the definitions made in connection with (5.2)-(5.3), along with 
some substitions of the form $y_{2k}=-v_k$\ $\blacksquare$

\bigskip

\goodbreak\proclaim Lemma~5.2. Let $n\in{\Bbb N}$. Then, for $m=1,\ldots ,n$, there is a 
linear operator ${\cal L}_m$ with domain ${\cal S}\bigl( {\Bbb C}^n\bigr)$, and 
range contained in ${\cal S}\bigl( {\Bbb C}^n\bigr)$, that is given by 
$${\cal L}_m f({\bf z}) 
=-\left({\partial^2\over\partial x_m^2}+{\partial^2\over\partial y_m^2}\right) f({\bf z})\qquad\qquad  
\hbox{($f\in{\cal S}\bigl( {\Bbb C}^n\bigr)$, ${\bf z}\in{\Bbb C}^n$),}\eqno(5.12)$$
where $x_m,y_m$ denote (respectively) the real and imaginary parts of $z_m$. 
Let $f$ lie in the space ${\cal S}\bigl( {\Bbb C}^n\bigr)$. Then, for $m=1,\ldots ,n$, 
the functions $f$ and ${\cal L}_m f$ have 
Fourier transforms $\hat{f},\widehat{{\cal L}_m f}\in{\cal S}\bigl( {\Bbb C}^n\bigr)$
that are related to one another by: 
$$\widehat{{\cal L}_m f}({\bf w})
=\bigl| 2\pi w_m\bigr|^2 \hat f({\bf w})\qquad\qquad  
\hbox{(${\bf w}\in{\Bbb C}^n$).}\eqno(5.13)$$
For all ${\bf w}\in{\Bbb C}^n$, and all ${\bf j}\in({\Bbb N}\cup\{ 0\})^n$, one has 
$$\bigl|\hat f({\bf w})\bigr|\prod_{m=1}^n \bigl( 2\pi\bigl| w_m\bigr|\bigr)^{2 j_m} 
=\left|\widehat{{\cal L}^{\bf j} f}({\bf w})\right| 
\leq\widehat{\,|{\cal L}^{\bf j}f|\,}({\bf 0}) 
=\int_{\Bbb C}\cdots\int_{\Bbb C} 
\bigl| {\cal L}^{\bf j} f({\bf z})\bigr|
{\rm d}_{+}z_1\,\cdots\,{\rm d}_{+}z_n\;,\eqno(5.14)$$
where ${\cal L}^{\bf j}$ denotes the operator ${\cal L}_1^{j_1}\cdots {\cal L}_n^{j_n}$. 

\medskip

\goodbreak 
\noindent{\bf Proof.}\quad 
Let $f\in{\cal S}\bigl( {\Bbb C}^n\bigr)$ and $m\in\{ 1,\ldots ,n\}$. 
Then, as an immediate corollary of the definitions of the spaces 
${\cal S}\bigl( {\Bbb C}^n\bigr)$ and ${\cal S}\bigl( {\Bbb R}^{2n}\bigr)$,  
one has also $({\partial/\partial x_m})^2 f\in{\cal S}\bigl( {\Bbb C}^n\bigr)$ and 
$({\partial/\partial y_m})^2 f\in{\cal S}\bigl( {\Bbb C}^n\bigr)$. 
Therefore, and since ${\cal S}\bigl( {\Bbb C}^n\bigr)$ is a complex vector space, 
we have ${\cal L}_m f\in{\cal S}\bigl( {\Bbb C}^n\bigr)$ when 
${\cal L}_m f : {\Bbb C}^n\rightarrow{\Bbb C}$ is the function given by (5.12). 

The case $n=1$ of (5.12) is [23, Lemma~4.2, Equation~(4.6)]: a short proof 
of that result is supplied in [23] (we do not repeat it here).  

Supposing now that $n>1$, it will suffice to prove (5.13) for $m=n$:  
the other cases may be proved similarly. By an appeal to the relevant definitions 
(of Fourier transforms) we may first express $\widehat{{\cal L}_n f}({\bf w})$ 
as an integral over ${\Bbb R}^{2n}$. Then, by Fubini's reduction theorem for 
higher dimensional real integrals [1, Section~15.7], we find that 
$$\widehat{{\cal L}_n f}({\bf w}) 
=\hat G\bigl( w_1 , \ldots , w_{n-1}\bigr)\;,\eqno(5.15)$$
where, for ${\bf s}\in{\Bbb C}^{n-1}$, 
$$G({\bf s})=\widehat{{\cal L}_1 f_{\bf s}}\bigl( w_n\bigr)\;,\eqno(5.16)$$
with $f_{\bf s} : {\Bbb C}\rightarrow{\Bbb C}$ being given by 
$$f_{\bf s}\bigl( z_1\bigr) 
=f\bigl( s_1,s_2,\ldots ,s_{n-1},z_1\bigr)\qquad\qquad\hbox{($z_1\in{\Bbb C}$).}\eqno(5.17)$$
Since it is a corollary of the relevant definitions that 
$f\in{\cal S}\bigl( {\Bbb C}^n\bigr)$ implies $f_{\bf s}\in{\cal S}({\Bbb C})\,$ 
(for all ${\bf s}\in{\Bbb C}^{n-1}$), we have here, by the case $m=n=1$ of (5.13), 
$$\widehat{{\cal L}_1 f_{\bf s}}\bigl( w_n\bigr) 
=\bigl| 2\pi w_n\bigr|^2 \hat f_{\bf s}\bigl( w_n\bigr)\;,$$ 
and so it follows by (5.16), (5.15) and the linearity of Fourier transforms 
that 
$$\widehat{{\cal L}_n f}({\bf w})=|2\pi w_n|^2\hat H(w_1,\ldots ,w_{n-1})\;,$$  
where, for ${\bf s}\in{\Bbb C}^{n-1}$, one has $H({\bf s})=H({\bf s};w_n)=\hat f_{\bf s}(w_n)$. 
By (5.17) and the definition of the Fourier transform (as it applies to 
${\cal S}({\Bbb C})$, 
${\cal S}\bigl( {\Bbb C}^{n-1}\bigr)$ and ${\cal S}\bigl( {\Bbb C}^n\bigr)$), 
the equation just obtained is the case $m=n$ of the result (5.13) of the lemma: 
as noted above, the proofs for 
$m=1,\ldots ,n-1$ are similar.

From (5.13) it follows by induction on ${\Bbb N}^n$ that, for ${\bf w}\in{\Bbb C}^n$ and 
${\bf j}\in({\Bbb N}\cup\{ 0\})^n$, one has 
$$\hat f({\bf w})\prod_{m=1}^n \bigl( 2\pi\bigl| w_m\bigr|\bigr)^{2 j_m} 
=\widehat{{\cal L}^{\bf j} f}({\bf w})\;,$$
where, as stated below (5.14), ${\cal L}^{\bf j}={\cal L}_1^{j_1}\cdots {\cal L}_n^{j_n}$. 
The first equality in (5.14) follows trivially; the last equality there is simply  
a statement of the relevant definition; and the 
inequality in (5.14) is also trivial, given that 
${\rm e}({\rm Re}({\bf w}\cdot{\bf z}))$ is, for 
${\bf w},{\bf z}\in{\Bbb C}^n$, a complex number of unit modulus.  
\ $\blacksquare$

\bigskip

\goodbreak 
\noindent{\bf Remarks.}\quad 
The operator ${\cal L}_m$ defined by the equation~(5.12), in Lemma~5.2, is (apart from 
the factor $-1$) the Euclidean Laplacian operator. By (5.12) one has 
$${\cal L}_m=-4\,{\partial\over\partial z_m}\,{\partial\over\partial\,\overline{z_m}}\;,
\eqno(5.18)$$
where 
$${\partial\over\partial z_m}={1\over 2}\left( {\partial\over\partial x_m}
-i\,{\partial\over\partial y_m}\right)\qquad\ \hbox{and}\qquad\  
{\partial\over\partial\,\overline{z_m}}
={1\over 2}\left( {\partial\over\partial x_m}
+i\,{\partial\over\partial y_m}\right)\eqno(5.19)$$
(with $x_m$ and $y_m$ denoting the real and imaginary parts of $z_m$).
Although Lemma~5.2 assigns ${\cal L}_m$ the domain ${\cal S}\bigl( {\Bbb C}^n\bigl)$, 
it is helpful not to be so restrictive when assigning the domains of 
the above operators $\partial /\partial z_m$ and $\partial /\partial\,\overline{z_m}$.  
Indeed, given any $n$ non-empty open regions $D_1,\ldots ,D_n\subset{\Bbb C}$,  
we may (assuming $n\geq m$) apply these operators to any function 
$f : D_1\times\cdots \times D_n\rightarrow{\Bbb C}$ 
such that, for all 
${\bf s}\in D_1\times\cdots\times D_n$, 
the function $z_m\mapsto f\bigl( s_1,\ldots ,s_{m-1},z_m,s_{m+1},\ldots ,s_n\bigr)$ 
is smooth on $D_m$ (in the sense defined at the start of Subsection~1.2); 
for such $f$, one has, when ${\bf d}\in{\Bbb C}^n$ and 
$g_{\bf d}$ maps ${\bf z}\in{\Bbb C}^n$ to $( d_1 z_1,\ldots , d_n z_n)\in{\Bbb C}^n$,  
the elementary identities  
$${\partial\over\partial z_m}\left( f\circ g_{\bf d}\right) 
=d_m\,\left( {\partial\over\partial z_m}\,f\right)\circ g_{\bf d}\qquad\  
\hbox{and}\qquad\  
{\partial\over\,\overline{\partial z_m}}\left( f\circ g_{\bf d}\right) 
=\overline{d_m}\,\left( {\partial\over\partial\,\overline{z_m}}\,f\right)\circ g_{\bf d}\;,$$
which, by (5.18), together imply that  
$${\cal L}_m\left( f\circ g_{\bf d}\right) 
=\bigl| d_m\bigr|^2\left( {\cal L}_m f\right)\circ g_{\bf d}\;.\eqno(5.20)$$

It is also worth noting, for use later, that if $q$ is a holomorphic complex function on some 
non-empty open region $D\subset{\Bbb C}$, then for $z_m\in D$ one has (as a consequence of 
the Cauchy-Riemann equations):  
$${\partial\over\partial z_m}\left( q\bigl( z_m\bigr)\right) =q'\bigl( z_m\bigr) ,\quad 
{\partial\over\partial z_m}\left( q\bigl( \overline{z_m}\bigr)\right) =0 ,\quad 
{\partial\over\partial\,\overline{z_m}}\left( q\bigl( z_m\bigr)\right) =0 ,\quad 
{\partial\over\partial\,\overline{z_m}}\left( q\bigl(\overline{z_m}\bigr)\right) 
=q'\bigl(\overline{z_m}\bigr) .\eqno(5.21)$$

\bigskip

\goodbreak\proclaim Lemma~5.3 (Poisson Summation over ${\Bbb R}^n$ and over ${\Bbb C}$). 
For $n\in{\Bbb N}$,  one has 
$$\sum_{{\bf x}\in{\Bbb Z}^n} F({\bf x})
=\sum_{{\bf y}\in{\Bbb Z}^n}\hat F ({\bf y})\qquad\qquad 
\hbox{($F\in{\cal S}\bigl( {\Bbb R}^n\bigr)$).}\eqno(5.22)$$
For $\tau\in{\Bbb C}$, and $f\in{\cal S}\bigl( {\Bbb C}\bigr)$, one has: 
$$\sum_{\nu\in{\frak O}}f(\nu){\rm e}\left( {\rm Re}(\tau\nu)\right)
=\sum_{\xi\in{\frak O}}\hat f(\xi-\tau)\;.\eqno(5.23)$$

\medskip

\goodbreak 
\noindent{\bf Proof.}\quad  This lemma is essentially a reproduction 
of part of [23, Lemma~4.1]; the proof given there sketches how to deduce (5.23) 
from the case $n=2$ of (5.22). 
For a proof of (5.22), see [17, Chapter13, Section~6]\ $\blacksquare$

\bigskip

\goodbreak
\proclaim Lemma~5.4. Let $n\in{\Bbb N}$, $\Delta\in(0,\infty)$, ${\bf\Omega}\in(0,\infty)^n$,  
$C>1$ and $f\in{\cal S}\bigl( {\Bbb C}^n\bigr)$.  
Suppose that 
$${\cal L}_1^{j_1}\cdots{\cal L}_n^{j_n} f({\bf z}) 
\ll_{\bf j}\prod_{m=1}^n\left(\Delta\bigl| z_m\bigr|^2\right)^{-j_m}\qquad\qquad  
\hbox{($\,{\bf z}\in\bigl({\Bbb C}^{*}\bigr)^n$, $\ {\bf j}\in ({\Bbb N}\cup\{ 0\})^n$),}\eqno(5.24)$$
where the operators ${\cal L}_1,\ldots ,{\cal L}_n$ are those defined by 
the equation~(5.12) of Lemma~5.2. 
Suppose moreover that, for 
all ${\bf z}\in{\Bbb C}^n$, 
$$f({\bf z})=0\qquad\ {\rm unless}\qquad\   
C^{-1}\Omega_m<\bigl| z_m\bigr|^2 <C\Omega_m\quad {\rm for}\quad  
m=1,\ldots ,n\;.\eqno(5.25)$$ 
Then, for $j=0,1,2,\ldots\ $, one has 
$$\hat f({\bf w})\ll_{j,n}\prod_{m=1}^n 
{C\Omega_m\over\bigl( 1+C^{-1}\Delta\Omega_m\bigl| w_m\bigr|^2\bigr)^j}\qquad\qquad   
\hbox{(${\bf w}\in{\Bbb C}^n$).}\eqno(5.26)$$

\medskip

\goodbreak 
\noindent{\bf Proof.}\quad  
Let ${\bf w}\in{\Bbb C}^n$; and let ${\bf j}\in({\Bbb N}\cup\{ 0\})^n$. 
Since $f\in{\cal S}\bigl( {\Bbb C}^n\bigr)$, it follows by the result (5.14) of 
Lemma~5.2, and by the hypotheses (5.24) and (5.25), that 
$$\bigl|\hat f({\bf w})\bigr|\prod_{m=1}^n \bigl( 2\pi\bigl| w_m\bigr|\bigr)^{2 j_m} 
\ll_{\bf j}\ \int\limits_{{\bf z}\in{\cal A}_1\times\cdots\times{\cal A}_n}
\!\!\left(\prod_{m=1}^n\left(\Delta\bigl| z_m\bigr|^2\right)^{\!-j_m}\!\!\right) 
\!{\rm d}_{+}z_1\,\cdots\,{\rm d}_{+}z_n
=\prod_{m=1}^n\Biggl(\!\Delta^{-j_m}
\!\int\limits_{{\cal A}_m}\!|z|^{-2 j_m}\,{\rm d}_{+}z\!\Biggr) ,$$
where, for $m=1,\ldots ,n$, 
$${\cal A}_m =\bigl\{ z\in{\Bbb C} : C^{-1}\Omega_m < |z|^2 < C\Omega_m\bigr\}\;.$$
Hence, by using the upper bounds 
$$\int\limits_{{\cal A}_m}\!|z|^{-2 j_m}\,{\rm d}_{+}z
\leq{1\over\bigl( C^{-1}\Omega_m\bigr)^{j_m}}\int\limits_{{\cal A}_m}{\rm d}_{+}z 
<{\pi C\Omega_m\over\bigl( C^{-1}\Omega_m\bigr)^{j_m}}\qquad\qquad  
\hbox{($m=1,\ldots ,n$),}$$
one finds that 
$$\bigl|\hat f({\bf w})\bigr|\prod_{m=1}^n \bigl( 2\pi\bigl| w_m\bigr|\bigr)^{2 j_m} 
\ll_{\bf j}\prod_{m=1}^n {C\Omega_m\over\bigl( C^{-1}\Delta\Omega_m\bigr)^{j_m}}\;.$$

Let $j$ be a non-negative integer. We apply 
apply the last bound above, for the unique   
${\bf j}\in({\Bbb N}\cup\{ 0\})^n$ 
having, for $m=1,\ldots ,n$,  
$$j_m=\cases{j &if $(2\pi |w_m|)^2>\bigl( C^{-1}\Delta\Omega_m\bigr)^{-1}$, \cr 
0 &otherwise.}$$ 
This yields the upper bound 
$$\bigl|\hat f({\bf w})\bigr| 
\ll_{j,n}\,\prod_{m=1}^n {C\,\Omega_m\over\left( 
\max\bigl\{ 4\pi^2 |w_m|^2C^{-1}\Delta\Omega_m\ ,\,1\bigr\}\right)^j}\;.$$
Since $\max\{ a,b\}\geq (a+b)/2>0$ for all positive $a,b$, the result (5.26) follows
\ $\blacksquare$

\bigskip

The next three lemmas are perfect Gaussian integer analogues of 
results contained in [21, Lemmas~2.2,~2.3 and~2.4]: 
although the relevant proofs are also 
analogous, we have nevertheless chosen (for the sake of completeness) 
to include sketched proofs of these lemmas. 

\bigskip

\goodbreak\proclaim Lemma~5.5. Let $\Delta,\Omega_1\in(0,\infty)$, $C>1$ and $f\in{\cal S}({\Bbb C})$ 
satisfy the case $n=1$ of the hypotheses (5.24), (5.25) stated in the previous lemma. 
Then, for $\tau\in{\Bbb C}-{\frak O}$ and $j=2,3,4,\ldots\ $, one has
$$\sum_{\xi\in{\frak O}}\bigl|\hat f(\xi -\tau)\bigr|\ll_{j,C} 
\bigl(\Delta\Omega_1\|\tau\|^2\bigr)^{-j}\Omega_1\;,\eqno(5.27)$$ 
where $\beta\mapsto\|\beta\|$ is the `distance to the nearest Gaussian integer' function, 
defined just above (5.7). 

\medskip

\goodbreak 
\noindent{\bf Proof.}\quad 
Let $\tau\in{\Bbb C}$; and let $j$ be an integer with $j\geq 2$. 
As noted below (5.7), one must have $\|\tau\|^2\leq 1/2$; there therefore 
exists some 
$\nu\in{\frak O}$ such that $1/\sqrt{2}\geq |\tau -\nu|=\|\tau\|>0$ 
(the last inequality  holding by virtue of the hypothesis that $\tau\not\in{\frak O}$). 
For such a $\nu\in{\frak O}$, the sum on the left-hand side of (5.27) 
may be rewritten as $\sum_{\xi'\in{\frak O}}|\hat f(\xi'+\nu -\tau)|$. 
Consequently it follows by the case $n=1$ of the bound (5.26), which Lemma~5.4 provides, 
that we have here 
$$\sum_{\xi\in{\frak O}}\bigl|\hat f(\xi -\tau)\bigr| 
=\sum_{\xi'\in{\frak O}}\bigl|\hat f(\xi'+\nu -\tau)\bigr| 
\ll_j {C\Omega_1\over\bigl( 1+C^{-1}\Delta\Omega_1\bigl|\nu-\tau\bigr|^2\bigr)^j}
+\sum_{0\neq\xi'\in{\frak O}}
{C\Omega_1\over\bigl( 1+C^{-1}\Delta\Omega_1\bigl|\xi'-(\tau -\nu)\bigr|^2\bigr)^j}\;.$$
The result (5.27) follows, since   
$|\xi'-(\tau-\nu)|\geq |\xi'|-1/\sqrt{2}>(1-1/\sqrt{2})|\xi'|$, 
for $0\neq\xi'\in{\frak O}$, and since one has both 
$\infty>|\nu -\tau|^{-2j}=\|\tau\|^{-2j}\geq 2^j$ and, given that $j\geq 2$,  
$\,\sum_{0\neq\xi'\in{\frak O}} |\xi'|^{-2j}\ll 1$\ $\blacksquare$

\bigskip 

\goodbreak 
\goodbreak\proclaim Lemma~5.6.  Let $\Delta,\Omega_1\in(0,\infty)$, $C>1$ and $f\in{\cal S}({\Bbb C})$ 
satisfy the case $n=1$ of the conditions~(5.24), (5.25) stated in Lemma~5.4; let 
$d\in{\frak O}-\{ 0\}$ and $j\in\{ 2,3,4,\ldots\ \}$. Then, for $\tau\in{\Bbb C}-{\frak O}$, 
$$\sum_{\scriptstyle m\in{\frak O}\atop\scriptstyle m\equiv 0\bmod d{\frak O}} 
f(m)\,{\rm e}\!\left( {\rm Re}\left(\tau m/d\right)\right)
\ll_{j,C}\left(\Delta |d|^{-2}\Omega_1\|\tau\|^2\right)^{-j} |d|^{-2}\Omega_1\;.\eqno(5.28)$$ 
For $h,k\in{\frak O}$, $q\in{\frak O}-\{ 0\}$ and $B>0$, one has moreover: 
$$\eqalignno{
\sum_{\scriptstyle m\in{\frak O}\atop\scriptstyle m\equiv 0\bmod d{\frak O}} 
\!\!f(m)\,&{\rm e}\!\left( {\rm Re}\!\left( {h (m/d)^{*}\over q}\right)\right) = {} 
 &(5.29)\cr 
 &={1\over |q|^2}\sum_{\scriptstyle b\bmod q{\frak O}\atop\scriptstyle 
\left\| {b\over q}\right\|^2\leq {B |d|^2\over\Delta\Omega_1}}\!\!S(-h,b;q) 
\!\!\!\sum_{\scriptstyle m\in{\frak O}\atop\scriptstyle m\equiv 0\bmod d{\frak O}} 
\!\!f(m)\,{\rm e}\!\left( {\rm Re}\!\left({b\over q}\,{m\over d}\right)\right)  
 +O_{j,C}\left( {\Delta^{-1} |q|^2\over B^j}\right) }$$
and 
$$\eqalignno{
\sum_{\scriptstyle m\in{\frak O}\atop\scriptstyle m\equiv 0\bmod d{\frak O}} 
\!\!f(m)\,&S(hm/d , k ; q) = {} 
 &(5.30)\cr 
 &=\sum_{\scriptstyle b\bmod q{\frak O}\atop\scriptstyle 
\left\| {b\over q}\right\|^2\leq {B |d|^2\over\Delta\Omega_1}}\!\!c_q(b,h;k) 
\!\!\!\sum_{\scriptstyle m\in{\frak O}\atop\scriptstyle m\equiv 0\bmod d{\frak O}} 
\!\!f(m)\,{\rm e}\!\left( {\rm Re}\!\left({b\over q}\,{m\over d}\right)\right)  
 +O_{j,C}\left( {\Delta^{-1} |q|^4\over B^j}\right) , }$$
where $(m/d)^{*}$ denotes an element of ${\frak O}$ satisfying 
$(m/d) (m/d)^{*}\equiv 1\bmod q{\frak O}$ 
(it thereby being an implicit condition of the summation on the 
left-hand side of (5.29) that $m$ must satisfy $(m/d , q)\sim 1$), 
while $S(u,v;w)$ is the `simple' Kloosterman sum defined in (1.3.6), 
and the term $c_q(b,h;k)$ in (5.30) is given by  
$$c_q(b,h;k)
=\sum_{\scriptstyle a\bmod q{\frak O}\atop{\scriptstyle (a,q)\sim 1\atop\scriptstyle 
ab\equiv h\bmod q{\frak O}}} {\rm e}\!\left( {\rm Re}\!\left( {ak\over q}\right)\right) . 
\eqno(5.31)$$ 
 
\goodbreak 
\noindent{\bf Proof.} The results (5.28)-(5.30) can be rewritten by 
expressing $f(m)$ as $g(n)$, where $n=m/d\in{\frak O}$ (in the relevant summations), 
and where $g(z)=f(dz)$ for $z\in{\Bbb C}$. One has therefore only to prove the 
case $d=1$ of the lemma: for it follows by the remarks subsequent to Lemma~5.2 
that the conditions on $f$ in (5.24) and (5.25) imply that the function $g$ 
satisfies similar conditions 
(differing only in that $\Omega_m$ is replaced by $|d|^{-2}\Omega_m$).

Since the case $d=1$ of (5.28) is an immediate corollary of the identity (5.23) of Lemma~5.3 
and the bound (5.27) of Lemma~5.5,
it therefore only remains to consider, for $d=1$, the results (5.29) and (5.30)-(5.31). 
If one lets $B\rightarrow +\infty$, then the sums over $b$ appearing in (5.29) and 
(5.30) become sums over all $b\bmod q{\frak O}$; 
the case $d=1$, $B\rightarrow +\infty$ of (5.29) is 
therefore a direct consequence of (5.8) (i.e. of the orthogonality of the characters 
of the additive group ${\frak O}/q{\frak O}$); while the case 
$d=1$, $B\rightarrow +\infty$ of (5.30)-(5.31) is an immediate consequence 
of the definition (1.3.6) of the `simple' Kloosterman sum. The $O$-terms in 
in (5.30) and (5.31) therefore need only serve 
as upper bounds for the sums 
$$E_1=\sum_{\scriptstyle b\bmod q{\frak O}\atop\scriptstyle 
\left\| {b\over q}\right\|^2 > {B\over\Delta\Omega_1}}\Biggl|\ {S(-h,b;q)\over q^2}  
\sum_{m\in{\frak O}} 
f(m)\,{\rm e}\!\left( {\rm Re}\!\left({bm\over q}\right)\right)\Biggr|$$
and 
$$E_2=\sum_{\scriptstyle b\bmod q{\frak O}\atop\scriptstyle 
\left\| {b\over q}\right\|^2 > {B\over\Delta\Omega_1}}\Biggl|\ c_q(b,h;k) 
\sum_{m\in{\frak O}} 
f(m)\,{\rm e}\!\left( {\rm Re}\!\left({bm\over q}\right)\right)\Biggr|\;,$$
respectively (note that the relevant definitions trivially imply 
that $\| b/q\|$ is determined by the residue class of $b\bmod q{\frak O}$). To verify that 
those $O$-terms do serve in this capacity, we note firstly that  one has 
$$\eqalign{ 
\sum_{\scriptstyle b\bmod q{\frak O}\atop\scriptstyle 
\left\| {b\over q}\right\|^2 > {B\over\Delta\Omega_1}}\Biggl|  
\sum_{m\in{\frak O}} 
f(m)\,{\rm e}\!\left( {\rm Re}\!\left({bm\over q}\right)\right)\Biggr|
 &=\sum_{\scriptstyle b\bmod q{\frak O}\atop\scriptstyle 
\left\| {b\over q}\right\|^2 > {B\over\Delta\Omega_1}}
O_{j,C}\left( \biggl(\Delta\Omega_1\left\| 
{b\over q}\right\|^2\,\biggr)^{\!\!-(j+1)}\Omega_1\right) \ll_{j,C} {}\cr 
 &\ll_{j,C} \left(\Delta\Omega_1\right)^{-(j+1)}\Omega_1 |q|^{2(j+1)} 
\sum_{\scriptstyle\beta\in{\frak O}\atop\scriptstyle |\beta|^2\geq {B |q|^2\over\Delta\Omega_1}} 
|\beta|^{-2(j+1)} = {}\cr 
 &=\left(\Delta\Omega_1\right)^{-(j+1)}\Omega_1 |q|^{2(j+1)} 
O\left( \Bigl( {B |q|^2\over\Delta\Omega_1}\Bigr)^{-j}\right) \ll\Delta^{-1} |q|^2 B^{-j}
}$$
(the first line being an application of the case $d=1$, $\tau =b/q$ of 
(5.28), with `$j+1$' substituted for `$j$'$\,$). The verification is 
completed by noting that both $|S(-h,b;q)|$ and $|c_q(b,h;k)|$
are (given (5.9), (1.3.6) and (5.31)) bounded above by $|q|^2$, so that, by the above,  
one has 
$E_1\ll_{j,C} \Delta^{-1} |q|^2 B^{-j}$ and 
$E_2\ll_{j,C} \Delta^{-1} |q|^4 B^{-j}\ \,\blacksquare$

\bigskip 

\goodbreak\proclaim Lemma~5.7. Let $b,h,k\in{\frak O}$ and $q\in{\frak O}-\{ 0\}$. 
For $u,v\in{\frak O}$, let $c_q(u,v;k)$ be given by the equation~(5.31), in the 
previous lemma. Then one has 
$$c_q(b,h;k)
=0\qquad\ \hbox{if $\quad (b,q)\not\sim (h,q)$.}\eqno(5.32)$$ 
If it is, however, the case that 
$$(b,q)\sim(h,q)\sim c\in{\frak O}\qquad\ \hbox{(say),}\eqno(5.33)$$
then $c\neq 0$, and 
$$c_q(b,h;k)
={1\over 4}\sum_{\scriptstyle t\in{\frak O}\atop\scriptstyle t\mid (c,k)}
\mu_{\frak O}\!\left( {c\over t}\right) 
\,|t|^2\,{\rm e}\!\left( {\rm Re}\!\left( {h\over c}\,{k\over t}\,{(b/t)^{*}\over q/c}\right)\right) , 
\eqno(5.34)$$ 
where it is implicit in the last summation that $(c/t,q/c)\sim 1$ 
(the factor $(b/t)^{*}/(q/c)=((b/c)(c/t))^{*}/(q/c)$ here having the 
meaning explained in Subsection~1.5, under the heading `Number-Theoretic Notation'). 
One has, in particular, 
$$S(k,0;q)=S(0,k;q)=c_q(0,0;k)
=\mu_{\frak O}\!\!\left( {q\over (q,k)}\right) |(q,k)|^2 
\!\!\prod_{\scriptstyle{\rm prime\ ideals}\ \varpi{\frak O}\subset{\frak O}\atop\scriptstyle 
\varpi{\frak O}\ni (q,k)\,,\,\varpi{\frak O}\not\ni q/(q,k)}
\!\!\left( 1-{1\over |\varpi|^2}\right) ,\eqno(5.35)$$
where $S(u,v;w)$ is the `simple Kloosterman sum' defined in (1.3.6).

\medskip 

\goodbreak 
\noindent{\bf Proof.}\quad 
The conditions of summation in (5.31) imply the congruences 
$h\equiv ab\bmod q{\frak O}$ 
and $b\equiv a^{*}h\bmod q{\frak O}\,$ (which are equivalent when one has $(a,q)\sim 1$). 
The sum defining $c_q(b,h;k)$ is therefore empty unless one has both 
$(b,q)\mid (h,q)$ and $(h,q)\mid (b,q)$; and so the result (5.32) follows. 

Suppose now that $c\in{\frak O}$, and that (5.33) holds. Then $c\mid q$, 
and $q\neq 0$, so $c\neq 0$. By (5.33), we have 
$$b=Bc,\qquad h=Hc\quad\hbox{and}\quad q=Qc\;,\eqno(5.36)$$
for a unique $Q\in{\frak O}-\{ 0\}$, and a unique pair $B,H\in{\frak O}$ satisfying 
$$(B,Q)\sim 1\sim (H,Q)\;.\eqno(5.37)$$
Hence we can find an $A\in{\frak O}$ satisfying $AB\equiv 1\bmod Q{\frak O}$. This $A$ 
necessarily satisfies $(A,Q)\sim 1$.  Then, by (5.36) and (5.37), we may 
deduce from the definition (5.31) of $c_q(b,h;k)$ that 
$$c_q(b,h;k) 
=\sum_{\scriptstyle a\bmod Qc{\frak O}\atop{\scriptstyle (a,c)\sim 1\atop\scriptstyle 
a\equiv AH\bmod Q{\frak O}}} {\rm e}\!\left( {\rm Re}\!\left( {ka\over Qc}\right)\right) 
=\sum_{\scriptstyle d\bmod c{\frak O}\atop\scriptstyle (AH+dQ,c)\sim 1} 
{\rm e}\!\left( {\rm Re}\!\left( {kAH\over Qc}+{kd\over c}\right)\right) .$$
By the property (5.6) of the M\"{o}bius function $\mu_{\frak O}$, 
we therefore obtain: 
$$c_q(b,h;k)
={\rm e}\!\left( {\rm Re}\!\left( {kAH\over Qc}\right)\right) 
\sum_{\scriptstyle s\in{\frak O}\atop\scriptstyle s\mid c}\mu_{\frak O}(s) 
\sum_{\scriptstyle d\bmod c{\frak O}\atop\scriptstyle dQ\equiv -AH\bmod s{\frak O}} 
{\rm e}\!\left( {\rm Re}\!\left( {kd\over c}\right)\right) ,$$
where, since $(AH,Q)\sim 1$, the innermost sum on the right is empty unless 
$(s,Q)\sim 1$. Consequently, 
$$\eqalign{
c_q(b,h;k)
 &={\rm e}\!\left( {\rm Re}\!\left( {kAH\over Qc}\right)\right) 
\sum_{\scriptstyle s\mid c\atop\scriptstyle (s,Q)\sim 1}\mu_{\frak O}(s) 
\sum_{\scriptstyle d\bmod c{\frak O}\atop\scriptstyle d\equiv -AHP_s\bmod s{\frak O}} 
{\rm e}\!\left( {\rm Re}\!\left( {kd\over c}\right)\right)  = {}\cr
 &={\rm e}\!\left( {\rm Re}\!\left( {kAH\over Qc}\right)\right) 
\sum_{\scriptstyle s\mid c\atop\scriptstyle (s,Q)\sim 1}\mu_{\frak O}(s) 
{\rm e}\!\left( -{\rm Re}\!\left( {kAHP_s\over c}\right)\right)
\sum_{g\bmod (c/s){\frak O}} 
{\rm e}\!\left( {\rm Re}\!\left( {kg\over (c/s)}\right)\right) ,}$$
where, for $(s,Q)\sim 1$, we take $P_s$ to be a Gaussian integer satisfying 
$Q P_s\equiv 1\bmod s{\frak O}$. By the orthogonality relation (5.8), 
it follows from the above that 
$$c_q(b,h;k)
=\sum_{\scriptstyle s\mid c\atop{\scriptstyle (s,Q)\sim 1\atop\scriptstyle (c/s)\mid k}}
\mu_{\frak O}(s)\left| {c\over s}\right|^2 
{\rm e}\!\left( {\rm Re}\!\left( {kAH\bigl( 1-QP_s\bigr)\over Qc}\right)\right) 
=\sum_{\scriptstyle st=c\atop{\scriptstyle (s,Q)\sim 1\atop\scriptstyle t\mid k}}
\mu_{\frak O}(s) | t|^2 
{\rm e}\!\left( {\rm Re}\!\left( {(k/t)HAR_s\over Q}\right)\right) ,$$
where $R_s=\bigl( 1-Q P_s\bigr) /s\in{\frak O}$. 
Now we have $AB\equiv 1\bmod Q{\frak O}$, and (by construction) $s R_s\equiv 1\bmod Q{\frak O}$ 
also; so it follows that $(A R_s)(Bs)\equiv (AB)(s R_s)\equiv (1)(1)\equiv 1\bmod Q{\frak O}$, 
and we are therefore able to rewrite the last expression obtained for the value of $c_q(b,h;k)$ as: 
$$c_q(b,h;k)
=\sum_{\scriptstyle st=c\atop{\scriptstyle (s,Q)\sim 1\atop\scriptstyle t\mid k}}
\mu_{\frak O}(s) | t|^2 
{\rm e}\!\left( {\rm Re}\!\left( {kH\over t}\,{(Bs)^{*}\over Q}\right)\right) ,$$
where, for $(s,Q)\sim 1$, one has $(Bs)^{*}\in{\frak O}$ and 
$(Bs)(Bs)^{*}\equiv 1\bmod Q{\frak O}$ (this defining $(Bs)^{*}\bmod Q{\frak O}$). 

In light of the point noted below it, the equation~(5.34) follows from the result 
just obtained (by the substitution of $c/t$ for $s$, and the reversal 
of the substitutions recorded in (5.36)). 

Since the mapping $d\bmod q{\frak O}\mapsto d^{*}\bmod q{\frak O}$ 
is a permutation on the multiplicative group $({\frak O}/q{\frak O})^{*}$, 
the first equality in (5.35) follows  
immediately from the definition (1.3.6). The second equality in 
(5.35) becomes obvious when one compares (1.3.6) and (5.31) 
for $u=0$, $v=k$, $w=q$, $b=h=0$. To prove the final equality in (5.35), 
we begin by noting that if $b=h=0$ then (5.33) holds with $c=q$, 
so that by (5.34) one has  
$$c_q(0,0;k)={1\over 4}\sum_{t\mid (q,k)}\mu_{\frak O}\!\left( {q\over t}\right) 
|t|^2 {\rm e}(0)
={1\over 4}\sum_{t\mid (q,k)}\mu_{\frak O}\!\left( {q\over t}\right) 
|t|^2\;.$$
One may rewrite the last sum over $t$ 
by means of the substitution  $t=(q,k)/s\,$ (where $s\mid (q,k)$);    
given the properties (5.4), (5.5) of the function $\mu_{\frak O}$, 
this substitution shows that  
$$c_q(0,0;k)={1\over 4}\sum_{s\mid (q,k)}\mu_{\frak O}\!\left( {qs\over (q,k)}\right) 
{|(q,k)|^2 \over |s|^2} 
=|(q,k)|^2\,\mu_{\frak O}\!\left( {q\over (q,k)}\right) 
\!\!\sum_{\scriptstyle s\mid (q,k)\atop\scriptstyle\bigl( s\,,\,{q\over(q,k)}\bigr)\sim 1}
{\mu_{\frak O}(s)\over 4|s|^2}\;,$$
where, since ${\frak O}$ is the principal ideal domain ${\Bbb Z}[i]$,   
the last sum over $s$ is (given (5.4)) just what one 
obtains on multiplying out the product
over prime ideals $\varpi{\frak O}\subset{\frak O}$ 
appearing in the equation~(5.35)\ $\blacksquare$ 

\bigskip 

\goodbreak\proclaim Lemma~5.8 (a general analytic large sieve for ${\Bbb Z}[i]$). 
Let $c_n\in{\Bbb C}$ for all non-zero $n\in{\frak O}={\Bbb Z}[i]$; 
let  
$$S(\alpha,N)=\sum_{\scriptstyle n\in{\frak O}\atop\scriptstyle 0<|n|^2\leq N} 
c_n\,{\rm e}\left( {\rm Re}(\alpha n)\right)\qquad\qquad\hbox{($\alpha\in{\Bbb C}$, 
$N>0$);}\eqno(5.38)$$
and let $\alpha_r\in{\Bbb C}$ for all $r\in{\Bbb N}$. 
Then, for $R\in{\Bbb N}$, $N\geq 1$ and $1/2\geq\delta >0$, one has 
$$\sum_{1\leq r\leq R}\left| S\bigl(\alpha_r ,N\bigr)\right|^2 
\leq 16\,M(\delta,R)\left( 2N+{\delta}^{-1}\right) 
\left\|{\bf c}_N\right\|_2^2\;,\eqno(5.39)$$
where 
$$M(\delta,R)=\max_{1\leq r\leq R}
\,\left|\left\{ p\in{\Bbb N} : 1\leq p\leq R\ {\rm and}\ 
\bigl\|\alpha_p - \alpha_r\bigr\|^2<\delta\right\}\right|\;,\eqno(5.40)$$
and where the definition of $\|{\bf c}_N\|_2$ is as 
indicated by the equation~(1.2.11), in Theorem~2.
\medskip

\goodbreak 
\noindent{\bf Proof.}\quad 
Let $R\in{\Bbb N}$ and $N\geq 1$.  
For $r=1,\ldots ,R$, we have (by (5.38)):
$$S\left(\alpha_r,N\right) 
=\,\sum\sum_{\!\!\!\!\!\!\!\!\!\!\!\!n_1,n_2\in{\Bbb Z}}
\,c\bigl(n_1,n_2\bigr)\,{\rm e}\!\left( n_1 x_1^{(r)}+n_2 x_2^{(r)}\right)\;,\eqno(5.41)$$
where 
$$\left( x_1^{(r)}\,,\,x_2^{(r)}\right) 
=\Bigl( {\rm Re}\bigl(\alpha_r\bigr)\,,\,{\rm Im}\bigl(\alpha_r\bigr)\Bigr) 
={\bf x}^{(r)}\in{\Bbb R}^2\qquad\ \hbox{(say)}\eqno(5.42)$$
and 
$$c\bigl( n_1,n_2\bigr) =\cases{c_{n_1-i n_2} &if $0<n_1^2+n_2^2\leq N$ \cr 
0 &otherwise.}\eqno(5.43)$$

Let $\delta\in (0,1/2]$; and let $M(\delta ,R)$ be given by (5.40). 
Supposing firstly that $M(\delta,R)=1$, we have  
$$\left\|\alpha_p -\alpha_r\right\|^2\geq\delta\qquad\quad\hbox{for $\quad 1\leq p<r\leq R$.}$$
In this case it follows, by (5.42) and (5.7), that we have
$$\max_{j=1,2}\,\bigl\| x_j^{(p)}-x_j^{(r)}\bigr\|\geq \left( \delta /2\right)^{1/2}\qquad\quad  
\hbox{for $\quad 1\leq p<r\leq R$,}$$
which is a special case 
($k=2$, $\delta_1=\delta_2=(\delta /2)^{1/2}\in (0,1/2]$)
of [9, Theorem~1, Condition~(3)]. Moreover, in the sum 
on the right-hand side of (5.41)$\,$ (which is similar in 
form to that in the case $k=2$ of [9, Theorem~1, Definition (1)]$\,$), 
the variables $n_1,n_2$ are effectively constrained (given (5.43)) 
to range within the interval $[-N^{1/2},N^{1/2}]$; and so, for arbitrary 
$N_1,N_2>2 N^{1/2}$, the sums in (5.41) satisfy the special case  
$k=2$, $M_1=-N_1 /2$, $M_2=-N_2 /2$ of [9, Theorem~1, Condition~(2)].
As a consequence of these observations, it follows by [9, Theorem~1] 
that, when $M(\delta , R)=1$, one has 
$$\sum_{1\leq r\leq R}\left| S\bigl(\alpha_r,N\bigr)\right|^2 
\leq B_1 B_2\,\sum\sum_{\!\!\!\!\!\!\!\!\!\!\!\!n_1,n_2\in{\Bbb Z}}
\,\left| c\bigl(n_1,n_2\bigr)\right|^2 
=B_1 B_2 \left\|{\bf c}_N\right\|_2^2\;,\eqno(5.44)$$
where 
$$B_1=B_2=\left( \bigl( 2 N^{1/2}\bigr)^{1/2}+(\delta /2)^{-1/4}\right)^2 
\leq 2\left( 2 N^{1/2}+(\delta /2)^{-1/2}\right)\;.$$
In this case, since $(2N^{1/2})^2+(\delta /2)^{-1}=2(2N+\delta^{-1})$, the result (5.39) 
of the lemma follows from (5.44) by way of the same inequality, $(a+b)^2\leq 2 a^2 +2 b^2$,  
that justifies the bound just given for $B_1$ and $B_2$. 

The above completes the proof in the case where $M(\delta ,R)=1$. 
Suppose now that $M(\delta , R)>1$. Then certainly we have $R>1$ also. 
Without loss of generality we may assume that 
$$\left| S\bigl(\alpha_R,N\bigr)\right|\geq\left| S\bigl(\alpha_p,N\bigr)\right|\qquad\quad  
\hbox{for $\quad 1\leq p\leq R$.}\eqno(5.45)$$
We may also renumber $\alpha_1,\ldots ,\alpha_{R-1}$ so that, for some non-negative 
integer $R(2)<R(1)=R$, we have: 
$$\left\|\alpha_r -\alpha_{R(1)}\right\|<\delta^{1/2}\quad\ \hbox{for $\ \,R(1)\geq r >R(2)$;}
\qquad\qquad   
\left\|\alpha_r -\alpha_{R(1)}\right\|\geq\delta^{1/2}\quad\ \hbox{for $\ \,R(2)\geq r\geq 1$.}$$
By (5.40), the above number $R(2)$ must satisfy $R-R(2)\leq M(\delta,R)$. 
Hence, and by (5.45), 
$$\sum_{1\leq r\leq R}\left| S\bigl(\alpha_r ,N\bigr)\right|^2 
\leq M(\delta ,R)\left| S\bigl(\alpha_R,N\bigr)\right|^2 
+\sum_{1\leq r\leq R(2)}\left| S\bigl(\alpha_r ,N\bigr)\right|^2\;.$$
Moreover, if $R(2)\neq 0$ then a similar upper bound 
can be obtained for the last sum here (the sum over $r=1,\ldots ,R(2)$). 
Hence, by iteration of the same procedure, one arrives at a bound of the form  
$$\sum_{1\leq r\leq R}\left| S\bigl(\alpha_r ,N\bigr)\right|^2 
\leq\sum_{j=1}^J M_j\left| S\bigl(\alpha_{R(j)},N\bigr)\right|^2 ,$$
where $J\geq 1$, and the sequence of integers  
$R(1),\ldots ,R(J)$ is strictly decreasing, with $R=R(1)\geq R(J)>0$ (and $R(J+1)=0$);  
while, given the nature of our iterative procedure (and the definition (5.40)),
$$\left\|\alpha_{R(p)}-\alpha_{R(r)}\right\|\geq\delta^{1/2}\qquad\quad  
\hbox{for $\quad 1\leq p<r\leq J$,}\eqno(5.46)$$
and the sequence of integers $M_1,\ldots M_j$ is non-increasing, 
with $M_1=M(\delta ,R(1))=M(\delta ,R)\,$ 
(subject to a suitable initial renumbering of $\alpha_1,\ldots ,\alpha_R$, 
prior to the start of our 
iterative procedure, one would have here $M_j=M(\delta,R(j))$ for $j=1,\ldots ,J$).
We deduce that 
$$\sum_{1\leq r\leq R}\left| S\bigl(\alpha_r ,N\bigr)\right|^2 
\leq M(\delta ,R)\sum_{j=1}^J \left| S\bigl(\alpha_{R(j)},N\bigr)\right|^2\;.\eqno(5.47)$$

In (5.46)-(5.47), we may put $J=R'$ (say) and may also, for $r=1,\ldots ,R'$, put  
$\alpha_{r}'=\alpha_{R(r)}$. Hence, given the bound (5.46), the previously established  
case $M(\delta ,R)=1$ of (5.39)-(5.40) shows that   
$$\sum_{j=1}^J \left| S\bigl(\alpha_{R(j)},N\bigr)\right|^2
\leq 16\left( 2N+\delta^{-1}\right)\left\|{\bf c}_N\right\|_2^2\;.$$
By this bound and that in (5.47), the proof of the results (5.39)-(5.40) is complete
\ $\blacksquare$

\bigskip

\goodbreak 
\noindent{\bf Remark.}\quad 
The above lemma is slightly more elaborate than we actually require: for in 
this paper it is used only to establish the next lemma, and we could do as much  
with just the case $M(\delta ,R)=1$ of Lemma~5.8. 

\bigskip

\goodbreak 
\proclaim Lemma~5.9 (a special analytic large sieve for ${\Bbb Z}[i]$). Let $c_n\in{\Bbb C}$ 
for all non-zero $n\in{\frak O}={\Bbb Z}[i]$. Then, for $Q,N\geq 1$ and 
$d\in{\frak O}-\{ 0\}$, one has 
$$\sum_{\scriptstyle 0<|q|^2\leq Q\atop\scriptstyle q\equiv 0\bmod d{\frak O}}
\sum_{\scriptstyle a\bmod q{\frak O}\atop\scriptstyle (a,q)\sim 1} 
\Biggl|\sum_{0<|n|^2\leq N} 
c_n\,{\rm e}\left( {\rm Re}\left( {an\over q}\right)\right)\Biggr|^2 
\leq 64\left( 2N+{Q^2\over |d|^2}\right)\left\|{\bf c}_N\right\|_2^2\;,\eqno(5.48)$$
where $q,a,n$ are Gaussian integer variables of summation, and where 
$\|{\bf c}_N\|_2$ is as (1.2.11) indicates.

\medskip

\goodbreak 
\noindent{\bf Proof.}\quad  The sum on the left-hand side of (5.48) may 
be written as $\sum_{1\leq r\leq R} |S(\alpha_r,N)|^2$, where 
$S(\alpha)$ is given by the equation~(5.38) of Lemma 5.8,
and where, to each $r\in\{1,\ldots ,R\}$ there corresponds a pair 
$(q_r, a_r\bmod q_r{\frak O})$ with $q=q_r$, $a=a_r$
satisfying the conditions of summation in (5.48), 
and with $a_r\equiv\alpha_r q_r\bmod q_r{\frak O}\,$   
(this correspondence $r\mapsto (q_r,a_r\bmod q_r{\frak O})$ being one-to-one).
Therefore, assuming that $R>0$ and $|d|^2/Q^2\leq 1/2$, it will suffice 
to show that the relevant sequence $\alpha_1\ldots ,\alpha_R$ is such 
that, when $M(\delta ,R)$ is as defined in Lemma~5.8, one has 
$M(|d|^2/Q^2,R)\leq 4$: for the bound (5.48) will, in that case, be 
implied by the result (5.39) of Lemma~5.8. 

In order to show that $M(|d|^2/Q^2,R)\leq 4$, we first note that if 
$1\leq p,r\leq R$ and $\|\alpha_p-\alpha_r\|\neq 0$ then, since 
$d\mid q_p,q_r$, one will have 
$$0<\left\|\alpha_p-\alpha_r\right\|
=\left\| {\alpha_p q_p\over q_p}-{\alpha_r q_r\over q_r}\right\|
=\left\| {a_p\over q_p}-{a_r\over q_r}\right\|
=\left\| {a_p q_r-a_r q_p\over q_p q_r}\right\|={|kd|\over\bigl| q_p q_r\bigr|}\eqno(5.49)$$
for some $k\in{\frak O}-\{ 0\}$, 
and hence $\|\alpha_p-\alpha_r\|^2\geq |d|^2 /|q_p q_r|^2\geq |d|^2 /Q^2$. 
By this and (5.40) it follows that, for some $r\in\{ 1,\ldots ,R\}$, one has 
$$M(|d|^2/Q^2,R)=\bigl|\bigl\{ p\in{\Bbb N} : 1\leq p\leq R\ {\rm and}\ 
\|\alpha_p -\alpha_r\|=0\}|\;.$$   
Moreover, the first three equalities in (5.49) show that 
$\|\alpha_p -\alpha_r\|=0$ if and only if $a_p q_r\equiv a_r q_p\bmod q_p q_r{\frak O}$, 
and so only if $a_p q_r\equiv 0\bmod q_p{\frak O}$ 
and $a_r q_p\equiv 0\bmod q_r{\frak O}$.  Since the conditions of summation in 
(5.48) are satisfied when either $q=q_p$ and $a=a_p$, or $q=q_r$ and $a=a_r$, 
we have $(a_p,q_p)\sim 1$ and $(a_r,q_r)\sim 1$.  The 
simultaneous congruences $a_p q_r\equiv 0\bmod q_p{\frak O}$ 
and $a_r q_p\equiv 0\bmod q_r{\frak O}$ therefore imply that we have 
both $q_p\mid q_r$ and $q_r\mid q_p$, and so $q_p\sim q_r$. 
It follows that $\|\alpha_p -\alpha_r\|=0$ if and only if, for some unit $\epsilon\in{\frak O}^{*}$, 
one has $q_p=\epsilon q_r$ and $a_p q_r\equiv a_r \epsilon q_r\bmod q_p q_r{\frak O}$. 
Since the last congruence implies $a_p\equiv a_r \epsilon\bmod q_p{\frak O}$, 
we may conclude that $M(|d|^2/Q^2,R)=|{\frak O}^{*}|=4$: as noted above, this proves  
the lemma in cases where $R>0$ and $|d|^2/Q^2\leq 1/2$. 

To complete the proof we observe firstly that (5.48) is essentially trivial in 
cases where $R\leq 16$: 
for in such cases the sum on the left-hand side of (5.48) 
is either empty (and hence equal to zero) or,  
for some $r\in\{1,2,\ldots ,16\}$, 
is less than or equal to $16|S(\alpha_r,N)|^2$, where, 
by the Cauchy-Schwarz inequality, one has 
$$|S(\alpha_r,N)|^2\leq |\{ n\in{\frak O} : 0<|n|^2\leq N\}|\,\|{\bf c}_N\|_2^2\leq 
\bigl( (2 N^{1/2}+1)^2-1\bigr)\,\|{\bf c}_N\|_2^2\leq 8N\,\|{\bf c}_N\|_2^2$$
(when $N\geq 1$). To complete the proof we note that the conditions of 
summation in (5.48) imply $|q|^4\leq Q^2$ and $2|d|^2\leq 2|q|^2$, so 
that if $|d|^2/Q^2>1/2$ (implying $Q^2<2 |d|^2$) then one will have 
$|q|^4<2 |q|^2$, for all $q$ in the sum, and hence $R\leq 4<16$
(the summation over $q$ being  
restricted to $q\in{\frak O}^{*}$)\ $\blacksquare$

\bigskip 
\bigskip 

\goodbreak\centerline{\bf \S 6. An Elementary Bound for a Sum of Kloosterman Sums.}

\bigskip 

In this section we consider a sum 
$${\cal R}
=\sum_{p\neq 0}\theta_p|p|^{-2}\sum_{q\neq 0}|q|^{-2}\sum_h\phi_h\sum_k\sum_{\ell}
\,S(hk,\ell;pq)\,\varphi(h,k,\ell,p,q)\,\Upsilon_{\ell}\;,\eqno(6.1)$$
where $S(u,v;w)$ is the `simple Kloosterman sum' defined in (1.3.6); 
and where the  summation is over the points $(p,q,h,k,\ell)\in{\frak O}^5$ 
with $pq\neq 0$. 

We suppose that the function $\varphi$ has domain ${\Bbb C}^5$, and is complex valued;   
and we assume that the function 
$\Phi : {\Bbb R}^{10}\rightarrow{\Bbb C}$  
given by $\Phi\bigl( x_1,\ldots ,x_{10}\bigr)=\varphi(x_1+ix_2, x_3+ix_4,\ldots ,x_9+ix_{10}\bigr)$ 
(${\bf x}\in{\Bbb R}^{10}$) is such that all its partial derivatives 
(of any given order) are defined and continuous at all points of ${\Bbb R}^{10}$. 
The function $\varphi : {\Bbb C}^5\rightarrow{\Bbb C}$ might therefore be termed `smooth'. 
We suppose moreover that, for some given $H,K,L,P,Q\geq 1$ and some given $\delta >0$,  
one has 
$$\varphi\bigl( z_1 , z_2 , z_3 , z_4 , z_5\bigr) =0\qquad{\rm unless}\qquad 
\left({\bigl| z_1\bigr|^2\over H} ,  {\bigl| z_2\bigr|^2\over K} , 
{\bigl| z_3\bigr|^2\over L} , {\bigl| z_4\bigr|^2\over P} , 
{\bigl| z_5\bigr|^2\over Q}\right)\in\left({1\over 2}\,,\,1\right)^5 ,\eqno(6.2)$$
and, for ${\bf j},{\bf k}\in({\Bbb N}\cup\{ 0\})^5$ and 
all ${\bf x},{\bf y}\in{\Bbb R}^5$ such that 
$x_h+iy_h\neq 0$ for $h=1,\ldots ,5$, 
$${\partial^{j_1+\cdots +j_5+k_1+\cdots +k_5}\over 
\partial x_1^{j_1}\cdots\partial x_5^{j_5} 
\partial y_1^{k_1}\cdots\partial y_5^{k_5}}\,\varphi\bigl( x_1+iy_1 , \ldots , x_5+iy_5\bigr) 
\ll_{{\bf j},{\bf k}} 
\prod_{h=1}^5\left(\delta\bigl| x_h + i y_h\bigr|\right)^{-(j_h+k_h)}\;.\eqno(6.3)$$

\bigskip 

\goodbreak 
\noindent{\bf Remark.} By (6.2) and the hypothesis of `smoothness', 
the function $\varphi$ lies in the Schwartz space 
${\cal S}\bigl( {\Bbb C}^5\bigr)$.  

\bigskip 

As for the coefficients 
$\theta_p$, $\phi_h$ and $\Upsilon_{\ell}$ in (6.1), 
we suppose that, for $p\in{\frak O}-\{ 0\}$ and $h,\ell\in{\frak O}$,  
these coefficients satisfy   
$$\theta_p,\phi_h,\Upsilon_{\ell}\in{\Bbb C}\;,\qquad 
\left|\theta_p\right|\ll 1\quad\ {\rm and}\quad\ \left|\phi_h\right|\ll 1\;;\eqno(6.4)$$
$$\theta_p=0\ {\rm unless}\ {P\over 2}<|p|^2\leq P\;;\quad 
\phi_h=0\ {\rm unless}\ {H\over 2}<|h|^2\leq H\;;\quad{\rm and}\quad   
\Upsilon_{\ell}=0\ {\rm unless}\ {L\over 2}<|\ell|^2\leq L\;.\eqno(6.5)$$ 

Given the above hypotheses, and given that (1.3.6) and (5.9) imply the bounds      
$$|S(u,v;w)|\leq\left|({\frak O}/w{\frak O})^{*}\right|
\leq\left|{\frak O}/w{\frak O}\right|=|w|^2\qquad\qquad  
\hbox{($u,v\in{\frak O}$ and $w\in{\frak O}-\{ 0\}$),}\eqno(6.6)$$ 
it is trivially the case that 
$$|{\cal R}|\leq\sum_p\sum_q\sum_h\sum_k\sum_{\ell}|\theta_p\phi_h\Upsilon_{\ell} 
\,\varphi(h,k,\ell,p,q)|\;.$$
Therefore, and by (6.2)-(6.4) and the Cauchy-Schwarz inequality, we have: 
$${\cal R}\ll PQHK
\Biggl( L\sum_{{L\over 2}<|\ell|^2\leq L} \left|\Upsilon_{\ell}\right|^2\Biggr)^{\!\!1/2}\;.
\eqno(6.7)$$
Our goal in this section (realised in Lemma~6.3) is the 
proof of a particular improvement of this last, essentially trivial, preliminary 
upper bound for $|{\cal R}|$. 
The implicit constants in both (6.7) and the result (6.61) of Lemma~6.3 
do of course depend on the implicit constants in the bounds of (6.3) and (6.4).

\bigskip 

\goodbreak 
\noindent{\bf Remark.}\quad 
In some of the proofs which follow (both in this section, and subsequently) 
we make use of the bounds 
$$\sum_{\scriptstyle m\neq 0} |m|^{-(2+\varepsilon)}\ll_{\varepsilon} 1\qquad 
{\rm and}\qquad 
\sum_{d\mid n} 1\ll_{\varepsilon} |n|^{2\varepsilon}\qquad\qquad   
\hbox{($0\neq n\in{\frak O}$),}$$
where, as usual, $\varepsilon$ denotes an arbitrary positive constant, 
and $m$ and $d$ are Gaussian integer valued variables of summation. 
Since these elementary bounds should be well known, we make no comment 
when using them. 

\bigskip

\goodbreak 
\proclaim Lemma~6.1. Let $H,K,L,P,Q\geq 1$, $\varphi\in{\cal S}\bigl({\Bbb C}^5\bigr)$ 
and $\delta >0$ be such that the conditions (6.2) and (6.3) are satisfied. Suppose that 
$$\max\{ HK\,,\,L\}\ll Q\ll (HKL)^{2/3}\;,\qquad 
H\ll K\qquad\hbox{and}\qquad HL\ll PQ\;.\eqno(6.8)$$
For $p\in{\frak O}-\{ 0\}$ and $h,\ell\in{\frak O}$, 
let $\theta_p,\phi_h,\Upsilon_{\ell}\in{\Bbb C}$ 
satisfy (6.4) and (6.5). Let ${\cal R}$ be given by (6.1). 
Let $\varepsilon >0$; and let 
$$E=(PQ)^{\varepsilon}\bigl( 1+\delta^{-1}\bigr)\;.$$  
Then either 
$${\cal R}\ll_{\varepsilon} 
\left( E^8 K^{-1} PQ+(PQ)^{\varepsilon} (HKL)^{1/2}\right) 
\biggl( HK\sum_{\ell}\bigl|\Upsilon_{\ell}\bigr|^2\biggr)^{1/2}\;,\eqno(6.9)$$ 
or else: 
$$\delta^2 K>16(PQ)^{\varepsilon}\;,\qquad\  
Q>4 E^2 H\eqno(6.10)$$ 
and, for some non-zero Gaussian integers $w,r,s,c,t,k,q$ satisfying  
$$0<|w|^2\leq H\;,\quad\ c\mid w\;,\quad\ t\mid c\;,\quad\ 
w=rs\;,\quad\ s\mid q\quad\ {\rm and}\quad\ {Q\over 2}<|q|^2\leq Q\;,\eqno(6.11)$$
one has 
$${\cal R}\ll_{\varepsilon} 
(PQ)^{\varepsilon} K |t|^2 \left| {\cal E}(w,c,t,r,s;k,q)\right|\;,\eqno(6.12)$$
with  
$$\eqalign{ 
{\cal E}(w,c,t,r,s;k,q) 
 &=\sum_{\scriptstyle 0<|a|^2\leq A\atop\scriptstyle (a,w)\sim c} 
{1\over |a/t|^2} \sum_{\scriptstyle h\in{\frak O}\atop\scriptstyle w\mid h}\phi_h 
\sum_{\scriptstyle \ell\in{\frak O}\atop\scriptstyle t\mid\ell}\Upsilon_{\ell} 
\sum_{\scriptstyle 0\neq p\in{\frak O}\atop\scriptstyle (p,w)\sim r} 
{\theta_p\over |p/r|^2}\,{\rm e}\!\left( {\rm Re}\!\left( 
{ak\over pq}+{h\ell\over apq}\right)\right)\times \cr 
 &\qquad\times \varphi(h,k,\ell,p,q)
\sum_{\scriptstyle 0\neq b\in{\frak O}\atop\scriptstyle 
|b|^2\leq V(sa/t)} e\!\left( {\rm Re}\!\left( {(q/s)b\over (a/t)}\right)\right) 
S\!\left( {h\over w}\,{\ell\over t}\,\left( {p\over r}\right)^{*},\,b\,;\,{a\over t}\right) 
,} \eqno(6.13)$$ 
where (with $S(u,v;w)$ being given by (1.3.6)) 
the factor $(p/r)^{*}$ has the meaning explained under the heading `Number-Theoretic Notation' 
in Subsection~1.5, while 
$$A={(PQ)^{1+\varepsilon}\over\delta^2 K}<{PQ\over 16}\qquad\quad{\rm and}\qquad\quad    
 V(z)={E^2 |z|^2\over Q}<{|z|^2\over 4H}\quad\ \hbox{for $\ \,z\in{\Bbb C}^{*}$.}\eqno(6.14)$$

\medskip

\goodbreak 
\noindent{\bf Proof.} 
By the trivial bound (6.7),  
$${\cal R}\ll \left( H K^3 L\right)^{1/2} K^{-1} PQ  
\Biggl( HK\sum_{{L\over 2}<|\ell|^2\leq L} \left|\Upsilon_{\ell}\right|^2\Biggr)^{\!\!1/2}\;,$$
where, by (6.8), one has $H K^3 L\ll H K^3 (HK)^2\ll K^8$. 
Consequently, subject to the hypotheses of the lemma, 
the bound (6.9) is obtained whenever $K\ll E^2$. 
Moreover, one has $K\ll E^2$ if at least one of the inequalities in 
(6.10) is false: for if $Q\leq 4 E^2 H$ then, by (6.8), one has $HK\ll E^2 H$, and so 
$K\ll E^2$; while if $\delta^2 K\leq 16 (PQ)^{\varepsilon}$ then 
$K\ll (PQ)^{\varepsilon}\delta^{-2}\leq (PQ)^{2\varepsilon}\delta^{-2}<E^2$. 
We may therefore suppose henceforth that the inequalities in (6.10) are satisfied:  
for otherwise we have $K\ll E^2$, and so (as observed above) obtain the bound (6.9). 

To complete this proof it will suffice to deduce 
(assuming the conditions in (6.10)) that either the bound 
(6.12) holds, for some $w,r,s,c,t,k,q\in{\frak O}-\{ 0\}$ satisfying (6.11),  
or else one has (6.9). This will be achieved in two steps, by applying the 
results (5.30) and (5.29) of Lemma~5.6. 

For our first application of Lemma~5.6, we suppose that 
$h,\ell,p,q\in{\frak O}-\{ 0\}$ are given and take 
$f : {\Bbb C}\rightarrow{\Bbb C}$ to be the function 
$z\mapsto\varphi(h,z,\ell,p,q)$. Since $\varphi\in{\cal S}\bigl({\Bbb C}^5\bigr)$, 
we have $f\in{\cal S}({\Bbb C})$. In order that Lemma~5.6 may be applied 
it will suffice that the function $f$ satisfies, 
for some $\Delta,\Omega_1>0$, and some $C>1$,  
the case $n=1$ of conditions~(5.24), (5.25) in Lemma~5.4. 
Taking $\Omega_1 =K$ and $C=2$, the case $n=1$ of (5.25) follows immediately 
from (6.2). Moreover, by (6.2) and (6.3), one has 
$${\partial^{j+k}\over\partial x^j\partial y^k}\,f(x+iy) 
=O_{j,k}\left( (\delta |x+iy|)^{-(j+k)}\right) 
\ll_{j+k} (\delta |x+iy|)^{-(j+k)}\;,$$ 
for all $j,k\in{\Bbb N}\cup\{ 0\}$, and all $x,y\in{\Bbb R}$ such that $x+iy\neq 0$. 
Since $f\in{\cal S}({\Bbb C})$, we therefore have 
(with ${\cal L}_1$ defined as in (5.12), $z\in{\Bbb C}$, 
$x={\rm Re}(z)$ and $y={\rm Im}(z)$): 
$${\cal L}_1^j f(z)=(-1)^j\left( {\partial^2\over\partial x^2}
+{\partial^2\over\partial y^2}\right)^{\!j}\!f(z) 
=\pm\sum_{r=0}^j \pmatrix{j\cr r} 
{\partial^{2j}\over\partial x^{2r}\partial y^{2j-2r}}\,f(x+iy)
\ll_{2j}\ \sum_{r=0}^j \pmatrix{j\cr r} (\delta |x+iy|)^{-(2j)}\;,$$
which implies that $f$ satisfies the case $n=1$ of (5.24) if 
one takes there $\Delta =\delta^2$. 
Therefore, by the case $d=1$, $B=(PQ)^{1+\varepsilon}/|pq|^2$ 
of the result (5.30) of Lemma~5.6, we have, for $j\geq 2$, 
$$\sum_k\varphi(h,k,\ell,p,q) S(hk,\ell;pq) 
=\sum_{\scriptstyle a\bmod pq{\frak O}\atop\scriptstyle 
\left\| {\textstyle{a\over pq}}\right\|^2\leq {\textstyle{A\over |pq|^2}}}
\!\!c_{pq}(a,h;\ell)\sum_k\varphi(h,k,\ell,p,q)\,{\rm e}\!\left( {\rm Re}\!\left( 
{ak\over pq}\right)\right) 
+O_j\left( {\delta^{-2}|pq|^{2j+4}\over (PQ)^{(1+\varepsilon)j}}\right) ,$$
with $A$ as in (6.14). By (6.2), the last $O$-term is zero unless 
$|pq|^2\in(PQ/4,PQ)$, and so may be replaced by 
$O_j\bigl(\delta^{-2} (PQ)^{2-j\varepsilon}\bigr)$.  
In the above, both $\| a/(pq)\|$ and the factor 
${\rm e}\bigl( {\rm Re}(ak/(pq))\bigr)$ are periodic, 
$\bmod\ pq{\frak O}$, as functions of the variable~$a$; by (5.31),  
so too is the factor $c_{pq}(a,h;\ell)$.  
Hence, and since  
$\| a/(pq)\|=|(a-pqm)/(pq)|=|a-pqm|/|pq|$ for some $m\in{\frak O}$, 
it may be assumed that $\| a/(pq)\|=|a|/|pq|$ in the above sum. 
Moreover, by (6.2) and (6.14), we have $\varphi(h,k,\ell,p,q)\neq 0$ 
only if $|pq|>2A^{1/2}$, and so only if no two distinct elements of 
the set $\bigl\{ a\in{\frak O} : |a|^2\leq A\bigr\}$ are congruent 
to one another, $\bmod\ pq{\frak O}$. Therefore the conditions of 
summation on the right-hand side of the above equation may be 
simplified to just: $a\in{\frak O}$ and $|a|^2\leq A$. Hence, by 
taking $j=[ 2/\varepsilon ]+1$, we obtain 
$$\sum_k\varphi(h,k,\ell,p,q) S(hk,\ell;pq) 
=\sum_{|a|^2\leq A}
\!\!c_{pq}(a,h;\ell)\sum_k\varphi(h,k,\ell,p,q)\,{\rm e}\!\left( {\rm Re}\!\left( 
{ak\over pq}\right)\right) 
+O_{\varepsilon}\left(\delta^{-2}\right) ,\eqno(6.15)$$
where, by the result (5.32) of Lemma~5.7, one has $c_{pq}(a,h;\ell)=0$ unless 
$(a,pq)\sim (h,pq)$. Consequently, given (6.2) and the second inequality in (6.10), 
one has 
$c_{pq}(a,h;\ell)\varphi(h,k,\ell,p,q)\neq 0$ only if 
$$|(a,pq)|^2=|(h,pq)|^2\leq |h|^2<H<{Q\over 4 E^2}<{Q\over 4}\leq {PQ\over 4}<|pq|^2\;,$$
and so only if 
$a\not\equiv 0\bmod pq{\frak O}$. It is therefore effectively an implicit 
condition of the summation over $a$ in (6.15) that $a\neq 0$. 

By (6.15) and the observations subsequent to it, and by (6.2), (6.3) (for 
${\bf j}={\bf k}={\bf 0}$), (6.4), (6.5) and the Cauchy-Schwarz inequality, 
we deduce that 
$$\eqalign{
{\cal R} 
&={\cal R}'+O_{\varepsilon}\Biggl( 
\delta^{-2} (PQ)^{-1}\sum_h\sum_{\ell}\sum_p 
\left|\phi_h\Upsilon_{\ell}\theta_p\right|\sum_{{\textstyle{Q\over 2}}<|q|^2<Q} 1
\Biggr) = {}\cr 
&={\cal R}'+O_{\varepsilon}\Biggl( 
\delta^{-2} P^{-1}\sum_{{\textstyle{H\over 2}}<|h|^2<H}
\ \sum_{{\textstyle{P\over 2}}<|p|^2<P}\ \sum_{{\textstyle{L\over 2}}<|\ell|^2<L}
\left|\Upsilon_{\ell}\right|\Biggr)  
={\cal R}'+O_{\varepsilon}\left( \delta^{-2}H L^{1/2}\biggl( 
\sum_{\ell}\left|\Upsilon_{\ell}\right|^2\biggr)^{1/2}\right) ,}$$
where 
$${\cal R}'
=\sum_{p\neq 0}\theta_p |p|^{-2}\sum_{q\neq 0} |q|^{-2}\sum_{\ell} \Upsilon_{\ell} 
\sum_h \phi_h 
\!\!\!\!\sum_{\scriptstyle 0<|a|^2\leq A\atop\scriptstyle (a,pq)\sim (h,pq)} 
\!\!\!\!c_{pq}(a,h;\ell) 
\sum_k \varphi(h,k,\ell,p,q) 
{\rm e}\!\left( {\rm Re}\!\left( {ak\over pq}\right)\right) .\eqno(6.16)$$
Moreover, by (6.10) we have $\delta^{-2} H L^{1/2}<HKL^{1/2}=(HKL)^{1/2}(HK)^{1/2}\leq 
(PQ)^{\varepsilon}(HKL)^{1/2}(HK)^{1/2}$ 
in the above; so it follows that either 
$$|{\cal R}|\leq 2\left|{\cal R}'\right|\;,\eqno(6.17)$$
or else the bound (6.9) holds. In the latter case we have nothing more to prove: 
we may therefore assume henceforth that the inequality (6.17) is satisfied. 

In the sum on the right-hand side of the equation (6.16) one has 
$(a,pq)\sim (h,pq)\sim c\,$ (say), with $c\in{\frak O}-\{ 0\}$ 
dependent upon $p$, $q$ and $h$.    
After grouping together summands in (6.16) corresponding to the same `$c$', 
we may apply the result (5.33)-(5.34) of Lemma~5.7, so as to obtain: 
$${\cal R}'={1\over 16}\sum_{c\neq 0}\sum_{t\mid c} 
\mu_{\frak O}\left( {c\over t}\right) |t|^2 {\cal R}'(c,t)\;,\eqno(6.18)$$ 
with   
$$\eqalignno{
{\cal R}'(c,t) 
 &=\sum_{p\neq 0}\sum_{\scriptstyle q\neq 0\atop\!\!\!\!\!\!\!\!\!\!\!\!{\scriptstyle 
c\mid (pq)}} \theta_p |pq|^{-2}
\sum_{\scriptstyle \ell\atop\scriptstyle t\mid\ell} \Upsilon_{\ell} 
\sum_{\scriptstyle h\atop\scriptstyle (h,pq)\sim c} \phi_h 
\!\!\!\!\sum_{\scriptstyle 0<|a|^2\leq A\atop\scriptstyle c\mid a} 
\!\!\!\!{\rm e}\!\left( {\rm Re}\!\left( {h\over c}\,{\ell\over t}\,{(a/t)^{*}\over (pq/c)} 
\right)\right)  
\sum_k \varphi(h,k,\ell,p,q) 
{\rm e}\!\left( {\rm Re}\!\left( {ak\over pq}\right)\right) = {}\cr 
 &=\sum_k 
\sum_{\scriptstyle 0<|a|^2\leq A\atop\scriptstyle c\mid a} 
\sum_{\scriptstyle h\atop\scriptstyle c\mid h} \phi_h 
\sum_{\scriptstyle \ell\atop\scriptstyle t\mid\ell} \Upsilon_{\ell} 
\sum_{p\neq 0}\sum_{\scriptstyle q\neq 0\atop\!\!\!\!\!\!\!\!\!\!\!\!\!\!\!{\scriptstyle 
(pq,h)\sim c}} \theta_p |pq|^{-2}
\varphi(h,k,\ell,p,q) 
{\rm e}\!\left( {\rm Re}\!\left( {ak\over pq} 
+{h\over c}\,{\ell\over t}\,{(a/t)^{*}\over (pq/c)}\right)\right) , &(6.19)}$$
where, by the definition of $(a/t)^{*}\bmod (pq/c){\frak O}$,   
it is an implicit condition of summation that    
$(pq/c,a/t)\sim 1$.
 
The final step in this proof is essentially Poisson summation with respect to the 
variable $q$ (by which the innermost sum in (6.19) is indexed). 
As things stand, in (6.19), the explicit condition $(pq,h)\sim c$ is 
an obstacle to the efficient implementation of Poisson summation with 
respect to $q$. Our (quite standard) solution for this difficulty is 
to note that, by (5.6), the restriction of summation to pairs 
$p,q$ satisfying $(pq,h)\sim c$ is identical in effect 
to the multiplication of all terms by the supplementary `weight' factor: 
$${1\over 4}\sum_{\scriptstyle d\atop\scriptstyle (cd)\mid (pq,h)} \mu_{\frak O}(d) 
=\cases{1 &if $(pq,h)\sim c$, \cr 0 &otherwise.}$$
This enables us to deduce from (6.19) that, for $0\neq c\in{\frak O}$ and $t\mid c$, one has 
$$\eqalignno{
{\cal R}'(c,t) 
 &={1\over 4}\sum_k 
\!\sum_{\scriptstyle 0<|a|^2\leq A\atop\scriptstyle c\mid a} 
\sum_{\scriptstyle h\atop\scriptstyle c\mid h} \phi_h 
\sum_{\scriptstyle \ell\atop\scriptstyle t\mid\ell} \Upsilon_{\ell} 
\sum_{p\neq 0}\sum_{q\neq 0} {\theta_p\over |pq|^2}
\,\varphi(h,k,\ell,p,q) 
{\rm e}\!\left( {\rm Re}\!\left( {ak\over pq} 
+{h\over c}\,{\ell\over t}\,{(a/t)^{*}\over (pq/c)}\right)\right) 
\!\!\!\!\!\sum_{\scriptstyle d\atop\scriptstyle (cd)\mid (pq,h)}\!\!\!\!\mu_{\frak O}(d) ={}\cr 
 &={1\over 4}\sum_{\scriptstyle d\neq 0\atop\scriptstyle (d,c/t)\sim 1} 
\mu_{\frak O}(d) {\cal R}''(c,t,d)\;, &(6.20)}$$ 
where 
$${\cal R}''(c,t,d)
=\sum_k 
\sum_{\scriptstyle 0<|a|^2\leq A\atop\scriptstyle (a,cd)\sim c} 
\sum_{\scriptstyle h\atop\scriptstyle (cd)\mid h} \phi_h 
\sum_{\scriptstyle \ell\atop\scriptstyle t\mid\ell} \Upsilon_{\ell} 
\sum_{p\neq 0}\sum_{\scriptstyle 
q\neq 0\atop\!\!\!\!\!\!\!\!\!\!\!\!\!\!\!\!{\scriptstyle (cd)\mid pq}} 
\theta_p |pq|^{-2}\varphi(h,k,\ell,p,q) 
{\rm e}\!\left( {\rm Re}\!\left( {ak\over pq} 
+{h\over c}\,{\ell\over t}\,{(a/t)^{*}\over (pq/c)}\right)\right)\eqno(6.21)$$
(in which it is implicit that $(d,a/t)\sim 1$, so that one has both  
$(d,a/c)\sim 1$ and $(d,c/t)\sim 1$). 

Subject to the explicit conditions of summation in (6.21), the congruence 
$(a/t)(a/t)^{*}\equiv 1\bmod (pq/c){\frak O}$ implies 
$(a/t)(a/t)^{*}\equiv 1\bmod (pq/(cd)){\frak O}$, 
so that one has 
$${\rm e}\!\left( {\rm Re}\!\left( 
{h\over c}\,{\ell\over t}\,{(a/t)^{*}\over (pq/c)}\right)\right)
={\rm e}\!\left( {\rm Re}\!\left( 
{h\over cd}\,{\ell\over t}\,{(a/t)^{*}\over (pq/cd)}\right)\right)\eqno(6.22)$$
when the left-hand side of this equation is defined 
(i.e. when $(a/t,pq/c)\sim 1$). 
Moreover, given that one assumes $t\mid c$, $\,(c/t,d)\sim 1$ and $(a,cd)\sim c$, 
the equation~(6.22) is 
effectively an identity: for, if the right-hand side of (6.22) is defined,  
then $(a/t,pq/(cd))\sim 1$ and, by assumption, $(a/t,d)\sim((a/c)(c/t),d)\sim 1$, 
so that one has $(a/t,pq/c)\sim 1$, which is
sufficient to ensure that both sides of (6.22) are defined and equal. 
Therefore it follows by (6.18), (6.20)-(6.22) and the definition (5.4) and 
property (5.5) of the M\"{o}bius function $\mu_{\frak O}$ that 
$$\eqalignno{
{\cal R}' 
 &={1\over 64}\sum_{c\neq 0}\sum_{t\mid c} 
\mu_{\frak O}\left( {c\over t}\right) |t|^2 
\sum_{\scriptstyle d\neq 0\atop\scriptstyle (d,c/t)\sim 1} 
\mu_{\frak O}(d) {\cal R}''(c,t,d) = {}\cr 
 &={1\over 64}\sum_{d\neq 0}
\sum_{c\neq 0}\sum_{t\mid c} 
\mu_{\frak O}\left( {cd\over t}\right) |t|^2 
{\cal R}''(c,t,d) 
={1\over 64}\sum_{w\neq 0}
\sum_{c\mid w}\sum_{t\mid c} 
\mu_{\frak O}\left( {w\over t}\right) |t|^2 
{\cal R}^{*}(w,c,t)\;, &(6.23)}$$
where 
$$\eqalign{
{\cal R}^{*}(w,c,t)
 &={\cal R}''(c,t,w/c) = {}\cr 
 &=\sum_k 
\sum_{\scriptstyle 0<|a|^2\leq A\atop\scriptstyle (a,w)\sim c} 
\sum_{\scriptstyle h\atop\scriptstyle w\mid h} \phi_h 
\sum_{\scriptstyle \ell\atop\scriptstyle t\mid\ell} \Upsilon_{\ell} 
\sum_{p\neq 0}\sum_{\scriptstyle 
q\neq 0\atop\!\!\!\!\!\!\!\!\!\!\!\!\!\!\!\!{\scriptstyle w\mid pq}} 
\theta_p |pq|^{-2}\varphi(h,k,\ell,p,q) 
{\rm e}\!\left( {\rm Re}\!\left( {ak\over pq} 
+{h\over w}\,{\ell\over t}\,{(a/t)^{*}\over (pq/w)}\right)\right) .
}$$
Now, in the last sum over $p$, we group together terms according to 
the highest common factor $(p,w)$. 
When $(p,w)\sim r$ (say), one has $w\mid (pq)$ if and only if 
$(w/r)\mid q$. Therefore, by this grouping of terms, we find that 
$${\cal R}^{*}(w,c,t) 
={1\over 4}\sum_{\scriptstyle r,s\in{\frak O}\atop\scriptstyle rs=w} 
{\cal R}^{*}(w,c,t,r,s)\;,\eqno(6.24)$$ 
where 
$$\eqalign{
{\cal R}^{*}(w,c,t,r,s) 
 &=\sum_k \sum_{\scriptstyle 0<|a|^2\leq A\atop\scriptstyle (a,w)\sim c} 
\sum_{\scriptstyle h\atop\scriptstyle w\mid h} \phi_h 
\sum_{\scriptstyle \ell\atop\scriptstyle t\mid\ell} \Upsilon_{\ell} 
\!\!\!\sum_{\scriptstyle p\neq 0\atop\scriptstyle (p,w)\sim r}
\!\!{\theta_p\over |p|^2} \times {}\cr 
 &\qquad\times \sum_{\scriptstyle q\neq 0\atop\scriptstyle s\mid q} 
{\varphi(h,k,\ell,p,q)\over |q|^2} 
\,{\rm e}\!\left( {\rm Re}\!\left( {ak\over pq} 
+{(h/w)(\ell/t)(a/t)^{*}\over (p/r)(q/s)}\right)\right) .}\eqno(6.25)$$

By (6.2) and (6.25) one has ${\cal R}^{*}(w,c,t,r,s)=0$ unless 
$K/2<|k|^2<K$. Moreover, since  
$w\mid h$ is a condition of summation on the right-hand side of (6.25), 
it is implied by the constraints in (6.5) on the 
coefficients $\phi_h\ $ ($h\in{\frak O}-\{ 0\}$) that one has 
${\cal R}^{*}(w,c,t,r,s)=0$ whenever $|w|^2>H$. Therefore, and since 
$$\eqalign{ 
\sum_{0<|w|^2\leq H}\sum_{c\mid w}\sum_{t\mid c}\left|\mu_{\frak O}
\left( {w\over t}\right)\right| 
\sum_{\scriptstyle r,s\atop\scriptstyle rs=w}\sum_{0<|k|^2<K}{1\over |r|^2 |s|^2} 
 &\leq \sum_{0<|w|^2\leq H}\biggl(\,\sum_{d\mid w} 1\biggr)^{\!\!3} 
O\!\left( {K\over |w|^2}\right) \ll_{\varepsilon} {}\cr 
 &\ll_{\varepsilon}\sum_{0<|w|^2\leq H} {K\over |w|^{2-\varepsilon}}  
\leq\sum_{w\neq 0} 
{H^{\varepsilon} K\over |w|^{2+\varepsilon}}
\ll_{\varepsilon} H^{\varepsilon} K\;,}$$
it follows from (6.23)-(6.25), (6.2) and (6.5) that, for some 
$w,r,s,c,t,k\in{\frak O}-\{ 0\}\,$ satisfying 
$${K\over 2}<|k|^2<K,\quad\  0<|w|^2\leq H,\quad\  c\mid w,\quad\  t\mid c\quad\  
\hbox{and}\quad\  rs=w\;,\eqno(6.26)$$ 
one has: 
$${\cal R}'\ll_{\varepsilon} H^{\varepsilon} K |t|^2 
\left| {\cal D}(w,c,t,r,s;k)\right|\;,\eqno(6.27)$$
with 
$${\cal D}(w,c,t,r,s;k) 
=\sum_{\scriptstyle 0<|a|^2\leq A\atop\scriptstyle (a,w)\sim c} 
\sum_{\scriptstyle h\atop\scriptstyle w\mid h} \phi_h 
\sum_{\scriptstyle \ell\atop\scriptstyle t\mid\ell} \Upsilon_{\ell} 
\!\!\!\sum_{\scriptstyle p\neq 0\atop\scriptstyle (p,w)\sim r}
\!\!\theta_p |p/r|^{-2}\,{\cal U}(h,\ell,p;a/t)\;,\eqno(6.28)$$
where ${\cal U}(h,\ell,p;a/t)={\cal U}(h,\ell,p;a/t;w,t,r,s,k)$ is given by
$${\cal U}(h,\ell,p;a/t) 
=\sum_{\scriptstyle q\neq 0\atop\scriptstyle s\mid q} 
|q/s|^{-2}\varphi(h,k,\ell,p,q) 
\,{\rm e}\!\left( {\rm Re}\!\left( 
{ak\over pq} 
+{(h/w)(\ell/t)(a/t)^{*}\over (p/r)(q/s)}\right)\right) .\eqno(6.29)$$

In the last summation it is implicit that one sums only over $q\in{\frak O}$ 
such that $((p/r)(q/s),a/t)\sim 1\,$ (the sum is therefore void unless 
$(p/r,a/t)\sim 1$). When this condition is satisfied one can find 
(by the Euclidean algorithm for ${\Bbb Z}[i]$) Gaussian integers 
$((p/r)(q/s))^{*}$ and $(a/t)^{*}$ such that 
$$(p/r)(q/s)((p/r)(q/s))^{*}+(a/t)(a/t)^{*}=1\;.\eqno(6.30)$$ 
One then has, by (6.30), 
$${(a/t)^{*}\over (p/r)(q/s)}
={1\over (p/r)(q/s)(a/t)}-{((p/r)(q/s))^*\over (a/t)}\;,\eqno(6.31)$$ 
where the use of `$*$' accords with the 
convention set down in Subsection~1.5, under the heading `Number-Theoretic Notation': it (for example) being implied (6.30)  
that $(a/t)(a/t)^{*}\equiv 1\bmod (p/r)(q/s){\frak O}$. 
Moreover, since one has $((p/r)(q/s),a/t)\sim 1$ if and only if 
$(p/r,a/t)\sim 1$ and $(q/s,a/t)\sim 1$, and since the relations 
$(p/r)(p/r)^{*}\equiv 1\bmod (a/t){\frak O}$ and 
$(q/s)(q/s)^{*}\equiv 1\bmod (a/t){\frak O}$ imply 
$(p/r)(q/s)(p/r)^{*}(q/s)^{*}\equiv 1\bmod (a/t){\frak O}$,  
one will have  
$((p/r)(q/s))^{*}\equiv (p/r)^{*}(q/s)^{*}\bmod (a/t){\frak O}$ whenever 
either one of the residue classes $((p/r)(q/s))^{*}\bmod (a/t){\frak O}$,  
$(p/r)^{*}(q/s)^{*}\bmod (a/t){\frak O}$ is defined. 
Hence, and by (6.26), (6.29) and (6.31), we find that, for 
$h,\ell,p,a$ satisfying the conditions of summation in (6.28), one has: 
$${\cal U}(h,\ell,p;a/t)
={|s|^2\over Q}\sum_{\scriptstyle q\atop\scriptstyle s\mid q} f(q) 
\,{\rm e}\!\left( {\rm Re}\!\left( {m(q/s)^{*}\over (a/t)}\right)\right) ,\eqno(6.32)$$
where 
$$m=m(h,\ell,p;a/t)\in{\frak O},\qquad m\equiv -(h/w)(\ell /t)(p/r)^{*}\bmod (a/t){\frak O}\;,\eqno(6.33)$$ 
while, for $z\in{\Bbb C}$, 
$$f(z)=f_{h,\ell,p,\nu}(z) =\cases{0 &if $z=0$, \cr 
Q|z|^{-2}\varphi(h,k,\ell,p,z) 
\,{\rm e}\!\left( {\rm Re}\bigl(\nu z^{-1}\bigr)\right) &otherwise,}\eqno(6.34)$$
with 
$$\nu =\nu(h,\ell,p;a/t)={ak\over p}+{h\ell\over ap}\;,\eqno(6.35)$$
so that if $\phi_h\Upsilon_{\ell}\theta_p\neq 0$ then, by (6.5), (6.8) 
and the first part of (6.14), 
$$|\nu|^2\ll {AK\over P}+{HL\over P}\ll {(PQ)^{1+\varepsilon}\over\delta^2 P}+Q 
\leq (PQ)^{\varepsilon}\left(\delta^{-2}+1\right) Q\;.\eqno(6.36)$$ 

Since $\varphi\in{\cal S}\bigl( {\Bbb C}^5\bigr)$,  and since  (by (6.2))   
$\varphi(h,k,\ell,p,z)=0$ for $|z|^2\leq Q/2\,$, the definition 
(6.34) ensures that we have $f\in{\cal S}({\Bbb C})$. 
Taking now 
$$\Delta =(PQ)^{-\varepsilon}\left( 1+\delta^{-1}\right)^{-2}\;,\qquad 
\Omega_1 =Q\;,\qquad{\rm and}\qquad C=2\;,\eqno(6.37)$$ 
we seek to verify that 
$\Delta$, $\Omega_1$, $C$ and $f$ satisfy the case $n=1$ 
of the conditions~(5.24), (5.25) of Lemma~5.4.   
This will enable us to obtain, by means of Lemma~5.6, 
an alternative expression for the sum over $q$ in (6.32). 

The verification that (5.25) is satisfied requires no work, since 
the condition (6.2) immediately implies the case $n=1$ of (5.25) 
(when $f$ is as in (6.34), and $\Omega_1,C$ as in (6.37)). 

Our verification of (5.24) (for $n=1$) begins with the observation that, by (6.34), 
$$f(z)=\varphi(h,k,\ell,p,z)
\left( {Q^{1/2}\over z}\,\exp\left( {\pi i\nu\over z}\right)\right) 
\left( {Q^{1/2}\over \overline{z}}
\,\exp\left( {\pi i\overline{\nu}\over\overline{z}}\right)\right)\qquad\quad  
\hbox{for $\quad z\in{\Bbb C}^{*}$.}$$ 
Hence, and by (5.18)-(5.19), (5.21) and Leibniz's rule for the higher order derivatives 
of a product, a short calculation suffices to show that, for $j\in{\Bbb N}\cup\{ 0\}$ and 
$z\in{\Bbb C}^{*}$, one has: 
$$\eqalign{{\cal L}_1^j f(z)
 &=(-4)^j\sum_{\lambda =0}^j\sum_{\mu =0}^j\pmatrix{j\cr\lambda}\pmatrix{j\cr\mu} 
\left( {\partial^{\lambda}\over\partial\overline{z}^{\lambda}}
\,{Q^{1/2}\over \overline{z}}
\,\exp\left( {\pi i\overline{\nu}\over\overline{z}}\right)\right) \times \cr
 &\qquad\qquad\qquad\quad\times
\left( {\partial^{j-\mu}\over\partial z^{j-\mu}}
\,{\partial^{j-\lambda}\over\partial\overline{z}^{j-\lambda}}
\ \varphi(h,k,\ell,p,z)\right) 
\left( {\partial^{\mu}\over\partial z^{\mu}} 
\,{Q^{1/2}\over z}\,\exp\left( {\pi i\nu\over z}\right)\right) ,}\eqno(6.38)$$
where, by (5.19) and (6.2)-(6.3), 
$${\partial^{j-\mu}\over\partial z^{j-\mu}}
\,{\partial^{j-\lambda}\over\partial\overline{z}^{j-\lambda}}
\ \varphi(h,k,\ell,p,z)
\ll_j (\delta |z|)^{\lambda +\mu-2j}\qquad\quad  
\hbox{for $\quad z\in{\Bbb C}^{*}\ $ and $\ \mu,\lambda\in\{ 0,1,\ldots ,j\}$.}\eqno(6.39)$$ 

In considering the other derivatives in (6.38), we may note that 
if $\nu\neq 0$, then 
$${\partial^{\mu}\over\partial z^{\mu}} 
\,{Q^{1/2}\over z}\,\exp\left( {\pi i\nu\over z}\right)
={Q^{1/2}\over\pi i\nu}\,{\partial^{\mu}\over\partial z^{\mu}}
\,g\!\left( {z\over \pi i\nu}\right)\qquad\qquad  
\hbox{($\mu\in{\Bbb N}\cup\{ 0\}$, $z\in{\Bbb C}^{*}$),}$$
where 
$$g(\tau)=\tau^{-1}\exp\left(\tau^{-1}\right)\qquad\qquad 
\hbox{($\tau\in{\Bbb C}^{*}$).}$$
By induction it may be established that, for each $\mu\in{\Bbb N}\cup\{ 0\}$, 
one has 
$$g^{(\mu)}(\tau)
=\biggl(\ \sum_{\kappa =1+\mu}^{1+2\mu}\alpha(\mu,\kappa)\tau^{-\kappa}\biggr) 
\exp\left(\tau^{-1}\right)\qquad\qquad 
\hbox{($\tau\in{\Bbb C}^{*}$),}$$
where the coefficients $\alpha(\mu,\kappa)$ are certain integer 
valued constants. Hence, for $\mu\in{\Bbb N}\cup\{ 0\}$ and 
$z\in{\Bbb C}^{*}$, one obtains: 
$$\eqalign{ 
{\partial^{\mu}\over\partial z^{\mu}} 
\,{Q^{1/2}\over z}\,\exp\left( {\pi i\nu\over z}\right)
 &=Q^{1/2}(\pi i\nu)^{-(\mu +1)} g^{(\mu)}\!\left( {z\over \pi i\nu}\right) \ll {}\cr 
 &\ll Q^{1/2} |\pi\nu|^{-(\mu +1)} O_{\mu}\!\left(\left|{\nu\over z}\right|^{1+\mu} 
+\left|{\nu\over z}\right|^{1+2\mu}\right) 
\left|\exp\left( {\pi i\nu\over z}\right)\right| \ll_{\mu} {}\cr 
 &\ll_{\mu} Q^{1/2} |z|^{-(\mu +1)} \left( 1 + \left|{\nu\over z}\right|\right)^{\mu} 
\left|\exp\left( {\pi i\nu\over z}\right)\right|\;.}$$
One obtains the same bound (more easily) when $\nu =0$. 
Similarly, one has  
$${\partial^{\lambda}\over\partial\overline{z}^{\lambda}} 
\,{Q^{1/2}\over\overline{z}}\,\exp\left( {\pi i\overline{\nu}\over\overline{z}}\right)
\ll_{\lambda} Q^{1/2} |z|^{-(\lambda +1)} \left( 1 + \left|{\nu\over z}\right|\right)^{\lambda} 
\left|\exp\left( {\pi i\overline{\nu}\over\overline{z}}\right)\right|\qquad\quad  
\hbox{for $\quad \lambda\in{\Bbb N}\cup\{ 0\}\ $ and $\ z\in{\Bbb C}^{*}$.}$$
Since 
$$\left|\exp\left( {\pi i\nu\over z}\right)\right|
\,\left|\exp\left( {\pi i\overline{\nu}\over\overline{z}}\right)\right|
=\left|\exp\left( 2\pi i{\rm Re}\left( {\nu\over z}\right)\right)\right| =1\;,$$
it follows by the last two upper bounds, and by (6.38) and (6.39), that 
$$\eqalign{ 
{\cal L}_1^j f(z) 
 &\ll_j\ \sum_{\lambda =0}^j \sum_{\mu =0}^j 
(\delta |z|)^{\lambda +\mu -2j} Q |z|^{-(\lambda +\mu +2)}
\left( 1 +\left|{\nu\over z}\right|\right)^{\lambda +\mu} = {}\cr 
 &={Q\over |z|^2}\,(\delta |z|)^{-2j}\biggl( 
\ \sum_{\lambda =0}^j \left( 1+\left|{\nu\over z}\right|\right)^{\lambda}\delta^{\lambda}
\biggr)^{\!2} \ll_j {}\cr 
 &\ll_j Q |z|^{-2} (\delta |z|)^{-2j} 
\left( 1 + \left( 1+\left|{\nu\over z}\right|\right)\delta\right)^{2j}
=Q |z|^{-(2j+2)}\left(\delta^{-1} +1 + \left|{\nu\over z}\right|\right)^{2j}\;.}$$
Since we have already verified that 
the function $f$ satisfies the case $n=1$ of (5.25), 
with $\Omega_1$ and $C$ as in (6.37),     
it may therefore be assumed in the above that 
$Q/2<|z|^2<2Q\,$ (for it is otherwise trivially the case that ${\cal L}_1^j f(z)=0$). 
Hence, and by (6.36), we obtain (for $j\in{\Bbb N}\cup\{ 0\}$ and $z\in{\Bbb C}^{*}$): 
$${\cal L}_1^j f(z) 
\ll_j |z|^{-2j}\Bigl(\left(\delta^{-1}+1\right)^2 (PQ)^{\varepsilon}\Bigr)^j 
=\left(\Delta |z|^2\right)^{-j}\;,$$
with $\Delta$ as in (6.37). 

Since the above completes the verification of (5.24),  
and since (5.25) has also been verified, we may now apply Lemma~5.6, 
with $f\in{\cal S}\bigl( {\Bbb C}\bigr)$ given by (6.34), and 
$\Delta,\Omega_1,C$ as in (6.37). By the case $d=s$, $B=(PQ)^{\varepsilon}$ 
of the result (5.29) of Lemma~5.6,
it follows that, for $m\in{\frak O}$, $a_1=a/t\in{\frak O}-\{ 0\}$ and $j\geq 1$, 
we have 
$$\eqalign{
\sum_{\scriptstyle q\atop\scriptstyle s\mid q} f(q)
\,{\rm e}\!\left( {\rm Re}\!\left( 
{m (q/s)^{*}\over a_1}\right)\right)  
 &={1\over\left| a_1\right|^2} 
\sum_{\scriptstyle b\bmod a_1{\frak O}\atop\scriptstyle 
\bigl\|{\textstyle{b\over a_1}}\bigr\|^2\leq V(s)} 
S\left( -m , b ; a_1\right) 
\sum_{\scriptstyle q\atop\scriptstyle s\mid q} f(q)
\,{\rm e}\!\left( {\rm Re}\!\left( 
{b\over a_1}\,{q\over s}\right)\right) + {}\cr 
 &\quad\ +O_j\!\left( \left| a_1\right|^2\left( 1+\delta^{-2}\right) 
(PQ)^{-j\varepsilon}\right) ,}$$
where (given (6.37)) $ V(z)=(PQ)^{\varepsilon}|z|^2 /(\Delta\Omega_1)
=(E^2 /Q)|z|^2$, with $E=(PQ)^{\varepsilon}(1+\delta^{-1})\,$   
(as stated in the lemma), so that $ V(z)$ is the function 
defined in (6.14). Hence, and by (6.32)-(6.35), we have, in (6.28), 
$$\eqalignno{
{\cal U}(h,\ell,p;a/t) 
 &=O_j\!\left( (PQ)^{-j\varepsilon}  V(sa/t)\right) + {} 
 &(6.40)\cr 
 &\quad\ +{1\over |a/t|^2}
\sum_{\scriptstyle b\bmod (a/t){\frak O}\atop\scriptstyle 
\bigl\|{\textstyle{b\over (a/t)}}\bigr\|^2\leq V(s)} 
S\left( (h/w)(\ell /t)(p/r)^{*} , b ; a/t\right) \times {}\cr 
 &\qquad\qquad\qquad\qquad\qquad  
\times\sum_{\scriptstyle q\neq 0\atop\scriptstyle s\mid q} 
{\varphi( h,k,\ell,p,q)\over |q/s|^2}
\,{\rm e}\!\left( {\rm Re}\!\left( {ak\over pq}+{h\ell\over apq}
+{(q/s)b\over (a/t)}\right)\right) ,}$$
for $j\geq 3$. 
Moreover, by the second inequality in (6.10), we have, as recorded in (6.14),   
$\,V(z)<|z|^2 /4H$; 
given the conditions~(6.26) which the Gaussian integers $r$, $s$ and $w$ satisfy, 
it therefore follows that 
$$ V(s)<{|s|^2\over 4H}\leq {|w|^2\over 4H}\leq {1\over 4}\;,\eqno(6.41)$$ 
and so (by reasoning similar to that which justified the 
simple condition `$|a|^2\leq A$' in (6.15)) we are able to 
simplify the conditions for summation over $b$, in (6.40), 
to just: $b\in{\frak O}$ and 
$|b|^2\leq V(s) |a/t|^2 =  V(sa /t)$.  

Amongst the terms of the sum over $b$ in (6.40), the term in the case 
$b=0$ is special: for by the result (5.35) of Lemma~5.7 one has 
$$\left| S\left( (h/w)(\ell /t)(p/r)^{*} , 0 ; a/t\right)\right| 
\leq\left|\bigl( (h/w)(\ell/t)(p/r)^{*} , a/t\bigr)\right|^2 
=\left|\bigl( (h/w)(\ell/t) , a/t\bigr)\right|^2\eqno(6.42)$$ 
(the $*$-notation implying, in this context,  
$(p/r)(p/r)^{*}\equiv 1\bmod (a/t){\frak O}$, so that 
$(p/r)^{*}$ and $a/t$ are coprime). 

For $b\neq 0$ the best 
available estimate for the Kloosterman sum $S(-m,b;a/t)$ is 
[3, Theorem~10], which shows that one has 
$$|S(-m,b;a/t)|\leq 2^{\omega(a/t)+7/2}
|(-m,b,a/t)a/t|\;,$$
where $\omega(n)$ denotes the number of 
prime ideals of ${\frak O}$ containing $n$; we do not use this upper bound, since 
(as our result in the final lemma of this section shows) there is an 
advantage to be gained in doing otherwise: we shall, in effect, exploit 
cancellations between different Kloosterman sums.  

In conjunction with (6.42) we shall use the bound 
$$\sum_{\scriptstyle q\neq 0\atop\scriptstyle s\mid q} {\varphi( h,k,\ell,p,q)\over |q/s|^2}
\,{\rm e}\!\left( {\rm Re}\!\left( {ak\over pq}+{h\ell\over apq}\right)\right)
\ll \sum_{\scriptstyle q\neq 0\atop\scriptstyle s\mid q} 
{\left|\varphi( h,k,\ell,p,q)\right|\over |q/s|^2} 
=\sum_{0<|q_1|^2<{\textstyle{Q\over |s|^2}}} O\!\left( {|s|^2\over Q}\right)\ll 1\;,\eqno(6.43)$$
which is implied by the hypotheses (6.2), (6.3). 

By (6.41), (6.26), (6.10) and the definitions of $A$ and $V(z)$, in (6.14), 
it follows that if $j=[2/\varepsilon]+1$ (where by hypothesis, $\varepsilon >0$) 
then 
$$(PQ)^{-j\varepsilon} V(sa/t)\leq (PQ)^{-2} |a/t|^2  V(s) 
< (16 A)^{-2} |a|^2 /4 
\leq 2^{-10}|a|^{-2}< |a/t|^{-2}\qquad\quad\hbox{for $\quad 0\neq a\in{\frak O}$.}$$
By this, together with (6.40), (6.42), (6.43), the observation following 
(6.41), and the equation~(6.28), we obtain: 
$${\cal D}(w,c,t,r,s;k) 
={\cal D}_0(w,c,t,r,s;k) +{\cal D}_1(w,c,t,r,s;k)\;,$$ 
where 
$${\cal D}_0(w,c,t,r,s;k)
=\sum_{\scriptstyle 0<|a|^2\leq A\atop\scriptstyle (a,w)\sim c} 
\sum_{\scriptstyle h\atop\scriptstyle w\mid h} \phi_h 
\sum_{\scriptstyle \ell\atop\scriptstyle t\mid\ell} \Upsilon_{\ell} 
\!\!\!\sum_{\scriptstyle p\neq 0\atop\scriptstyle (p,w)\sim r}
\!\!\theta_p |p/r|^{-2} 
\,O_{\varepsilon}\!\!\left( {\left|\bigl( (h/w)(\ell/t) , a/t\bigr)\right|^2\over 
|a/t|^2}\right)\eqno(6.44)$$
and 
$${\cal D}_1(w,c,t,r,s;k)
=\sum_{\scriptstyle q\neq 0\atop\scriptstyle s\mid q} 
|q/s|^{-2} {\cal E}(w,c,t,r,s;k,q)\;,$$
with ${\cal E}(w,c,t,r,s;k,q)$ as given by~(6.13).  

By the result just obtained, we either have 
$${\cal D}(w,c,t,r,s;k)\ll |{\cal D}_0(w,c,t,r,s;k)|\;,$$ 
or else have 
$$0<|{\cal D}(w,c,t,r,s;k)|\ll |{\cal D}_1(w,c,t,r,s;k)|\;.$$ 

In the latter case it follows by (6.13), (6.2) and the rightmost bound in (6.43) 
that, for some $q\in s{\frak O}-\{ 0\}$ satisfying $Q/2<|q|^2<Q$, one will have 
the upper bound ${\cal D}(w,c,t,r,s;k)\ll |{\cal E}(w,c,t,r,s;k,q)|$, which,   
by (6.17), (6.26)-(6.27) and (6.8), implies the result (6.11)-(6.13) 
of the lemma. 

In the former case, where 
${\cal D}(w,c,t,r,s;k)\ll |{\cal D}_0(w,c,t,r,s;k)|$, 
one obtains the bound (6.9): for, by (6.44), (6.26), (6.4), (6.5), (6.8), (6.14) and 
the Cauchy-Schwarz inequality, 
\goodbreak 
$$\eqalign{ 
{\cal D}_0(w,c,t,r,s;k) 
 &\ll_{\varepsilon} \sum_{\scriptstyle 0<|a|^2\leq A\atop\scriptstyle c\mid a} 
|a/t|^{-2}\sum_{\scriptstyle{\textstyle{H\over 2}}<|h|^2\leq H\atop\scriptstyle w\mid h} 
\sum_{\scriptstyle{\textstyle{L\over 2}}<|\ell|^2\leq L\atop\scriptstyle t\mid\ell} 
\left|\Upsilon_{\ell}\right|\,\left|\left( {h\over w}\,{\ell\over t}\,,\,{a\over t}\right)\right|^2 \leq {}\cr 
 &\leq 
\sum_{\scriptstyle{\textstyle{H\over 2}}<|h|^2\leq H\atop\scriptstyle w\mid h} 
\sum_{\scriptstyle{\textstyle{L\over 2}}<|\ell|^2\leq L\atop\scriptstyle t\mid\ell} 
\left|\Upsilon_{\ell}\right| 
\sum_{\scriptstyle 0<|a|^2\leq A\atop\scriptstyle c\mid a} 
|a/c|^{-2} \left|\left( h\ell\,,\,{a\over c}\right)\right|^2 \leq {}\cr 
 &\leq 
\sum_{\scriptstyle{\textstyle{H\over 2}}<|h|^2\leq H\atop\scriptstyle w\mid h} 
\sum_{\scriptstyle{\textstyle{L\over 2}}<|\ell|^2\leq L\atop\scriptstyle t\mid\ell} 
\left|\Upsilon_{\ell}\right| 
\sum_{d\mid h\ell} |d|^2 
\sum_{\scriptstyle 0<|a_1|^2\leq{\textstyle{A\over |c|^2}}\atop\scriptstyle d\mid a_1} 
\left| a_1\right|^{-2}  = {}\cr 
 &=\sum_{\scriptstyle{\textstyle{H\over 2}}<|h|^2\leq H\atop\scriptstyle w\mid h} 
\sum_{\scriptstyle{\textstyle{L\over 2}}<|\ell|^2\leq L\atop\scriptstyle t\mid\ell} 
\left|\Upsilon_{\ell}\right| 
\sum_{d\mid h\ell}  
\ \sum_{0<|a_2|^2\leq{\textstyle{A\over |cd|^2}}} 
\left| a_2\right|^{-2} = {}\cr 
 &=\sum_{\scriptstyle{\textstyle{H\over 2}}<|h|^2\leq H\atop\scriptstyle w\mid h} 
\sum_{\scriptstyle{\textstyle{L\over 2}}<|\ell|^2\leq L\atop\scriptstyle t\mid\ell} 
\left|\Upsilon_{\ell}\right| O_{\varepsilon}\!\!\left( 
|h\ell|^{\varepsilon /2}\log(A+1)\right) 
\ll_{\varepsilon} (PQ)^{\varepsilon /2}
\,{H\over |w|^2}\left( {L\over |t|^2}
\sum_{\ell}\left|\Upsilon_{\ell}\right|^2\right)^{\!\!1/2}\;; 
}$$
and so, when 
${\cal D}(w,c,t,r,s;k)\ll |{\cal D}_0(w,c,t,r,s;k)|$, 
it follows by (6.17), (6.26)-(6.27) and (6.8), that 
$$\eqalign{ 
{\cal R} &\ll_{\varepsilon} 
H^{\varepsilon} (PQ)^{\varepsilon /2} |t| |w|^{-2} K H L^{1/2} 
\biggl(\sum_{\ell}\left|\Upsilon_{\ell}\right|^2\biggr)^{\!1/2} \leq {}\cr 
 &\leq H^{\varepsilon} (PQ)^{\varepsilon /2} (HKL)^{1/2} 
\biggl(HK\sum_{\ell}\left|\Upsilon_{\ell}\right|^2\biggr)^{\!1/2}
\ll_{\varepsilon}  (PQ)^{\varepsilon} (HKL)^{1/2} 
\biggl(HK\sum_{\ell}\left|\Upsilon_{\ell}\right|^2\biggr)^{\!1/2}\;.
}$$
This completes the proof of the lemma\ $\blacksquare$

\bigskip

\goodbreak\proclaim Lemma~6.2. Let 
$A_1,H_1,L_1,P_1>0$; and let $0\neq c_1 \in{\frak O}$. For 
$h,\ell\in{\frak O}$, $a\in c_1{\frak O}-\{ 0\}$ and 
$b,p\in{\frak O}-\{ 0\}$, let $\Phi_h$, $B_{\ell}$, $\xi_a$, $\psi(a;b)$ 
and $\Theta(a;p)$ be complex numbers such that  
$$\Phi_h=0\ {\rm unless}\ {H_1\over 2}<|h|^2\leq H_1\;,\ \,  
B_{\ell}=0\ {\rm unless}\ {L_1\over 2}<|\ell|^2\leq L_1\;,\ \,  
\Theta(a;p)=0\ {\rm unless}\ {P_1\over 2}<|p|^2\leq P_1\;,\eqno(6.45)$$
$$\Phi_h\ll 1\;,\qquad\xi_a\ll |a|^{-2}\;,\qquad\psi(a;b)\ll 1\quad\ 
{\rm and}\quad\ \Theta(a;p)\ll 1\;.\eqno(6.46)$$
Let $\Delta ,\varepsilon >0$; let $0<\rho<1/4$; and let  
$f\in{\cal S}\bigl( {\Bbb C}^3\bigr)$. Suppose moreover that 
$$f\left( z_1 , z_2 , z_3\right) =0\qquad\quad{\rm unless}\quad\  
\left( {\left| z_1\right|^2\over H_1}\,,\,{\left| z_2\right|^2\over L_1}\,,
\,{\left| z_3\right|^2\over P_1}\right)\in\left( {1\over 2}\,,\,1\right)^3 ,\eqno(6.47)$$
and that 
$${\cal L}_1^{j_1}{\cal L}_2^{j_2}{\cal L}_3^{j_3} f({\bf z}) 
\ll_{\bf j}\,\prod_{k=1}^3\left(\Delta\bigl| z_m\bigr|^2\right)^{-j_m}\qquad\qquad  
\hbox{(${\bf z}\in\left({\Bbb C}^{*}\right)^{\!3}$, ${\bf j}\in ({\Bbb N}\cup\{ 0\})^3$),} 
\eqno(6.48)$$
where ${\cal L}_m$ is the linear operator on ${\cal S}\bigl( {\Bbb C}^3\bigr)$ 
defined in the equation~(5.12), in Lemma~5.2. 
Put 
$${\cal E}^{*} 
=\sum_{\scriptstyle 0<|a|^2\leq A_1\atop\scriptstyle c_1 \mid a} 
\xi_a \sum_h \Phi_h \sum_{\ell} B_{\ell} \sum_p \Theta(a;p) 
f(h,\ell,p) \sum_{0<|b|^2\leq\rho |a|^2}\psi(a;b) S\left( h \ell p^{*} , b ; a\right)\;, 
\eqno(6.49)$$
where $S(u,v;w)$ is given by (1.3.6), and where the superscript notation `$*$' 
has the meaning explained in Subsection~1.5, under `Number-Theoretic Notation' (it therefore being an implicit condition of 
the summation on the right-hand side of the equation (6.49) that $(p,a)\sim 1$). 
Then 
$${\cal E}^{*} 
\ll_{\varepsilon} \Delta^{-3} \left( H_1 L_1 P_1 A_1\right)^{\varepsilon} 
\left( 1+{H_1 L_1 \left| c_1\right|^2\over A_1^2}\right)^{1/2} 
\left( 1 +\rho P_1\right)^{1/2} 
\left( \left| c_1\right|^{-4} 
\rho A_1^3 P_1 H_1 \sum_{\ell}\left| B_{\ell}\right|^2\right)^{1/2} , 
\eqno(6.50)$$
where the implicit constant is determined by those in (6.46) and (6.48),  
and by the positive constant $\varepsilon$. 

\medskip 

\goodbreak 
\noindent{\bf Proof.}\quad 
We may suppose that $A_1,H_1,L_1,P_1\geq 1$: for it is otherwise trivially  
implied by (6.47) and (6.49) that ${\cal E}^{*}=0$. 
Since $f\in{\cal S}\bigl( {\Bbb C}^3\bigr)$, it follows by (6.49) and 
Fourier's inversion formula (the case $n=3$ of Lemma~5.1, Equation (5.11)) that we 
have 
$${\cal E}^{*}
=\int_{\Bbb C}\int_{\Bbb C}\int_{\Bbb C}\hat{f}({\bf w}) {\cal E}({\bf w}) 
\,{\rm d}_{+}w_1\,{\rm d}_{+}w_2\,{\rm d}_{+}w_3\;,\eqno(6.51)$$
where $\hat{f}({\bf w})$ is the Fourier transform of $f$ defined in (5.2)-(5.3), 
while, for ${\bf w}\in{\Bbb C}^3$, 
$$\eqalignno{
{\cal E}({\bf w})  
 &=\sum_{\scriptstyle 0<|a|^2\leq A_1\atop\scriptstyle c_1\mid a} 
\xi_a \sum_h \Phi_h\,{\rm e}\bigl({\rm Re}\bigl( h w_1\bigr)\bigr) 
\sum_{\ell} B_{\ell}\,{\rm e}\bigl({\rm Re}\bigl( \ell w_2\bigr)\bigr)\times {} 
 &(6.52)\cr  
 &\qquad\qquad\quad\times\sum_p \Theta(a;p)\,{\rm e}\bigl({\rm Re}\bigl( p w_3\bigr)\bigr) 
\sum_{0<|b|^2\leq\rho |a|^2}\psi(a;b) S\left( h \ell p^{*} , b ; a\right)\;. 
}$$
Let ${\bf\Omega_1}=(H_1 , L_1 , P_1)\in(0,\infty)^3$ and $C=2$. 
Then it follows by (6.47) and (6.48) that 
$\Delta$, ${\bf\Omega}$, $C$ and $f$ satisfy 
the case $n=3$ of the conditions  
(5.24) and (5.25) of Lemma~5.4.
That lemma therefore applies, giving: 
$$\hat{f}({\bf w}) 
\ll\prod_{k=1}^3 {\Omega_k\over\bigl( 1+\Delta\Omega_k\left| w_k\right|^2\,\bigr)^{\!2}}\qquad\qquad 
\hbox{for $\quad {\bf w}\in{\Bbb C}^3$}$$
(this being the case $j=2$ of the result (5.26)). 
By this bound for $\hat{f}({\bf w})$, one has: 
$$\int_{\Bbb C}\int_{\Bbb C}\int_{\Bbb C}\bigl| \hat{f}({\bf w})\bigr|  
\,{\rm d}_{+}w_1\,{\rm d}_{+}w_2\,{\rm d}_{+}w_3 
\ll\prod_{k=1}^3 
\int_{\Bbb C} 
{\Omega_k\over\bigl( 1+\Delta\Omega_k\left| w_k\right|^2\,\bigr)^{\!2}}
\,{\rm d}_{+}w_k =\left( {\pi\over\Delta}\right)^3 .$$
Therefore, given (6.51),  it must be the case that 
$${\cal E}^{*}\ll \Delta^{-3} \bigl| {\cal E}({\bf w})\bigr|\qquad\qquad  
\hbox{for some $\quad {\bf w}\in{\Bbb C}^3$.}\eqno(6.53)$$

Let ${\bf w}\in{\Bbb C}^3$. Then, by (6.52), 
$${\cal E}({\bf w})  
=\sum_{\scriptstyle 0<|a|^2\leq A_1\atop\scriptstyle c_1\mid a} 
\xi_a \sum_m \tau_m \sum_p \vartheta(a;p)  
\sum_{0<|b|^2\leq\rho |a|^2}\psi(a;b) S\left( m p^{*} , b ; a\right)\;,\eqno(6.54)$$
where, for $m,a,p\in{\frak O}-\{ 0\}$,  
$$\tau_m =\sum_{\scriptstyle h , \ell\atop\scriptstyle h\ell =m} 
\Phi_h B_{\ell}\,{\rm e}\bigl({\rm Re}\bigl( h w_1 + \ell w_2\bigr)\bigr)\qquad\ 
{\rm and}\qquad\  
\vartheta(a;p)=\Theta(a;p)\,{\rm e}\bigl({\rm Re}\bigl( p w_3\bigr)\bigr)\;.\eqno(6.55)$$
In the above it is implicit in the $*$-notation that $(p,a)\sim 1$ and 
$p^{*} p\equiv p p^{*}\equiv 1\bmod a{\frak O}$, 
so that by (1.3.6) one has $S(m p^{*} , b ; a)=S( m p^{*}, b p^{*} p ; a)$. 
Moreover, when $a\in{\frak O}$ and $(p,a)\sim 1$, one has 
$p^{*}d^{*}\equiv (pd)^{*}\bmod a{\frak O}$ for all $d\in{\frak O}$ such that 
$(d,a)\sim 1$; and if one restricts $d$ here to lie in some complete 
set of reduced residues $\bmod\,a{\frak O}$ then the mapping 
$d\bmod a{\frak O}\mapsto p d\bmod a{\frak O}$ is 
a permutation of the multiplicative group $({\frak O}/a{\frak O})^{*}$.  
It is therefore an immediate consequence of the definition (1.3.6) of the 
`simple Kloosterman sum' that 
$S( m p^{*}, b p^{*} p ; a) =S( m , b p^{*} ; a)$, so that 
in (6.54) one has $S(m p^{*} , b ; a)=S( m , b p^{*} ; a)$. 
Hence, by applying the definition (1.3.6), for $u=m$, 
$v\equiv b p^{*}\bmod a{\frak O}$ and $w=a$, we are able to 
rewrite (6.54) as: 
$${\cal E}({\bf w}) 
=\sum_{\scriptstyle 0<|a|^2\leq A_1\atop\scriptstyle c_1\mid a}  
\sum_{\scriptstyle d\bmod a{\frak O}\atop\scriptstyle (d,a)\sim 1} 
T\!\left( {d^{*}\over a}\right) 
U(a;d)\;,$$
where 
$$T(z)=\sum_m \tau_m\,{\rm e}({\rm Re}(mz))\qquad\qquad\hbox{($z\in{\Bbb C}$)}\eqno(6.56)$$
and 
$$U(a;d) 
=\xi_a \sum_p \vartheta(a;p)  
\sum_{0<|b|^2\leq\rho |a|^2}\psi(a;b)
\,{\rm e}\!\left({\rm Re}\left( {b p^{*} d\over a}\right)\right)\qquad\qquad  
\hbox{($a\in{\frak O}-\{ 0\}$, $d\in{\frak O}$).}\eqno(6.57)$$

By the Cauchy-Schwarz inequality, we deduce from the result just obtained that 
$$\bigl| {\cal E}({\bf w})\bigr|^2 
\leq {\cal X} {\cal Y}\;,\eqno(6.58)$$
where 
$${\cal X}
=\sum_{\scriptstyle 0<|a|^2\leq A_1\atop\scriptstyle c_1\mid a}  
\sum_{\scriptstyle d\bmod a{\frak O}\atop\scriptstyle (d,a)\sim 1} 
\biggl| T\!\left( {d^{*}\over a}\right)\biggr|^2 
=\sum_{\scriptstyle 0<|a|^2\leq A_1\atop\scriptstyle c_1\mid a}  
\sum_{\scriptstyle d\bmod a{\frak O}\atop\scriptstyle (d,a)\sim 1} 
\biggl| T\!\left( {d\over a}\right)\biggr|^2 \geq 0$$
and 
$${\cal Y} 
=\sum_{\scriptstyle 0<|a|^2\leq A_1\atop\scriptstyle c_1\mid a}  
\sum_{\scriptstyle d\bmod a{\frak O}\atop\scriptstyle (d,a)\sim 1} 
|U(a;d)|^2\geq 0\;.$$
By (6.55) and the hypothesis (6.45), we have $\tau_m =0$ unless 
$0<|m|^2\leq H_1 L_1$. Hence, given (6.56), the application of 
Lemma~5.9 (a special analytic large sieve for ${\Bbb Z}[i]$), 
yields the bound 
$${\cal X}\ll\left( H_1 L_1 +{A_1^2\over \left| c_1\right|^2}\right) \sum_m \left|\tau_m\right|^2\;.
\eqno(6.59)$$

The structure of the sum in (6.57) prevents us from obtaining a bound analogous 
to (6.59) for ${\cal Y}$. We fall back on the  observation that one has  
(trivially) the upper bound 
$${\cal Y} 
\leq\sum_{\scriptstyle 0<|a|^2\leq A_1\atop\scriptstyle c_1\mid a}  
\sum_{d\bmod a{\frak O}} 
|U(a;d)|^2\;,$$
which, by (6.57) and the orthogonality relations (5.8), implies: 
$$\eqalignno{ 
{\cal Y}
 &\leq\sum_{\scriptstyle 0<|a|^2\leq A_1\atop\scriptstyle c_1\mid a}  
\left|\xi_a\right|^2 \sum_{d\bmod a{\frak O}} 
\left| \sum_p \vartheta(a;p)  
\sum_{0<|b|^2\leq\rho |a|^2}\psi(a;b)
\,{\rm e}\!\left({\rm Re}\left( {b p^{*} d\over a}\right)\right)\right|^2 = {}\cr 
 &=\sum_{\scriptstyle 0<|a|^2\leq A_1\atop\scriptstyle c_1\mid a}  
|a|^{2}\left|\xi_a\right|^2 
\sum_{p_1} \sum_{p_2} 
\vartheta\bigl( a;p_1\bigr)\,\overline{\vartheta\bigl( a;p_2\bigr)}
\quad\ \sum\!\!\!\!\!\!\!\!
\sum_{\!\!\!\!\!\!\!\!\!\!\!\!\!
{\scriptstyle 0<|b_1|^2,|b_2|^2\leq\rho |a|^2\atop\scriptstyle 
b_1 p_1^{*}\equiv b_2 p_2^{*}\bmod a{\frak O}}}
\psi\bigl( a;b_1\bigr)\,\overline{\psi\bigl( a;b_2\bigr)} = {}\cr 
 &=\sum_{\scriptstyle 0<|a|^2\leq A_1\atop\scriptstyle c_1\mid a}  
|a|^{2}\left|\xi_a\right|^2 
\ \sum_{p_1}\!\!\!\sum_{\scriptstyle p_2\atop
\!\!\!\!\!\!\!\!\!\!\!\!\!{\scriptstyle (p_1 p_2 , a)\sim 1}} 
\quad\ \ \sum\!\!\!\!\!\!\!\!
\sum_{\!\!\!\!\!\!\!\!\!\!\!\!\!
{\scriptstyle 0<|b_1|^2,|b_2|^2\leq\rho |a|^2\atop\scriptstyle 
b_1 p_2\equiv b_2 p_1\bmod a{\frak O}}}
\overline{\psi\bigl( a;b_2\bigr)}\,\vartheta\bigl( a;p_1\bigr)
\,\psi\bigl( a;b_1\bigr)\,\overline{\vartheta\bigl( a;p_2\bigr)}\;. 
 &(6.60)}$$

By (6.55), (6.45), (6.46) and the Cauchy-Schwarz inequality, we have 
$$\eqalign{ 
\sum_m\left|\tau_m\right|^2 
 &\leq\sum_{0<|m|^2\leq H_1 L_1}\Biggl(\,\sum_{\ell\mid m} 
\left|\Phi_{m/\ell}\right|\,\left| B_{\ell}\right|\Biggr)^2 \leq {}\cr 
 &\leq \sum_{0<|m|^2\leq H_1 L_1}\Biggl(\,\sum_{\ell\mid m} 
\left| O(1)\right|^2\Biggr)\Biggl(
\,\sum_{\ell\mid m}\left| B_{\ell}\right|^2\Biggr) \ll_{\varepsilon} {}\cr 
 &\ll_{\varepsilon} \left( H_1 L_1\right)^{\varepsilon} 
\sum_{\ell}\left| B_{\ell}\right|^2 
\sum_{\scriptstyle 0<|m|^2\leq H_1 L_1\atop\scriptstyle \ell\mid m} 1 
=\left( H_1 L_1\right)^{\varepsilon} 
\sum_{\ell}\left| B_{\ell}\right|^2 O\!\left( {H_1 L_1\over |\ell|^2}\right) 
=\left( H_1 L_1\right)^{\varepsilon} 
\sum_{\ell}\left| B_{\ell}\right|^2 O\left( H_1\right)\;,
}$$
and, given that $0<\rho <1/4$, 
$$\eqalign{
\ \sum\!\!\!\sum_{\!\!\!\!\!\!\!\!\!\!\!\!\!
{\scriptstyle p_1\quad p_2\atop\scriptstyle (p_1 p_2 , a)\sim 1}} 
\quad\ \ \sum\!\!\!\!\!\!\!\!
 &\sum_{\!\!\!\!\!\!\!\!\!\!\!\!\!
{\scriptstyle 0<|b_1|^2,|b_2|^2\leq\rho |a|^2\atop\scriptstyle 
b_1 p_2\equiv b_2 p_1\bmod a{\frak O}}}
 \overline{\psi\bigl( a;b_2\bigr)}\,\vartheta\bigl( a;p_1\bigr)
\,\psi\bigl( a;b_1\bigr)\,\overline{\vartheta\bigl( a;p_2\bigr)} = {}\cr
 &\qquad\quad\ =  
\qquad\,\sum\!\!\!\!\!\!\!\!\!
\sum_{\!\!\!\!\!\!\!\!\!\!\!\!\!\!\!
{\scriptstyle 0<|n_1|^2,|n_2|^2\leq\rho |a|^2 P_1\atop\scriptstyle 
n_2\equiv n_1\bmod a{\frak O}}}
\biggl(\sum_{p_1\mid n_1} O(1)\biggr) \biggl(\sum_{p_2\mid n_2} O(1)\biggr) = {}\cr 
 &\qquad\quad\ =\sum_{0<|n_1|^2\leq \rho |a|^2 P_1}
\sum_{\scriptstyle k\in{\frak O}\atop\scriptstyle |ak+n_1|^2\leq \rho |a|^2 P_1} 
O_{\varepsilon}\!\left(\bigl( P_1 |a|^2\bigr)^{\varepsilon}\right) \ll_{\varepsilon} {}\cr 
 &\qquad\quad\ \ll_{\varepsilon} 
\bigl( P_1 |a|^2\bigr)^{\varepsilon}\Biggl(\,\sum_{0<|n_1|^2\leq \rho |a|^2 P_1} 1\Biggr)  
\Biggl(\,\sum_{|k|^2\leq 4\rho P_1} 1\Biggr)  
=\bigl( P_1 |a|^2\bigr)^{\varepsilon} O\left(\rho |a|^2 P_1\right) 
\left( 1 + O\bigl(\rho P_1\bigr)\right) .
}$$
It therefore  follows by (6.59), (6.60) and (6.46) that 
$${\cal X}
\ll_{\varepsilon} \left( H_1 L_1\right)^{\varepsilon} 
\left( H_1 L_1 +{A_1^2\over \left| c_1\right|^2}\right) H_1 \sum_{\ell}\left| B_{\ell}\right|^2 
$$
and 
$${\cal Y}
=\sum_{\scriptstyle 0<|a|^2\leq A_1\atop\scriptstyle c_1\mid a} 
O\left( |a|^{-2}\right)\,O_{\varepsilon}\!\Bigl( 
\bigl( P_1 |a|^2\bigr)^{\varepsilon} \rho |a|^2 P_1\left( 1+\rho P_1\right)\Bigr) 
\ll_{\varepsilon} \bigl( P_1 A_1\bigr)^{\varepsilon} \left| c_1\right|^{-2} A_1 
\rho P_1\left( 1+\rho P_1\right) .
$$
By these bounds for 
${\cal X},{\cal Y}\in [0,\infty)$, and by (6.53) and (6.58), it follows that 
$$\eqalign{ 
{\cal E}^{*} 
 &\ll \Delta^{-3} {\cal X}^{1/2} {\cal Y}^{1/2} \ll_{\varepsilon} {}\cr 
 &\ll_{\varepsilon} \Delta^{-3} 
\bigl( H_1 L_1 P_1 A_1\bigr)^{\varepsilon /2}
\left( H_1 L_1 +{A_1^2\over \left| c_1\right|^2}\right)^{1/2} H_1^{1/2} 
\biggl(\sum_{\ell}\left| B_{\ell}\right|^2\biggr)^{1/2}  
\left| c_1\right|^{-1} A_1^{1/2} 
\rho^{1/2} P_1^{1/2}\left( 1+\rho P_1\right)^{1/2} , 
}$$
and so (given that $H_1 L_1 P_1 A_1\geq 1$) the result (6.50) is obtained
\ $\blacksquare$

\bigskip

\goodbreak\proclaim Lemma~6.3. Let the hypotheses of Lemma~6.1 be satisfied. Then  
$${\cal R}\ll_{\varepsilon} 
E^{11}\left( (HKL)^{1/2} +(P/K)Q+(P/K)^{3/2}(QHK)^{1/2}\right) 
\biggl( HK\sum_{\ell}\left|\Upsilon_{\ell}\right|^2\biggr)^{\!\!1/2} .\eqno(6.61)$$

\medskip 

\goodbreak 
\noindent{\bf Proof.}\quad 
By Lemma 6.1, it may be supposed that the inequalities in (6.10) hold, and that, for 
certain non-zero Gaussian integers $w,c,t,r,s,k,q$ satisfying the 
conditions in (6.11), one has the upper bound (6.12), with 
${\cal E}(w,c,t,r,s;k,q)$ as defined in (6.13)-(6.14): for otherwise 
one obtains the result in (6.61) as a consequence of the 
stronger bound (6.9) that Lemma~6.1 implies. 
Let ${\cal E}_1={\cal E}(w,c,t,r,s;k,q)$. 
Then, by (6.11), (6.13) and (6.14),   
$$\eqalignno{ 
{\cal E}_1  
 &=\sum_{\scriptstyle 0<|a_1|^2\leq A_1\atop\scriptstyle (a_1,w_1)\sim c_1} 
{1\over |a_1|^2} \sum_{h_1}\phi_{wh_1}  
\sum_{\ell_1}\Upsilon_{t\ell_1} 
\sum_{\scriptstyle p_1\neq 0\atop\scriptstyle (p_1,s)\sim 1} 
\theta_{rp_1} \left| p_1\right|^{-2}\,{\rm e}\!\left( {\rm Re}\!\left( 
{k a_1\over w_1 p_1 q_1}\right)\right)\times {} 
 &(6.62)\cr 
 &\qquad\quad\times {\rm e}\!\left( {\rm Re}\!\left( {h_1\ell_1\over a_1 p_1 q_1}\right)\right)  
f\left( h_1 , \ell_1 , p_1\right)
\sum_{0<|b|^2\leq\rho |a_1|^2} 
e\!\left( {\rm Re}\!\left( {q_1 b\over a_1}\right)\right) 
S\!\left( h_1\ell_1 p_1^{*} , b ; a_1\right)\;,}$$
where 
$$w_1,c_1,q_1\in{\frak O}-\{ 0\}\;,\qquad w_1=w/t\;,\qquad c_1=c/t\;,\quad\  
{\rm and}\quad\  q_1=q/s\;,\eqno(6.63)$$
$$A_1=|t|^{-2} A\quad\ {\rm and}\quad\ \rho =E^2 Q_1^{-1}\;,\eqno(6.64)$$
with 
$$A=(PQ)^{1+\varepsilon}\delta^{-2} K^{-1}<PQ    
\;,\qquad E=\bigl( 1+\delta^{-1}\bigr)(PQ)^{\varepsilon}\quad\  
{\rm and}\quad\ Q_1=|s|^{-2} Q\;,\eqno(6.65)$$
while 
$$f({\bf z})=\varphi\left( w z_1 , k , t z_2 , r z_3 , q\right)\qquad\quad  
\hbox{for $\quad {\bf z}\in{\Bbb C}^3$.}\eqno(6.66)$$
Given that we have (6.12), the completion of this proof requires only 
a sufficiently strong upper bound for $|{\cal E}_1|=|{\cal E}(w,c,t,r,s;k,q)|$. 
We show next that such a bound may be deduced from Lemma~6.2. 

Let $w_1,c_1,q_1\in{\frak O}-\{ 0\}$, $A_1,\rho,A,E,Q_1>0$ and 
$f : {\Bbb C}^3\rightarrow{\Bbb C}$ be as stated in (6.63)-(6.66). 
Then, by (6.66) and (6.2), the condition~(6.47) in Lemma~6.2 is satisfied 
by $f$ when one has there:  
$$H_1 = |w|^{-2} H\;,\qquad L_1 = |t|^{-2} L\quad\ 
{\rm and}\quad\  P_1 = |r|^{-2} P\;.\eqno(6.67)$$ 
Since $\varphi\in{\cal S}\bigl({\Bbb C}^5\bigr)$, and since one has, for $m=1,2,3$, 
$${\partial\over\partial z_m}
+{\partial\over\partial \overline{z_m}}={\partial\over\partial x_m}\qquad\   
{\rm and}\qquad\ {\partial\over\partial z_m}-{\partial\over\partial \overline{z_m}}
=-i{\partial\over\partial y_m}$$ 
(with $x_m={\rm Re}(z_m)$, $y_m={\rm Im}(z_m)$ 
and $\partial/\partial z_m$, $\partial/\partial \overline{z_m}$ 
the linear operators defined in (5.19)), it 
follows by (6.66) and the pair of equations from which (5.20) is deduced  
that $f$ lies in the Schwartz space ${\cal S}\bigl( {\Bbb C}^3\bigr)$. 
Moreover, by (6.66), (5.20), (5.12) and (6.3) (and since 
${\cal L}_m$, as defined in (5.12), is a linear 
operator on ${\cal S}\bigl( {\Bbb C}^n\bigr)$), one has 
$$\eqalign{
\bigl( {\cal L}_1^{j_1}
{\cal L}_2^{j_2}
{\cal L}_3^{j_3} f\bigr) ({\bf z}) 
 &=|w|^{2 j_1} |t|^{2 j_2} |r|^{2 j_3} 
\bigl( {\cal L}_1^{j_1}
{\cal L}_3^{j_2}
{\cal L}_4^{j_3} \varphi\bigr)\bigl( w z_1 , k , t z_2 , r z_3 , q\bigr) \cr 
 &=|w|^{2 j_1} |t|^{2 j_2} |r|^{2 j_3}  
\,O_{\bf j}\left( \bigl(\delta\bigl| w z_1\bigr|\bigr)^{-2 j_1} 
\bigl(\delta\bigl| t z_2\bigr|\bigr)^{-2 j_2} 
\bigl(\delta\bigl| r z_3\bigr|\bigr)^{-2 j_3} \right) 
}$$
for ${\bf j}\in\bigl( {\Bbb N}\cup\{ 0\}\bigr)^3$,  
$\,{\bf z}\in\bigl( {\Bbb C}^{*}\bigr)^3$; 
and so $f$ satisfies the condition~(6.48) of Lemma~6.2 when 
one has there 
$$\Delta =\delta^2\;.\eqno(6.68)$$ 

Given the conclusions reached in the last paragraph, 
it is only the factor ${\rm e}\bigl( {\rm Re}(h_1\ell_1 /a_1 p_1 q_1)\bigr)$,  
occurring in the sum on the right-hand side of the equation~(6.62),  
that prevents us from bounding ${\cal E}_1$  
by the direct application of Lemma~6.2. This factor is, however, 
an essentially trivial obstacle to the application of Lemma~6.2. For, 
by (6.2), (6.8), (6.11) and (6.63)-(6.67), 
the summation on the right-hand side 
of the equation~(6.62) is effectively constrained to points 
$(a_1,h_1,\ell_1,p_1,b)\in\bigl( {\frak O}-\{ 0\}\bigr)^5$ such that 
$${\left| h_1 \ell_1\right|^2\over A_1 P_1 Q_1} 
< \left|{h_1\ell_1\over a_1 p_1 q_1}\right|^2 
<{H_1 L_1\over\left| a_1\right|^2\left( P_1 /2\right) \left( Q_1 /2\right)}
={4HL\over\left| t a_1\right|^2 PQ}\leq {4HL\over PQ}\ll 1\;,\eqno(6.69)$$
and so the factor  ${\rm e}\bigl( {\rm Re}(h_1\ell_1 /a_1 p_1 q_1)\bigr)$  
may be very well approximated by a partial sum of just 
$O(1)$ terms from the product of Taylor series: 
$${\rm e}\!\left( {\rm Re}\!\left( {h_1\ell_1\over a_1 p_1 q_1}\right)\right) 
=\exp\!\left(  \pi i\,{h_1\ell_1\over a_1 p_1 q_1}\right) 
\exp\!\left(\pi i\,\overline{h_1\ell_1\over a_1 p_1 q_1}\ \right) 
=\sum_{m=0}^{\infty}\sum_{n=0}^{\infty} 
{(\pi i)^{m+n}\over (m!)(n!)}\left( {h_1\ell_1\over a_1 p_1 q_1}\right)^{\!\!m} 
\left( \,\overline{h_1\ell_1\over a_1 p_1 q_1}\ \right)^{\!\!n}\;.$$
By employing this last expansion of the offending factor 
in (6.62), and then making the trivial substitutions 
$a_1=a$, $h_1=h$, $\ell_1 =\ell$ and $p_1=p$, 
one obtains an absolutely convergent sum over 
$a\in{\frak O}-\{ 0\}$, $h,\ell\in{\frak O}$, $p\in{\frak O}-\{ 0\}$, 
$m,n\in{\Bbb N}\cup\{ 0\}$ and 
$b\in{\frak O}-\{ 0\}$ (in that order). Any change in the 
order of summation can be justified, so that one has, in particular, 
$${\cal E}_1 
=\sum_{m=0}^{\infty}\sum_{n=0}^{\infty} 
{(\pi i)^{m+n}\over (m!)(n!)}
\bigl( P_1 /2\bigr)^{-1}\left( {4 H_1 L_1\over P_1 Q_1}\right)^{\!\!(m+n)/2}
{\cal E}_1^{*}(m,n)\;,$$ 
with 
$${\cal E}_1^{*}(m,n)
=\sum_{\scriptstyle 0<|a|^2\leq A_1\atop\scriptstyle c_1\mid a} 
\xi_{a}^{m,n} \sum_{h} \Phi_{h}^{m,n} \sum_{\ell} B_{\ell}^{m,n} 
\sum_{p} \Theta_{m,n}\left( a ; p\right) 
f\left( h , \ell, p\right) \sum_{0<|b|^2\leq\rho \left| a\right|^2}
\psi_{m,n}\left( a ; b\right)  
S\left( h \ell p^{*} , b ; a\right)\;,$$
where, for $h,\ell\in{\frak O}$, $\,p,b\in{\frak O}-\{ 0\}$, 
$\,a\in c_1{\frak O}-\{ 0\}$, 
$$\eqalignno{
\psi_{m,n}(a;b) &={\rm e}\!\left({\rm Re}\!\left( {q_1 b\over a}\right)\!\right) ,\quad 
\Theta_{m,n}(a;p) 
=\cases{\!{\displaystyle{(P_1 /2)^{(m+n+2)/2}\over p^{m+1} (\overline{p})^{n+1}}}   
\,{\rm e}\!\left( {\rm Re}\!\left( {\displaystyle{k a\over w_1 q_1 p}}\right)\!\right)
\!\theta_{rp} &if $(p , s)\sim 1$, \cr \quad &\quad \cr 0 &otherwise,}\qquad\quad\ 
 &(6.70)\cr \quad & \cr}$$ 
$$B_{\ell}^{m,n} ={\ell^m (\overline{\ell})^n \Upsilon_{t\ell}\over L_1^{(m+n)/2}}\;,\ \,   
\Phi_h^{m,n} ={h^m (\overline{h})^n \phi_{w h}\over H_1^{(m+n)/2}}\ \      
{\rm and}\ \ \xi_{a}^{m,n} =\cases{{\displaystyle{\left( Q_1 /2\right)^{(m+n)/2}\over 
|a|^2 \left( a q_1\right)^m \left( \overline{a q_1}\right)^n}}   
 &if $\left( {\displaystyle {a\over c_1}\,,\,{w_1\over c_1}}\right)\sim 1$, \cr \quad &\quad \cr 
0 &otherwise.}
\eqno(6.71)$$

In the above definition of ${\cal E}_1^{*}(m,n)$ 
the summation over $a,h,\ell,p,b$ is subject to 
the same effective constraint (6.69) as applied to the 
summation in (6.62)$\,$ (the factor $f(h_1,\ell_1,p_1)$ being present 
in both cases). Since the number of points $(a_1,h_1,\ell_1,p_1,b)\in\bigl( {\frak O}-\{ 0\}\bigr)^5$ 
satisfying the first two inequalities in (6.69) is finite, and since, for all 
such points, one has 
$$\left( {h_1\ell_1\over a_1 p_1 q_1}\right)^{\!\!m} 
\left( \,\overline{h_1\ell_1\over a_1 p_1 q_1}\ \right)^{\!\!n}
\left( {\left( P_1 /2\right) \left( Q_1 /2\right)\over H_1 L_1}\right)^{\!\!(m+n)/2} 
\longrightarrow 0\qquad\quad\hbox{as $\quad (m+n)\longrightarrow +\infty\;$,}$$
it therefore must be the case that, for some pair $m,n\in{\Bbb N}\cup\{ 0\}$, one has:  
$$\left|{\cal E}_1\right|\leq \tau\left|{\cal E}_1^{*}(m,n)\right|\;,\eqno(6.72)$$ 
where  
$$\eqalignno{
\tau 
 &=\sum_{\mu =0}^{\infty}\sum_{\nu =0}^{\infty} 
{\pi^{\mu +\nu}\over (\mu !)(\nu !)}\left( P_1 /2\right)^{-1}\left( {H_1 L_1\over
\left( P_1 /2\right) \left( Q_1 /2\right)}\right)^{\!\!(\mu +\nu)/2} = {}  
\qquad\qquad\qquad\qquad\qquad\qquad\qquad\quad\cr 
 &={2\over P_1}\left( 
\exp\left( 2\pi\left( {H_1 L_1\over P_1 Q_1}\right)^{\!\!1/2}\right)\right)^{\!\!2} 
={2\over P_1}\exp\left( 
4\pi\left( {H_1 L_1\over P_1 Q_1}\right)^{\!\!1/2}\right) 
= {2\over P_1}\exp\left( O(1)\right)\ll {1\over P_1}
&(6.73)}$$
(with the upper bound $O(1)$ used here following by (6.69), since 
$0\neq t\in{\frak O}$ implies $|t|\geq 1$). 
 
Let $m,n\in{\Bbb N}\cup\{ 0\}$ be one of the pairs for which one has 
(6.72)-(6.73). Then, given that $m,n\geq 0$, and that 
$w,c,t,r,s,k,q$ are non-zero Gaussian 
integers satisfying the conditions in~(6.11), it follows 
by (6.4), (6.5), (6.63) and (6.65) that 
the conditions~(6.45) and (6.46) of 
Lemma~6.2 are satisfied
when $H_1,L_1,P_1$ are as in (6.67) and   
$\xi_{a}=\xi_{a}^{m,n}$, 
$\Phi_{h}=\Phi_{h}^{m,n}$, $B_{\ell}=B_{\ell}^{m,n}$,   
$\Theta(a;p)=\Theta_{m,n}(a;p)$ and  
$\psi(a;b)=\psi_{m,n}(a;b)\,$   
(for $a\in c_1{\frak O}-\{ 0\}$, $\,h,\ell\in{\frak O}$, 
$\,p\in{\frak O}-\{ 0\}$, $\,m,n\in{\Bbb N}\cup\{ 0\}$  
and $b\in{\frak O}-\{ 0\}$), with  
$\xi_{a}^{m,n}$, 
$\Phi_{h}^{m,n}$, $B_{\ell}^{m,n}$,   
$\Theta_{m,n}(a;p)$ and  
$\psi_{m,n}(a;b)$ as defined in (6.70) and (6.71). 
Since we already verified that (with $H_1,L_1,P_1$ as in (6.67), and 
$\Delta$ as in (6.68)) the conditions 
(6.47) and (6.48) of Lemma~6.2 are satisfied by $f$, 
we may therefore  apply Lemma~6.2 with 
the coefficients 
$\xi_{a}$, 
$\Phi_{h}$, $B_{\ell}$,   
$\Theta(a;p)$ and  
$\psi(a;b)$ as just indicated, and with 
$c_1$, $A_1$ and $\rho$ given by (6.63)-(6.65), and  
$f\in{\cal S}\bigl( {\Bbb C}^3\bigr)$ given by (6.66).  
Moreover, in respect of this particular application of 
Lemma~6.2, the term ${\cal E}^{*}$ is, by the equation~(6.49), 
evidently equal to the term ${\cal E}_1^{*}(m,n)$ that we 
defined earlier in this proof (i.e. the definitions 
of ${\cal E}^{*}$ and ${\cal E}_1^{*}(m,n)$ coincide). Consequently, 
by the upper bound for $|{\cal E}^{*}|$ in the result (6.50) of Lemma~6.2, 
one has 
$${\cal E}_1^{*}(m,n) 
\ll_{\varepsilon} \Delta^{-3} \left( H_1 L_1 P_1 A_1\right)^{\varepsilon} 
\left( 1+{H_1 L_1 \left| c_1\right|^2\over A_1^2}\right)^{1/2} 
\left( 1 +\rho P_1\right)^{1/2} 
\left( \left| c_1\right|^{-4} \rho A_1^3 P_1 H_1 
\sum_{\ell}\left| B_{\ell}^{m,n}\right|^2\right)^{1/2} . 
$$ 
Since $w,t,r$ are non-zero Gaussian integers, it follows 
by (6.68), (6.67), (6.64), (6.65) and (6.8) that we have 
here
$$\Delta^{-3}\left( H_1 L_1 P_1 A_1\right)^{\varepsilon} 
\leq\delta^{-6} (HLPA)^{\varepsilon} 
\ll\delta^{-6} (PQ)^{3\varepsilon}<(PQ)^{-\varepsilon}\delta^{-2}E^4\;.$$
Moreover, by using (in addition) (6.63), (6.11), (6.71) and (6.5), 
one finds that 
$${H_1 L_1 \left| c_1\right|^2\over A_1^2}
={HL |c|^2\over A^2 |w|^2} 
={\delta^4 K^2 H L |c|^2\over (PQ)^{2+2\varepsilon} |w|^2} 
\leq {\delta^4 H K^2 L |c|^2\over P^2 Q^2 |w|^2}\;,\qquad\qquad  
\rho P_1 =E^2 Q_1^{-1} P_1={E^2 P |s|^2\over Q |r|^2}\;,$$
$$\left| c_1\right|^{-4} \rho A_1^3 P_1 H_1 
={|t|^4 E^2 |s|^2\over |c|^4 Q}
\left( {(PQ)^{1+\varepsilon}\over\delta^2 K |t|^2}\right)^3 
{PH\over |r|^2 |w|^2}
={(PQ)^{3\varepsilon}\delta^{-6} E^2 P^4 Q^2 H\over K^3 |c|^4 |t|^2 |r|^4} 
\leq {E^8 P^4 Q^2 H\over K^3 |c|^4 |t|^2 |r|^4} $$
and
$$\sum_{\ell}\left| B_{\ell}^{m,n}\right|^2
\leq\sum_{{\textstyle{L_1\over 2}}<|\ell|^2\leq L_1}
\left( {|\ell|^{2}\over L_1}\right)^{\!\!m+n}\left|\Upsilon_{t\ell}\right|^2 
\leq\sum_{\ell}\left|\Upsilon_{t\ell}\right|^2\leq\sum_{\ell}\left|\Upsilon_{\ell}\right|^2\;.$$
Therefore the bound that we have obtained for ${\cal E}_1^{*}(m,n)$ 
implies: 
$${\cal E}_1^{*}(m,n) 
\ll_{\varepsilon} (PQ)^{-\varepsilon} \delta^{-2} E^{9} 
\left( 1+{\delta^2 H^{1/2} K L^{1/2} |c|\over PQ|w|}\right) 
\left( 1+{P^{1/2} |s|\over Q^{1/2} |r|}\right) 
{P^2 Q H^{1/2}\over K^{3/2} |c|^2 |t| |r|^2} 
\left(\sum_{\ell} \left|\Upsilon_{\ell}\right|^2\right)^{1/2}\;.$$

Since we have the 
bound (6.12) (in Lemma~6.1), where ${\cal E}(w,c,t,r,s;k,q)={\cal E}_1$, 
it follows that 
$${\cal R}\ll_{\varepsilon} 
(PQ)^{\varepsilon} K |t|^2 \left| {\cal E}_1\right|\;.$$
Therefore, using (6.72), (6.73), (6.67) and the last bound obtained 
for ${\cal E}_1^{*}(m,n)$, we deduce that 
$$\eqalign{
{\cal R} 
 &\ll_{\varepsilon} (PQ)^{\varepsilon}  K |t|^2 P^{-1} |r|^2 \left| {\cal E}_1^{*}(m,n)\right| 
 \ll_{\varepsilon} {}\cr 
 &\ll_{\varepsilon} \delta^{-2} E^{9} 
\left( 1+{\delta^2 H^{1/2} K L^{1/2} |c|\over PQ|w|}\right) 
\left( 1+{P^{1/2} |s|\over Q^{1/2} |r|}\right) 
{P Q H^{1/2} |t|\over K^{1/2} |c|^2} 
\left(\sum_{\ell} \left|\Upsilon_{\ell}\right|^2\right)^{1/2} .
}$$
Moreover, given (6.65), and 
the conditions~(6.11) satisfied by $w$, $c$, $t$, $r$ and $s$, 
one has  
$$\eqalign{
\delta^{-2}\left( 1+{\delta^2 H^{1/2} K L^{1/2} |c|\over PQ|w|}\right) 
\left( 1+{P^{1/2} |s|\over Q^{1/2} |r|}\right) 
{|t|\over |c|^2} 
 &\leq \left( \delta^{-2}+{H^{1/2} K L^{1/2} |c|\over PQ|w|}\right) 
\left( 1+{P^{1/2} |w|\over Q^{1/2}}\right) 
{1\over |c|} < {}\cr
 &<\left( {E^2\over |c|}+{H^{1/2} K L^{1/2}\over PQ|w|}\right) 
\left( 1+{P^{1/2} |w|\over Q^{1/2}}\right) \leq {}\cr 
 &\leq E^2 + {E^2 P^{1/2} |w|\over Q^{1/2}} + {H^{1/2} K L^{1/2}\over PQ}
+{H^{1/2} K L^{1/2}\over P^{1/2}Q^{3/2}} \leq {}\cr 
 &\leq E^2 + {E^2 P^{1/2} H^{1/2}\over Q^{1/2}} + {H^{1/2} K L^{1/2}\over PQ} 
+\left( {HK\over Q}\right)\left( {L\over Q}\right)^{1/2} .
}$$
where, by (6.8), one has $HK/Q=O(1)$ and $L/Q=O(1)$. Consequently we may deduce that 
$${\cal R}\ll_{\varepsilon} 
E^{11}\left( 1+{P^{1/2} H^{1/2}\over Q^{1/2}} + {H^{1/2} K L^{1/2}\over PQ}\right) 
{P Q H^{1/2}\over K^{1/2}} 
\left(\sum_{\ell} \left|\Upsilon_{\ell}\right|^2\right)^{1/2} .
$$
Since $PQ H^{1/2}/K^{1/2}=(PQ/K)(HK)^{1/2}$, while 
$$\left( 1+{P^{1/2} H^{1/2}\over Q^{1/2}} + {H^{1/2} K L^{1/2}\over PQ}\right) 
{PQ\over K}
={PQ\over K}+\left( {P\over K}\right)^{3/2} (QHK)^{1/2} +(HL)^{1/2}\;,$$
the result (6.61) therefore follows\ $\blacksquare$  

\bigskip 
\bigskip 

\goodbreak\centerline{\bf \S 7. Switching to Levels of Greater Modulus}

\bigskip 

Lemma~7.3 below shows that the mean value $S_t(Q,X,N)$ is, 
in a certain sense, `approximately' a monotonic non-decreasing 
function of the level related parameter $Q$.  
This result (the inequality (7.31) below) has an important application in the next section, 
where it enables us to work around the lower bound condition $Q\gg HK$ in (6.8);  
that condition would otherwise adversely limit our use of Lemma~6.3. 

We prove Lemma~7.3 by deducing it (via elementary arguments) from the 
simpler bound given by Lemma~7.2. For the proof of Lemma~7.2 
we need the results of Lemma~7.1, below. 

\bigskip 

\goodbreak\proclaim Lemma~7.1. Let $0\neq q\in{\frak O}$ and $0\neq r\in q{\frak O}$; 
let $\Gamma=\Gamma_0(q)\leq SL(2,{\frak O})$ 
and $\tilde\Gamma=\Gamma_0(r)\leq SL(2,{\frak O})$; 
let $V'$ and $V''$ be amongst the cuspidal subspaces $V$ 
occurring as factors in the decomposition (1.1.3) of ${}^{0}L^2(\Gamma\backslash G)$; 
and let $f_{0,0}^{V'}$ and $f_{0,0}^{V''}$ be (as in (1.1.6)) generators of 
the corresponding spaces $V_{K,0,0}'$ and $V_{K,0,0}''$, normalised so as to satisfy 
(1.1.9). Suppose, moreover, that 
$p_{V''}=p_{V'}=0$ and $\nu_{V''}=\nu_{V'}=\nu$ (say). 
Then $\Gamma\geq\tilde\Gamma$, $\bigl[\Gamma :\tilde\Gamma\bigr]<\infty$,  
and 
$f_{0,0}^{V'},f_{0,0}^{V''}\in L^2(\Gamma\backslash G)\cap L^2(\tilde\Gamma\backslash G;0,0)$ 
(where the latter space is that given by the case $\Gamma=\tilde\Gamma$ of (1.1.20)); 
the functions $f_{0,0}^{V'}$ and $f_{0,0}^{V''}$ are bounded and continuous on $G$, 
and one has 
$$\bigl\| f_{0,0}^{V'}\bigr\|_{\tilde\Gamma\backslash G}^2
=\bigl[\Gamma :\tilde\Gamma\bigr]
\bigl\| f_{0,0}^{V'}\bigr\|_{\Gamma\backslash G}^2\eqno(7.1)$$ 
and 
$$\sum_{\scriptstyle W\atop\scriptstyle (\nu_W,p_W)=(\nu,0)}^{(\tilde\Gamma)} 
{1\over \bigl\| f_{0,0}^{W}\bigr\|_{\tilde\Gamma\backslash G}^2}
\left\langle f_{0,0}^{V'},f_{0,0}^W\right\rangle_{\tilde\Gamma\backslash G} 
\left\langle f_{0,0}^W,f_{0,0}^{V''}\right\rangle_{\tilde\Gamma\backslash G}
=\cases{\bigl\| f_{0,0}^{V'}\bigr\|_{\tilde\Gamma\backslash G}^2 &if $V''=V'$,\cr 
0 &otherwise,}\eqno(7.2)$$
where the meaning of the bracketed `$\,\tilde\Gamma$' is that the summation is 
restricted to irreducible cuspidal subspaces $W$ of 
${}^{0}L^2(\tilde\Gamma\backslash G)$ (with the equations~(1.1.6) and~(1.1.9), 
as they apply when $\Gamma =\tilde\Gamma$, 
determining to within a constant factor of unit modulus 
the $\tilde\Gamma$-automorphic function  
$f_{0,0}^W : G\rightarrow{\Bbb C}$). 

\medskip

\goodbreak 
\noindent{\bf Proof.}\quad 
Let $q$, $r$, $\Gamma$, $\tilde\Gamma$, $V'$, $V''$, $\nu$, $f_{0,0}^{V'}$ and 
$f_{0,0}^{V''}$ satisfy the hypotheses of the lemma. 
Then, since $r\in q{\frak O}$, 
the congruence $c\equiv 0\bmod r{\frak O}$ implies $c\equiv 0\bmod q{\frak O}$. Therefore, 
it is a trivial corollary of the definition of $\Gamma_0(q)$ given in Subsection~1.1 
that we have here $\tilde\Gamma\leq\Gamma$. 
It follows that the $\Gamma$-automorphic functions $f_{0,0}^{V'},f_{0,0}^{V''} : G\rightarrow{\Bbb C}$ 
are {\it a fortiori\/} also $\tilde\Gamma$-automorphic. 
Furthermore, since any set of right-coset representatives 
$\gamma_1,\gamma_2,\ldots ,\gamma_{[SL(2,{\frak O}):\tilde\Gamma]}\in SL(2,{\frak O})$ 
for the quotient $\tilde\Gamma\backslash SL(2,{\frak O})$ is a union 
of $[\Gamma : \tilde\Gamma]$ sets of coset representatives for 
$\Gamma\backslash SL(2,{\frak O})$, one has 
$$\int_{\tilde\Gamma\backslash G} f(g) {\rm d}g
=[\Gamma : \tilde\Gamma]\int_{\Gamma\backslash G} f(g) {\rm d}g\;,\eqno(7.3)$$
for any measurable $\Gamma$-automorphic function $f : G\rightarrow{\Bbb C}$ 
such that the latter integral exists; 
by the pairwise orthogonality of the irreducible cuspidal subspaces 
$V\subset{}^{0}L^2(\Gamma\backslash G)\subset L^2(\Gamma\backslash G)$ 
occurring as factors in the 
decomposition (1.1.3),
one has, in particular: 
$$\left\langle f_{0,0}^{V'},f_{0,0}^{V''}\right\rangle_{\tilde\Gamma\backslash G} 
=[\Gamma : \tilde\Gamma]\left\langle f_{0,0}^{V'},f_{0,0}^{V''}\right\rangle_{\Gamma\backslash G} 
=\cases{[\Gamma : \tilde\Gamma]\bigl\| f_{0,0}^{V'}\bigr\|_{\Gamma\backslash G}^2 &if $V''=V'$, \cr 
0 &otherwise,}\eqno(7.4)$$
which contains the result (7.1).

Since $\tilde\Gamma\leq\Gamma\leq SL(2,{\frak O})$, and since 
one has (see (7.39) below) 
$\,[SL(2,{\frak O}) : \tilde\Gamma]
=[SL(2,{\frak O}) : \Gamma_0(r)]<\infty$, 
the index $[\Gamma : \tilde\Gamma]$ is certainly finite. 
We have, moreover, 
$f_{0,0}^{V'}\in V_{K,0,0}'\subset V'\subset 
{}^{0}L^2(\Gamma\backslash G)\subset L^2(\Gamma\backslash G)\,$  
(by (1.1.2), (1.1.3), (1.1.5) and (1.1.6)), so it follows 
from what has so far been established that $f_{0,0}^{V'}\in L^2(\tilde\Gamma\backslash G)$; 
by the observations between (1.1.5) and (1.1.6), 
we have also $\Omega_K f_{0,0}^{V'}=0$ and $({\partial /\partial\psi})f_{0,0}^{V'}=0$, 
and so may deduce that 
$f_{0,0}^{V'}$ lies in the space 
$L^2(\tilde\Gamma\backslash G ; 0,0)$ defined by (1.1.20).

Turning now to the proof of (7.3) we 
seek to apply the Parseval identity [22, Theorem~A], 
taking there $\ell =q=0$, $\Gamma =\tilde\Gamma$ 
and $f_1=f_{0,0}^{V'}$, $f_2=f_{0,0}^{V''}$: 
in this case the hypotheses of 
[22, Theorem~A] require only that 
$f_1$ and $f_2$ lie in the space $L^2(\tilde\Gamma\backslash G ;0,0)$ 
and are both bounded and smooth 
(possessing continuous partial derivatives of all orders, 
with respect to $x,y,r,\theta,\varphi,\psi$, where 
$(x+iy,r,\theta,\varphi,\psi)$ are the Iwasawa coordinates for $G$ 
described in Subsection~1.1). This need only be verified for $f_1=f_{0,0}^{V'}$, 
since similar conclusions will apply to $f_2=f_{0,0}^{V''}$ 
(given the symmetry in our hypotheses concerning $f_{0,0}^{V'}$ and $f_{0,0}^{V''}$). 

For the smoothness property see, for example, [22, Relation (1.7.10),  
Definitions (1.2.2), (1.4.5)-(1.4.7)] (and the accompanying justification).   
Were $\tilde\Gamma =\Gamma_0(r)$ a cocompact subgroup of $G$, the smoothness 
would imply the boundedness; since, however, the fundamental domain 
${\cal F}_{\tilde\Gamma\backslash\Gamma}\subset G$ is non-compact, we  
need the growth condition (1.1.10) in order to prove the boundedness 
of $f_{0,0}^{V'}$. A short calculation shows that, for each $\gamma_k\in SL(2,{\frak O})$ 
featuring in the description of ${\cal F}_{\Gamma\backslash G}$ in Subsection~1.1, 
there exists a $\delta=\delta(\gamma_k)>0$ such that if 
${\frak c}={\frak c}(\gamma_k)$ is the cusp $\gamma_k\infty$ then one has 
$$g_{\frak c}^{-1}\gamma_k N a[r] K\subseteq N a[\delta r] K\qquad\qquad 
\hbox{($r>0$).}\eqno(7.5)$$
Given that the fundamental domain ${\cal F}_{\tilde\Gamma\backslash G}$  
is similar in description to ${\cal F}_{\Gamma\backslash G}$, 
it therefore follows   
(since we have $[SL(2,{\frak O}) : \tilde\Gamma]<\infty$ and  
$r^{1/2}e^{-\pi r}\leq (2\pi e)^{-1/2}$, for $r>0\,$)  
that the application of (1.1.10) for a finite number of cusps ${\frak c}$ 
suffices to show that $f_{0,0}^{V'}$ is bounded on a set 
${\cal F}_{\tilde\Gamma\backslash G}^{\infty}\subseteq {\cal F}_{\tilde\Gamma\backslash G}$ 
such that ${\cal F}_{\tilde\Gamma\backslash G}-{\cal F}_{\tilde\Gamma\backslash G}^{\infty}$ 
is compact. By the smoothness of $f_{0,0}^{V'}$, 
the function $f_{0,0}^{V'}$ is also bounded on the latter (compact) set, and so 
is bounded on the set 
${\cal F}_{\tilde\Gamma\backslash G}^{\infty}\cup\bigl( 
{\cal F}_{\tilde\Gamma\backslash G}-{\cal F}_{\tilde\Gamma\backslash G}^{\infty}\bigr)
={\cal F}_{\tilde\Gamma\backslash G}\,$. Therefore, 
with $f_{0,0}^{V'}$ being $\tilde\Gamma$-automorphic, and 
${\cal F}_{\tilde\Gamma\backslash G}\subset G$ a fundamental domain for 
the action of $\tilde\Gamma$ on~$G$, we may conclude that $f_{0,0}^{V'}$ is bounded on $G$.

By the above we have verified that the case $\Gamma =\tilde\Gamma$, $\ell =q=0$ of 
[22, Theorem~A] may be applied with $f_1=f_{0,0}^{V'}$, $f_2=f_{0,0}^{V''}$: 
note that the transform `$T_{V}\varphi_{\ell,q}(\nu_V,p_V)$' which appears 
on the right-hand side of [22, Equation~(1.8.7)] is  
that function (or `generator') which we  refer to in (1.1.6) and (1.1.8)-(1.1.10)  
as `$f_{\ell,q}^V$'. 
In the case that concerns us, 
that theorem shows firstly that, for all cusps ${\frak c}$ of $\tilde\Gamma$,  
the function 
$t\mapsto\langle f_{0,0}^{V'} , E_{0,0}^{\frak c}(it,0)\rangle_{\tilde\Gamma\backslash G}$ 
(with $E_{0,0}^{\frak c}(it,0)$ given, for all real $t$, by the case 
case $\Gamma=\tilde\Gamma$ of the definition (1.1.12)) 
is in the space $L^2(-\infty,\infty)$; 
secondly it shows (given (1.1.6), (1.1.9) and [22, (1.7.8), (1.7.14), 
(1.6.7) and (1.6.8)]) that 
$$\eqalignno{
\left\langle f_{0,0}^{V'} , f_{0,0}^{V''}\right\rangle_{\tilde\Gamma\backslash G} 
 &={\left\langle f_{0,0}^{V'} , \vec 1\right\rangle_{\tilde\Gamma\backslash G}
\left\langle \vec 1, f_{0,0}^{V''}\right\rangle_{\tilde\Gamma\backslash G}\over 
\| \vec 1\|_{\tilde\Gamma\backslash G}^2}  
+\sum_{\scriptstyle W\atop\scriptstyle p_W=0}^{(\tilde\Gamma)} 
{\left\langle f_{0,0}^{V'} , f_{0,0}^W\right\rangle_{\tilde\Gamma\backslash G}
\left\langle f_{0,0}^W , f_{0,0}^{V''}\right\rangle_{\tilde\Gamma\backslash G}\over 
\left\| f_{0,0}^W\right\|_{\tilde\Gamma\backslash G}^2} + {} 
 &(7.6)\cr
 &\quad\   
+\sum_{{\frak c}\in {\frak C}}^{(\tilde\Gamma)} 
{\bigl[\tilde\Gamma_{\frak c} : \tilde\Gamma_{\frak c}'\bigr]\over 4\pi i} 
\int\limits_{(0)} 
\left\langle f_{0,0}^{V'} , E_{0,0}^{\frak c}(\nu ,0)\right\rangle_{\tilde\Gamma\backslash G}
\left\langle E_{0,0}^{\frak c}(\nu ,0) , f_{0,0}^{V''}\right\rangle_{\tilde\Gamma\backslash G}\, 
{\rm d}\nu\;,}$$
where `$\vec 1$' denotes the constant function 
$\varphi_{0,0}(-1,0) : G\rightarrow\{ 1\}\subset{\Bbb R}$. 
Here, since $f_{0,0}^{V'}\in{}^{0}L^2(\Gamma\backslash G)\subset L^2(\Gamma\backslash G)$, 
it follows by (7.3) and the orthogonality of the 
subspaces ${}^{0}L^2(\Gamma\backslash G)$ and ${\Bbb C}={\Bbb C}\vec 1$ 
in (1.1.2) that  we have 
$$\bigl\langle f_{0,0}^{V'} , \vec 1\,\bigr\rangle_{\tilde\Gamma\backslash G}=0\;.\eqno(7.7)$$ 
Hence the first term on the right-hand side of the equation~(7.6) equals zero. 
By a somewhat more roundabout argument we shall next show that the 
terms of the sum over ${\frak c}$ in (7.6) also vanish. 

Suppose that $f_{0,0}^{V'}\in{}^{0}L^2(\tilde\Gamma\backslash G)$. Then it follows 
by the orthogonality of the subspaces ${}^{0}L^2(\tilde\Gamma\backslash G)$
and ${}^{\rm e}L^2(\tilde\Gamma\backslash G)$, and by the case $\Gamma =\tilde\Gamma$ of 
(1.1.19), combined with the square integrability (mentioned before (7.6)) of 
the function 
$t\mapsto\langle f_{0,0}^{V'} , E_{0,0}^{\frak c}(it,0)\rangle_{\tilde\Gamma\backslash G}$, 
that one will have also, for all cusps ${\frak c}$ of $\tilde\Gamma$, 
$$\eqalignno{
0 &=\biggl\langle\ \int\limits_{(0)} 
\left\langle f_{0,0}^{V'} , E_{0,0}^{\frak c}(\nu ,0)\right\rangle_{\tilde\Gamma\backslash G} 
E_{0,0}^{\frak c}(\nu ,0)\,{\rm d}\nu\, ,\, 
f_{0,0}^{V''}\biggr\rangle_{\!\tilde\Gamma\backslash G} = {}\cr 
 &=\int\limits_{(0)} 
\left\langle f_{0,0}^{V'} , E_{0,0}^{\frak c}(\nu ,0)\right\rangle_{\tilde\Gamma\backslash G}
\left\langle E_{0,0}^{\frak c}(\nu ,0) , f_{0,0}^{V''}\right\rangle_{\tilde\Gamma\backslash G}\, 
{\rm d}\nu\;. &(7.8)}$$
We therefore now seek to establish the validity of the premise here 
(that $f_{0,0}^{V'}\in{}^{0}L^2(\tilde\Gamma\backslash G)$). Firstly, we may 
note that by (1.1.6) and [22, (1.7.10)], 
$$f_{0,0}^{V'}\in A_{\Gamma}^0(\Upsilon_{\nu,0};0,0)
\subseteq A_{\Gamma}^{\rm pol}(\Upsilon_{\nu,0};0,0)
\subseteq A_{\Gamma}(\Upsilon_{\nu,0};0,0)\;,\eqno(7.9)$$ 
where the latter three sets are the subspaces of ${}^{0}L^2(\Gamma\backslash G)$ 
defined in [22, (1.4.1)-(1.4.7)] 
($A_{\Gamma}^0(\Upsilon_{\nu,0};0,0)$ being a space of cusp forms);  
and where $\Upsilon_{\nu,p}$ is the character of 
${\Bbb C}\bigl[\Omega_{+},\Omega_{-}\bigr]$ given by 
[22, (1.3.3)]$\,$ (so that one may write the equation in~(1.1.4) as 
$\Omega_{\pm}f=\Upsilon_{\nu_V,p_V}\bigl(\Omega_{\pm}\bigr) f\,$). 
Since $\Gamma$-automorphicity implies $\tilde\Gamma$-automorphicity, 
it is immediate from the relevant definitions in 
[22, Sections~1.2-1.4] that  
$A_{\Gamma}(\Upsilon_{\nu,0};0,0)$ is a subspace of 
$A_{\tilde\Gamma}(\Upsilon_{\nu,0};0,0)$; so, by (7.9) we obtain: 

$$f_{0,0}^{V'}\in A_{\tilde\Gamma}(\Upsilon_{\nu,0};0,0)\;.\eqno(7.10)$$
If we can furthermore show it to be the case that 
$$f_{0,0}^{V'}\in A_{\tilde\Gamma}^0(\Upsilon_{\nu,0};0,0)
\eqno(7.11)$$
then the sought for conclusion, that $f_{0,0}^{V'}\in{}^{0}L^2(\tilde\Gamma\backslash G)$, 
will follow: for ${}^{0}L^2(\tilde\Gamma\backslash G)$ is 
(see [22, Subsection~1.7]) defined to be the closure of the 
subspace of $L^2(\tilde\Gamma\backslash G)$ generated by the 
set of all $\tilde\Gamma$-automorphic cusp forms,    
and so, since each non-zero element of the set 
$A_{\tilde\Gamma}^0(\Upsilon_{\nu,0};0,0)$ 
is (by definition) a cusp form, it is trivially the case that 
${}^{0}L^2(\tilde\Gamma\backslash G)\supseteq 
A_{\tilde\Gamma}^0(\Upsilon_{\nu,0};0,0)$.

Given the relevant definitions in [22, Subsection~1.4], and given (7.10), 
the verification of (7.11) may be achieved in two steps: the first of these being 
to show that $f_{0,0}^{V'}$ having `polynomial growth' as a $\Gamma$-automorphic 
implies that $f_{0,0}^{V'}$ also has `polynomial growth' as a $\tilde\Gamma$-automorphic function; 
the second step being to demonstrate the like implication in respect of 
the `cuspidality' critereon 
$$\left( F_{0}^{\frak c}f_{0,0}^{V'}\right) (g) =0\qquad\qquad   
\hbox{$\left(\, g\in G,\ {\frak c}\in{\Bbb P}^{1}\bigl({\Bbb Q}(i)\bigr)
={\Bbb Q}(i)\cup\{\infty\}\,\right)$,}\eqno(7.12)$$
where, as indicated in Subsection~1.1,  $F_{0}^{\frak c} f$ is 
the $0$-th order term in the Fourier expansion of $f\,$ (as a 
$\Gamma$-automorphic function)  at the cusp~${\frak c}$. 

We address first the question of `polynomial growth' (the reader 
may refer to [22, Subsection~1.4] for the meaning of this terminology). 
Since the parabolic stabiliser subgroups $\Gamma_{\frak c}'$ and 
$\tilde\Gamma_{\frak c}'$ may differ, our insistence that all scaling 
matrices satisfy the condition (1.1.1) necessitates that we indicate 
when the scaling matrix should be one appropriate for the Fourier expansion 
of $\tilde\Gamma$-automorphic functions: we do this by marking 
the relevant `$g_{\frak c}$' with a tilde. Similarly we write 
$\tilde F_{0}^{\frak c} f$ for the $0$-th order term in the Fourier 
expansion of $f$ as a $\tilde\Gamma$-automorphic function. 
 
By (1.1.1), one has 
$g_{\frak c}^{-1}\Gamma_{\frak c}' g_{\frak c}
=\tilde g_{\frak c}^{-1}\tilde\Gamma_{\frak c}' \tilde g_{\frak c}
=\{ n[\alpha] : \alpha\in{\frak O}\}\leq G$. 
A calculation then enables one to deduce that, for some $z_{\frak c}\in{\Bbb C}$, 
and some $u_{\frak c}\in{\Bbb C}^{*}$ with $u_{\frak c}^{2}\in{\frak O}$ 
and $|u_{\frak c}|^{4}=\bigl[\Gamma_{\frak c}' : \tilde\Gamma_{\frak c}'\bigr]$, 
one has the equation 
$$g_{\frak c}^{-1} \tilde g_{\frak c}
=n\bigl[ z_{\frak c}\bigr] h\bigl[ u_{\frak c}\bigr]\;.\eqno(7.13)$$
Then, through a result similar to that in (7.5), one finds that,  
since the function $g\mapsto f_{0,0}^{V'}\bigl( g_{\frak c} g\bigr)\ $ ($g\in G$) 
has (by virtue of (7.9)) polynomial growth along $A$ 
(in the sense defined in [22, (1.4.1)]), so too does the 
function $g\mapsto f_{0,0}^{V'}\bigl( \tilde g_{\frak c} g\bigr)\ $ ($g\in G$). 
This applies for all cusps ${\frak c}\in{\Bbb Q}(i)\cup\{\infty\}$, 
and so $f_{0,0}^{V'}$ meets the criteria stated in [22] for being a 
$\tilde\Gamma$-automorphic function of polynomial growth; given (7.10), 
we therefore have   
$$f_{0,0}^{V'}\in A_{\tilde\Gamma}^{\rm pol}(\Upsilon_{\nu ,0};0,0)\;.\eqno(7.14)$$ 

The first step in our verification of (7.11) is now complete. 
For the second step, relating to the cuspidality criteron (7.12), we may note 
that, from (1.1.1), (7.13) and the definition of 
$F_{\omega}^{\frak c}f : G\rightarrow{\Bbb C}$ in [22, (1.4.3), (1.4.4)], 
one can work out that
$$\bigl( \tilde F_0^{\frak c} f_{0,0}^{V'}\bigr) (g) 
=\left( F_0^{\frak c} f_{0,0}^{V'}\right) \left( h\bigl[ u_{\frak c}\bigr] g\right)\qquad\qquad  
\hbox{for $\quad g\in G$, $\ {\frak c}\in{\Bbb Q}(i)\cup\{\infty\}$.}$$
Therefore, given that $h\bigl[ u_{\frak c}\bigr]$ is an element of the group $G$, 
it follows from (7.12) that 
$$\bigl( \tilde F_{0}^{\frak c}f_{0,0}^{V'}\bigl) (g) =0\qquad\qquad 
\hbox{$\left(\, g\in G,\ {\frak c}\in{\Bbb P}^{1}\bigl({\Bbb Q}(i)\bigr)
={\Bbb Q}(i)\cup\{\infty\}\,\right)$.}$$
This, together with (7.14), makes the verification of (7.11) complete; 
by (7.11) and the discussion around it, we have 
$f_{0,0}^{V'}\in {}^{0}L^2(\tilde\Gamma\backslash G)$, the premise on 
which our deduction of (7.8) depended. 

Now we may apply (7.7) and (7.8), so that the equation~(7.6) is simplified to:
$$\sum_{\scriptstyle W\atop\scriptstyle p_W=0}^{(\tilde\Gamma)} 
{\left\langle f_{0,0}^{V'} , f_{0,0}^W\right\rangle_{\tilde\Gamma\backslash G}
\left\langle f_{0,0}^W , f_{0,0}^{V''}\right\rangle_{\tilde\Gamma\backslash G}\over 
\left\| f_{0,0}^W\right\|_{\tilde\Gamma\backslash G}^2}
=\left\langle f_{0,0}^{V'} , f_{0,0}^{V''}\right\rangle_{\tilde\Gamma\backslash G}\;.\eqno(7.15)$$
Let $W$ be one of the cuspidal irreducible spaces 
of ${}^{0}L^2(\tilde\Gamma\backslash G)$ by which the summation in (7.15) is indexed. 
We already have $p_{V''}=p_{V'}=0$: 
suppose also that  $p_W=0$. Then, by the points noted in the paragraph containing 
(1.1.11) (understood as applying to the case $\Gamma=\tilde\Gamma$), 
the functions $f_{0,0}^{V'}$, $f_{0,0}^{V''}$ and 
$f_{0,0}^W$ are elements of the space $C^{\infty}(G/K)$ 
(defined above (1.1.11)), and one has 
$\nu_X\in i[0,\infty)\cup(0,2/9)$ and   
$\bigl( 1-\nu_X^2\bigr) f_{0,0}^X =-{\bf \Delta} f_{0,0}^{X}$ 
for $X=W,V',V''$, where the operator $-{\bf \Delta} =-4\bigl(\Omega_{+}+\Omega_{-}\bigr)$  
is symmetric on a subspace of $L^2(\tilde\Gamma\backslash G)$ containing 
${\Bbb C}f_{0,0}^W\oplus {\Bbb C}f_{0,0}^{V'}\oplus {\Bbb C}f_{0,0}^{V''}$. 
Since two eigenspaces corresponding to distinct eigenvalues of 
the same symmetric operator are necessarily orthogonal to one another, 
it follows that the term on the left-hand side 
of (7.15) indexed by $W$ is non-zero only if   
$\nu_{V'}^2=\nu_W^2=\nu_{V''}^2$; in which case, given that   
$\nu_{V'}=\nu_{V''}=\nu$ and $\nu,\nu_W\in i[0,\infty)\cup (0,\infty)$, 
one would have $\nu_W=\nu$. Therefore the sums on the left-hand sides of 
the equations~(7.15) and~(7.2) are equal. 
Since it is,  by 
(7.4) (and its corollary (7.1)),
also the case that the terms on the right-hand sides of 
the equations~(7.15) and~(7.2) are equal, 
the proof of (7.2) (and hence of the lemma) is complete
\ $\blacksquare$

\bigskip

\goodbreak\proclaim Lemma~7.2.
Let $\nu\in i[0,\infty)\cup (0,1)$; let $N\geq 1$; and let $a_{\omega}\in{\Bbb C}$ 
for all $\omega\in{\frak O}$ satisfying $0<|\omega|^2\leq N$. For $0\neq q\in{\frak O}$ 
and $\Gamma =\Gamma_0(q)\leq SL(2,{\frak O})$, put 
$$S(\Gamma)=S_N(\Gamma , \nu)
=\sum_{\scriptstyle V\atop\scriptstyle (\nu_V,p_V)=(\nu,0)}^{(\Gamma)} 
\Biggl|\sum_{0<|n|^2\leq N} 
a_n c_V^{\infty}(n;\nu,0)\Biggr|^2\eqno(7.16)$$
(where the relevant scaling matrix $g_{\infty}$ is as in (1.3.3)).  
Then, for $q\in{\frak O}-\{ 0\}$, $r\in q{\frak O}-\{ 0\}$, 
$\Gamma =\Gamma_0(q)$ and $\tilde\Gamma =\Gamma_0(r)$, one has: 
$$S(\Gamma)\leq\bigl[\Gamma : \tilde\Gamma\,\bigr] S\bigl(\tilde\Gamma\bigr)\;.\eqno(7.17)$$

\medskip

\goodbreak 
\noindent{\bf Proof.}\quad  
Let $\nu$, $N$ and the coefficients $a_{\omega}$ 
satisfy the stated hypotheses. Suppose moreover that $0\neq q,r\in{\frak O}$; 
that $q$ is a factor of $r$; and that $\Gamma$ and $\tilde\Gamma$ are (respectively) 
the subgroups $\Gamma_0(q)$ and $\Gamma_0(r)$ of $SL(2,{\frak O})$. 
By the final point noted in the paragraph of (1.1.11), all the summations 
on the right-hand side of (7.16) are finite. Hence both $S(\Gamma)$ and $S(\tilde\Gamma)$ 
are well-defined sums, and we have $|S(\Gamma)|,|S(\tilde\Gamma)|<\infty$. 
By (1.1.21), the modified Fourier coefficients $c_V^{\infty}(n ;\nu,0)$ 
occurring in the sum (7.16) satisfy 
$$c_V^{\infty}(\omega ;\nu,0)
=(\pi |\omega|)^{\nu} c_V^{\infty}(\omega)\;,\eqno(7.18)$$
where (for $0\neq\omega\in{\frak O}$) the factor $c_V^{\infty}(\omega)$ 
is the same coefficient as appears in the Fourier expansion (1.1.8) of the 
chosen generator $f_{0,0}^V$ for 
the one-dimensional subspace $V_{K,0,0}\,$ (occurring 
in the orthogonal decomposition (1.1.5) of $V\,$). Recall that $c_V(\omega)$, in (1.1.8), 
is independent of the parameters $q$ and $\ell$ there. With 
$g_{\infty}$ given by (1.3.3), it follows by [4, Lemma~5.1]  
that the case $\nu_V=\nu$, $p_V=q=\ell=0$, ${\frak c}=\infty$ of the equation~(1.1.8)  
may, for $z\in{\Bbb C}$, $r>0$, $k\in K$, and $g=n[z]a[r]k\in G$, be cast in more 
classical terminology as: 
$$f_{0,0}^V(g)
=\sum_{0\neq\omega\in{\frak O}} c_V^{\infty}(\omega) 
\,{2\pi^{\nu +1}|\omega|^{\nu}\over\Gamma(\nu +1)}\,{\rm e}\bigl( {\rm Re}(\omega z)\bigr) 
r K_{\nu}(2\pi |\omega| r)\;,\eqno(7.19)$$
where the Bessel function $K_{\nu} : {\Bbb C}-(-\infty ,0]\rightarrow{\Bbb C}$ 
(differing from that `$K_{\nu}(z)$' defined in [24] by the omission of a factor 
$\cos(\pi\nu)$\/) is, by virtue of the relevant asymptotic expansion 
given in [24, Section~17.7], non-zero for all  positive values of 
the argument that are sufficiently large (in terms of $\nu$). Hence, when $r>0$, one has 
$${2\pi^{\nu +1}|\omega|^{\nu}\over\Gamma(\nu +1)}\,r K_{\nu}(2\pi |\omega| r)
c_V^{\infty}(\omega) 
=\bigl( F_{\omega}^{\infty} f_{0,0}^V\bigr) (a[r])\qquad\qquad  
\hbox{($0\neq\omega\in{\frak O}$),}\eqno(7.20)$$
where $F_{\omega}^{\frak c}f$ is the same term seen in the 
Fourier expansion displayed just below the equation~(1.1.1), and is uniquely determined 
by virtue of the classical integral representation of Fourier coefficients.  
When $r>0$ is sufficiently large, the equations in (7.20) determine 
the coefficients $c_V^{\infty}(\omega)\,$   
(for all $\omega\in{\frak O}-\{ 0\}$). 

For our proof of (7.17) we shall need to express  
Fourier coefficients $c_V^{\infty}(\omega)$, associated (through (1.1.6) and (1.1.8))  
with the irreducible cuspidal 
subspaces $V\subseteq {}^{0}L^2(\Gamma\backslash G)$ having  
$\bigl(\nu_V,p_V)=(\nu,0)$, in terms of the corresponding Fourier coefficients, $c_W^{\infty}(\omega)$, 
associated with irreducible cuspidal subspaces $W\subseteq{}^{0}L^2(\tilde\Gamma\backslash G)$ 
having the same pair of spectral parameters, $(\nu ,0)$. To this end, we shall first determine an expression for  
the function $f_{0,0}^V\in V\subset{}^{0}L^2(\Gamma\backslash G)$ in terms of the 
corresponding functions, $f_{0,0}^W$, lying in relevant cuspidal irreducible subspaces 
$W$ of ${}^{0}L^2(\tilde\Gamma\backslash G)$: the required 
relations 
between Fourier coefficients will then be seen to follow through an appeal to the 
final remark of the previous paragraph. 

Let ${\cal V}_{\nu,0}(\Gamma)$ be the set of all of those of the irreducible cuspidal 
subspaces $V\subset {}^{0}L^2(\Gamma\backslash G)$ occurring in 
the orthogonal decomposition (1.1.3) that have their spectral parameters 
$\bigl(\nu_{V},p_{V}\bigr)$ equal to $\bigl(\nu,0\bigr)$ (so that ${\cal V}_{\nu,0}(\Gamma)$ 
is the range of the variable of summation, $V$, in the sum on the right-hand side of 
the equation~(7.16)). If the set ${\cal V}_{\nu,0}(\Gamma)$ is empty, 
then since sums with no terms are (by definition) equal to zero, 
it follows from the definition in (7.16) that $S(\Gamma)=0$ and $S(\tilde\Gamma)\geq 0$: 
the result (7.17) of the lemma is, in that case, a trivial 
consequence of the lower bound $[\Gamma :\tilde\Gamma\,]\geq 1$ implied by Lemma~7.1. 
We may therefore assume that  the set ${\cal V}_{\nu,0}(\Gamma)$ contains 
at least one element. 

Suppose that $V\in {\cal V}_{\nu,0}(\Gamma)$; 
and let $f_{0,0}^{V}$ be a generator of the 
subspace $V_{K,0,0}\subset V$, 
normalised so as to satisfy (1.1.9).   
Then, by the case $V'=V''=V$ of Lemma~7.1 (the equation (7.2), in particular), 
one has  
$$\sum_{\scriptstyle W\atop\scriptstyle (\nu_W,p_W)=(\nu,0)}^{\left(\tilde\Gamma\right)} 
\Biggl|\Biggl\langle f_{0,0}^V\,,\,{1\over
\bigl\| f_{0,0}^W\bigr\|_{\tilde\Gamma\backslash G}}
\,f_{0,0}^W\Biggr\rangle_{\tilde\Gamma\backslash G}\Biggr|^2 
=\bigl\| f_{0,0}^V\bigr\|_{\tilde\Gamma\backslash G}^2\;,\eqno(7.21)$$
which is an example of Bessel's inequality holding with equality:  
for on the left-hand side of this equation the variable of summation $W$ 
indexes a set of functions $f_{0,0}^W$ that are 
(by the discussion around (1.1.3)-(1.1.6), as it applies for $\Gamma=\tilde\Gamma$) 
pairwise orthogonal elements of the space $L^2(\tilde\Gamma\backslash G)$. 
The same case of 
Lemma 7.1 shows also that the function $f_{0,0}^V : G\rightarrow{\Bbb C}$ 
is bounded and continuous. 
Similarly, for each $W\in{\cal V}_{\nu,0}(\tilde\Gamma)$ (i.e. each space 
$W$ indexing a summand on the left hand side of the equation (7.21)), 
the corresponding normalised generator $f_{0,0}^W$ of the subspace 
$W_{K,0,0}\subset W$ is a bounded and continuous function on $G$.  
It therefore follows from (7.21) that 
$$f_{0,0}^V=\sum_{\scriptstyle W\atop\scriptstyle (\nu_W,p_W)=(\nu,0)}^{\left(\tilde\Gamma\right)} 
\beta_{V,W} f_{0,0}^W\;,\eqno(7.22)$$
where 
$$\beta_{V,W}
=\Biggl\langle f_{0,0}^V\,,\,{1\over
\bigl\| f_{0,0}^W\bigr\|_{\tilde\Gamma\backslash G}}
\,f_{0,0}^W\Biggr\rangle_{\!\!\tilde\Gamma\backslash G}\,{1\over
\bigl\| f_{0,0}^W\bigr\|_{\tilde\Gamma\backslash G}}
={\bigl\langle f_{0,0}^V\,,\,f_{0,0}^W\bigr\rangle_{\tilde\Gamma\backslash G}\over 
\bigl\| f_{0,0}^W\bigr\|_{\tilde\Gamma\backslash G}^2}\;.\eqno(7.23)$$
Indeed, the equation~(7.21) implies that the 
$L^2(\tilde\Gamma\backslash G)$-norm of the difference between the two sides 
of equation~(7.22) is equal to zero. Since 
that difference is a continuous function on $G$, and has $L^2(\tilde\Gamma\backslash G)$-norm 
equal to zero, it must therefore have range $\{ 0\}$ and domain $G$.  

In the equation~(7.19) one may substitute, in place of $V$, any of the spaces $W$ by which 
the summation in (7.22) is indexed: for (7.19) would not fail to apply 
if we had $\tilde\Gamma=\Gamma$ and $V=W$. Hence, and by (7.22), it may be deduced that, 
for $g=n[z]a[r]k\in G$ with 
$z\in{\Bbb C}$, $r>0$ and $k\in K$, one has 
$$\eqalignno{
f_{0,0}^V
 &=\sum_{\scriptstyle W\atop\scriptstyle (\nu_W,p_W)=(\nu,0)}^{\left(\tilde\Gamma\right)} 
\beta_{V,W} \sum_{0\neq\omega\in{\frak O}} c_W^{\infty}(\omega) 
\,{2\pi^{\nu +1}|\omega|^{\nu}\over\Gamma(\nu +1)}\,{\rm e}\bigl( {\rm Re}(\omega z)\bigr) 
r K_{\nu}(2\pi |\omega| r) = {}\cr 
 &=\sum_{0\neq\omega\in{\frak O}}\tilde c_V^{\infty}(\omega) 
{2\pi^{\nu +1}|\omega|^{\nu}\over\Gamma(\nu +1)}\,{\rm e}\bigl( {\rm Re}(\omega z)\bigr) 
r K_{\nu}(2\pi |\omega| r)\;, &(7.24)}$$
where, for $0\neq\omega\in{\frak O}$, 
$$\tilde c_V^{\infty}(\omega)
=\sum_{\scriptstyle W\atop\scriptstyle (\nu_W,p_W)=(\nu,0)}^{\left(\tilde\Gamma\right)} 
\beta_{V,W} c_W^{\infty}(\omega)\;.\eqno(7.25)$$

Similarly to how (7.20) was deduced from (7.19), one may deduce from (7.24) 
that (7.20) continues to hold for all $r>0$ if, for all $\omega\in{\frak O}-\{ 0\}$,  
one substitutes for the Fourier coefficient $c_V^{\infty}(\omega)$ in (7.20)  
the number $\tilde c_V^{\infty}(\omega)$ just defined: given the 
point noted below (7.20), it must therefore be the case that 
$\tilde c_V^{\infty}(\omega)=c_V^{\infty}(\omega)$  
for all $\omega\in{\frak O}-\{ 0\}$. Hence, given (7.25) and (7.18) 
(which remains valid when 
$V$ is replaced by any one of the subspaces $W$ by which the summation in (7.25) is indexed), 
we are able to deduce that 
$$c_V^{\infty}(\omega;\nu,0) 
=\sum_{\scriptstyle W\atop\scriptstyle (\nu_W,p_W)=(\nu,0)}^{\left(\tilde\Gamma\right)} 
\beta_{V,W}\,c_W^{\infty}(\omega ;\nu,0)\qquad\qquad  
\hbox{($0\neq\omega\in{\frak O}$)}.\eqno(7.26)$$

Since (7.26) has been shown to hold for an arbitrary  
member $V$ of the set of spaces ${\cal V}_{\nu,0}(\Gamma)$ (defined 
earlier in this proof), we may apply (7.26) to expand every one of the modified 
Fourier coefficients 
$c_V^{\infty}(n;\nu,0)$ occurring in the definition (7.16). We consequently 
find that 
$$S(\Gamma)
=\sum_{\scriptstyle V\atop\scriptstyle (\nu_V,p_V)=(\nu,0)}^{(\Gamma)} 
\left|\sigma_V\right|^2
\;,\eqno(7.27)$$
\goodbreak 
\noindent 
where, for $V\in{\cal V}_{\nu,0}(\Gamma)$,  
$$\eqalignno{
\sigma_V
 &=\sum_{0<|n|^2\leq N} 
a_n c_V^{\infty}(n ;\nu,0) = {}\cr
 &=\sum_{0<|n|^2\leq N} 
a_n\sum_{\scriptstyle W\atop\scriptstyle (\nu_W,p_W)=(\nu,0)}^{\left(\tilde\Gamma\right)} 
\beta_{V,W} c_W^{\infty}(n ;\nu,0) = {}\cr
&=\sum_{\scriptstyle W\atop\scriptstyle (\nu_W,p_W)=(\nu,0)}^{\left(\tilde\Gamma\right)} 
\beta_{V,W}\sum_{0<|n|^2\leq N}
a_n c_W^{\infty}(n ;\nu,0)
=\sum_{\scriptstyle W\atop\scriptstyle (\nu_W,p_W)=(\nu,0)}^{\left(\tilde\Gamma\right)} 
\beta_{V,W} \sigma_W^{*}\qquad\quad\hbox{(say).} &(7.28)}$$
Hence 
$$\eqalign{
S(\Gamma) 
 &=\sum_{\scriptstyle V\atop\scriptstyle (\nu_V,p_V)=(\nu,0)}^{(\Gamma)} 
\overline{\sigma_V}\,\sigma_V = {}\cr
 &=\sum_{\scriptstyle V\atop\scriptstyle (\nu_V,p_V)=(\nu,0)}^{(\Gamma)} 
\overline{\sigma_V}\,\sum_{\scriptstyle W\atop\scriptstyle (\nu_W,p_W)=(\nu,0)}^{\left(\tilde\Gamma\right)} 
\beta_{V,W} \sigma_W^{*} 
=\sum_{\scriptstyle W\atop\scriptstyle (\nu_W,p_W)=(\nu,0)}^{\left(\tilde\Gamma\right)} 
\sigma_W^{*} \sum_{\scriptstyle V\atop\scriptstyle (\nu_V,p_V)=(\nu,0)}^{(\Gamma)} 
\overline{\sigma_V}\,\beta_{V,W}\;,}$$
and so it follows by the Cauchy-Schwarz inequality that 
$$|S(\Gamma)|^2\leq 
\Biggl( \sum_{\scriptstyle W\atop\scriptstyle (\nu_W,p_W)=(\nu,0)}^{\left(\tilde\Gamma\right)} 
\left|\sigma_W^{*}\right|^2\Biggr)
\Biggl( \sum_{\scriptstyle W\atop\scriptstyle (\nu_W,p_W)=(\nu,0)}^{\left(\tilde\Gamma\right)} 
\Biggl|\sum_{\scriptstyle V\atop\scriptstyle (\nu_V,p_V)=(\nu,0)}^{(\Gamma)} 
\overline{\sigma_V}\,\beta_{V,W}\Biggr|^2\,\Biggr) 
=S\bigl(\tilde\Gamma\bigr)\,T\bigl(\Gamma\,,\,\tilde\Gamma\bigr)\;,\eqno(7.29)$$
where (given that we take $\sigma_W^{*}$, on the right-hand side of (7.28), 
to equal the sum over $n$ 
on the same line)  $S(\tilde\Gamma)$ is given by 
the case $\Gamma=\tilde\Gamma$ of (7.16), while 
$$\eqalignno{ 
T\bigl(\Gamma\,,\,\tilde\Gamma\bigr)
 &=\sum_{\scriptstyle W\atop\scriptstyle (\nu_W,p_W)=(\nu,0)}^{\left(\tilde\Gamma\right)} 
\Biggl|\sum_{\scriptstyle V\atop\scriptstyle (\nu_V,p_V)=(\nu,0)}^{(\Gamma)} 
\overline{\sigma_V}\,\beta_{V,W}\Biggr|^2 = {}\cr
 &=\sum_{\scriptstyle W_{\ }\atop\scriptstyle (\nu_W,p_W)=(\nu,0)}^{\left(\tilde\Gamma\right)} 
\sum_{\scriptstyle V''\atop\scriptstyle\!\!\!=(\nu_{V''},p_{V''})}^{(\Gamma)} 
\sigma_{V''}\,\overline{\beta_{V'',W}}  
\sum_{\scriptstyle V'\atop\scriptstyle (\nu_{V'},p_{V'})=(\nu,0)}^{(\Gamma)} 
\overline{\sigma_{V'}}\,\beta_{V',W} = {}\cr 
 &=\quad
\sum_{\scriptstyle V''}^{(\Gamma)}
\!\!\!\sum_{\scriptstyle V'_{\hbox{\ }}\atop
{\!\!\!\!\!\!\!\!\!\!\!{\scriptstyle \nu_{V''}=\nu=\nu_{V'}\atop 
\scriptstyle p_{V''}=0=p_{V'}}}}^{(\Gamma)} 
\sigma_{V''}\,\overline{\sigma_{V'}} 
\sum_{\scriptstyle W\atop\scriptstyle (\nu_W,p_W)=(\nu,0)}^{\left(\tilde\Gamma\right)}
\overline{\beta_{V'',W}}\,\beta_{V',W}\;. &(7.30)}$$
By the equation (7.23) (for $V=V'$, and for $V=V''$), 
the normalisation (1.1.9) (for $V=W$, and for $V=V'\,$) 
and the result (7.2) of Lemma~7.1, one finds that  
the inner sum on the right-hand side of (7.30) is 
$$\eqalign{
\sum_{\scriptstyle W\atop\scriptstyle (\nu_W,p_W)=(\nu,0)}^{(\tilde\Gamma)} 
{\overline{\langle f_{0,0}^{V''},f_{0,0}^W\rangle_{\tilde\Gamma\backslash G}} 
\ \langle f_{0,0}^{V'},f_{0,0}^{W}\rangle_{\tilde\Gamma\backslash G}\over 
\bigl\| f_{0,0}^{W}\bigr\|_{\tilde\Gamma\backslash G}^4}
 &=\sum_{\scriptstyle W\atop\scriptstyle (\nu_W,p_W)=(\nu,0)}^{(\tilde\Gamma)} 
{\overline{\langle f_{0,0}^{V''},f_{0,0}^W\rangle_{\tilde\Gamma\backslash G}} 
\,\langle f_{0,0}^{V'},f_{0,0}^{W}\rangle_{\tilde\Gamma\backslash G}\over 
\bigl\| f_{0,0}^{V'}\bigr\|_{\Gamma\backslash G}^2 
\bigl\| f_{0,0}^{W}\bigr\|_{\tilde\Gamma\backslash G}^2} = {}\cr 
 &={1\over \bigl\| f_{0,0}^{V'}\bigr\|_{\Gamma\backslash G}^2}
\sum_{\scriptstyle W\atop\scriptstyle (\nu_W,p_W)=(\nu,0)}^{(\tilde\Gamma)} 
{\langle f_{0,0}^{V'},f_{0,0}^{W}\rangle_{\tilde\Gamma\backslash G}   
\,\langle f_{0,0}^W,f_{0,0}^{V''}\rangle_{\tilde\Gamma\backslash G}\over 
\bigl\| f_{0,0}^{W}\bigr\|_{\tilde\Gamma\backslash G}^2} = {}\cr 
 &=\cases{
\bigl\| f_{0,0}^{V'}\bigr\|_{\Gamma\backslash G}^{-2}  
\bigl\| f_{0,0}^{V'}\bigr\|_{\tilde\Gamma\backslash G}^2
 &if $V''=V'$,\cr 0 &otherwise,}
}$$
where, by the result (7.1) of Lemma~7.1, one has 
$$\bigl\| f_{0,0}^{V'}\bigr\|_{\Gamma\backslash G}^{-2}  
\bigl\| f_{0,0}^{V'}\bigr\|_{\tilde\Gamma\backslash G}^2 =
\bigl[\Gamma : \tilde\Gamma\bigr]\;.$$    
By this result, and (7.30) and (7.27), one obtains:  
$$T\bigl(\Gamma\,,\,\tilde\Gamma\bigr)
=\sum_{\scriptstyle V'\atop\scriptstyle (\nu_{V'},p_{V'})=(\nu,0)}^{(\Gamma)} 
\sigma_{V'}\,\overline{\sigma_{V'}}\,\bigl[\Gamma : \tilde\Gamma\bigr]\ \,
=\ \,\bigl[\Gamma : \tilde\Gamma\bigr]
\!\!\!\sum_{\scriptstyle V'\atop\scriptstyle (\nu_{V'},p_{V'})=(\nu,0)}^{(\Gamma)}
\left| \sigma_{V'}\right|^2\  
=\,\bigl[\Gamma : \tilde\Gamma\bigr] S(\Gamma)\;.$$
Hence and by (7.29), it follows that 
$$|S(\Gamma)|^2
\leq S\bigl(\tilde\Gamma\bigr)\bigl[\Gamma : \tilde\Gamma\bigr] S(\Gamma)\;,$$
where $\bigl[\Gamma : \tilde\Gamma\bigr]\geq 1$,   
and where (by the definition (7.16)) $S(\Gamma),S(\tilde\Gamma)\geq 0$. 
One  therefore must have (7.17)   
\ $\blacksquare$

\bigskip

\goodbreak\proclaim Lemma~7.3.
Let $t\in{\Bbb R}$ and $Q,X,N\geq 1$. Then, for $Q_1\geq{5\over 2}\,Q$, one has 
$$S_t(Q,X,N)\leq{\textstyle {36\over 5}\!\left( 1-{\log 2\over\log 5}\right)^{\!\!-2}} 
\left(\log Q_1\right) 
\bigl( S_t\left( Q_1,X,N\right) +S_t\left( 2Q_1,X,N\right)\bigr) .\eqno(7.31)$$

\medskip

\goodbreak 
\noindent{\bf Proof.}\quad 
By the definition (1.3.2) of $S_t(Q,X,N)$, it will suffice to prove the case 
$t=0$ of this lemma (application of that case with $a_n|n|^{2it}$ substituted 
for $a_n$, for all $n\in{\frak O}-\{ 0\}$, will imply the cases where $t\neq 0$). 
By (1.3.2) and (1.1.11), one has  
$$S_0(Q,X,N)
=\sum_{Q/2<|q|^2\leq Q}\ \sum_{0<\nu<1}  
X^{\nu} S_N\bigl(\Gamma_0(q)\,,\,\nu\bigr)\;,\eqno(7.32)$$ 
with $S_N(\Gamma , \nu)$ as given by the equation~(7.16) in Lemma~7.2 
(note that, by the remarks below (1.1.11), the above summation summation over $\nu$ 
is effectively finite).  

Let $Q_1\geq {5\over 2}Q$; and let $\varpi$ be a Gaussian prime satisfying 
$${Q_1\over Q}<|\varpi|^2\leq{2Q_1\over Q}\;.\eqno(7.33)$$ 
Then, by the case $r=\varpi q$ of the result (7.17) of Lemma~7.2, one has 
$$S_N\bigl(\Gamma_0(q)\,,\,\nu\bigr)
\leq\left[\Gamma_0(q) : \Gamma_0(\varpi q)\right] 
S_N\bigl(\Gamma_0(\varpi q)\,,\,\nu\bigr)\qquad\qquad 
\hbox{($0\neq q\in{\frak O}$ and $0<\nu<1$),}$$
and so it follows by (7.32), the case $\Gamma =\Gamma_0(\varpi q)$ of the 
definition (7.16) in Lemma~7.2, and (1.1.11), that 
$$\eqalign{
S_0(Q,X,N) &\leq
\sum_{Q/2<|q|^2\leq Q}\left[\Gamma_0(q) : \Gamma_0(\varpi q)\right]
\sum_{0<\nu<1} X^{\nu} S_N\bigl(\Gamma_0(\varpi q)\,,\,\nu\bigr) = {}\cr 
 &=\sum_{Q/2<|q|^2\leq Q}\left[\Gamma_0(q) : \Gamma_0(\varpi q)\right]
\sum_{0<\nu<1} X^{\nu}
\sum_{\scriptstyle V\atop\scriptstyle 
(\nu_V,p_V)=(\nu,0)}^{\left(\Gamma_0(\varpi q)\right)} 
\Biggl|\sum_{0<|n|^2\leq N} 
a_n c_V^{\infty}(n;\nu,0)\Biggr|^2 ={}\cr 
 &=\sum_{Q/2<|q|^2\leq Q}\left[\Gamma_0(q) : \Gamma_0(\varpi q)\right]
\sum_{\scriptstyle V\atop\scriptstyle \nu_V>0}^{\left(\Gamma_0(\varpi q)\right)} 
X^{\nu_V} 
\Biggl|\sum_{0<|n|^2\leq N} 
a_n c_V^{\infty}(n;\nu,0)\Biggr|^2\;.}$$
Since this holds for all Gaussian primes $\varpi$ satisfying (7.33), we may 
therefore sum the above bound over all such $\varpi$ so as to obtain:
$$P\bigl( 2Q_1 /Q\bigr) S_0(Q,X,N) 
\leq\sum_{Q_1/2<|q_1|^2\leq 2Q_1} M\bigl( Q , Q_1 ; q_1\bigr) 
\sum_{\scriptstyle V\atop\scriptstyle \nu_V>0}^{\left(\Gamma_0(q_1)\right)} 
X^{\nu_V} 
\Biggl|\sum_{0<|n|^2\leq N} 
a_n c_V^{\infty}(n;\nu,0)\Biggr|^2\;,\eqno(7.34)$$
where 
$$P(x)=\left|\left\{\varpi\in{\Bbb Z}[i] : 
\varpi\ {\rm is\ prime}\ {\rm and}
\ x\geq |\varpi|^2>{x\over 2}\right\}\right|\eqno(7.35)$$
and 
$$M\bigl( Q , Q_1 ; q_1\bigr) 
=\sum_{\scriptstyle{Q_1\over Q}<|\varpi|^2\leq{2Q_1\over Q}\atop\scriptstyle
\varpi\ {\rm is\ prime}}
\sum_{\scriptstyle{Q\over 2}<|q|^2\leq Q\atop\scriptstyle\varpi q=q_1} 
\bigl[\Gamma_0(q) : \Gamma_0\bigl( q_1\bigr)\bigr]\eqno(7.36)$$
(the term `prime' here signifying a `Gaussian prime' of the ring ${\Bbb Z}[i]$, 
as distinct from a `rational prime' of the ring ${\Bbb Z}$).
The result (7.31) will be seen to follow from (7.34), once we have 
obtained a suitable lower bound for $P\bigl( 2Q_1/Q\bigr)$, and a suitable 
upper bound for $M\bigl( Q , Q_1 ; q_1\bigr)$. 

To bound $P\bigl( 2Q_1/Q\bigr)$ from below, we observe firstly that, 
by (7.35) and [8, Theorems~251~and~252],   
$$P(x)\geq 8\left(\pi(x;4,1)-\pi\bigl({x\over 2};4,1\bigr)\right) 
\geq{8\over\log x}\left(\theta(x;4,1)-\theta\bigl({x\over 2};4,1\bigr)\right)\qquad\qquad  
\hbox{($x\geq 5$),}\eqno(7.37)$$
where $\pi(x;k,h)$ denotes the number of rational primes $p\equiv h\bmod k{\Bbb Z}$ 
satisfying $x\geq p>0$, while 
$$\theta(x;k,h)=\sum_{\scriptstyle {\rm rational\ primes}\ p\atop 
{\scriptstyle p\equiv h\bmod k{\Bbb Z}\atop\scriptstyle x\geq p>0}}\log p\;.$$
By [20, Theorem~1, Table 1 and Theorem~5.2.1],  
$$\left|\theta(y;4,1)-{y\over 2}\right|\leq\delta y+\sqrt{y}\qquad\quad
\hbox{for $\quad y\geq 14$,}$$
with the constant $\delta =(0.002238)/2=0.001119<1/36$. 
Using this result one finds that if $x\geq 14^2$ then 
$$\left|\left(\theta(x;4,1)-\theta\bigl({x\over 2};4,1\bigr)\right)-{x\over 4}\right| 
<{x\over 6}\;,$$
so that one has, by (7.37), the lower bound   
$8\left(\pi(x;4,1)-\pi\bigl(x/2;4,1\bigr)\right)\geq {2\over 3}\,x(\log x)^{-1}$. 
The latter conclusion can be shown (by means of some elementary numerical computation) 
also to hold good when $14^2>x\geq 14$. 

The lower bounds in (7.37) are no help at all when $13>x\geq 10$, 
though we do have $x\geq |3|^2>x/2$   
for such $x$. By taking account of $3$ and its associates, we are able to 
deduce from the results of the previous paragraph that 
$$P(x)\geq{2x\over 3\log x}\qquad\quad\hbox{for $\quad x\geq 5$.}$$
Since $Q_1\geq{5\over 2}Q$ and $5>e$, this lower bound on $P(x)$ implies, in particular, 
that we have: 
$$P\bigl( 2Q_1/Q\bigr) 
\geq{4Q_1\over 3Q\log\bigl( 2Q_1 /Q\bigr)}
>{2e\over 3}>1\;.\eqno(7.38)$$
 
To obtain a suitable upper bound on $M\bigl( Q , Q_1 ; q_1\bigr)$, in (7.34) and (7.36), 
we note firstly that, by [22, Equation~(1.1.5)] 
(the $SL(2,{\Bbb Z})$-analogue of which is proved in [12, Section~2.4]), 
one has 
$$\bigl[ SL(2,{\frak O}) : \Gamma_0\bigl( r\bigr)\bigr] 
=\left| r\right|^2\prod_{\scriptstyle{\rm prime\ ideals}\ \varpi_1{\frak O}\subset{\frak O}\atop\scriptstyle 
\varpi_1{\frak O}\ni r} 
\left( 1+{1\over\left|\varpi_1\right|^2}\right)\qquad\quad  
\hbox{for $\quad 0\neq r\in{\frak O}$}\eqno(7.39)$$
(i.e. with, in the last product, only one factor, not four, per prime ideal of 
the ring ${\frak O}={\Bbb Z}[i]$). Hence, in the sum on the right-hand side of 
(7.36) (where $q_1 /q$ equals the Gaussian prime $\varpi$), one has:
$$\eqalign{ 
\bigl[\Gamma_0\bigl( q\bigr) : \Gamma_0\bigl( q_1\bigr)\bigr]
 &={\bigl[ SL(2,{\frak O}) : \Gamma_0\bigl( q_1\bigr)\bigr]\over 
\bigl[ SL(2,{\frak O}) : \Gamma_0\bigl( q\bigr)\bigr]} = {}\cr
 &=\left| {q_1\over q}\right|^2 
\prod_{\scriptstyle{\rm prime\ ideals}\ \varpi_1{\frak O}\subset{\frak O}\atop\scriptstyle 
\varpi_1{\frak O}\ni q_1\ {\rm and}\ \varpi_1{\frak O}\not\ni q} 
\left( 1+{1\over\left|\varpi_1\right|^2}\right)  
=\cases{|\varpi|^2 &if $\varpi\mid q$, \cr 
|\varpi|^2 +1 &otherwise.}}$$
Given this evaluation of the index $\bigl[ SL(2,{\frak O}) : \Gamma_0\bigl( q_1\bigr)\bigr]$,  
it follows trivially from the definition (7.36) that 
$$M\bigl( Q , Q_1 ;q_1\bigr) 
\leq \left( {2Q_1\over Q} +1\right) F\!\left( q_1\,,\,{Q_1\over Q}\right) ,\eqno(7.40)$$
where 
$$F\bigl( q_1 , z\bigr)
=\left|\left\{\varpi\in{\Bbb Z}[i] : 
\varpi\ {\rm is\ prime,}\ 
\ \varpi\mid q_1\ {\rm and}\ |\varpi|^2>z\right\}\right| .$$
Here we have $F\bigl( q_1 , z\bigr)=4k$ (say), where $k$ is, in all cases, a 
non-negative integer, and where, when $k\neq 0$,  there is a set of $k$ pairwise non-associated 
Gaussian primes $\varpi_1,\varpi_2,\ldots ,\varpi_k$ such that, 
for $j=1,2,\ldots ,k$,  
one has both $|\varpi_j|^2>z$ and $\varpi_j\mid q_1$; 
so if one has also $q_1\neq 0$ and $z\geq 0$, then it 
must follow that $|\varpi_1\varpi_2\cdots \varpi_k|^2>z^k$ 
and $q_1\in \varpi_1\varpi_2\cdots \varpi_k{\frak O}-\{0\}$, and hence that    
$|q_1|^2\geq |\varpi_1\varpi_2\cdots \varpi_k|^2>z^k$.
Since the last two inequalities imply that $k\log z<\log\bigl( |q_1|^2\bigr)$, 
we may therefore deduce that 
$$F\bigl( q_1 , z\bigr)\leq {4\log\bigl( |q_1|^2\bigr)\over\log z}\qquad\quad   
\hbox{for $\quad 0\neq q_1\in{\frak O}\ $ and $\ z>1$.}$$
Hence and by (7.40), we have (given that $Q_1 /Q\geq 5/2>1$):   
$$M\bigl( Q , Q_1 ;q_1\bigr) 
\leq {4\log\bigl( 2Q_1\bigr)\over\log\bigl( Q_1 /Q\bigr)}
\left( {2Q_1\over Q} +1\right)\qquad\quad  
\hbox{for $\quad q_1\in{\frak O}\ \,$ with   
$\ \,{\displaystyle{Q_1\over Q}}<\bigl| q_1\bigr|^2\leq 2Q_1$.}\eqno(7.41)$$

By (7.34), (7.38), (7.41) and the definition (1.3.2) of the sum $S_t(Q,X,N)$, 
we have now  
$$\eqalign{
{4Q_1\over 3Q\log\bigl( 2Q_1 /Q\bigr)}\,S_0(Q,X,N) 
 &\leq {4\log\bigl( 2Q_1\bigr)\over\log\bigl( Q_1 /Q\bigr)}
\left( {2Q_1\over Q} +1\right)\!\!\sum_{{\textstyle{Q_1\over 2}}<|q_1|^2\leq 2Q_1} 
\sum_{\scriptstyle V\atop\scriptstyle \nu_V>0}^{\left(\Gamma_0(q_1)\right)}\!X^{\nu_V} 
\Biggl|\sum_{0<|n|^2\leq N} 
a_n c_V^{\infty}(n;\nu,0)\Biggr|^2 = {}\cr 
 &={4\log\bigl( 2Q_1\bigr)\over\log\bigl( Q_1 /Q\bigr)}
\left( {2Q_1\over Q} +1\right) 
\left( S_0\bigl( Q_1 , X , N\bigr) + S_0\bigl( 2Q_1 , X , N\bigr)\right) ,}$$
which gives: 
$$S_0(Q,X,N)\leq {3\log\bigl( 2Q_1 /Q\bigr)\over\log\bigl( Q_1 /Q\bigr)}\, 
\log\bigl( 2Q_1\bigr)\left( 2+{Q\over Q_1}\right) 
\left( S_0\bigl( Q_1 , X , N\bigr) + S_0\bigl( 2Q_1 , X , N\bigr)\right)\;.$$
The case $t=0$ of the result (7.31) follows from this last bound, 
as an elementary consequence of the 
hypotheses that $Q_1\geq{5\over 2}\,Q$ and $Q\geq 1$\ $\blacksquare$ 

\bigskip 
\bigskip

\goodbreak\centerline{\bf \S 8. The Proof of Theorem~9}

\bigskip

We begin this section with a sequence of six lemmas. 
Lemmas~8.1-8.3 are concerned with   
elementary points of analysis. Lemmas~8.4-8.6, 
concerning the sum $S_t(Q,X,N)$ defined in (1.3.2),  
are deduced from Lemma~6.3 by means of Theorem~6 and Lemma~7.3. 
Lemma~8.6 enables the `proof by induction' of Theorem~9, which follows it.  
In the statements and proofs of 
Lemmas~8.2 and 8.3, and in the proof 
of Lemma~8.4, it is 
to be supposed that $x$, $y$ and $z$ are a system of dependent variables 
such that $z\in{\Bbb C}$, $\,x={\rm Re}(z)\in{\Bbb R}$ and 
$y={\rm Im}(z)\in{\Bbb R}$; the same is to be understood when subsripts are used 
(i.e. one will have $x_{\ell}={\rm Re}(z_{\ell})$ and $y_{\ell}={\rm Im}(z_{\ell})$ 
for any given $\ell$). In the statement of Lemma~8.2, and 
in the proofs of Lemmas~8.2-8.4, we furthermore take ${\cal A}^{\Bbb R}_m I$ (when  
$m\in{\Bbb N}$ and $I\subseteq[0,\infty)$)  to 
signify the subset of ${\Bbb R}^{2m}$ given by  
$${\cal A}^{\Bbb R}_m I=\left\{ \left( x_1,y_1,\ldots ,x_m,y_m\right)\in{\Bbb R}^{2m} : 
\left| x_d +i y_d\right|\in I\ {\rm for}\ d=1,\ldots ,m\right\}\;.\eqno(8.1)$$
One example of this notation is ${\cal A}^{\Bbb R}_{1}[0,1)$, which denotes the open disc, in 
${\Bbb R}^2$, with radius $1$ and centre $(0,0)$. Another example 
is ${\cal A}^{\Bbb R}_{1}(0,\infty)$, which denotes the set ${\Bbb R}^2 -\{ (0,0)\}$. 

\bigskip 

\goodbreak\proclaim Lemma~8.1. Let $m\in{\Bbb N}$ and $\,N\in{\Bbb N}\cup\{ 0\}$. 
Let $U$ be a non-empty open subset of 
$\,{\Bbb R}^m$; let $V$ be a non-empty open subset of $\,{\Bbb R}$;  
and let $f : U\rightarrow V$ and $g : V\rightarrow{\Bbb R}$. 
Let the function $f$ be such that, for 
all $\,n\in{\Bbb N}\cup\{ 0\}$, each one of its $m^n$ partial derivatives of order $n$  
is a continuous real-valued function on $U$.
Suppose, moreover, that $g$ is infinitely differentiable on $V$.   
Then the function $g\circ f : U\rightarrow{\Bbb R}$  
is such that every one of its $m^N$ partial derivatives 
of order $N$ is a continuous real-valued function on $U$. 

\goodbreak 
\noindent{\bf Proof.}\quad 
For ${\bf u}\in U\subset{\Bbb R}^m$, one has $f({\bf u})\in V$, so that 
$g(f({\bf u}))$ is defined. Moreover, the set $U$ is the domain of $f$; and 
the function $g$ is real-valued. Therefore 
the function ${\bf u}\mapsto g(f({\bf u}))$ 
is a real-valued function with domain $U$. In other words, we 
have $g\circ f : U\rightarrow{\Bbb R}$. 

Since $g$ is infinitely differentiable on $V$ it is, in particular, 
continuous on $V$. Moreover, since the (unique) partial derivative of 
$f$ of order $0$ is $f$ itself, the hypotheses of the lemma 
imply the continuity of $f$ on $U$. Since $f$ and $g$ are continuous, so
too is their composition $g\circ f$; so the case $N=0$ if the lemma follows. 

Suppose now that $N'\in{\Bbb N}$, and that 
the lemma is true in all cases where $N<N'$. 
By the chain-rule, 
$${\partial\over\partial u_j}\left( g\circ f\right)({\bf u})=
\left( g'\circ f\right)({\bf u})\,{\partial\over\partial u_j}\,f({\bf u})
\qquad\qquad\hbox{(${\bf u}\in U$, $\,j=1,\ldots ,m$)}\eqno(8.2)$$
where, by hypothesis, every partial derivative of the 
function $(\partial /\partial u_j)f$ is a continuous real-valued function on $U$ 
(all partial derivatives of this function being also partial derivatives of $f$). 
Furthermore, since $g$ is an infinitely differentiable real-valued function 
on $V$, so too is its derivative, $g'$: consequently 
it follows by 
the cases $N=0,\ldots ,N'-1$ of the lemma that $g'\circ f : U\rightarrow{\Bbb R}$, 
and that every partial derivative of this function of order not greater than 
$N'-1$ is a continuous real-valued function on $U$. Therefore, either by  (8.2) alone
(if $N'=1$), or  by (8.2) and the product rule of differential calculus 
(when $N'>1$), it follows that any partial derivative 
of the function $({\partial /\partial u_j})(g\circ f) : U\rightarrow{\Bbb R}$ 
of order $N'-1$ is a continuous real-valued function on $U$ (products and 
sums of continuous functions being continuous also). Since this conclusion 
holds for $j=1,\ldots ,m$, it has therefore been established that, 
when the cases $N=0,\ldots ,N'-1$ of the lemma are true, so too is the 
case $N=N'$. This holds for all $N'\in{\Bbb N}$, so that 
(with the case $N=0$ of the lemma having been proved in the preceding paragraph)  
it follows by induction 
that the lemma is true in all cases\ $\blacksquare$

\bigskip 

\goodbreak\proclaim Lemma~8.2. Let $m\in{\Bbb N}$; 
let ${\cal A}^{\Bbb R}_m (0,\infty)\subset{\Bbb R}^{2m}$ be given by (8.1); 
let ${\bf c}\in{\Bbb R}^{m+1}$ 
and let $\Psi : {\Bbb R}\rightarrow{\Bbb C}$ be infinitely differentiable.  
For ${\bf z}\in\left( {\Bbb C}^{*}\right)^{m}$, let 
$$F\bigl( x_1,y_1,\ldots ,x_m,y_m\bigr)
=f\bigl( x_1+iy_1,\ldots ,x_m+iy_m\bigr) 
=f({\bf z})=\Psi\Biggl( c_{m+1}+\sum_{\ell =1}^m c_{\ell}
\log\bigl(\,\left| z_{\ell}\right|^2\,\bigr)\Biggr) .\eqno(8.3)$$
Then the function $F : {\cal A}^{\Bbb R}_m (0,\infty)\rightarrow{\Bbb C}$ 
so defined is such that, for all $n\in{\Bbb N}$, every one of its   
$(2m)^{n}$ partial derivatives of order $n$ is a continuous 
complex-valued function on $U={\cal A}^{\Bbb R}_m (0,\infty)$; and  
if ${\bf j},{\bf k}\in({\Bbb N}\cup\{ 0\})^m$ and 
${\bf r},{\bf R},{\bf R}-{\bf r}\in(0,\infty)^m$ then 
$${\partial^{j_1+k_1+\ldots +j_m+k_m}\over 
\partial x_1^{j_1}\partial y_1^{k_1}\cdots\partial x_m^{j_m}\partial y_m^{k_m}} 
\, f({\bf z}) =O_{\Psi,{\bf c},{\bf r},{\bf R},{\bf j},{\bf k}}(1)\quad\   
\hbox{for all $\ {\bf z}\in{\Bbb C}^m\ $
with 
$\ \left(\left| z_{1}\right| ,\ldots ,\left| z_{m}\right|\right) 
\in{m\atop{\displaystyle\times\atop\scriptstyle\ell =1}} 
\left[ r_{\ell} , R_{\ell}\right]$.}\eqno(8.4)$$
Moreover, if $m=1$, $c_1\neq 0$ and $\,{\rm Supp}(\Psi)\subseteq[a,b]\subset(-\infty,\infty)$,  
then the Schwartz space ${\cal S}({\Bbb C})$ contains a unique function $f$ satisfying 
(8.3) for all ${\bf z}\in {\Bbb C}^{*}$; 
and this function $f$ is such that, when $z\in{\Bbb C}$,  one has:  
$${\partial^{j+k}\over\partial x^j\partial y^k}\,f(z)=O_{\Psi,{\bf c},j,k}(1)\qquad\qquad\quad   
\hbox{($j,k\in{\Bbb N}\cup\{ 0\}$)}\eqno(8.5)$$
and 
$$f(z)\neq 0\quad\ \,\hbox{only if}\quad\, 
\exp\!\left( \min\left\{ {a\over c_1}\,,\,{b\over c_1}\right\} - {c_2\over c_1}\right) 
< |z|^2 < 
\exp\!\left( \max\left\{ {a\over c_1}\,,\,{b\over c_1}\right\} - {c_2\over c_1}\right)
 .\eqno(8.6)$$

\goodbreak 
\noindent{\bf Proof.}\quad 
Since $x^2+y^2$ is a polynomial, and since the function $\log(v)$ is 
real-valued and infinitely differentiable for $v>0$, it is easily verified 
that the hypotheses of Lemma~8.1 are satisfied when one takes there 
$U={\cal A}^{\Bbb R}_1 (0,\infty)$, $V=(0,\infty)$, $f(x,y)=x^2+y^2\ $   
($(x,y)\in U$) and $g(v)=\log(v)\ $ ($v\in V$). Lemma~8.1 therefore shows 
that all partial derivatives of the function $(x,y)\mapsto\log\bigl( |x+iy|^2\bigr)$ 
are continuous real-valued functions on ${\cal A}^{\Bbb R}_1 (0,\infty)$. 
This trivially implies that all partial derivatives of the 
$m$ distinct functions, 
$(x_1,y_1,\ldots ,x_m,y_m)\mapsto\log\bigl( |x_{\ell}+iy_{\ell}|^2\bigr)\ $   
($\ell =1,\ldots ,m$), are continuous real-valued functions on 
${\cal A}^{\Bbb R}_m (0,\infty)$. 
Since partial derivatives are linear operators, and since sums and 
products of continuous real-valued functions are themselves 
continuous real-valued functions, it follows that all partial derivatives 
of the function 
$\bigl( x_1,y_1,\ldots ,x_m,y_m\bigr)\mapsto 
c_{m+1}+\sum_{\ell =1}^m c_{\ell}\log\bigl( | x_{\ell}+iy_{\ell}|^2\bigr)$
are continuous real-valued functions on ${\cal A}^{\Bbb R}_m (0,\infty)$. 
Hence, by appropriate applications of Lemma~8.1,  
with either $g(v)={\rm Re}(\Psi(v))$, or $g(v)={\rm Im}(\Psi(v))$, 
one finds that all partial derivatives of the two functions 
${\bf u}\mapsto {\rm Re}(F({\bf u}))$, ${\bf u}\mapsto {\rm Im}(F({\bf u}))\,$  
(where $F({\bf u})$ is given by~(8.3)) are continuous and real-valued 
on ${\cal A}^{\Bbb R}_m (0,\infty)$. This proves the first result of 
the lemma. The second result, in (8.4), follows almost immediately. 
Indeed, since the set 
$$W=\left\{ \bigl( x_1,y_1,\ldots ,x_m,y_m\bigr)\in{\Bbb R}^{2m} : 
\,r_{\ell}\leq\left| x_{\ell}+iy_{\ell}\right|\leq R_{\ell}\ \,{\rm for}\ \, 
\ell =1,\ldots ,m\right\}$$ 
is a compact subset of ${\cal A}^{\Bbb R}_m (0,\infty)$, 
the continuous real-valued function 
$$\bigl( x_1,y_1,\ldots ,x_m,y_m\bigr)\mapsto 
\left| {\partial^{j_1+k_1+\ldots +j_m+k_m}\over 
\partial x_1^{j_1}\partial y_1^{k_1}\cdots\partial x_m^{j_m}\partial y_m^{k_m}} 
\,F\bigl( x_1,y_1,\ldots ,x_m,y_m\bigr)\right| 
=\left| {\partial^{j_1+k_1+\ldots +j_m+k_m}\over 
\partial x_1^{j_1}\partial y_1^{k_1}\cdots\partial x_m^{j_m}\partial y_m^{k_m}} 
\,f({\bf z})\right|$$ 
must therefore attain its supremum on $W$; and, since ${\bf r}$ and ${\bf R}$ 
determine $W$, while 
the function concerned is determined by $\Psi$, ${\bf c}$, ${\bf j}$ 
and ${\bf k}$, that supremum 
is therefore determined by $\Psi$, ${\bf c}$, ${\bf j}$, ${\bf k}$, ${\bf r}$ and ${\bf R}$. 

Suppose now that $m=1$, $c_1\neq 0$, and that the support of $\Psi$ is contained 
in the bounded closed interval $[a,b]$. Let $f : {\Bbb C}\rightarrow{\Bbb C}$ 
be given by 
$$f(z)=\cases{ 0 &if $z=0$,\cr \Psi\!\left( c_2+c_1\log\left( |z|^2\right)\right)
 &otherwise.}$$
Since the equation~(8.3) is satisfied for all $z\in{\Bbb C}^{*}$, it 
follows by the first result of the lemma (proved above) that all partial 
derivatives of the functions $(x,y)\mapsto f(x+iy)$ are continuous 
on ${\cal A}^{\Bbb R}_1 (0,\infty)$. Moreover, since 
${\rm Supp}(\Psi)\subseteq [a,b]$, it follows by the definition of $f$ that 
$f(z)=0$ unless $a\leq c_2+c_1\log( |z|^2)\leq b$, and so 
(by elementary properties of the exponential and logarithm functions) 
we obtain the result (8.6). By~(8.6), one has $f(x+iy)=0$ for all $(x,y)$ in 
a neighbourhood of the point $(0,0)\in{\Bbb R}^2$. Therefore,  
in that neighbourhood, all partial derivatives of the function 
$(x,y)\mapsto f(x+iy)$ are defined and equal to zero, and so  are (in particular) 
continuous at the point $(0,0)$. This proves that the function 
$f : {\Bbb C}\rightarrow{\Bbb C}$ is smooth: for we showed already 
that all partial derivatives of the function 
$(x,y)\mapsto f(x+iy)$ are continuous on ${\cal A}^{\Bbb R}_1 (0,\infty)$. 
Since the Schwartz space ${\cal S}({\Bbb C})$ contains 
all smooth and compactly supported complex functions, it follows  
(given (8.6)) that we have $f\in{\cal S}({\Bbb C})$.   
The bound in (8.5) is an immediate corollary of (8.6) and the case 
$r=\exp\bigl(\min\{ a/c_1,b/c_1\}-(c_2/c_1)\bigr)$, 
$R=\exp\bigl(\max\{ a/c_1,b/c_1\}-(c_2/c_1)\bigr)$ of (8.4): 
no dependence on $a$ or $b$ is shown in (8.5), since 
(8.6) holds when $a$ and $b$ are, repectively, 
the infimum and supremum of ${\rm Supp}(\Psi)\ \blacksquare$ 

\bigskip 

\goodbreak\proclaim Lemma~8.3. Let $1\geq\delta >0$ and $t\in{\Bbb R}$. 
Let $f : {\Bbb C}\rightarrow{\Bbb C}$ be smooth. Suppose moreover that 
$${\partial^{j+k}\over\partial x^j\partial y^k}\,f(z)
\ll_{j,k} (\delta |z|)^{-(j+k)}\qquad\qquad  
\hbox{($j,k\in{\Bbb N}\cup\{ 0\}$, $\,z\in{\Bbb C}^{*}$),}\eqno(8.7)$$
and that ${\rm Supp}(f)\subset{\Bbb C}^{*}$.  
Let $g : {\Bbb C}\rightarrow{\Bbb C}$ be given by: 
$$g(z)=\cases{0 &if $z=0$,\cr f(z)|z|^{2it} &otherwise.}\eqno(8.8)$$
Then the function $g$ is smooth, has the same support as $f$, and is such that 
$${\partial^{j+k}\over\partial x^j\partial y^k}\,g(z)
\ll_{j,k} \left(\bigl(\delta^{-1}+|t|\bigr)^{-1} |z|\right)^{-(j+k)}\qquad\qquad  
\hbox{($j,k\in{\Bbb N}\cup\{ 0\}$, $\,z\in{\Bbb C}^{*}$).}\eqno(8.9)$$

\goodbreak 
\noindent{\bf Proof.}\quad 
The hypotheses of Lemma~8.2 are satisfied when $m=1$, ${\bf c}=(t,0)$ and 
$\Psi(v)=\exp(iv)$ ($v\in{\Bbb R}$). Hence it follows 
by Lemma~8.2 that the function $\chi : {\Bbb C}^{*}\rightarrow{\Bbb C}$ given by 
$$\chi(z)=|z|^{2it}=\exp\bigl( it\log( |z|^2 )\bigr)\qquad\qquad 
\hbox{($z\in{\Bbb C}^{*}$)}$$
is smooth. Since $f$ is (by hypothesis) a smooth function with domain ${\Bbb C}$, 
it follows by the product rule of differential calculus that the 
function $z\mapsto f(z)\chi(z)\ $ ($z\in{\Bbb C}^{*}$) is smooth. 
It can furthermore be deduced that the function $g : {\Bbb C}\rightarrow{\Bbb C}$, 
given by (8.8), is smooth. Indeed, since $0\not\in{\rm Supp}(f)$, one has 
$f(z)=0$ for all complex numbers $z$ lying in some neighbourhood of $0$; 
and so the function $(x,y)\mapsto g(x+iy)$, and all its partial derivatives, 
are defined and equal to zero on some neighbourhood of the point  
$(0,0)$ in ${\Bbb R}^2$; and those partial derivatives are therefore 
certainly continuous at $(0,0)$. This suffices to establish the 
smoothness of $g$, given that, for $(x,y)\in{\cal A}^{\Bbb R}_1 (0,\infty)$, we 
have $g(x+iy)=f(x+iy)\chi(x+iy)$, where the function 
$z\mapsto f(z)\chi(z)$ is smooth on ${\Bbb C}^{*}$. 
Since $|r^{2it}| =1\ $ ($r>0$), and since 
$g(0)=f(0)=0$, it follows by the definition (8.8) that 
all zeros of $f$ are zeros of $g$, and vice versa. 
Therefore ${\rm Supp}(g)={\rm Supp}(f)$. 

Suppose now that $j,k\in{\Bbb N}\cup\{ 0\}$, and that $z\in{\Bbb C}^{*}$. 
In order to obtain the bound (8.9), we note firstly that, by (8.7) and 
Leibniz's rule for higher order derivatives of products, one has: 
$$\left|{\partial^{j+k}\over\partial x^j\partial y^k}\,f(z)\chi(z)\right| 
\leq 2^{j+k} 
\left|{\partial^{\ell +m}\over\partial x^{\ell}\partial y^m}\,\chi(z)\right| 
O_{j-\ell,k-m}\left(\left(\delta |z|\right)^{-\left( (j-\ell)+(k-m)\right)}\right) , 
\eqno(8.10)$$
for some non-negative integers $\ell ,m$ with $\ell\leq j$ and $m\leq k$. 
Then we observe that, by (5.19), one has 
$$\eqalignno{
\left|{\partial^{\ell +m}\over\partial x^{\ell}\partial y^m}\,\chi(z)\right| 
 &=\left| i^m\left( {\partial\over\partial z}+{\partial\over\partial\overline{z}}\right)^{\ell} 
\left( {\partial\over\partial z}-{\partial\over\partial\overline{z}}\right)^m \chi(z)\right| \leq {}\cr  
 &\leq 2^{\ell +m}\left| {\partial^{r+s}\over\partial z^r\partial\overline{z}^{\,s}}\,\chi(z)\right|  
=2^{\ell +m}\left| {\partial^{r+s}\over\partial z^r\partial\overline{z}^{\,s}}
\,|z|^{2it}\right|\;,&(8.11)}$$ 
for some non-negative integers $r,s$ satisfying $r+s=\ell+m\leq j+k$. 

Real-variable calculus shows that 
$${\partial\over\partial z}\,|z|^{2it}
=it|z|^{2it} z^{-1}\qquad\ \hbox{and}\qquad\  
{\partial\over\partial\overline{z}}\,|z|^{2it}
=it|z|^{2it} (\overline{z})^{-1}$$
(this also follows, by (5.21), from the fact that $|z|^{2it}=z^{it}(\overline{z})^{it}$). 
By induction (and with the aid of (5.21) and the product rule of differential calculus), 
it may therefore be deduced that, for $n\in{\Bbb N}\cup\{ 0\}$, one has: 
$${\partial^n\over\partial z^n}\,|z|^{2it}
=|z|^{2it} z^{-n}(it-n+1)_n\qquad\ \hbox{and}\qquad\  
{\partial^n\over\partial\overline{z}^{\,n}}\,|z|^{2it}
=|z|^{2it} (\overline{z})^{-n}(it-n+1)_n\;,$$
where, as in (2.1), $(\alpha)_m=\alpha(\alpha +1)\cdots (\alpha +m-1)=\Gamma(\alpha +m)/\Gamma(\alpha)$. 
Hence, and by (5.21),  
$$\eqalign{
{\partial^{r+s}\over\partial z^r\partial\overline{z}^{\,s}}
\,|z|^{2it} 
 &=(it-s+1)_s\left(\overline{z}\right)^{-s}
{\partial^r\over\partial z^r}\,|z|^{2it} = {}\cr
 &=(it-s+1)_s (it-r+1)_r 
\left(\overline{z}\right)^{-s} |z|^{2it} z^{-r} \ll_{r+s} {}\cr 
 &\ll_{r+s} \left( 1+|t|\right)^{r+s} |z|^{-(r+s)} 
=\left( 1+|t|\right)^{\ell +m} |z|^{-(\ell +m)}\;.}$$
By this last bound, in combination with (8.10) and (8.11), 
we find that 
$$\eqalign{
\left|{\partial^{j+k}\over\partial x^j\partial y^k}\,f(z)\chi(z)\right|  
 &\ll_{j,k} \left( 1+|t|\right)^{\ell +m} |z|^{-(\ell +m)}
 \left(\delta |z|\right)^{-\left( (j-\ell)+(k-m)\right)} = {}\cr 
 &= \left( 1+|t|\right)^{\ell +m} \delta^{\ell+m -(j+k)} |z|^{-(j+k)} \leq {}\cr 
 &\leq\left( 1+\left( 1+|t|\right)^{j+k}\delta^{j+k}\right)\delta^{-(j+k)} |z|^{-(j+k)} = {}\cr 
 &=\left(\delta^{-(j+k)}+ \left( 1+|t|\right)^{j+k}\right) |z|^{-(j+k)}
\leq\left(\delta^{-1}+1+|t|\right)^{j+k} |z|^{-(j+k)}\;.
}$$
Consequently (since $\delta^{-1}\geq 1$, and since $g(w)=f(w)\chi(w)$ for 
$w\in{\Bbb C}^{*}$), we have the bound (8.9)\ $\blacksquare$ 

\bigskip 

\goodbreak 
\proclaim Lemma~8.4. Let $H,K,N,\delta >0$ and the functions 
$\alpha,\beta : {\Bbb C}\rightarrow{\Bbb C}$ be such as to 
satisfy the hypotheses 
of Theorem~9 (so that one has, in particular $HK=N\geq 1$); 
and, for $n\in{\frak O}-\{ 0\}$, let $a_n$ be given by 
the equation~(1.3.14). Suppose moreover that 
$$K\geq H\;,\eqno(8.12)$$ 
and that $\varepsilon >0$, $t\in{\Bbb R}$, $X\geq 1$ and 
$$R\in\bigl[ N\,,\,8 N^{4/3}\bigr]\;.\eqno(8.13)$$ 
Then 
$$S_t(R,X,N)
\ll_{\varepsilon} (XN)^{\varepsilon} 
\left(\delta^{-1}+|t|\right)^{11} 
\!\left( R+\left( {X\over K}\right) N 
+\left( {X\over R H^{-1}}\right)^{\!\!3/2} R^{1/2} N^{1/2}\right) \!N\;.\eqno(8.14)$$ 

\goodbreak 
\noindent{\bf Proof.}\quad 
Let $\delta_t =(\delta^{-1}+|t|)^{-1}$; and let 
$\alpha_t,\beta_t : {\Bbb C}\rightarrow{\Bbb C}$ be the complex functions 
satisfying both $\alpha_t(0)=\beta_t(0)=0$ and, 
for all $z\in{\Bbb C}^{*}$, 
$\,\alpha_t(z)=\alpha(z)|z|^{2it}$ and $\beta_t(z)=\beta(z)|z|^{2it}$. 
Then, by Lemma~8.3, the hypotheses of Theorem~9 concerning $H,K,N,\delta$ and the 
functions $\alpha$ and $\beta$ will continue to be satisfied if 
we substitute $\alpha_t$, $\beta_t$ and $\delta_t\in(0,1]$ for $\alpha$, $\beta$ 
and $\delta$, respectively (while making no change to $H$, $K$ and $N$). 
Therefore, given the definition (1.3.2) of $S_t(Q,X,N)$, and given (1.3.14),  
the cases of Lemma~8.4 in which $t\neq 0$ are a corollary of the  
particular case $t=0\,$ (i.e. a corollary obtained by applying that case of Lemma~8.4 
with $\alpha_t$, $\beta_t$ and $\delta_t$ 
substituted for $\alpha$, $\beta$ and $\delta$, respectively): 
for, by (1.3.14), one has, for $n\in{\frak O}-\{ 0\}$, 
$$\sum_{h\mid n}\alpha_t(h) \beta_t\left( {n\over h}\right)  
=\sum_{h\mid n}\alpha(h) |h|^{2it}\beta\left( {n\over h}\right)
\left| {n\over h}\right|^{2it} 
=\sum_{h\mid n}\alpha(h)\beta\left( {n\over h}\right)|n|^{2it}
=a_n |n|^{2it}$$
(which is the coefficient independent of $V$ in (1.3.2));  
and, with regard to the factor involving $\delta$ in (8.14),  
one has also $(\delta_t^{-1}+|0|)^{11}=\delta_t^{-11}=(\delta^{-1}+|t|)^{11}$. 

By the foregoing observations,  
we now have only to prove the case $t=0$ of the lemma. 
Moreover, in doing so we may suppose that 
$$X\geq 64 K\;.\eqno(8.15)$$ 
For if $1\leq X<64K$ then it follows by Theorem~6 that 
$S_0(R,X,N)\leq S_0(R,64K,N)$; and so, if the bound (8.14) holds  
for $t=0$, $X=64K$ and all $\varepsilon >0$, then, when $1\leq X<64K$, one has:  
$$\eqalign{ 
S_0(R,X,N) &\ll_{\varepsilon /2} 
(KN)^{\varepsilon /2}\delta^{-11}\left( R+N+(HK/R)^{3/2}R^{1/2}N^{1/2}\right) N = {}\cr 
 &=\left( H^{-1} N^2\right)^{\varepsilon /2}\delta^{-11}\left( R+N+R^{-1}N^2\right) N \ll {}\cr 
 &\ll N^{\varepsilon}\delta^{-11}RN\leq (XN)^{\varepsilon}\delta^{-11}RN}$$
(given (8.13)), which implies the case $t=0$ of the bound in (8.14).

For $z\in{\Bbb C}$ and $Z>0$, let 
$$\omega(Z;z)=\Omega\left( Z^{-1/2} z\right) ,\eqno(8.16)$$
where the function $\Omega : {\Bbb C}\rightarrow [0,1/e]$ is given by 
$$\Omega(w)=\cases{0 &if $w=0$, \cr 
\ &\ \cr
\Phi\!\left( 1+\displaystyle{\log\bigl( |w|^2\bigr)\over\log\bigl( 2^{1/2}\bigr)}\right) &otherwise,}
\eqno(8.17)$$
with $\Phi : {\Bbb R}\rightarrow [0,1/e]$ being the specific infinitely 
differentiable function defined, below the equation~(3.5), 
in the proof of Theorem~4.  Then, since ${\rm Supp}(\Phi)=[-1,1]$, and since 
$\Phi(v)\geq\Phi(1/2)=\exp(-4/3)>1/4$ for $-1/2\leq v\leq 1/2$, it 
follows that, for 
$w\in{\Bbb C}$, $z\in{\Bbb C}$ and $Z>0$, one has: 
$$\Omega(w)=0\quad{\rm unless}\quad{1\over 2}<|w|^2<1\;;\qquad\quad  
\omega(Z;z)=0\quad{\rm unless}\quad{Z\over 2}<|z|^2<Z\;;\eqno(8.18)$$
$$\omega(Z;z)\geq 0\;;\qquad\quad{\rm and}\qquad\quad  
\omega(Z;z)>{1\over 4}\quad{\rm if}\quad 2^{-3/4}Z\leq |z|^2\leq 2^{-1/4}Z\;.\eqno(8.19)$$ 
We furthermore define, for $Q>0$, 
$$T(Q,X,N) 
=\sum_{0\neq q\in{\frak O}}\omega(Q;q)
\sum_{\scriptstyle V\atop\scriptstyle\nu_V>0}^{(\Gamma_0(q))} 
{\bf K}f\bigl(\nu_V,0\bigr)
\Biggl|\sum_{{\textstyle{N\over 4}}<|n|^2\leq N} 
a_n c_V^{\infty}\bigl( n;\nu_V,0\bigr)\Biggr|^2 , 
\eqno(8.20)$$ 
with $f : {\Bbb C}^{*}\rightarrow[0,1/e]$ as defined 
below the equation~(3.5) in the proof of 
Theorem~4 (so that $f$ depends on, and is determined by, $X$); and with the 
transform ${\bf K}f(\nu,p)$ as defined in the statement of Theorem~1. 
Given the inequalities in (8.19) and the lower bound on $X$ in (8.15), 
which implies that $X\geq 64$, it follows therefore (similarly to 
how, in the proof of Theorem~7, the result (4.3) was obtained) 
that, for $Q>0$, we have 
$$\eqalign{ T(Q,X,N)
 &\gg \sum_{0\neq q\in{\frak O}}\omega(Q;q)
\sum_{\scriptstyle V\atop\scriptstyle\nu_V>0}^{(\Gamma_0(q))} 
X^{\nu_V} 
\Biggl|\sum_{{\textstyle{N\over 4}}<|n|^2\leq N} 
a_n c_V^{\infty}\bigl( n;\nu_V,0\bigr)\Biggr|^2 \gg {}\cr 
 &\gg\sum_{2^{-3/4}Q<|q|^2\leq 2^{-1/4}Q} 
\,\sum_{\scriptstyle V\atop\scriptstyle\nu_V>0}^{(\Gamma_0(q))} 
X^{\nu_V} 
\Biggl|\sum_{{\textstyle{N\over 4}}<|n|^2\leq N} 
a_n c_V^{\infty}\bigl( n;\nu_V,0\bigr)\Biggr|^2 .}$$
Hence, and since $(R/2,R]=(R/2,2^{-1/2}R]\cup(2^{-1/2}R,R]$, 
it follows by (1.3.2) that 
$$ T(Q,X,N)\gg S_0(R,X,N)\;,\eqno(8.21)$$
for some $Q$ satisfying 
$$Q\in\left\{ 2^{-1/4}R\,,\,2^{1/4}R\right\}\qquad\quad{\rm and}\qquad\quad 
Q\geq 2^{1/4}\;.\eqno(8.22)$$
We may therefore assume, in what follows, that which is stated in (8.22) and (8.21).  

Given the similar forms of the sums over $q$ and $V$ occurring in (8.20) and 
in the result (4.3)$\,$ (within the proof of Theorem~7), it follows by 
steps differing in only one minor respect from the steps 
taken in passing from (4.3) to (4.5)-(4.6) that we must either have 
$$ T(Q,X,N)
\ll (\log X)\left( Q+O_{\varepsilon}\left( N^{1+\varepsilon}\right)\right) 
\left\|{\bf a}_N\right\|_2^2\;,\eqno(8.23)$$ 
or else 
$$ T(Q,X,N)
\ll\Biggl|\sum_{2^{-3/4}P<|p|^2\leq 2^{-1/4}P} 
\quad\ \,\sum\!\!\!\!\!\!\sum_{
\!\!\!\!\!\!\!\!\!\!\!\!\!{{\textstyle{N\over 4}}<|\ell|^2,|m|^2\leq N}} 
\overline{a_{\ell}}\ a_m
\sum_{0\neq q\in{\frak O}}{S(\ell,m;pq)\over |pq|^2}
\,f\!\!\left( {2\pi\sqrt{\ell m}\over pq}\right)\omega(Q;q)\Biggr|\;,\eqno(8.24)$$
for some $P$ satisfying 
$$P\in\left\{ 
2^{-3/4}Q^{*}\,,\,2^{-5/4}Q^{*}\,,\,\ldots\,,\,2^{-27/4}Q^{*}\,,\,2^{-29/4}Q^{*}\right\} 
 ,\eqno(8.25)$$
where 
$$Q^{*}=64\pi^2 XN/Q\eqno(8.26)$$
(so that the relationship between $Q$ and $Q^{*}$ is the same as  
it is in Theorem~7). 
In either case we find  by (1.3.14), (1.3.15) and (1.2.11) that  
$$\left\|{\bf a}_N\right\|_2^2 
=\sum_{0<|n|^2\leq N}\biggl|\sum_{h\mid n}O(1)\biggr|^2 
=\sum_{0<|n|^2\leq N}\left(O_{\varepsilon}\!\left( |n|^{\varepsilon /4}\right)\right)^2 
\ll_{\varepsilon} N^{1+(\varepsilon /4)}\;.\eqno(8.27)$$
Since $HK=N\geq 1$, it follows by (8.15), (8.12) and (8.13)  
that $XN\geq 64KN\geq 64 N^{3/2}\geq 64 N^{4/3}\geq 8 R$. Consequently it is 
implied by the conditions in
(8.22), (8.25) and (8.26) that we have $P\geq 2^{3/2}\pi^2$. 

Note that in obtaining (8.24) we employ a division of the sum over $p$ 
which is, in a sense, `twice as fine' as the corresponding division of a sum  
used to obtain the bound (4.6), in the proof of Theorem~7. 
The conditional conclusion (8.24)-(8.26) is justified, 
since for $L/2<M<L$, and any coefficients 
$c_p\in{\Bbb C}$ ($p\in{\frak O}-\{ 0\}$), one has 
$${\cal C}(L/2,L)\leq {\cal C}(L/2,M)+{\cal C}(M,L)
\leq 2\max\!\left\{\,{\cal C}(L/2,M)\,,\,{\cal C}(M,L)\,\right\}\;,$$ 
where ${\cal C}(a,b)=|\sum_{a<|p|^2\leq b} c_p\,|$.

If the bound (8.23) holds then, by (8.21)-(8.23), (8.13) and (8.27), one has 
$$S_0(R,X,N)\ll (\log X)\,O_{\varepsilon}\!\left( R N^{\varepsilon}\right) 
\left\|{\bf a}_N\right\|_2^2 
\ll_{\varepsilon}X^{\varepsilon /2} R N^{1+(5\varepsilon /4)}<(XN)^{\varepsilon}RN
\eqno(8.28)$$
(the last inequality following since, by (8.15) and (8.12), one has 
$X^2>K^2\geq HK$ where, by hypothesis, $HK=N$).       
This means that we obtain the case $t=0$ 
of the bound in (8.14) when (8.23) holds 
(given that we have, by hypothesis, $1\geq\delta >0$). Therefore, bearing in mind  
what we concluded in (8.23)-(8.26), we may complete  
this proof by showing that the case $t=0$ of the bound in (8.14)  
holds if the hypotheses of the lemma and the 
conditions in (8.15) and (8.24)-(8.26) are satisfied. Accordingly, we 
now add to our hypotheses by supposing 
that the conditions in (8.24)-(8.26) are satisfied. 

By (8.24), (8.18) and (8.19), 
$$\eqalignno{ T(Q,X,N)
 &\ll\sum_{2^{-3/4}P<|p|^2\leq 2^{-1/4}P} 
\Biggl|\quad\ \,\sum\!\!\!\!\!\!\sum_{
\!\!\!\!\!\!\!\!\!\!\!\!\!{{\textstyle{N\over 4}}<|\ell|^2,|m|^2\leq N}} 
\overline{a_{\ell}}\ a_m
\sum_{0\neq q\in{\frak O}}{S(\ell,m;pq)\over |pq|^2}
\,f\!\!\left( {2\pi\sqrt{\ell m}\over pq}\right)\omega(Q;q)\Biggr| \ll {}\qquad\quad\cr 
 &\ll\sum_{0\neq p\in{\frak O}} 
\omega(P;p)\Biggl|\quad\ \,\sum\!\!\!\!\!\!\sum_{
\!\!\!\!\!\!\!\!\!\!\!\!\!{{\textstyle{N\over 4}}<|\ell|^2,|m|^2\leq N}} 
\overline{a_{\ell}}\ a_m
\sum_{0\neq q\in{\frak O}}{S(\ell,m;pq)\over |pq|^2}
\,f\!\!\left( {2\pi\sqrt{\ell m}\over pq}\right)\omega(Q;q)\Biggr| = {}\cr 
 &=\sum_{0\neq p\in{\frak O}} 
\omega(P;p)\theta_p\quad\ \,\sum\!\!\!\!\!\!\sum_{
\!\!\!\!\!\!\!\!\!\!\!\!\!{{\textstyle{N\over 4}}<|\ell|^2,|m|^2\leq N}} 
\overline{a_{\ell}}\ a_m
\sum_{0\neq q\in{\frak O}}{S(\ell,m;pq)\over |pq|^2}
\,f\!\!\left( {2\pi\sqrt{\ell m}\over pq}\right)\omega(Q;q)\;, &(8.29)}$$
where $\omega(Z;z)$ is as defined in (and below) (8.16) and (8.17);     
while, for $p\in{\frak O}-\{ 0\}$, the coefficient $\theta_p$ is a 
complex number determined by $Q,X,N,P$, the coefficients $a_n\ $    
($0\neq n\in{\frak O}$) and $p$, and moreover satisfies: 
$$\left|\theta_p\right| =\cases{1 &if $P/2<|p|^2<P$, \cr 0 &otherwise.}\eqno(8.30)$$

By a division of the range of summation for the variable $\ell$ in (8.29), one may 
deduce that 
$$ T(Q,X,N)
\ll\Biggl|\sum_{2^{-3/4}L<|\ell|^2\leq 2^{-1/4}L}\overline{a_{\ell}}
\sum_{p\neq 0} 
\omega(P;p)\theta_p\sum_{{\textstyle{N\over 4}}<|m|^2\leq N} a_m 
\sum_{q\neq 0}{S(\ell,m;pq)\over |pq|^2}
\,f\!\!\left( {2\pi\sqrt{\ell m}\over pq}\right)\omega(Q;q)\Biggr|\;,$$
for some $L$ satisfying 
$$L\in\left\{ 2^{1/4}N\,,\,2^{-1/4}N\,,\,2^{-3/4}N\,,\,2^{-5/4}N\right\}\qquad\quad 
{\rm and}\qquad\quad  
L\geq 2^{1/4}\;.\eqno(8.31)$$
Hence, by steps similar to those by which (8.29) was obtained from (8.24), 
we find that 
$$ T(Q,X,N)
\ll\sum_{\ell\neq 0}\omega(L;\ell) \Upsilon_{\ell} 
\sum_{p\neq 0} 
\omega(P;p)\theta_p\sum_{{\textstyle{N\over 4}}<|m|^2\leq N} a_m 
\sum_{q\neq 0}{S(\ell,m;pq)\over |pq|^2}
\,f\!\!\left( {2\pi\sqrt{\ell m}\over pq}\right)\omega(Q;q)\;,\eqno(8.32)$$ 
where $\omega(Z;z)$ is as defined in (and below) (8.16) and (8.17), 
while, for $\ell\in{\frak O}-\{ 0\}$, the coefficient $\Upsilon_{\ell}$ is a 
complex number determined by $Q,X,N,P,L$, the coefficients $a_n\ $  
($0\neq n\in{\frak O}$) and the variable $\ell$, and is such that 
$$\left|\Upsilon_{\ell}\right| =\cases{\left| a_{\ell}\right| &if $L/2<|\ell|^2<L$, \cr 
0 &otherwise.}\eqno(8.33)$$

If $w\in{\frak O}-\{ 0\}$ then the mapping 
$d\bmod w{\frak O}\mapsto d^{*}\bmod w{\frak O}$ is 
(as is evident from our definition of the meaning of $d^{*}$ in this context) an 
involution on the set of elements of the multiplicative group 
$({\frak O}/w{\frak O})^{*}$; and so it follows by the definition 
(1.3.6) that, in (8.32), one has 
$$S(\ell , m ; pq)=S(m , \ell; pq)\;.$$
Given this elementary fact, and given the hypotheses of Theorem~9 regarding 
${\rm Supp}(\alpha)$ and ${\rm Supp}(\beta)$, it follows by (8.32) and the 
definition (1.3.14) of the coefficient $a_n$, that we have now 
$$ T(Q,X,N)
\ll\sum_{p\neq 0}\theta_p |p|^{-2} 
\sum_{q\neq 0} |q|^{-2}\sum_h \phi_h\sum_k 
\sum_{\ell\neq 0}  S(hk,\ell ;pq)\,\varphi(h,k,\ell,p,q)\,\Upsilon_{\ell}\;,\eqno(8.34)$$ 
where, for ${\bf z}\in{\Bbb C}^5$,  
$$\varphi\bigl( z_1 , z_2 , z_3 , z_4 , z_5\bigr) 
=\cases{
f\!\!\left( \displaystyle{2\pi\sqrt{z_1 z_2 z_3}\over z_4 z_5}\right)
\!\alpha\bigl( z_1\bigr)\beta\bigl( z_2\bigr)\omega\bigl( L; z_3\bigr)
\omega\bigl( P; z_4\bigr)\omega\bigl( Q; z_5\bigr) &if $z_1,\ldots,z_5\neq 0$, \cr 
\ &\ \cr
0 &otherwise,}\eqno(8.35)$$
while, for $h\in{\frak O}-\{ 0\}$,  the (effectively 
redundant) factor $\phi_h$ is given by: 
$$\phi_h=\cases{1 &if $H/2<|h|^2\leq H$, \cr 0 &otherwise.}\eqno(8.36)$$

We shall complete this proof by applying Lemma~6.3 to obtain an upper bound 
for the sum on the right-hand side of (8.34) (that sum being similar in 
form to the sum ${\cal R}$ defined by the equation~(6.1)). 
In order to justify this it is necessary to first verify that the 
function $\varphi : {\Bbb C}^5\rightarrow{\Bbb C}$ given by (8.35) 
satisfies all the relevant hypotheses stated in the first paragraph of 
Section~6. To simplify matters we first reformulate those hypotheses 
in terms of the function $\xi : {\Bbb C}^5\rightarrow{\Bbb C}$ given by:  
$$\xi({\bf z})
=\varphi\bigl( H^{1/2}z_1\,,\,K^{1/2}z_2\,,\,L^{1/2}z_3\,,\,P^{1/2}z_4\,,\,Q^{1/2}z_5\bigr)\qquad\quad  
\hbox{for $\quad{\bf z}\in{\Bbb C}^5$.}\eqno(8.37)$$ 
The relevant hypotheses concerning $\varphi$ are that 
all partial derivatives of the function 
$(x_1,y_1,\ldots ,x_5,y_5)\mapsto\varphi(x_1+iy_1,\ldots ,x_5+y_5)$ 
are defined and continuous on ${\Bbb R}^{10}$;  
that $\varphi : {\Bbb C}^5\rightarrow{\Bbb C}$ has the property (6.2); 
and that, when 
${\bf x},{\bf y}\in{\Bbb R}^5$ are such that 
$x_d+iy_d\neq 0$ for $d=1,\ldots ,5$,  
the bound (6.3) holds for   
all ${\bf j},{\bf k}\in({\Bbb N}\cup\{ 0\})^5$. 
Since $H,K,L,P,Q>0$, it follows by 
the chain rule of differential 
calculus that these hypotheses concerning $\varphi$ 
are satisfied if and only if the function $\xi : {\Bbb C}^5\rightarrow{\Bbb C}$ 
given by (8.37) satisfies three particular conditions. 
The first of these conditions is that all partial derivatives of the function 
$(x_1,y_1,\ldots ,x_5,y_5)\mapsto\xi(x_1,y_1,\ldots ,x_5,y_5)$ 
be defined and continuous on ${\Bbb R}^{10}$. The second condition is that 
$\xi : {\Bbb C}^5\rightarrow{\Bbb C}$ satisfy  
$$\xi({\bf z})=0\qquad\hbox{unless}\qquad 
{1\over 2}<\bigl| z_1\bigr|^2,\ldots ,\bigl| z_5\bigr|^2<1\;.\eqno(8.38)$$
The third (and final) condition is that one have, 
for ${\bf j},{\bf k}\in({\Bbb N}\cup\{ 0\})^5$ and ${\bf x},{\bf y}\in{\Bbb R}^5$, 
the bound: 
$${\partial^{j_1+\cdots +j_5+k_1+\cdots +k_5}\over 
\partial x_1^{j_1}\cdots\partial x_5^{j_5} 
\partial y_1^{k_1}\cdots\partial y_5^{k_5}}\,\xi\bigl( x_1+iy_1 , \ldots , x_5+iy_5\bigr) 
\ll_{{\bf j},{\bf k}} 
\prod_{h=1}^5 \delta^{-(j_h+k_h)}\;.\eqno(8.39)$$
We do not claim that (8.39) is, by itself, equivalent to (6.3): it is 
the combination of (8.38) and (8.39) which is equivalent to the  
combination of (6.2) and (6.3). 

In preparation for the application of Lemma~6.3, we show now that 
the function $\xi({\bf z})$ satisfies the three conditions just mentioned: 
that being sufficient (given the observations of the preceding paragraph) 
to establish that the function $\varphi({\bf z})$, 
in (8.34) and (8.35), satisfies all of the requisite hypotheses.
By (8.16), (8.35) and (8.37), we have, for ${\bf z}\in{\Bbb C}^5$, 
$$\xi({\bf z})
=\cases{ g({\bf z}) A\bigl( z_1\bigr) B\bigl( z_2\bigr) \prod\limits_{d=3}^5 \Omega\bigl( z_d\bigr) 
 &if $z_1,\ldots ,z_5\neq 0$, \cr 0 &otherwise,}\eqno(8.40)$$
where $\Omega(z)$ is given by (8.17), while 
$$A(z)=\alpha\bigl( H^{1/2} z\bigr)\qquad{\rm and}\qquad    
B(z)=\beta\bigl( K^{1/2} z\bigr)\qquad\qquad\quad\hbox{($z\in{\Bbb C}$),}\eqno(8.41)$$
and 
$$g\bigl( z_1,\ldots ,z_5\bigr) 
=f\left( {Y\sqrt{z_1 z_2 z_3}\over z_4 z_5}\right)\qquad\qquad\quad   
\hbox{($\,{\bf z}\in\left( {\Bbb C}^{*}\right)^5\,$),}\eqno(8.42)$$
with (given (8.25), (8.26), (8.31) and the hypothesis that $HK=N$): 
$$Y={2\pi (HKL)^{1/4}\over (PQ)^{1/2}}=2^{\eta /16} X^{-1/2}\;,\eqno(8.43)$$
for some odd integer $\eta\in[-31,27]$, independent of the variables 
$z_1,\ldots ,z_5$.  
In order that we may reach the desired conclusions concerning the function $\xi({\bf z})$, 
we must first establish certain related facts about the above  
functions $A(z)$, $B(z)$, $\Omega(z)$ and 
$g({\bf z})$. 

By hypothesis, the functions $\alpha , \beta : {\Bbb C}\rightarrow{\Bbb C}$ 
are smooth. Hence, and by the chain-rule of differential calculus, the 
functions $A , B : {\Bbb C}\rightarrow{\Bbb C}$ given by (8.41) are also smooth. 
Moreover, the hypotheses of Theorem~9 concerning ${\rm Supp}(\alpha)$ 
and ${\rm Supp}(\beta)$ imply that, when $z_1,z_2\in{\Bbb C}$, one has: 
$$A\left( z_1\right) =0\quad{\rm unless}\quad 
{1\over 2}<\left| z_1\right|^2<1\;;\qquad\quad  
B\left( z_2\right) =0\quad{\rm unless}\quad 
{1\over 2}<\left| z_2\right|^2<1\;.\eqno(8.44)$$
By (8.41), (8.44) and the hypothesis (1.3.15), we may deduce that, 
for $j,k\in{\Bbb N}\cup\{ 0\}$ and $x,y\in{\Bbb R}$,  
$$\max\left\{ \left| {\partial^{j+k}\over\partial x^j\partial y^k}\,A(x+iy)\right|\,, 
\,\left| {\partial^{j+k}\over\partial x^j\partial y^k}\,B(x+iy)\right|\right\} 
\ll_{j,k}\delta^{-(j+k)}\eqno(8.45)$$
(note that we are here using (8.44) for the cases where $|x+iy|^2<1/2$). 

We turn next to the function $\Omega(z)$, which is defined by (8.17)$\,$  
(with $\Phi : {\Bbb R}\rightarrow [0,1/e]$ there as given below the equation~(3.5)). 
By the case $m=1$, ${\bf c}=(2/\log 2 , 1)$, $\Psi=\Phi$ of Lemma~8.2, the function 
$\Omega : {\Bbb C}\rightarrow{\Bbb C}$ is smooth (indeed it lies in 
the Schwartz space ${\cal S}({\Bbb C})$); the results (8.5), (8.6) of that lemma 
imply the result already noted in the first part of (8.18), 
and show also that     
$${\partial^{j+k}\over\partial x^j\partial y^k}\,\Omega(x+iy) 
=O_{\Phi,j,k}(1)=O_{j,k}(1)\qquad\qquad  
\hbox{($j,k\in{\Bbb N}\cup\{ 0\}$, $z=x+iy\in{\Bbb C}$).}\eqno(8.46)$$ 

The last two paragraphs contain all that we need concerning the functions 
$A(z)$, $B(z)$ and $\Omega(z)$. As for the function 
$g : \left( {\Bbb C}^{*}\right)^5\rightarrow{\Bbb C}$, it follows by 
(8.42), (8.43), and the definition of $f\,$ (below (3.5)), that 
$$g({\bf z})
=\Phi\!\left( (\log 2)^{-1}\log\!\left( X^{1/2}\left| 
{Y\sqrt{z_1 z_2 z_3}\over z_4 z_5}\right|\right)\right) 
=\Phi\!\left( c_6+\sum_{d=1}^5 c_d \log\bigl(\left| z_d\right|^2\bigr)\right) 
\qquad\qquad\hbox{(${\bf z}\in\left(\,{\Bbb C}^{*}\right)^5\,$),}$$
with ${\bf c}\in{\Bbb R}^6$ given by  
$c_1=c_2=c_3=(4\log 2)^{-1}$, $\,c_4=c_5=-(2\log 2)^{-1}$, $\,c_6=\eta/16$ 
(where $\eta\in[-31,27]$ is the integer constant in (8.43)),   
and with $\Phi : {\Bbb R}\rightarrow [0,1/e]$ as defined below the equation~(3.5). 
Therefore, by the case $m=5$, $\Psi =\Phi$ of Lemma~8.2, it 
follows that all partial derivatives of the function 
$(x_1,y_1,\ldots ,x_5,y_5)\mapsto g(x_1+iy_1,\ldots ,x_5+iy_5)$ are defined 
and continuous on the set ${\cal A}^{\Bbb R}_5 (0,\infty)\subset{\Bbb R}^{10}$ 
(defined as in (8.1)); given the 
result (8.4) of Lemma~8.2, it is moreover the case that, 
for ${\bf j},{\bf k}\in({\Bbb N}\cup\{ 0\})^5$ and 
${\bf z}={\bf x}+i{\bf y}\in\bigl({\Bbb C}^{*}\bigr)^5$, one has 
$${\partial^{j_1+\cdots +j_5+k_1+\cdots +k_5}\over 
\partial x_1^{j_1}\cdots\partial x_5^{j_5} 
\partial y_1^{k_1}\cdots\partial y_5^{k_5}}\,g\bigl( {\bf z}\bigr) 
=O_{\Phi,\eta,{\bf j},{\bf k}}(1)
=O_{{\bf j},{\bf k}}(1)\qquad\      
{\rm when}\quad\ {1\over 2}\leq\left| z_1\right|^2\!,\ldots ,\left| z_5\right|^2\leq 1\;.\eqno(8.47)$$

Given (8.40) and (8.42), it is an immediate consequence of the 
properties of the functions $\Omega,A,B : {\Bbb C}\rightarrow{\Bbb C}$ 
noted in (8.18) and (8.44) that $\xi({\bf z})$ satisfies 
the condition (8.38). Moreover, given the conclusions of the 
preceding paragraph (and since we showed earlier that 
the functions $A,B,\Omega : {\Bbb C}\rightarrow{\Bbb C}$ are smooth), 
it follows by (8.40) and the product rule that all partial derivatives 
of the function $(x_1,y_1,\ldots ,x_5,y_5)\mapsto\xi(x_1+iy_1,\ldots ,x_5+iy_5)$ 
are defined and continuous on ${\cal A}^{\Bbb R}_5 (0,\infty)$. 
All those partial derivatives are, furthermore, defined and continuous 
on ${\Bbb R}^{10}$: for if 
${\bf p}=(x_1',y_1',\ldots ,x_5',y_5')$ is a point of ${\Bbb R}^{10}$ 
not included in the set ${\cal A}^{\Bbb R}_5 (0,\infty)$ 
then, by (8.1) and (8.38), one has $\xi(x_1+iy_1,\ldots ,x_5+iy_5)=0$ for all points 
$(x_1,y_1,\ldots ,x_5,y_5)$ in the open Euclidean ball in ${\Bbb R}^{10}$ with 
centre ${\bf p}$ and radius $1/\sqrt{2}$, and so any partial derivative  
of the function $(x_1,y_1,\ldots ,x_5,y_5)\mapsto\xi(x_1+iy_1,\ldots ,x_5+iy_5)$ 
is defined and equal to zero on that open Euclidean ball, and is therefore 
certainly continuous at the point ${\bf p}$. 

The above shows that the function $\xi$ satisfies the first two of the 
three conditions stated below (8.37). In order to verify that the final condition 
there is also satisfied, we begin by observing that, since the set 
${\cal A}^{\Bbb R}_5 \bigl[ 2^{-1/2} , 1\bigr]$
is a closed subset of ${\Bbb R}^{10}$, it therefore follows from (8.38) 
that the bound (8.39) holds when  one has 
${\bf j},{\bf k}\in\bigl({\Bbb N}\cup\{ 0\}\bigr)^5$ and 
$(x_1,y_1,\ldots ,x_5,y_5)\in{\Bbb R}^{10}-{\cal A}^{\Bbb R}_5 \bigl[ 2^{-1/2} , 1\bigr]$: 
for in that case the 
partial derivative which appears in (8.39) is equal to zero. 
On the other hand, when 
${\bf j},{\bf k}\in\left({\Bbb N}\cup\{ 0\}\right)^5$ and 
${\bf z}={\bf x}+i{\bf y}\in\left( {\Bbb C}^{*}\right)^5$, it follows by (8.40) and 
Leibniz's rule for 
higher order derivatives of a product that,  
for some 
${\bf r},{\bf s},{\bf t},{\bf u}$ with 
$${\bf r},{\bf s},{\bf t},{\bf u}\in({\Bbb N}\cup\{ 0\})^5\;,\qquad\quad  
{\bf r}+{\bf t}={\bf j}\qquad\hbox{and}\qquad 
{\bf s}+{\bf u}={\bf k}\;,\eqno(8.48)$$
one has  
$$\eqalignno{ 
\left| {\partial^{j_1+\cdots +j_5+k_1+\cdots +k_5}\over 
\partial x_1^{j_1}\cdots\partial x_5^{j_5} 
\partial y_1^{k_1}\cdots\partial y_5^{k_5}}\,\xi\bigl( {\bf z}\bigr)\right| 
 &\leq 2^J 
\biggl|\left({\partial^{r_1+\cdots +r_5+s_1+\cdots +s_5}\over 
\partial x_1^{r_1}\cdots\partial x_5^{r_5} 
\partial y_1^{s_1}\cdots\partial y_5^{s_5}}\,g\bigl( {\bf z}\bigr)\right)\times \cr 
 &\qquad\qquad\qquad\qquad\times\biggl({\partial^{t_1+\cdots +t_5+u_1+\cdots +u_5}\over 
\partial x_1^{t_1}\cdots\partial x_5^{t_5} 
\partial y_1^{u_1}\cdots\partial y_5^{u_5}}\,
A\bigl( z_1\bigr) B\bigl( z_2\bigr) \prod\limits_{d=3}^5 \Omega\bigl( z_d\bigr)\biggr) 
\biggr| = {}\quad \cr 
 &= 2^J \biggl|\left({\partial^{r_1+\cdots +r_5+s_1+\cdots +s_5}\over 
\partial x_1^{r_1}\cdots\partial x_5^{r_5} 
\partial y_1^{s_1}\cdots\partial y_5^{s_5}}\,g\bigl( {\bf z}\bigr)\right)  
\!\left({\partial^{t_1+u_1}\over 
\partial x_1^{t_1}\partial y_1^{u_1}}\,
A\bigl( z_1\bigr)\right)\times \cr 
 &\qquad\qquad\qquad\qquad\times\left({\partial^{t_2+u_2}\over 
\partial x_2^{t_2}\partial y_2^{u_2}}
\,B\bigl( z_2\bigr)\right)\prod\limits_{d=3}^5 
\left({\partial^{t_d+u_d}\over 
\partial x_d^{t_d}\partial y_d^{u_d}}\,\Omega\bigl( z_d\bigr)\right)\biggr|\;, 
 &(8.49)}$$
where $J=\sum_{d=1}^5 (j_d+k_d)$. 
Given that $\delta\in(0,1]$, it follows by (8.45)-(8.49) that the  
condition (8.39) is satisfied 
at all points 
$(x_1,y_1,\ldots ,x_5,y_5)\in {\cal A}^{\Bbb R}_5 \bigl[ 2^{-1/2} , 1\bigr]$;  since 
we have already seen that the same is true when 
$(x_1,y_1,\ldots ,x_5,y_5)\in{\Bbb R}^{10}-{\cal A}^{\Bbb R}_5 \bigl[ 2^{-1/2} , 1\bigr]$, 
we may therefore conclude that the condition (8.39) is satisfied 
whenever ${\bf x},{\bf y}\in{\Bbb R}^5$ and 
${\bf j},{\bf k}\in\bigl({\Bbb N}\cup\{ 0\}\bigr)^5$. 

We have now shown that the function $\xi({\bf z})$ satisfies all three of the 
conditions stated towards the end of the paragraph containing (8.37). 
As noted in that paragraph, it follows (therefore) that 
the function $\varphi({\bf z})$ and parameters $H,K,L,P,Q$ and $\delta$ 
satisfy 
the initial hypotheses of Section~6 (up to, and including (6.3)): by the remark 
below (6.3), we have in particular $\varphi\in{\cal S}\bigl( {\Bbb C}^5\bigr)$. 
By (8.30), (8.33) and (8.36), the coefficients 
$\theta_p,\phi_h,\Upsilon_{\ell}$ ($p\in{\frak O}-\{ 0\}$, $h,\ell\in{\frak O}$)  
satisfy the hypotheses (6.4) and (6.5) of Section~6. 
Moreover, By (8.22), (8.13) and (8.12), we have $N\ll Q\ll N^{4/3}$ and $H\ll K$; 
and so it follows by (8.25), (8.26) and  (8.15) that 
$$PQ\asymp Q^{*} Q\asymp XN\gg KN\gg HN\;.\eqno(8.50)$$
Since we have also $L\asymp N$ (by (8.31)) and $N=HK$ 
(by hypothesis), 
the parameters $H,K,L,P,Q$ 
therefore satisfy the conditions~(6.8) of 
Lemma~6.1. This, combined with the preceding observations, 
shows that, if $\varepsilon_1>0$, and if $E$ is given by 
$E=(PQ)^{\varepsilon_1}(1+\delta^{-1})\asymp (PQ)^{\varepsilon_1} \delta^{-1}$, then 
the hypotheses of the case $\varepsilon =\varepsilon_1$ of 
Lemma~6.1 are satisfied. 
Consequently, and since (8.34) shows that we have 
$T(Q,X,N)\ll {\cal R}$, where ${\cal R}$ is the sum defined by the equation~(6.1), 
it follows by Lemma~6.3 that, for $\varepsilon_1>0$,  
$$T(Q,X,N)
\ll_{\varepsilon} (PQ)^{11\varepsilon_1} \delta^{-11}   
\left( (NL)^{1/2}+(P/K)Q+(P/K)^{3/2}(QN)^{1/2}\right) 
\left( N\sum_{\ell}\left| \Upsilon_{\ell}\right|^2\right)^{\!\!1/2} .$$
Hence and by (8.31), (8.50), (8.33), (8.27), (8.15) and (8.22), we find that  
$$\eqalign{
T(Q,X,N)
 &\ll_{\varepsilon} (XN)^{11\varepsilon_1}\delta^{-11} 
\left( N+K^{-1}XN+\left( Q^{-1}XH\right)^{3/2} Q^{1/2} N^{1/2}\right) 
N^{1+(\varepsilon /8)} \asymp {}\cr 
 &\asymp (XN)^{11\varepsilon_1} N^{\varepsilon /8} \delta^{-11} 
\left( K^{-1}XN +\left( R^{-1} XH\right)^{3/2} R^{1/2} N^{1/2}\right) N\qquad\qquad\quad\   
\hbox{($\varepsilon_1>0$).}}$$ 
Taking $\varepsilon_1=\varepsilon /22$, we have 
$(XN)^{11\varepsilon_1}N^{\varepsilon /8}=X^{\varepsilon/2} N^{5\varepsilon /8}\leq 
(XN)^{\varepsilon}\,$ (since $X,N\geq 1$).   
Therefore the last bound for $T(Q,X,N)$, combined with (8.21), 
completes our proof of the case $t=0$ of the lemma.
As noted at the start of this proof, the remaining cases of the lemma follow 
\ $\blacksquare$ 

\bigskip 

\goodbreak 
\proclaim Lemma~8.5. Suppose that, with the exception of (8.13), the hypotheses of 
Lemma~8.4 are satisfied. Suppose moreover that 
$$1\leq Q\leq 8 N^{4/3}\;.\eqno(8.51)$$ 
Then 
$$S_t(Q,X,N)\ll_{\varepsilon} 
(QN)^{\varepsilon}\left(\delta^{-1}+|t|\right)^{11}
\left( Q+N+(Q+N)^{1-\vartheta} (HX)^{\vartheta}\right) N\;,\eqno(8.52)$$ 
where $\vartheta$ is the absolute constant given by (1.2.20) and (1.2.21). 
The implicit constant in (8.52) is determined by $\varepsilon$ and 
the matrix $( c_{jk} )_{j,k\geq 0}$, 
where $c_{jk}$ is the implicit constant in the term $O_{j,k}(1)$ in (1.3.15).

\goodbreak 
\noindent{\bf Proof.}\quad 
We consider firstly the cases where 
$$8N^{4/3}\geq Q\geq N\;.\eqno(8.53)$$
In these cases $Q H^{-1}=Q N^{-1} K\geq K\geq 1$. 
Therefore, when (8.53) holds, it follows by Theorem~6 and 
Lemma~8.4 (with $\varepsilon /4$ substituted for $\varepsilon$) 
that 
$$\eqalignno{ 
S_t(Q,X,N) 
 &\leq\max\left\{ 1\,,\,\left( {X\over Q H^{-1}}\right)^{\!\vartheta}\right\} 
S_t\left( Q , Q H^{-1} , N\right) \ll_{\varepsilon} {}\cr 
 &\ll_{\varepsilon} 
\max\left\{ 1\,,\,\left( {X\over Q H^{-1}}\right)^{\!\vartheta}\right\} 
\left( Q H^{-1} N\right)^{\varepsilon /4}
\left(\delta^{-1} +|t|\right)^{11} 
\!\!\left( Q+ \left( {Q H^{-1}\over K}\right)\!N + 1^{3/2} Q^{1/2} N^{1/2}\right)\!N \asymp {}\cr 
 &\asymp \left( 1+{HX\over Q}\right)^{\vartheta}  
(QK)^{\varepsilon /4} \left(\delta^{-1} +|t|\right)^{11}  QN \ll {}\cr 
 &\ll (QN)^{\varepsilon /4} \left(\delta^{-1} +|t|\right)^{11} 
\left( Q + (HX)^{\vartheta} Q^{1-\vartheta}\right) N\;, &(8.54)}$$
so that the bound (8.52) is obtained. 

When (8.53) fails to hold, one has (by (8.51)) $\,1\leq Q<N$; an 
application of Lemma~7.3 then shows that 
$$S_t(Q,X,N)\ll \bigl( S_t(5N/2,X,N)+S_t(5N,X,N)\bigr)\log(5N/2)\;.$$
Moreover, since $N<5N/2,5N<8 N^{4/3}$, we may here apply the 
bound (8.54), with either $5N/2$ or $5N$ substituted for $Q$. 
As a result, we find that if (8.53) does not hold then 
$$\eqalign{ 
S_t(Q,X,N) 
 &\ll O_{\varepsilon}\left( N^{\varepsilon /2}
\left(\delta^{-1} +|t|\right)^{11} 
\left( N + (HX)^{\vartheta} N^{1-\vartheta}\right) N \right) 
\log(5N/2) \ll_{\varepsilon} {}\cr 
 &\ll_{\varepsilon} 
N^{\varepsilon} 
\left(\delta^{-1} +|t|\right)^{11} 
\left( N + (HX)^{\vartheta} N^{1-\vartheta}\right) N\;.}$$
This bound (valid when $Q<N$) implies that in (8.52), and so completes 
this proof\ $\blacksquare$ 

\bigskip 

\goodbreak\proclaim Lemma~8.6. Let $1/4\geq\varepsilon >0$; let $Q_0\geq 1$; and 
let $\vartheta$ be given by (1.2.20) and (1.2.21). Suppose moreover  
that, with the exception of the condition (8.13),  
the hypotheses of Lemma~8.4 are satisfied. 
Let $Q\geq 1$ satisfy either $Q^{1-\varepsilon}\leq N$, or $Q^{1+\varepsilon}\geq XN$, 
or $Q\leq Q_0$. 
Then one has 
$$S_t(Q,X,N)
\leq C_1^{*}\!\left(\varepsilon , Q_0\right) 
(QN)^{\varepsilon} \left(\delta^{-1} +|t|\right)^{11} 
\left( Q+N+(X/K)^{\vartheta} N\right) N\;,\eqno(8.55)$$
where $C_1^{*}\bigl(\varepsilon , Q_0\bigr)$ is a constant, greater than or equal to $1$, 
and depends only upon $\varepsilon$, $Q_0$ and the matrix 
$( c_{jk} )_{j,k\geq 0}$ of constants $c_{jk}$ implicit in the 
term $O_{j,k}(1)$ in the condition (1.3.15).

\goodbreak 
\noindent{\bf Proof.}\quad 
We shall deal firstly with the cases where $Q^{1+\varepsilon}\geq XN$. 
In the proof of Lemma~4.1 (where the hypotheses are more general than 
what is currently supposed) the bound (4.20) is shown to hold for 
$t\in{\Bbb R}$ and $Q,X,N\geq 1$ such that $Q^{1+3\varepsilon}\geq XN$. 
The same holds true if $\varepsilon /3$ is substituted for $\varepsilon$, 
so that if $Q^{1+\varepsilon}\geq XN$ then one has 
$$S_t(Q,X,N)\ll_{\varepsilon} (QN)^{\varepsilon /3}\left( Q+X^{\vartheta} N\right) 
\left\| {\bf a}_N\right\|_2^2\;,$$ 
which, by Theorem~3 and the bound (8.27) for $\left\| {\bf a}_N\right\|_2^2$, 
implies that 
$$S_t(Q,X,N)=O_{\varepsilon}\!\left( (QN)^{\varepsilon /3}
\left( Q+X^{2/9} N\right) N^{1+(\varepsilon /4)} \right)  
\ll_{\varepsilon} (QN)^{\varepsilon}
\left( Q+X^{2/9}N^{1-(5/12)\varepsilon}Q^{-(2/3)\varepsilon}\right) N\;.$$
Given that $0<\varepsilon\leq 1/4<7/3$ and $X,N\geq 1$, we have, in the above, 
$X^{2/9}N^{1-(5/12)\varepsilon}\leq (XN)^{1-(\varepsilon /3)}$. Therefore, and since    
$(1-(\varepsilon /3))(1+\varepsilon)<1+(2/3)\varepsilon$ 
(as $\varepsilon$ is positive), 
we find that 
$$S_t(Q,X,N)\ll_{\varepsilon} 
(QN)^{\varepsilon} QN\qquad\quad{\rm if}\quad\ Q^{1+\varepsilon}\geq XN\;.\eqno(8.56)$$

The above supplies all that we need in respect of the cases where 
$Q^{1-\varepsilon}\geq XN$. We now need only to obtain suitable 
bounds for $S_t(Q,X,N)$ in the cases where 
$Q\leq\max\left\{ N^{1/(1-\varepsilon)} , Q_0\right\}$. 
In each such case one has either 
$Q^{1-\varepsilon}> N$ and $Q\leq Q_0$, or else $Q^{1-\varepsilon}\leq N$. 
We shall consider, in turn, these two possibilities. 

If $Q^{1-\varepsilon}> N$ and $Q\leq Q_0$ then, since 
$\varepsilon >0$ and $N,Q\geq 1$, one has 
$N<Q^{1-\varepsilon}\leq Q\leq Q_0$. Moreover, given Theorem~3, 
one has $N^{\vartheta}<Q_0^{2/9}$ and $N^{\varepsilon /4}<Q_0^{\varepsilon /4}$ 
when $N<Q_0$; and so,  from  (8.27) and 
the bound (4.21)$\,$ (obtained within the proof of Lemma~4.1), 
one may deduce that  
$$S_t(Q,X,N)
\ll_{\varepsilon,Q_0} (QN)^{\varepsilon}
\left( Q+(X/N)^{\vartheta} N\right) N\qquad\quad{\rm if}\quad\    
N^{1/(1-\varepsilon)}<Q\leq Q_0\;.\eqno(8.57)$$  
In comparison to the implicit constant in (4.21), the 
implicit constant in (8.57) potentially accommodates 
an extra factor $D(\varepsilon)Q_0^{(2/9)+(\varepsilon /4)}$, where 
$D(\varepsilon)$, which is the implicit constant in (8.27), 
is determined by $\varepsilon$ and the relevant implicit 
constant in respect of the  case $j=k=0$ of the condition (1.3.15). 

If $Q^{1-\varepsilon}\leq N$ then, since $0<\varepsilon\leq 1/4$, 
one has $Q\leq N^{1/(1-\varepsilon)}\leq N^{4/3}$, 
so that the condition (8.51) in Lemma~8.5 is satisfied (as are 
all the other hypotheses of that lemma). Lemma~8.5 is valid for 
arbitrary $\varepsilon >0$. Therefore, by applying Lemma~8.5 with 
$\varepsilon^2\in(0,1/16]$ substituted for 
$\varepsilon$, we find that if $Q^{1-\varepsilon}\leq N$ then   
$$S_t(Q,X,N)\ll_{\varepsilon} 
(QN)^{(\varepsilon^2)}\left(\delta^{-1}+|t|\right)^{11}
\left( Q+N+Q^{1-\vartheta} (HX)^{\vartheta}+N^{1-\vartheta} (HX)^{\vartheta}\right) N\;.\eqno(8.58)$$
If $Q^{1-\varepsilon}\leq N$ then one  has  also  
$Q\leq N(QN)^{\varepsilon /(2-\varepsilon)}$,  
so that 
$$(QN)^{(\varepsilon^2)} Q^{\!{1-\vartheta}}
\leq (QN)^{(\varepsilon^2)}
\left( N (QN)^{\varepsilon /(2-\varepsilon)}\right)^{\!{1-\vartheta}} 
=(QN)^{\varepsilon F(\varepsilon)} N^{1-\vartheta}\;,$$ 
where 

$$F(\varepsilon)={1-\vartheta\over 2-\varepsilon}+\varepsilon 
\leq {1\over 2-\varepsilon}+{\varepsilon}<1$$
(with the last two inequalities following since $\vartheta\geq 0$ 
and $0<\varepsilon\leq 1/4$).    
Since the hypotheses of the lemma imply that $(QN)^{(\varepsilon^2)}\leq (QN)^{\varepsilon}$, 
it may therefore be deduced from 
the conditional bound (8.58) that  
$$S_t(Q,X,N) 
\ll_{\varepsilon} 
(QN)^{\varepsilon}\left(\delta^{-1}+|t|\right)^{11}
\left( Q+N+N^{1-\vartheta} (HX)^{\vartheta}\right) N\qquad\quad{\rm if}\quad\  
Q^{1-\varepsilon}\leq N\;.\eqno(8.59)$$ 

Let $A(\varepsilon),B(\varepsilon , Q_0),C(\varepsilon)>0$ be  
sufficiently large to serve as the 
implicit constants in 
(8.56), (8.57) and (8.59), respectively. Then, given that $0<\delta\leq 1$ and 
$N=HK\geq K\geq 1$ and $\vartheta\geq 0$, it follows by those conditional results, (8.56), (8.57) and (8.59), 
that if either $1\leq Q\leq\max\{ N^{1/(1-\varepsilon)}\,,\,Q_0\}$, or $Q^{1+\varepsilon}\geq XN$, 
then the bound (8.55) will hold with  
$ C_1^{*}(\varepsilon ,Q_0)
=\max\{ 1 , A(\varepsilon) , B(\varepsilon , Q_0) , C(\varepsilon)\}$ 
\ $\blacksquare$ 

\bigskip 
\bigskip 

\goodbreak 
\centerline{\it The Proof of Theorem~9.} 

\bigskip 

It will suffice to obtain the bound (1.3.16) in cases where $0<\varepsilon\leq 1/4$.  
Indeed, since $(QN)^{1/4}\leq (QN)^{\varepsilon}$ 
when $Q,N\geq 1$ and $\varepsilon >1/4$, 
all relevant cases where $\varepsilon >1/4$ follow 
from cases in which one has $\varepsilon =1/4$. 
Therefore we assume henceforth that 
$$0<\varepsilon\leq 1/4\;.$$
We shall assume also that 
$$K\geq H\;.$$
This latter assumption is justified, since the case $K<H$ of Theorem~9 
follows from the case $K\geq H$ 
by the substitution of $\beta$, $K$, $\alpha$ and $H$ for 
$\alpha$, $H$, $\beta$ and $K$, respectively. 

Given the hypotheses of Theorem~9, it follows (similarly to (8.27)) that 
$$\left\|{\bf a}_N\right\|_2^2\leq C_{\infty}^{*}(\varepsilon) N^{1+\varepsilon}\;,\eqno(8.60)$$
where $C_{\infty}^{*}=C_{\infty}^{*}(\varepsilon)\in[1,\infty)$ is 
a constant depending only upon $\varepsilon$ 
and the implicit constant in the case $j=k=0$ of the conditions (1.3.15). 
Taking now $C_{13}=C_{13}(\varepsilon)\in[1,\infty)$ to be one of those 
numbers whose existence is established by the case $j=13$ of Lemma~4.2, 
we put 
$$Q_0=Q_0^{*}(\varepsilon)=\left( 2^{15} C_{13}(\varepsilon)\right)^{(3\varepsilon^{-2})} 
\qquad{\rm and}\qquad 
C_0^{*}=C_0^{*}(\varepsilon)=\max\left\{ C_1^{*}\!\left(\varepsilon,Q_0\right)\,,
\,2^{-14} C_{\infty}^{*}(\varepsilon)\right\}\;,\eqno(8.61)$$
where $C_1^{*}(\varepsilon , Q_0)=C_1^{*}\in[1,\infty)$ is any one of those 
numbers whose existence is established in Lemma~8.6 
(note that we are certain to have here $Q_0\geq 2^{720}>1$). 
The functions $\alpha,\beta : {\Bbb C}\rightarrow {\Bbb C}$ 
(and hence also the associated parameters $H,K,N,\delta$) will remain 
fixed throughout this proof, as does $\varepsilon$.
The same is therefore true of all the numbers 
$C_{\infty}^{*}=C_{\infty}^{*}(\varepsilon)$, 
$C_{13}=C_{13}(\varepsilon)$, $Q_0=Q_0^{*}(\varepsilon)$, 
$C_1^{*}=C_1^{*}(\varepsilon , Q_0)$ and $C_0^{*}=C_0^{*}(\varepsilon)$, 
just described. 

For each $Q\in [1,\infty)$, let $A^{*}(Q)$ denote the proposition that, 
for all $X\geq 1$ and all $t\in{\Bbb R}$, one has 
$$S_t(Q,X,N)\leq C_0^{*}(\varepsilon) (QN)^{\varepsilon}
\left( Q+N+X^{\vartheta}N^{\vartheta}Q^{1-2\vartheta}
+X^{\vartheta}K^{-\vartheta} N\right)\left(\delta^{-1}+|t|\right)^{11} N\;.\eqno(8.62)$$
Given that we have $K\geq H$ and $0\leq\vartheta\leq 2/9$ (by Theorem~3), 
the inequality (8.62), if true, would imply the result (1.3.16) of Theorem~9. 
Therefore, in order to complete this proof of Theorem~9, it will 
suffice that we show that $A^{*}(Q)$ is true for all $Q\in[1,\infty)$.
Since the equalities noted in (4.30) (within the proof of Theorem~5) 
remain valid in the current context, we have, moreover: 
$$A^{*}(Q)\ \,{\rm implies}\ \,A^{*}(P)\qquad\ {\rm if}\quad\ Q+1>P\geq Q\in{\Bbb N}\;.$$
Consequently we may complete this proof of Theorem~9 simply by 
showing that $A^{*}(Q)$ is true for all $Q\in{\Bbb N}$: 
this we shall achieve through a `proof by contradiction'. 

Suppose that $A^{*}(Q)$ is false for some $Q\in{\Bbb N}$. 
Then, for reasons similar to those which justify (4.32) and (4.33) 
in the proof of Theorem~5, there must exist a unique $R\in{\Bbb N}$ 
such that 
$$A^{*}(R)\ \ {\rm is\ \,false}\eqno(8.63)$$
and 
$$A^{*}(Q)\ \ {\rm is\ \,true}\qquad{\rm for}\quad Q\in[1,R)\;.\eqno(8.64)$$
Lemma~8.6 implies that $A^{*}(Q)$ is true for all 
real $Q\geq 1$ satisfying either $Q\leq Q_0^{*}(\varepsilon)$, 
or $Q^{1-\varepsilon}\leq N$: 
for the definition (8.61) ensures that 
$C_0^{*}(\varepsilon)\geq C_1^{*}(\varepsilon ,Q_0)$, and so, in all the relevant cases, 
the result (8.55) of Lemma~8.6 implies  
the inequality (8.62). Hence, given (8.63), we must have:  
$$R>Q_0^{*}(\varepsilon)=Q_0>1\eqno(8.65)$$
and 
$$R^{1-\varepsilon}>N\;.\eqno(8.66)$$
By Lemma~8.6 (again), we moreover have 
$$S_t(R,X,N)\leq C_1^{*}\!\left(\varepsilon , Q_0\right) 
(RN)^{\varepsilon}\!\left( R+N+X^{\vartheta} K^{-\vartheta} N\right) 
\left(\delta^{-1}+|t|\right)^{11}\!N\quad 
\hbox{for $t\in{\Bbb R}$, $X\in\left[ 1 , R^{1+\varepsilon}/N\right]$.}\eqno(8.67)$$

The steps we shall now take in order to complete this proof 
are similar to those taken (after (4.34)) in completing the proof of Theorem~5. 
We shall deduce from (8.60), (8.61) and (8.64)-(8.67) that 
the proposition $A^{*}(R)$ is true: since that conclusion will directly 
contradict (8.63), we shall thereby have given a 
`proof by contradiction' that $A^{*}(Q)$ is true for all $Q\in{\Bbb N}$, 
and so (given the points noted below (8.62)) 
shall have completed the proof of Theorem~9. 

By (8.61) and (8.67), we obtain the bound (8.62) for $Q=R$, all $t\in{\Bbb R}$ 
and all $X\in[1,R^{1+\varepsilon}/N]$. Therefore, 
if it can be shown that 
(8.62) holds for $Q=R$, all $t\in{\Bbb R}$ and all $X>R^{1+\varepsilon}/N$, 
then we may  
deduce that $A^{*}(R)$ is true. Accordingly, let us suppose that 
$$t\in{\Bbb R}\qquad\ {\rm and}\qquad\ X>R^{1+\varepsilon}/N\eqno(8.68)$$
(by (8.66) this ensures that $X>R^{2\varepsilon}\geq 1$). 
Then, given the observations immediately preceding  our 
supposition of (8.68), the proposition $A^{*}(R)$ is true if it can 
now be deduced that 
$$S_t(R,X,N)
\leq C_0^{*}(\varepsilon) (RN)^{\varepsilon}
\left( R+N+X^{\vartheta}N^{\vartheta}R^{1-2\vartheta}
+X^{\vartheta}K^{-\vartheta} N\right)\left(\delta^{-1}+|t|\right)^{11}\!N\;.\eqno(8.69)$$

By (8.68) and (8.66), the case $Q=R$ of the condition (4.22) in Lemma~4.2 
is satisfied. Therefore, by applying Lemma~4.2 for $Q=R$ and $j=13$,  
we find that, for 
$Y=\min\{ X\,,\,R^{2-\varepsilon}/N\}$, some 
$v\in{\Bbb R}$ and some $L\in(YN/R , 2^{10}YN/R)$,  
one has 
$$S_t(R,X,N)
\leq C_{13}(\varepsilon)\left( {X\over Y}\right)^{\vartheta} 
\left( (1+|v-t|)^{-11} S_v(L,Y,N)
+R^{1+(2-\varepsilon)\varepsilon /3}\left\|{\bf a}_N\right\|_2^2\right) .\eqno(8.70)$$
Here, by (8.66), (8.68) and the same calculations as carried out below (4.36) 
(in the proof of Theorem~5), it follows that 
$$Y\geq 1\qquad{\rm and}\qquad 
1<L<2^{10}YN/R\leq 2^{10}R^{1-\varepsilon}\;.\eqno(8.71)$$
Moreover, given that $0<\varepsilon\leq 1/4$ 
and $C_{13}(\varepsilon)\geq 1$, it is implied by (8.65) and (8.61) that 
$$R^{\varepsilon}
>\left( 2^{15} C_{13}(\varepsilon)\right)^{(3\varepsilon^{-1})} 
\geq 2^{180}>2^{10}\;.\eqno(8.72)$$

The inequalities in (8.72) and the second part of (8.71) imply that 
$1<L<R$. It therefore follows, by (8.64), that the proposition $A^{*}(L)$ 
is true. In particular, a valid inequality is obtained in (8.62) 
when we there substitute $L$, $Y\in[1,\infty)$ and $v\in{\Bbb R}$ 
for $Q$, $X$ and $t$, respectively. It is therefore the case that 
$$S_v(L,Y,N)
\leq C_0^{*}(\varepsilon) 
(LN)^{\varepsilon}  
\left( L+N+Y^{\vartheta}N^{\vartheta}L^{1-2\vartheta}
+Y^{\vartheta}K^{-\vartheta} N\right)\left(\delta^{-1}+|v|\right)^{11}\!N\;.\eqno(8.73)$$
Since $0<\varepsilon\leq 1/4$ and $0<\delta\leq 1$, it follows by 
(8.66), (8.71) and Theorem~3 that we have, in the above, 
$$\eqalign{
L+N+Y^{\vartheta}N^{\vartheta}L^{1-2\vartheta}
 &\leq 2^{10}R^{1-\varepsilon}+R^{1-\varepsilon}
+R^{(2-\varepsilon)\vartheta}\left( 2^{10}R^{1-\varepsilon}\right)^{1-2\vartheta} = {}\cr 
 &=\left( 2^{10}+1\right) R^{1-\varepsilon} +2^{10(1-2\vartheta)}R^{1-(1-\vartheta)\varepsilon} \leq {}\cr 
 &\leq\left( 2^{11}+1\right) R^{1-(7/9)\varepsilon}<2^{23/2} R^{1-\varepsilon /3}\;,}$$
$$(LN)^{\varepsilon}\leq \left( 2^{10} R^{1-\varepsilon} N\right)^{\varepsilon} 
<2^{5/2} \left( R^{1-\varepsilon /3} N\right)^{\varepsilon}$$ 
and 
$$\delta^{-1}+|v|=\delta^{-1}+|t+(v-t)| 
\leq \delta^{-1}+|t|+|v-t|
\leq \left( \delta^{-1}+|t|\right)\left( 1+|v-t|\right)\;.$$
Therefore the bound (8.73) certainly implies that 
$$\left( 1+|v-t|\right)^{-11} S_v(L,Y,N) 
<2^{14} C_0^{*}(\varepsilon) \left( R^{1-\varepsilon /3} N\right)^{\varepsilon}
\left( R^{1-\varepsilon /3}
+Y^{\vartheta}K^{-\vartheta} N\right)\left(\delta^{-1}+|t|\right)^{11}\!N\;.\eqno(8.74)$$
Moreover, by (8.60) and (8.61), we have 
$$R^{1+(2-\varepsilon)\varepsilon /3}\left\|{\bf a}_N\right\|_2^2
\leq R^{1+(2-\varepsilon)\varepsilon /3}
C_{\infty}^{*}(\varepsilon) N^{1+\varepsilon} 
=C_{\infty}^{*}(\varepsilon) \left( R^{1-\varepsilon /3} N\right)^{\varepsilon +1}
\leq 2^{14}C_0^{*}(\varepsilon)
\left( R^{1-\varepsilon /3} N\right)^{\varepsilon +1} .\eqno(8.75)$$ 

Since $\delta^{-1}\geq 1$, it follows by (8.70), (8.74) and (8.75) 
that 
$$\eqalign{ 
S_t(R,X,N) 
 &\leq 2^{15} C_{13}(\varepsilon) C_0^{*}(\varepsilon)\left( {X\over Y}\right)^{\vartheta} 
\left( R^{1-\varepsilon /3} N\right)^{\varepsilon}
\left( R^{1-\varepsilon /3}
+Y^{\vartheta}K^{-\vartheta} N\right)\left(\delta^{-1}+|t|\right)^{11}\!N = {}\cr 
 &=2^{15} C_{13}(\varepsilon) R^{-\varepsilon^2 /3} C_0^{*}(\varepsilon) 
(RN)^{\varepsilon}
\left( \left( {X\over Y}\right)^{\vartheta} R^{1-\varepsilon /3}
+X^{\vartheta}K^{-\vartheta} N\right)\left(\delta^{-1}+|t|\right)^{11}\!N\;,
}$$
where, just as at the end of the proof of Theorem~5, one has: 
$$\left( {X\over Y}\right)^{\vartheta} R^{1-\varepsilon /3}
\leq\max\left\{ R^{1-\varepsilon /3}\,,
\,\left( {XN\over R^{2-\varepsilon}}\right)^{\vartheta} R^{1-\varepsilon /3}\right\} 
\leq R+X^{\vartheta} N^{\vartheta} R^{1-2\vartheta} 
<R+N+X^{\vartheta} N^{\vartheta} R^{1-2\vartheta}$$
(with the penulutimate inequality following by Theorem~3). 
By (8.72) we have, in the above,  
$$2^{15} C_{13}(\varepsilon) R^{-\varepsilon^2 /3}<1\;.$$
It may therefore be deduced that the inequality in (8.69) 
holds: this (as observed immediately above (8.69)) is sufficient 
to establish that the proposition $A^{*}(R)$ is true, so that 
the statement (8.63) is contradicted. Consequently, as explained 
in the paragraph below (8.67), the proof of Theorem~9 is now 
complete\ $\blacksquare$ 

\bigskip 
\bigskip 

\goodbreak\centerline{\bf \S 9. The Proofs of Theorems~10~and~11}

\bigskip 

In this section we prove first Theorem~11, and then Theorem~10. 
We begin with two lemmas required in the first of these proofs. 

\bigskip 

\goodbreak\proclaim Lemma~9.1. Let $q,r,s\in{\frak O}-\{ 0\}$ be such that 
$q=rs$ and $(r,s)\sim 1$; let $\Gamma =\Gamma_0(q)\leq SL(2,{\frak O})$;  
let $g_{\infty}\in SL(2,{\Bbb C})$ be as in (1.3.3); and let 
$u,t\in{\frak O}$ and $g_{1/s}\in SL(2,{\Bbb C})$ be such that 
the equations~(1.4.1) and~(1.4.2) hold. Then 
$$\eqalignno{
g_{1/s}^{-1}\Gamma g_{\infty}
 &=\left\{\pmatrix{A\sqrt{r} &B/\sqrt{r}\cr Cs\sqrt{r} &D\sqrt{r}} : 
A,B,C,D\in{\frak O}\ {\rm and}\ ADr=1+BCs\right\} , &(9.1)\cr 
g_{1/s}^{-1}\Gamma g_{1/s} 
 &=\left\{\pmatrix{A &B\cr Cq &D} : 
A,B,C,D\in{\frak O}\ {\rm and}\ AD=1+BCq\right\} 
=\Gamma =g_{\infty}^{-1}\Gamma g_{\infty} , &(9.2)}$$ 
and the condition (1.1.1) is satisfied when ${\frak c}\in\{\infty , 1/s\}$. 

\goodbreak 
\noindent{\bf Proof.}\quad 
Let $G(r,s)$ denote the set on the right-hand side of the equation~(9.1). 
Since $\Gamma=\Gamma_0(q)$, and since 
(1.3.3) makes $g_{\infty}$ the identity element of $SL(2,{\Bbb C})$, 
it follows that 
the final two equalities in (9.2) are trivial consequences of 
the definition of the Hecke congruence subgroup $\Gamma_0(q)$, 
and that $g_{1/s}^{-1}\Gamma g_{\infty}=g_{1/s}^{-1}\Gamma$.
Therefore, in order to complete the proofs of (9.1) and (9.2), 
it will suffice to show that 
$$g_{1/s}^{-1}\Gamma =G(r,s)\qquad\ {\rm and}\qquad\  
G(r,s) g_{1/s}=\Gamma\;.\eqno(9.3)$$

We begin with a proof of the first equation in~(9.3). Suppose, firstly, that 
$\gamma\in\Gamma$. 
Then, for some $a,b,c,d\in{\frak O}$ satisfying 
$ab-cd=1$ and $q\mid c$, one has 
$$\pmatrix{a &b\cr c &d}=\gamma\in\Gamma\leq SL(2,{\Bbb C})\;.\eqno(9.4)$$
By (1.4.1)-(1.4.2), the matrix $g_{1/s}$ is an element of 
the group $SL(2,{\Bbb C})$, and so it follows from (9.4) that 
$$g_{1/s}^{-1}\gamma\in SL(2,{\Bbb C})\;.\eqno(9.5)$$
Moreover, by (1.4.1)-(1.4.2) and (9.4) (again), 
$$\eqalignno{
g_{1/s}^{-1}\gamma 
 &=\pmatrix{u\sqrt{r} &-t/\sqrt{r}\cr -s\sqrt{r} &\sqrt{r}}
\pmatrix{a &b\cr c &d}  = {}\cr
 &=\pmatrix{(ua-tc/r)\sqrt{r} &(urb-td)/\sqrt{r}\cr (-a+c/s)s\sqrt{r} &(-sb+d)\sqrt{r}} 
=\pmatrix{A\sqrt{r} &B/\sqrt{r}\cr Cs\sqrt{r} &D\sqrt{r}}\qquad\quad\ \hbox{(say),} &(9.6)}$$
where, since $c\in q{\frak O}=rs{\frak O}$, one has $A,B,C,D\in{\frak O}$. 
By (9.5) and (9.6) the determinant of the last matrix is equal to $1$, so that 
$ADr=1+BCs$. It has therefore been shown that 
$$g_{1/s}^{-1}\Gamma\subseteq G(r,s)\;.\eqno(9.7)$$

Suppose now that 
$$h=\pmatrix{A'\sqrt{r} &B'/\sqrt{r}\cr C's\sqrt{r} & D'\sqrt{r}}\in G(r,s)\;.\eqno(9.8)$$
Then, since $g_{1/s}\in SL(2,{\Bbb C})$, and since 
$G(r,s)$ is a subset of the set of elements of $SL(2,{\Bbb C})$, 
one has 
$$g_{1/s}h\in SL(2,{\Bbb C})\;.$$
Hence, and by (1.4.1), 
$$g_{1/s}h
=\pmatrix{\sqrt{r} &t/\sqrt{r}\cr s\sqrt{r} &u\sqrt{r}} 
\pmatrix{A'\sqrt{r} &B'/\sqrt{r}\cr C's\sqrt{r} & D'\sqrt{r}} 
=\pmatrix{A'r+C'st &B'+D't\cr (A'+C'u)rs &B's+D'ru} 
\in\Gamma_0(q)=\Gamma$$
(for we have $r,s,t,u\in{\frak O}$, $rs=q$ and, by (9.8), $A',B',C',D'\in{\frak O}$).  
Since the relation $g_{1/s}h\in\Gamma$ implies that $h\in g_{1/s}^{-1}\Gamma$, the 
above therefore shows that $G(r,s)\subseteq g_{1/s}^{-1}\Gamma$. 
This, together with (9.7), completes the proof of the first equation in~(9.3). 

We may employ a similar strategy to prove next the second equation in~(9.3).
Suppose that $\gamma$ is as in~(9.4). Then, 
since $SL(2,{\Bbb C})\ni\gamma,g_{1/s}$, 
one has $\gamma g_{1/s}^{-1}\in SL(2,{\Bbb C})$.
Since it moreover follows from (9.4) and (1.4.1)-(1.4.2) that 
$$\gamma g_{1/s}^{-1}
=\pmatrix{a &b\cr c &d}\pmatrix{u\sqrt{r} &-t/\sqrt{r}\cr -s\sqrt{r} &\sqrt{r}} 
=\pmatrix{(au-bs)\sqrt{r} &(-at+br)/\sqrt{r}\cr ((c/s)u-d)s\sqrt{r} &(-(c/r)t+d)\sqrt{r}}$$
(where $c/s,c/r\in{\frak O}$, since $q\mid c$), we have therefore that 
$\gamma g_{1/s}^{-1}\in G(r,s)$. This proves that 
$\Gamma g_{1/s}^{-1}\subseteq G(r,s)$, and so enables us to deduce that 
$$\Gamma\subseteq G(r,s) g_{1/s}\;.\eqno(9.9)$$
On the other hand, for $h$ as in (9.8), one has 
$h g_{1/s}\in SL(2,{\Bbb C})$ (since $h,g_{1/s}\in SL(2,{\Bbb C})$) 
and so 
$$h g_{1/s}
=\pmatrix{A'\sqrt{r} &B'/\sqrt{r}\cr C's\sqrt{r} & D'\sqrt{r}}
\pmatrix{\sqrt{r} &t/\sqrt{r}\cr s\sqrt{r} &u\sqrt{r}} 
=\pmatrix{A'r+B's &A't+B'u\cr (C'+D')rs &C'st+D'ur}\in\Gamma_0(q)=\Gamma\;.$$
Therefore we have that $G(r,s)g_{1/s}\subseteq\Gamma$. This, together  
with (9.9), proves the second equation in~(9.3).

In order to complete the proof of the lemma we must show that 
(1.1.1) holds for all ${\frak c}\in\{\infty , 1/s\}$. 
Let ${\frak c}\in\{\infty , 1/s\}$. Then, by either (1.1.3) or (1.4.1)-(1.4.2) 
(whichever is appopriate), the scaling matrix $g_{\frak c}$ is an element of 
$SL(2,{\Bbb C})$ satisfying $g_{\frak c}\infty ={\frak c}$.  
Indeed, if ${\frak c}=\infty$ then, by (1.3.3), $g_{\infty}[1,0]=[1,0]\in{\Bbb P}^1({\Bbb C})$, 
while if instead ${\frak c}=1/s$ then, by (1.4.1)-(1.4.2), 
$g_{1/s}[1,0]=[\sqrt{r},s\sqrt{r}]=[1,s]\in{\Bbb P}^1({\Bbb C})$. 
It follows that when $\gamma\in\Gamma$ one has $\gamma{\frak c}={\frak c}$ 
if and only if $g_{\frak c}^{-1}\gamma g_{\frak c}\infty =\infty$, for the 
latter equation is equivalent to the equation $\gamma g_{\frak c}\infty =g_{\frak c}\infty$, 
and, as we have just seen, $g_{\frak c}\infty ={\frak c}$. 
Therefore, and since the trace of any $\gamma\in\Gamma$ is invariant under 
conjugation by an element of $SL(2,{\Bbb C})$, it follows from the definitions 
of $\Gamma_{\frak c}$ and $\Gamma_{\frak c}'$ preceding (1.1.1) that 
$$g_{\frak c}^{-1}\Gamma_{\frak c}' g_{\frak c} 
=\left\{ g\in g_{\frak c}^{-1}\Gamma g_{\frak c} : g\infty =\infty\ {\rm and}\ 
{\rm Tr}(g)=2\right\} .\eqno(9.10)$$
Recall now that, if  
$$g=\pmatrix{a &b\cr c &d}\in SL(2,{\Bbb C})\;,$$
then one has $g\infty =\infty$ if and only if $c=0$. Given this fact, it follows  
from (9.10) and (9.2)  that one has 
$$\eqalign{
g_{\frak c}^{-1}\Gamma_{\frak c}' g_{\frak c} 
 &=\left\{ g\in\Gamma : g\infty =\infty\ {\rm and}\ 
{\rm Tr}(g)=2\right\} = {}\cr
 &=\left\{\pmatrix{A &B\cr 0 &D} : A,B,D\in{\frak O},\ AD=1\ {\rm and}\ A+D=2\right\} 
=\left\{\pmatrix{1 &B\cr 0 &1} : B\in{\frak O}\right\} ,}$$
which is the required result (1.1.1)\quad$\blacksquare$ 

\bigskip 

\goodbreak 
\noindent{\bf Remark.}\quad The calculation below (9.9) (in the same paragraph), 
is somewhat superfluous. Indeed, by (9.9) and the first equation in~(9.3), one 
has $\Gamma\leq g_{1/s}^{-1}\Gamma g_{1/s}$; and it is not possible 
that $\Gamma$ be a proper subgroup of $g_{1/s}^{-1}\Gamma g_{1/s}$, for 
the covolumes of these two discrete and cofinite subgroups of 
$G=SL(2,{\Bbb C})$ are equal. 

\bigskip 

\goodbreak 
\proclaim Lemma~9.2. Let the hypotheses of Lemma~9.1 be satisfied. 
Let $\omega,\omega'\in{\frak O}$. For ${\frak a},{\frak b}\in\{\infty , 1/s\}$, 
put 
$$\delta_{\omega,\omega'}^{{\frak a},{\frak b}} 
=\sum_{\scriptstyle\Gamma_{\frak a}'\gamma\in\Gamma_{\frak a}'\backslash\Gamma\ 
:\ \gamma{\frak b}={\frak a}\atop\scriptstyle 
g_{\frak a}^{-1}\gamma g_{\frak b}=\pmatrix{u(\gamma) &\beta(\gamma)\cr 0 &1/u(\gamma)}} 
{\rm e}\left({\rm Re}\left(\beta(\gamma)u(\gamma)\omega\right)\right)  
\delta_{u(\gamma)\omega , \omega'/u(\gamma)}\;,\eqno(9.11)$$
where $\delta_{w,z}$ equals $1$ if $w=z$, and is otherwise zero; 
and let ${}^{\frak a}{\cal C}^{\frak b}$ and the generalised Kloosterman sums 
$S_{{\frak a},{\frak b}}\left(\omega,\omega';c\right)\ $ ($c\in{}^{\frak a}{\cal C}^{\frak b}$) 
be given by (1.1.13)-(1.1.15). 
Then one has what is stated in (1.4.3), (1.4.4), (1.4.16) and (1.4.17); 
and it is moreover the case that 
$$\delta_{\omega,\omega'}^{1/s,1/s}
=\delta_{\omega,\omega'}^{\infty,\infty} 
=\sum_{u\in{\frak O}^{*}}\delta_{u\omega , \omega'/u} 
=\cases{4 &if $\omega' =\omega = 0$,\cr 2 &if $\omega'=\pm\omega\neq 0$, \cr 
0 &otherwise.}\eqno(9.12)$$ 

\goodbreak 
\noindent{\bf Proof.}\quad 
By Lemma~9.1, we have (9.1), (9.2) and (1.1.1) for ${\frak c}\in\{\infty , 1/s\}$, 
which implies that 
$$\Gamma_{\frak c}'=g_{\frak c} B^{+} g_{\frak c}^{-1}\qquad\qquad\quad  
\hbox{(${\frak c}\in\{\infty , 1/s\}$),}\eqno(9.13)$$
where, recalling the notation of Subsection~1.1, one has $B^{+}=\{ n[\alpha] : \alpha\in{\frak O}\}$. 
Given the definitions (1.1.13) and (1.1.14), the result (1.4.3) 
concerning ${}^{1/s}{\cal C}^{\infty}$ will follow if it can be shown that, 
when $c\in{\Bbb C}-\{ 0\}$, the set $g_{1/s}^{-1}\Gamma g_{\infty}$ 
contains an element of the form 
$$\pmatrix{* &*\cr c &*}\eqno(9.14)$$
if and only if 
$$c/s\sqrt{r}\in{\frak O}\qquad\quad{\rm and}\qquad\quad\left( c/s\sqrt{r} , r\right)\sim 1\;.\eqno(9.15)$$
Accordingly, let $c\in{\Bbb C}-\{ 0\}$. By (9.1), the set 
$g_{1/s}^{-1}\Gamma g_{\infty}$ contains an element of the form (9.14) 
if and only if 
$c=Cs\sqrt{r}$ for some $C\in{\frak O}$ such that 
the congruence $rX\equiv 1\bmod Cs{\frak O}$ has a solution in ${\frak O}$. 
Moreover, since ${\Bbb Z}[i]$ is a principal ideal domain, and since 
$(r,s)\sim 1$, the conguence in question is soluble if and only if 
the Gaussian integer $C=c/s\sqrt{r}$ is coprime to $r$. 
Hence, in (9.15) we have necessary and sufficient conditions 
for $g_{1/s}^{-1}\Gamma g_{\infty}$ to contain at least one element of the 
form (9.14). This completes the proof of (1.4.3). 

The result (1.4.16) follows similarly (but even more easily) 
from (9.2): we omit the relevant details. 

As a first step towards the proof of (1.4.4) and (1.4.17), we observe that, 
as a consequence of (9.13) and the definitions (1.1.13)-(1.1.15), one has 
$$\eqalignno{
S_{{\frak a},{\frak b}}\left(\omega,\omega';c\right)  
 &=\sum_{B^{+}\pmatrix{a &*\cr c &d}B^{+}\in 
B^{+}\backslash g_{\frak a}^{-1}{}^{\frak a}\Gamma^{\frak b}g_{\frak b}/B^{+}} 
{\rm e}\!\left({\rm Re}\left(\omega\,{a\over c}+\omega'\,{d\over c}\right)\right) = {}\cr 
 &=\sum_{d+c{\frak O}\in{\Bbb C}/(c{\frak O})}\!\!\sum_{\scriptstyle a+c{\frak O}\in{\Bbb C}/(c{\frak O})\atop\scriptstyle 
\pmatrix{a &*\cr c &d}\in g_{\frak a}^{-1}\Gamma g_{\frak b}} 
{\rm e}\!\left({\rm Re}\left({\omega a +\omega' d\over c}\right)\right) &(9.16)}$$
whenever $({\frak a} , {\frak b})\in\{\infty , 1/s\}\times\{\infty , 1/s\}$ and 
$c\in{}^{\frak a}{\cal C}^{\frak b}$. Note that the final sum in (9.16) 
is completely determined by $\omega$, $\omega'$, $c$ and the set 
$g_{\frak a}^{-1}\Gamma g_{\frak b}\subset SL(2,{\Bbb C})$. 
Hence, and since we have $g_{1/s}^{-1}\Gamma g_{1/s}=g_{\infty}^{-1}\Gamma g_{\infty}$, 
by (9.2), and ${}^{1/s}{\cal C}^{1/s}={}^{\infty}{\cal C}^{\infty}=q{\frak O}-\{ 0\}$ 
(by the result (1.4.16), which follows from (9.2)), it is therefore the case that 
$$S_{{1/s},{1/s}}\left(\omega,\omega';Cq\right) 
=S_{{\infty},{\infty}}\left(\omega,\omega';Cq\right)\qquad\quad  
\hbox{for $\quad C\in{\frak O}-\{ 0\}$.}\eqno(9.17)$$
Moreover, by substituting into the case ${\frak a}={\frak b}=\infty$ of (9.16) 
the explicit description of $g_{\infty}^{-1}\Gamma g_{\infty}=\Gamma$ given in~(9.2),
we find that, for $0\neq C\in{\frak O}$, 
$$\eqalign{ 
S_{\infty,\infty}\left(\omega,\omega';Cq\right) 
 &=\sum_{D+Cq{\frak O}\in{\frak O}/(Cq{\frak O})}
\ \sum_{\scriptstyle A+Cq{\frak O}\in{\frak O}/(Cq{\frak O})\atop\scriptstyle 
\pmatrix{A &*\cr Cq &D}\in\Gamma_0(q)} 
{\rm e}\!\left({\rm Re}\left({\omega A +\omega' D\over Cq}\right)\right) = {}\cr 
 &=\qquad\sum\!\!\!\!\!\sum_{\!\!\!\!\!\!\!\!\!\!\!\!\!\!\!{\scriptstyle A,D\bmod Cq{\frak O}\atop\scriptstyle 
AD\equiv 1\bmod Cq{\frak O}}} 
{\rm e}\!\left({\rm Re}\left({\omega A +\omega' D\over Cq}\right)\right) .}$$
The conditions of summation here ensure that $(D,Cq)\sim 1$ and 
$A\equiv D^{*}\bmod Cq{\frak O}$. The above therefore shows that we have 
$$S_{\infty,\infty}\left(\omega,\omega';Cq\right) 
=S\left(\omega,\omega';Cq\right)\qquad\qquad\hbox{($C\in{\frak O}-\{ 0\}$),}$$
with $S(u,v;w)$ as defined in (1.3.6). By this and (9.17), the result (1.4.17) follows. 

To prove (1.4.4), we observe that, by (9.1), (9.16) and the (already proven) 
result (1.4.3), it follows that when $c=Cs\sqrt{r}$, with $C\in{\frak O}-\{ 0\}$ 
and $(C,r)\sim 1$, one has 
$$\eqalign{ 
S_{1/s,\infty}\left(\omega,\omega';c\right) 
 &=\sum_{D+Cs{\frak O}\in{\frak O}/(Cs{\frak O})}
\ \sum_{\scriptstyle A+Cs{\frak O}\in{\frak O}/(Cs{\frak O})\atop\scriptstyle 
ADr\equiv 1\bmod Cs{\frak O}}
{\rm e}\!\left({\rm Re}\left({\omega A\sqrt{r}+\omega' D\sqrt{r}\over Cs\sqrt{r}}\right)\right) = {}\cr 
 &=\sum_{\scriptstyle D\bmod Cs{\frak O}\atop\scriptstyle (D,Cs)\sim 1} 
\sum_{\scriptstyle A\bmod Cs{\frak O}\atop\scriptstyle A\equiv r^{*}D^{*}\bmod Cs{\frak O}} 
{\rm e}\!\left({\rm Re}\left({\omega A +\omega' D\over Cs}\right)\right) .}$$
The result (1.4.4) follows, since the last sum above is (by the definition (1.3.6)) 
equal to the simple Kloosterman sum $S(\omega r^{*},\omega';Cs)$. 

To obtain the result (9.12) (and so complete the proof of the lemma), we note firstly 
that, by (9.11) and (9.13), one has 
$$\delta_{\omega,\omega'}^{{\frak a},{\frak a}} 
=\sum_{B^{+}\pmatrix{u &\beta\cr 0 &1/u}\in B^{+}\backslash g_{\frak a}^{-1}\Gamma g_{\frak a}} 
{\rm e}\left( {\rm Re}\left(\beta u\omega\right)\right) 
\delta_{u\omega , \omega' /u}\eqno(9.18)$$
when ${\frak a}\in\{\infty , 1/s\}\,$ (it should be noted here that if $\gamma\in\Gamma$ 
is such that $g_{\frak a}^{-1}\gamma g_{\frak a}\infty =\infty$, then $\gamma{\frak a}={\frak a}$). 
As an immediate consequence of (9.18) and (9.2), we find that 
$$\eqalign{ 
\delta_{\omega,\omega'}^{1/s,1/s} 
=\delta_{\omega,\omega'}^{\infty,\infty} 
 &=\sum_{B^{+}\pmatrix{u &\beta\cr 0 &1/u}\in B^{+}\backslash\Gamma} 
{\rm e}\left( {\rm Re}\left(\beta u\omega\right)\right) 
\delta_{u\omega , \omega' /u} = {}\cr 
 &=\sum_{u\in{\frak O}^{*}}
\sum_{\beta+(1/u){\frak O}\in{\frak O}/((1/u){\frak O})} 
{\rm e}\left( {\rm Re}\left(\beta u\omega\right)\right) 
\delta_{u\omega , \omega' /u} = {}\cr 
 &=\sum_{u\in{\frak O}^{*}}
\sum_{\beta+{\frak O}\in{\frak O}/{\frak O}} 
{\rm e}\left( {\rm Re}\left(\beta u\omega\right)\right) 
\delta_{u\omega , \omega' /u}\;.}$$
Therefore, since ${\frak O}/{\frak O}=\{ 0+{\frak O}\}$, since  
${\rm e}(0u\omega)={\rm e}(0)=1\,$ (for $u\in{\Bbb C}$), and 
since 
$$\sum_{u\in{\frak O}^{*}}\delta_{u\omega , \omega' /u}
=\sum_{u\in{\frak O}^{*}}\delta_{u^2\omega , \omega'}
=2\delta_{\omega , \omega'}+2\delta_{-\omega , \omega'}\;,$$ 
we obtain all parts of (9.12)\quad$\blacksquare$ 

\bigskip 
\bigskip 

\goodbreak 
\centerline{\it The Proof of Theorem~11.} 

\bigskip 

This proof is an application of the preceding lemma, in conjunction with 
two results from [22]. The first of the latter two results [22, Theorem~B] 
is a `spectral to Kloosterman' summation formula  (inverse in effect to 
the `Kloosterman to spectral' summation formula in~(1.2.1)); 
the other is [22, Theorem~1], which has been reproduced in 
Section~1.2 of the present paper (it appears there as Theorem~2). 

Let the hypotheses of Theorem~11 be satisfied.  Put 
${\cal N}=\{\nu\in{\Bbb C} : -2/3\leq {\rm Re}(\nu)\leq 2/3\}$; and let the function 
$h : {\cal N}\times{\Bbb Z}\rightarrow{\Bbb C}$ be given by: 
$$h(\nu,p)=\cases{\left( X^{\nu}+X^{-\nu}\right)\exp\!\left(\nu^2\right) 
&if $\nu\in{\cal N}$ and $p=0$,\cr 0 &if $\nu\in{\cal N}$ and $p\in{\Bbb Z}-\{ 0\}$.} 
\eqno(9.19)$$

This function $h$ satisfies all of the relevant hypotheses of the 
case $\sigma =2/3$ of [22, Theorem~B], 
as summarised in [22, Theorem~B,~Conditions~(i)-(iii)]. 
Indeed,  by (9.19) we have, for $(\nu,p)\in{\cal N}\times{\Bbb Z}$, 
$$h(-\nu,-p)=h(\nu,p)$$ 
and 
$$\eqalign{ 
h(\nu,p) &\ll (1+|p|)^{-4}\exp\!\left( |(\log X){\rm Re}(\nu)|+|{\rm Re}(\nu)|^2 
-|{\rm Im}(\nu)|^2\right) \leq {}\cr 
 &\leq (1+|p|)^{-4}\exp\!\left( {2\over 3}\,|\log X| +{4\over 9}\right) 
\left( 1+{1\over 2!}\,|{\rm Im}(\nu)|^4\right)^{\!\!-1} 
\ll_X (1+|p|)^{-4}\left( 1+|{\rm Im}(\nu)|\right)^{-4}
}$$ 
(which takes care of [22, Theorem~B, Conditions~(i)~and~(iii)]); 
and, with regard to [22, Theorem~B, Condition~(ii)] 
(requiring that, for each $p\in{\Bbb Z}$, the function 
$\nu\mapsto h(\nu,p)$ have a holomorphic continuation into a 
neighbourhood of the strip ${\cal N}\,$), it 
suffices to note that, since $X$ is positive, both the functions $\nu\mapsto 0$ 
and $\nu\mapsto (X^{\nu}+X^{-\nu})\exp(\nu^2)$ are entire. 
Since we have, moreover, $\Gamma =\Gamma_0(q)$, where 
$q,r,s\in{\frak O}-\{ 0\}$ satisfy (1.4.18), and since the last part of 
Lemma~9.1 (and its proof) shows that the scaling matrices 
$g_{\infty}$ and $g_{1/s}$ satisfy the 
relevant hypotheses (including the condition (1.1.1)), it therefore follows by [22, Theorem~B] that, for 
${\frak a}\in\{\infty , 1/s\}$ and $m,n\in{\frak O}-\{ 0\}$, 
$$\eqalignno{ 
 &\sum_V^{(\Gamma)}\,\overline{c_V^{\frak a}\left(m;\nu_V,p_V\right)}\,
c_V^{\frak a}\left(n;\nu_V,p_V\right) h\left(\nu_V , p_V\right) + {}\cr 
 &\qquad\ +\sum_{{\frak c}\in{\frak C}}^{(\Gamma)}
{1\over 4\pi i\left[\Gamma_{\frak c} : \Gamma_{\frak c}'\right]} 
\sum_{p\in{1\over 2}\left[\Gamma_{\frak c} : \Gamma_{\frak c}'\right]{\Bbb Z}}
\ \,\int\limits_{(0)} \overline{B_{\frak c}^{\frak a}\left(m;\nu,p\right)}\,
B_{\frak c}^{\frak a}\left(n;\nu,p\right) h(\nu,p)\,{\rm d}\nu =\cr
 &\qquad\qquad\qquad\qquad\quad 
=\;{\delta^{{\frak a},{\frak a}}_{m,n}\over 4\pi^3 i}\,
\sum_{p\in{\Bbb Z}}\ \,\int\limits_{(0)} h(\nu,p) \left( p^2 -\nu^2\right)\,{\rm d}\nu  
+\sum_{c\in {}^{\frak a}{\cal C}^{\frak a}}^{(\Gamma)}  
\,{S_{{\frak a},{\frak a}}\left(m , n ; c\right)\over |c|^2}\,
{\bf B}h\left( {2\pi\sqrt{mn}\over c}\right) ,\qquad\quad   
 &(9.20)}$$
where $\delta_{m,n}^{{\frak a},{\frak a}}$ is as defined in the equation~(9.11) 
of Lemma~9.2, and where the ${\bf B}$-transform is that defined above (1.2.5), in the proof 
of Theorem~1; while the meaning of 
any other non-standard notation used is explained in  
Subsection~1.1.

By the results of Lemma~9.2 
(specifically (1.4.16), (1.4.17) and (9.12)), it follows that 
the right-hand side of the equation~(9.20) is independent of the choice of cusp 
(that choice being between having ${\frak a}=1/s$, or else ${\frak a}=\infty$): 
the same is therefore true 
(when $m$ and $n$ are given) 
of the numerical value of the left-hand side of the equation~(9.20). 
Therefore, and since our choice of test-function $h$ (in (9.19)) ensures that 
$h(\nu,p)\neq 0$ only if $p=0$, we may deduce that 
$$\eta_q^{1/s}(m,n;h)
=\eta_q^{\infty}(m,n;h)\qquad\qquad\hbox{($m,n\in{\frak O}-\{ 0\}$),}$$
where 
$$\eqalign{\ \eta_q^{\frak a}(m,n;h) 
 &=\sum_{V\,:\,p_V=0}^{(\Gamma)}\,\overline{c_V^{\frak a}\left(m;\nu_V,0\right)}\,
c_V^{\frak a}\left(n;\nu_V,0\right) h\left(\nu_V , 0\right) + {}\cr 
 &\qquad\quad\, {} +\sum_{{\frak c}\in{\frak C}}^{(\Gamma)}
{1\over 4\pi i\left[\Gamma_{\frak c} : \Gamma_{\frak c}'\right]} 
\,\int\limits_{(0)} \overline{B_{\frak c}^{\frak a}\left(m;\nu,0\right)}\,
B_{\frak c}^{\frak a}\left(n;\nu,0\right) h(\nu,0)\,{\rm d}\nu\;. 
}$$
Consequently one has 
$$H_q^{1/s}({\bf b},N;h)
=H_q^{\infty}({\bf b},N;h)\;,\eqno(9.21)$$
where, for ${\frak a}\in\{ \infty , 1/s\}$, 
$$\eqalign{
H_q^{\frak a}({\bf b},N;h)
 &=\quad\ \,\sum\!\!\!\!\sum_{\!\!\!\!\!\!\!\!\!\!\!\!\!\!\!{{\textstyle{N\over 4}}<|m|^2,|n|^2\leq N}}
\overline{b_m}\,b_n\,\eta_q^{\frak a}(m,n;h) =\cr 
 &=\ \sum_{V\,:\,p_V=0}^{(\Gamma)}\,\Biggl| 
\sum_{{\textstyle{N\over 4}}<|n|^2\leq N} b_n c_V^{\frak a}\left(n;\nu_V,0\right)
\Biggr|^{\,2} h\left(\nu_V , 0\right) + {}\cr 
 &\qquad\quad +\sum_{{\frak c}\in{\frak C}}^{(\Gamma)}
{1\over 4\pi i\left[\Gamma_{\frak c} : \Gamma_{\frak c}'\right]} 
\,\int\limits_{(0)}\,\Biggl| 
\sum_{{\textstyle{N\over 4}}<|n|^2\leq N} b_n 
B_{\frak c}^{\frak a}\left(n;\nu,0\right)\Biggr|^{\,2} h(\nu,0)\,{\rm d}\nu\;.   
}$$

Recall now that the index $V$ in the second last summation denotes a 
cuspidal subspace occurring in the orthogonal decomposition (1.1.3) 
of the space ${}^{0}L^{2}(\Gamma\backslash G)$. 
For each such $V$ the associated spectral parameter $\nu_V$ is either 
positive (and less than $2/9$), or else lies on the ray $i[0,\infty)$ 
in the complex plane (see (1.1.2)-(1.1.4) and the paragraph containing 
(1.1.11)). Hence, and by (9.19), we may rewrite the expression just 
obtained for $H_q^{\frak a}({\bf b},N;h)$ so as to obtain: 
$$H_q^{\frak a}({\bf b},N;h) 
=\rho_q^{\frak a}({\bf b},N;X) 
+2\sum_{j=0}^1 R_{q,j}^{\frak a}({\bf b},N;X)\qquad\qquad 
\hbox{(${\frak a}\in\{\infty , 1/s\}$),} 
\eqno(9.22)$$ 
where $\rho_q^{\frak a}({\bf b},N;X)$ is the sum defined in (1.4.19),  
$$R_{q,0}^{\frak a}({\bf b},N;X) 
=\sum_{\scriptstyle V\,:\,p_V=0\atop\scriptstyle\nu_V\in i[0,\infty)}^{(\Gamma)}
\,\Biggl| 
\sum_{{\textstyle{N\over 4}}<|n|^2\leq N} b_n c_V^{\frak a}\left(n;\nu_V,0\right)
\Biggr|^{\,2} \cos\left( (\log X)\,{\rm Im}\!\left(\nu_V\right)\right) 
\exp\!\left( -\left({\rm Im}\!\left(\nu_V\right)\right)^2\right)$$
and 
$$R_{q,1}^{\frak a}({\bf b},N;X) 
=\sum_{{\frak c}\in{\frak C}}^{(\Gamma)}
{1\over 4\pi\left[\Gamma_{\frak c} : \Gamma_{\frak c}'\right]} 
\,\int\limits_{-\infty}^{\infty}\,\Biggl| 
\sum_{{\textstyle{N\over 4}}<|n|^2\leq N} b_n 
B_{\frak c}^{\frak a}\left(n;it,0\right)\Biggr|^{\,2} 
\cos\left( (\log X)t\right) \exp\!\left( -t^2\right) {\rm d}t\;.$$ 
Moreover, since $-1\leq\cos\theta\leq 1$ for all real $\theta$, and since 
$$\int\limits_{|t|}^{\infty} T\exp\!\left( -T^2\right) {\rm d}T 
={1\over 2}\,\exp\!\left( -t^2\right)\qquad\qquad\quad\hbox{($t\in{\Bbb R}$),}$$
we have here 
$$R_{q,j}^{\frak a}({\bf b},N;X) 
\leq 4\sum_{h=0}^{1}\ \int\limits_0^{\infty} 
E_j^{\frak a}(q,1/2,T;2^{-h}N,{\bf b})\,T\exp\!\left( -T^2\right) {\rm d}T\qquad\quad  
\hbox{(${\frak a}\in\{\infty , 1/s\}$, $j=0,1$),}\eqno(9.23)$$ 
where the sums $E_j^{\frak a}(q,P,K;N,{\bf b})\ $ (${\frak a}\in{\Bbb Q}(i)\cup\{\infty\}$, 
$j=0,1$) are those defined by (1.2.7)-(1.2.8), in Theorem~2. 

In the above, each sum $E_j^{\frak a}(q,1/2,T;N,{\bf b})$ is, by its definition, 
a real-valued and monotonic increasing function of the real variable $T$, 
and satisfies 
$0\leq E_j^{\frak a}(q,1/2,K;N,{\bf b})\leq E_j^{\frak a}(q,1,K;N,{\bf b})$ 
for all real $K$. 
Hence, and by Theorem~2, it follows from (9.23) that, for 
${\frak a}\in\{ 0 , 1/s\}$, $j=0,1$ and any $\varepsilon >0$, one 
has 
$$\eqalignno{ 
R_{q,j}^{\frak a}({\bf b},N;X) 
 &\leq 4\sum_{h=0}^{1}\ \int\limits_1^{\infty} 
E_j^{\frak a}(q,1,K;2^{-h}N,{\bf b})\,(K-1)\exp\!\left( -(K-1)^2\right) {\rm d}K 
 \ll {}\qquad\qquad\qquad\qquad\qquad\quad\ \cr 
 &\ll\int\limits_1^{\infty}  
\left( 1+K^2\right) 
\left( K+O_{\varepsilon}\!\left( N^{1+\varepsilon}|\mu({\frak a})|^2 K^{-1/2}\right)\right) 
\left\|{\bf b}_N\right\|_2^2 (K-1) K^{-6}\,{\rm d}K \ll {}\cr 
 &\ll\left( 1+O_{\varepsilon}\!\left( N^{1+\varepsilon}|\mu({\frak a})|^2\right)\right) 
\left\|{\bf b}_N\right\|_2^2\;. &(9.24)}$$ 
Given the definition of $\mu({\frak a})$ in (1.2.10), and in 
light of Remark~3 below Theorem~2, we have here 
$${1\over\mu(1/s)}\sim{q\over\bigl( (s,q)\,,\,q/(s,q)\bigr)}\sim{q\over (s,r)}\sim q\sim 
{1\over\mu(\infty)}\,,$$
and so we may deduce from (9.22) and (9.24) that 
$$H_q^{\frak a}({\bf b},N;h) 
=\rho_q^{\frak a}({\bf b},N;X) 
+O\!\left( \left( 1+O_{\varepsilon}\!\left( N^{1+\varepsilon}|q|^{-2}\right)\right) 
\left\|{\bf b}_N\right\|_2^2\right)\qquad\qquad  
\hbox{(${\frak a}\in\{\infty , 1/s\}$, $\varepsilon >0$).}$$
By the substitution of these results into (9.21), one obtains the 
result seen in (1.4.20)\quad$\blacksquare$ 

\bigskip 

We end this section with the proof of Theorem~10. 
By way of preparation, we include here three more lemmas. 
The first of these is a corollary of Theorems~3, 4, 8, 9 and~11:  
the others are of a technical nature. 

\bigskip 

\goodbreak\proclaim Lemma~9.3. Let the hypotheses of Theorem~10, concerning 
$\vartheta,\varepsilon,N,L,\delta,P,Q,R,S,X\in{\Bbb R}$, 
the function $A : {\Bbb C}\rightarrow{\Bbb C}$, the set 
${\cal B}(R,S)\subset{\frak O}\times{\frak O}$ and the function 
$b$ be satisfied;  let $a_n\in{\Bbb C}$ for $n\in{\frak O}-\{ 0\}$;   
and, for $u,y\in{\Bbb R}$,  let  
$S_{u,y}^{\infty,*}=S_{u,y}^{\infty,*}(R,S;X;L,N)$ be given by: 
$$S_{u,y}^{\infty,*}
=\sum_{(r,s)\in{\cal B}(R,S)}\!\!\!\!|b(r,s)| 
\!\sum_{\scriptstyle V\atop\scriptstyle\nu_V>0}^{\left(\Gamma_0(rs)\right)}\!\!X^{\nu_V} 
\left|\sum_{{\textstyle{L\over 2}}<|\ell|^2\leq L}\!\!A(\ell) |\ell|^{2iu} 
c_V^{\infty}\!\left(\ell;\nu_V,0\right)
\!\!\sum_{{\textstyle{N\over 4}}<|n|^2\leq N}\!\!\!\overline{a_n}\,|n|^{2iy} 
c_V^{1/s}\left( n;\nu_V,0\right)\right| .\eqno(9.25)$$ 
Then the function  $(u,y)\mapsto S_{u,y}^{\infty,*}$ is continuous and 
bounded on ${\Bbb R}\times{\Bbb R}$; and, for 
$(u,y)\in{\Bbb R}\times{\Bbb R}$,  one has   
$$\left( S_{u,y}^{\infty,*}\right)^{\!2} 
\ll_{\varepsilon} 
Q^{\varepsilon -1} \| b\|_2^2 \left\|{\bf a}_N\right\|_2^2 L (L+Q)(N+Q) 
\!\left( 1+{X^2 LN\over (L+Q)^2 (N+Q)}\right)^{\!\!\vartheta} 
\!\left(\delta^{-1}+|u|\right)^{11}\;,\eqno(9.26)$$
where $\|{\bf a}_N\|_2$ and $\| b\|_2$ are as defined in (1.2.11) and (1.4.13), 
respectively. \hfill\break 
$\hbox{\qquad}$If it is moreover the case that the hypotheses of 
Theorem~9 concerning $H,K\in{\Bbb R}$, $N$, the functions 
$\alpha,\beta : {\Bbb C}\rightarrow{\Bbb C}$ and the coefficients $a_n$ 
($n\in{\frak O}-\{ 0\}$) are satisfied, then one has also 
$$\eqalignno{
\left( S_{u,y}^{\infty,*}\right)^2 
 &\ll_{\varepsilon} 
Q^{\varepsilon} \| b\|_{\infty}^2 N L (L+Q)\times {} 
 &(9.27)\cr 
 &\quad\ {} \times \left(\!\left( 1+{X^2 L\over (H+K)(L+Q)^2}\right)^{\!\!\vartheta}\!N
+\left( 1+{X^2 LN\over Q^2(L+Q)^2}\right)^{\!\!\vartheta}\!Q\!\right)  
\left(\delta^{-1}+|u|\right)^{11}\!\left(\delta^{-1}+|y|\right)^{11} ,}$$
where $\| b\|_{\infty}$ is as defined by (1.4.15). 

\goodbreak 
\noindent{\bf Proof.}\quad 
Let 
$$Z={Q^2\over L}+L\qquad\quad{\rm and}\qquad\quad Z'={X^2\over Z}\;.\eqno(9.28)$$
Then $Z>L\geq 1$ and $Z'>0$, and for $\nu>0$ one has 
$$1\leq X^{\nu}=\left(\sqrt{ Z Z'}\right)^{\nu}=Z^{\nu/2} \left( Z'\right)^{\nu/2}\;.$$
It therefore follows, by (9.25) and the Cauchy-Schwarz inequality, that 
$$0\leq\left( S_{u,y}^{\infty,*}(R,S;X;L,N)\right)^2 
\leq S_{u}^{\infty}(R,S;Z;L)\,S_{y}^{*}(R,S;Z';N)\;,\eqno(9.29)$$ 
where 
$$S_{u}^{\infty}(R,S;Z;L)
=\sum_{(r,s)\in{\cal B}(R,S)} 
\sum_{\scriptstyle V\atop\scriptstyle\nu_V>0}^{\left(\Gamma_0(rs)\right)} Z^{\nu_V}  
\Biggl|\sum_{{\textstyle{L\over 2}}<|\ell|^2\leq L}\!\!\!\!A(\ell) |\ell|^{2iu} 
c_V^{\infty}\!\left(\ell;\nu_V,0\right)\Biggr|^2\eqno(9.30)$$
and
$$S_{y}^{*}(R,S;Z';N)
=\sum_{(r,s)\in{\cal B}(R,S)}\!\!\!|b(r,s)|^2  
\sum_{\scriptstyle V\atop\scriptstyle\nu_V>0}^{\left(\Gamma_0(rs)\right)} \left( Z'\right)^{\nu_V} 
\Biggl|\sum_{{\textstyle{N\over 4}}<|n|^2\leq N}\!\!\!\!\overline{a_n}\,|n|^{2iy} 
c_V^{1/s}\left( n;\nu_V,0\right)\Biggr|^2\;.\eqno(9.31)$$ 

Since the relation $(r,s)\in{\cal B}(r,s)$ implies that 
$0\neq rs\in{\frak O}$ and $RS/4\leq |rs|^2\leq RS=Q$, and since 
$$\left|\left\{ (r,s)\in{\cal B}(R,S) : rs=q\right\}\right| 
\leq\sum_{r\mid q} 1
=O_{\varepsilon}\!\left( |q|^{\varepsilon /4}\right)\qquad\quad  
\hbox{for $\quad q\in{\frak O}-\{ 0\}$,}\eqno(9.32)$$
it consequently follows by (9.30) that 
$$S_{u}^{\infty}(R,S;Z;L)
\ll_{\varepsilon} Q^{\varepsilon /8}\sum_{\scriptstyle w\in\{ 0,1\}\atop\scriptstyle 2^w\leq Q}
\sum_{{\textstyle{Q\over 2^{w+1}}}<|q|^2\leq{\textstyle{Q\over 2^w}}}  
\sum_{\scriptstyle V\atop\scriptstyle\nu_V>0}^{\left(\Gamma_0(q)\right)} Z^{\nu_V}  
\Biggl|\sum_{{\textstyle{L\over 2}}<|\ell|^2\leq L}\!\!\!\!A(\ell) |\ell|^{2iu} 
c_V^{\infty}\!\left(\ell;\nu_V,0\right)\Biggr|^2\;.$$
Therefore, given that the parameters $L$ and $\delta$, and the function $A : {\Bbb C}\rightarrow{\Bbb C}$ 
satisfy the hypotheses stated in Theorem~10, and given that $\vartheta\geq 0$, 
it follows   
by Theorem~8 (applied with $\varepsilon /6$, $Z$, $L$ and $A$ 
substituted for $\varepsilon$, $X$, $H$ and $\alpha$, respectively,  
and, when it is appropriate, with $Q/2$ substituted for $Q$) 
that one has: 
$$S_{u}^{\infty}(R,S;Z;L)
\ll_{\varepsilon} Q^{\varepsilon /6}\left(\delta^{-1}+|u|\right)^{11} 5^{\vartheta} 
(Q+L)^{1+\varepsilon /6} L\;.$$
By this bound, together with Theorem~3 and the conditions in (1.4.6), we find that 
$$S_{u}^{\infty}(R,S;Z;L)
\ll_{\varepsilon} Q^{\varepsilon /2} L (L+Q) \left(\delta^{-1}+|u|\right)^{11}\;.\eqno(9.33)$$ 

In order to obtain a suitable bound for $S_y^{*}(R,S;Z';N)$, we note firstly that, 
since $Z'>0$, it follows that if $q=rs$, with $(r,s)\in{\cal B}(R,S)$, and if 
$$Y=Z'+1\eqno(9.34)$$
(so that $Y>Z'>0$ and $Y>1$), then 
$$\eqalign{
0\leq\sum_{\scriptstyle V\atop\scriptstyle\nu_V>0}^{\left(\Gamma_0(q)\right)} 
\!\!\left( Z'\right)^{\nu_V} 
\!\Biggl|\sum_{{\textstyle{N\over 4}}<|n|^2\leq N}\!\!\!\!\overline{a_n}\,|n|^{2iy} 
c_V^{1/s}\!\left( n;\nu_V,0\right)\Biggr|^{\,2} 
 &\leq 
\sum_{\scriptstyle V\atop\scriptstyle\nu_V>0}^{\left(\Gamma_0(q)\right)} Y^{\nu_V} 
\Biggl|\sum_{{\textstyle{N\over 4}}<|n|^2\leq N}\!\!\overline{a_n}\,|n|^{2iy} 
c_V^{1/s}\!\left( n;\nu_V,0\right)\Biggr|^{\,2} \leq {}\cr 
 &\leq 2\!\!\!\sum_{\scriptstyle w\in\{ 0,1\}\atop\scriptstyle 2^w\leq N} 
\!\!\sum_{\scriptstyle V\atop\scriptstyle\nu_V>0}^{\left(\Gamma_0(q)\right)} 
\!\!Y^{\nu_V} 
\Biggl|\sum_{{\textstyle{N\over 2^{w+1}}}<|n|^2\leq{\textstyle{N\over 2^w}}}\!\!\overline{a_n}\,|n|^{2iy} 
c_V^{1/s}\!\left( n;\nu_V,0\right)\Biggr|^{\,2} ,
}$$
and so, by Theorem~4 (applied with $\varepsilon /8$ and $Y$ substituted for 
$\varepsilon$ and $X$, respectively, with $N/2$ substituted for $N$, 
when appropriate, and with 
$b_n=\overline{a_n}\,|n|^{2iy}$, 
for $n\in{\frak O}-\{ 0\}$, and ${\frak a}=1/s$), one has: 
$$\eqalign{
\sum_{\scriptstyle V\atop\scriptstyle\nu_V>0}^{\left(\Gamma_0(q)\right)} 
\!\!\left( Z'\right)^{\nu_V} 
 &\Biggl|\sum_{{\textstyle{N\over 4}}<|n|^2\leq N}\!\!\!\!\overline{a_n}\,|n|^{2iy} 
c_V^{1/s}\!\left( n;\nu_V,0\right)\Biggr|^{\,2} \ll {}\cr  
 &\qquad\qquad\qquad\ll \left( 1+Y M_{1/s}\,N\right)^{\Theta(q)} 
\left( 1+O_{\varepsilon}\bigl( M_{1/s}\,N^{1+\varepsilon /8}\,\bigr)\right)^{\!1-\Theta(q)} 
\left\|{\bf a}_N\right\|_2^2 \log\bigl( 2+M_{1/s}^{-1}\,N^{-1}\,\bigr) ,}$$ 
where $\Theta(q)$ is as defined in (1.2.20), and where $M_{1/s}=|\mu(1/s)|^2=|q|^{-2}$ 
(for, as noted towards the end of the proof of Theorem~11,  
one has $1/\mu(1/s)\sim q\,$ when $\Gamma=\Gamma_0(q)$ and $s\in{\frak O}$ 
is a factor of $q$ 
such that $(s,q/s)\sim 1$). 
Hence, and by (9.31), Theorem~3, the definition (1.4.8) and 
the conditions in (1.4.6), we find that 
$$\eqalignno{ 
 &S_y^{*}(R,S;Z';N) = {}\cr 
 &\quad\ =\sum_{(r,s)\in{\cal B}(R,S)} |b(r,s)|^2  
\,O_{\varepsilon}\!\!\left( 
\!N^{\varepsilon /8}\!\left( 1+{Y N\over |rs|^2}\right)^{\!\!\Theta(rs)} 
\!\left( 1+{N\over |rs|^2}\right)^{\!\!{1-\Theta(rs)}} 
\left\|{\bf a}_N\right\|_2^2 \log\!\left( 2+{|rs|^2\over N}\right)\right) \ll_{\varepsilon} {}\cr 
 &\quad\ \ll_{\varepsilon} Q^{\varepsilon /2} 
\left( 1+{Y N\over Q}\right)^{\!\!\vartheta} 
\!\left( 1+{N\over Q}\right)^{\!\!{1-\vartheta}} 
\left\|{\bf a}_N\right\|_2^2 \| b\|_2^2\;, &(9.35)}$$
where $\| b\|_2$ is as defined in (1.4.13), while $Y$ is given by 
(9.34) and (9.28) (ensuring that we have $Y>1$, and so justifying the  
upper bound given here for the sum over $r$ and $s$).   

By the combination of results in (9.29), (9.33) and (9.35), we have the bound  
$$\left( S_{u,y}^{\infty,*}(R,S;X;L,N)\right)^2 
\ll_{\varepsilon} Q^{\varepsilon} L (L+Q) 
\left( 1+{Y N\over Q}\right)^{\!\!\vartheta} 
\!\left( 1+{N\over Q}\right)^{\!\!{1-\vartheta}} 
\left\|{\bf a}_N\right\|_2^2 \| b\|_2^2\left(\delta^{-1}+|u|\right)^{11}\;,\eqno(9.36)$$ 
where, by (9.28) and (9.34), 
$$1+{YN\over Q} 
=1+{N\over Q}+{Z' N\over Q} 
=1+{N\over Q}+{X^2 L N\over\left( Q^2+L^2\right) Q}
\leq\left( 1+{2 X^2 L N\over (Q+L)^2 (Q+N)}\right)\left( 1+{N\over Q}\right) .$$
Since $0\leq\vartheta\leq 2/9$, the relations on the last line imply that  
$$\left( 1+{Y N\over Q}\right)^{\!\!\vartheta} 
\!\left( 1+{N\over Q}\right)^{\!\!{1-\vartheta}} 
\leq\left( 1+{2X^2 L N\over (Q+L)^2 (Q+N)}\right)^{\vartheta}\left( 1+{N\over Q}\right) 
\ll \left( 1+{X^2 L N\over (Q+L)^2 (Q+N)}\right)^{\vartheta}\left( {Q+N\over Q}\right) ,
$$ 
and so, by the bound in (9.36), we obtain the 
result stated in (9.26). 

The statement preceding (9.26) merits some justification. 
Recall that, for each Hecke congruence subgroup $\Gamma\leq SL(2,{\frak O})$,  
there can be at most a finite number, $E(\Gamma)$ (say), 
of irreducible cuspidal subspaces 
$V\subset L^2(\Gamma\backslash G)$ that have $\nu_V>0\,$ (these $V$ corresponding 
to the `exceptional' eigenvalues $\lambda_V$ discussed below (1.1.11)). 
Hence, and by (1.4.6), (1.4.8), (9.25) and (9.32), one has 
$$S_{u,y}^{\infty,*}
=\sum_{0<|q|^2\leq Q}\!C_q 
\!\sum_{j=1}^{E(\Gamma_0(q))}X^{\nu(q,j)} 
\left|\sum_{{\textstyle{L\over 2}}<|\ell|^2\leq L} 
\ \sum_{{\textstyle{N\over 4}}<|n|^2\leq N} w_{q,j}(\ell) z_{q,j}(n) |\ell|^{2iu} |n|^{2iy} 
\right|\qquad\qquad\hbox{($u,y\in{\Bbb R}$),}$$ 
where $Q$, $L$, $N$, $C_q$, $\nu(q,j)$, $w_{q,j}(\ell)$ and $z_{q,j}(n)$ 
denote complex numbers (real and non-negative, in the 
case of $Q$, $L$, $N$, $C_q$ and $\nu(q,j)$) 
that are independent of the variables $u$ and $y$. 
Since one has here $E(\Gamma_0(q))<\infty\,$ 
and $|\{ m\in{\frak O} : 0<|m|^2\leq M\}|\leq 4M<\infty\,$ 
(for $0\neq q\in{\frak O}$, $M\geq 0$), 
and since all functions of the form
$(u,y)\mapsto |\ell|^{2iu}|n|^{2iy}$ are continuous, 
the continuity of the function $(u,y)\mapsto S_{u,y}^{\infty,*}$ on 
${\Bbb R}\times{\Bbb R}$ 
therefore follows.  We have, moreover, both
$\min\{ S_{u,y}^{\infty,*} : u,y\in{\Bbb R}\}\geq 0$, and 
$$\max\{ S_{u,y}^{\infty,*} : u,y\in{\Bbb R}\}
\leq\sum_{0<|q|^2\leq Q}\!C_q 
\!\sum_{j=1}^{E(\Gamma_0(q))}X^{\nu(q,j)} 
\sum_{{\textstyle{L\over 2}}<|\ell|^2\leq L} |w_{q,j}(\ell)|
\sum_{{\textstyle{N\over 4}}<|n|^2\leq N} |z_{q,j}(n)|<\infty ,$$ 
so that the function $(u,y)\mapsto S_{u,y}^{\infty,*}$ is 
bounded, as asserted in the statement of the lemma.

In proving the one remaining result of the lemma, which is the bound (9.27), we 
may of course assume the relevant premise, stated in the lemma. 
Accordingly, it is to be supposed that the hypotheses of Theorem~9 
concerning $H,K\in{\Bbb R}$, $\,N$, the functions 
$\alpha , \beta : {\Bbb C}\rightarrow{\Bbb C}$ and coefficients 
$a_n$ ($n\in{\frak O}-\{ 0\}$) are satisfied. Since we assume these 
hypotheses in addition to (and not as a substitute for) 
those that were previously assumed, it follows that all the results previously 
obtained in this proof remain valid: we refer, in particular, to the results (9.29), 
(9.33) and (9.35), through which (9.36) was obtained.  
We shall show that, given the additional hypotheses, one can obtain a stronger bound 
for the sum $S_y^{*}(R,S;Z';N)$ than that 
obtained in (9.35). By using this stronger bound, together with (9.29) and (9.33), 
we shall obtain the result in (9.27).  

Let $Y$ and $Z'$ be given (as previously) by (9.34) and (9.28), so that 
$Y>Z'>0$. Then, for $\nu >0$, 
$$\left( Z'\right)^{\nu}<Y^{\nu}<Y^{\nu}+Y^{-\nu}
<\left( Y^{\nu}+Y^{-\nu}\right)\exp\!\left( \nu^2\right) 
=h_{Y}\!(\nu)\qquad\ \,\hbox{(say).}\eqno(9.37)$$
Hence, and by Theorem~11 (with $b_n=\overline{a_n}\,|n|^{2iy}$ for $n\in{\frak O}-\{ 0\}$, 
and with $\varepsilon /8$ substituted for $\varepsilon$), 
we find that if $(r,s)\in{\cal B}(R,S)$ and $q=rs$, then  
$$\eqalignno{
\sum_{\scriptstyle V\atop\scriptstyle\nu_V>0}^{\left(\Gamma_0(q)\right)} \left( Z'\right)^{\nu_V} 
\Biggl|\sum_{{\textstyle{N\over 4}}<|n|^2\leq N}\!\!\!\!\overline{a_n}\,|n|^{2iy} 
c_V^{1/s}\left( n;\nu_V,0\right)\Biggr|^2 
 &\leq\sum_{\scriptstyle V\atop\scriptstyle\nu_V>0}^{\left(\Gamma_0(q)\right)} 
h_Y\!\left(\nu_V\right)  
\Biggl|\sum_{{\textstyle{N\over 4}}<|n|^2\leq N}\!\!\!\!\overline{a_n}\,|n|^{2iy} 
c_V^{1/s}\left( n;\nu_V,0\right)\Biggr|^2 = {}\qquad\quad\ \cr 
 &=\sum_{\scriptstyle V\atop\scriptstyle\nu_V>0}^{\left(\Gamma_0(q)\right)} 
h_Y\!\left(\nu_V\right) 
\Biggl|\sum_{{\textstyle{N\over 4}}<|n|^2\leq N}\!\!\!\!\overline{a_n}\,|n|^{2iy} 
c_V^{\infty}\left( n;\nu_V,0\right)\Biggr|^2 + \cr 
 &\quad\ +O\!\left(\left( 1+O_{\varepsilon}\!\left( |q|^{-2} N^{1+\varepsilon /8}\right)\right) 
\left\|{\bf a}_N\right\|_2^2\right) . &(9.38)}$$

Since $Y=Z'+1>1$, the function $h_Y:(0,\infty)\rightarrow(0,\infty)$ defined in (9.37) 
satisfies 
$$2<h_Y(\nu)<2Y^{\nu}\exp\!\left(\nu^2\right)\qquad\quad\hbox{for $\quad\nu >0$.}$$
Morever, it follows by 
(1.2.20), (1.2.21) and (1.1.11) that,  
for each cuspidal subspace 
$V\subset{}^0L^2\bigl(\Gamma_0(q)\backslash G\bigr)$ 
indexing a term of the sum on the right-hand side of the equation (9.38), 
one has $0<\nu_V<2/9$, and so 
$$1<\exp\!\left(\nu_V^2\right)<\exp(4/81)\;.$$ 
With the aid of these observations, one may deduce from (9.38) and (9.31) 
that 
$$\eqalign{
0\leq S_y^{*}(R,S;Z';N) 
 &\leq 2\exp(4/81) 
\sum_{(r,s)\in{\cal B}(R,S)}\!\!\!|b(r,s)|^2  
\sum_{\scriptstyle V\atop\scriptstyle\nu_V>0}^{\left(\Gamma_0(rs)\right)} Y^{\nu_V} 
\Biggl|\sum_{{\textstyle{N\over 4}}<|n|^2\leq N}\!\!\!\!\overline{a_n}\,|n|^{2iy} 
c_V^{\infty}\left( n;\nu_V,0\right)\Biggr|^2 +\cr 
 &\quad\ \, +O_{\varepsilon}\!\left( 
N^{\varepsilon/8}\left\|{\bf a}_N\right\|_2^2 
\sum_{(r,s)\in{\cal B}(R,S)}\!\!\!|b(r,s)|^2 \left( 1+{N\over |rs|^2}\right)\right) .}$$ 
Hence, and by (1.4.6), (1.4.8) and (9.32), we obtain: 
$$\eqalignno{ 
S_y^{*}(R,S;Z';N) 
 &\ll_{\varepsilon} Q^{\varepsilon /8} \| b\|_{\infty}^2 
\sum_{{\textstyle{Q\over 4}}<|q|^2\leq Q}  
\sum_{\scriptstyle V\atop\scriptstyle\nu_V>0}^{\left(\Gamma_0(q)\right)} Y^{\nu_V} 
\Biggl|\sum_{{\textstyle{N\over 4}}<|n|^2\leq N}\!\!\!\!\overline{a_n}\,|n|^{2iy} 
c_V^{\infty}\left( n;\nu_V,0\right)\Biggr|^2 + 
 &(9.39)\cr 
 &\qquad +N^{\varepsilon/8}\left\|{\bf a}_N\right\|_2^2 \| b\|_{\infty}^2 (Q+N)\;,}$$ 
where, as noted in (8.60), within the proof of Theorem~9, one has 
$$\left\|{\bf a}_N\right\|_2^2\ll_{\eta} N^{1+\eta}\qquad\qquad\hbox{($\eta >0$),}\eqno(9.40)$$
with an implicit constant that is determined by $\eta$ and the value of the 
implicit constant associated with the case $j=k=0$ of the condition (1.3.15). 

Given what we currently suppose concerning the coefficients $a_n\ $   
($n\in{\frak O}-\{ 0\}$), it follows by Theorem~9 
(applied with $\varepsilon /8$, $\overline{\alpha} : {\Bbb C}\rightarrow{\Bbb C}$, 
$\overline{\beta} : {\Bbb C}\rightarrow{\Bbb C}$, $y$, $Y$ and, when appropriate, 
$Q/2$ substituted for $\varepsilon$, $\alpha : {\Bbb C}\rightarrow{\Bbb C}$, 
$\beta : {\Bbb C}\rightarrow{\Bbb C}$, $t$, $X$ and $Q$, respectively) 
that one has, for $Q_1\in\{ Q/2 , Q\}$ with $Q_1\geq 1$, 
$$\eqalignno{
\!\!\!\!\!\sum_{{\textstyle{Q_1\over 2}}<|q|^2\leq Q_1}  
\sum_{\scriptstyle V\atop\scriptstyle\nu_V>0}^{\left(\Gamma_0(q)\right)} Y^{\nu_V} 
 &\Biggl|\sum_{{\textstyle{N\over 4}}<|n|^2\leq N}\!\!\!\!\overline{a_n}\,|n|^{2iy} 
c_V^{\infty}\left( n;\nu_V,0\right)\Biggr|^2 \ll_{\varepsilon} {} 
 &(9.41)\cr 
 &\qquad\ll_{\varepsilon} 
\left(\delta^{-1}+|y|\right)^{11}\!\left( 
\!\left( 1+{Y\over Q_1^2 N^{-1}}\right)^{\!\!\vartheta}\!Q 
+\left( 1+{Y\over H+K}\right)^{\!\!\vartheta}\!N\!\right) Q^{\varepsilon /8} 
N^{1+\varepsilon /8}\;.}$$
Since $0\leq\vartheta\leq 2/9$, since $0<\delta\leq 1$, 
since $N\leq Q^2\,$ (by (1.4.6)), 
and since $Y=Z'+1$, where $Z'$ is given by 
(9.28),  it follows by (9.39), (9.40) (for $\eta =\varepsilon /8$) and (9.41)   
that one has 
$$S_y^{*}(R,S;Z';N) 
\ll_{\varepsilon} 
Q^{\varepsilon /2} N \| b\|_{\infty}^2 
\left( \!\left( 1+{Y\over Q^2 N^{-1}}\right)^{\!\!\vartheta}\!Q 
+\left( 1+{Y\over H+K}\right)^{\!\!\vartheta}\!N\!\right)  
\left(\delta^{-1}+|y|\right)^{11}\;,\eqno(9.42)$$
where 
$${Y\over Q^2 N^{-1}}={1+Z'\over Q^2 N^{-1}} 
\leq 1+{Z' N\over Q^2}=1+{X^2 N\over\left( Q^2 L^{-1}+L\right) Q^2} 
\ll 1+{X^2 NL\over (Q+L)^2 Q^2}$$
and 
$${Y\over H+K}={1+Z'\over H+K}\leq {1\over 2}+{Z'\over H+K} 
={1\over 2}+{X^2\over\left( Q^2 L^{-1}+L\right) (H+K)} 
\ll 1+{X^2 L\over (Q+L)^2 (H+K)}\;.$$
By (9.29), (9.33), and (9.42) and the bounds just noted, 
the result in (9.27) follows\ $\blacksquare$ 

\bigskip 

\goodbreak\proclaim Lemma~9.4. Let $B>1$, $X>0$ and $t\in{\Bbb R}$; 
let the function $\Omega : (0,\infty)\rightarrow{\Bbb C}$ be infinitely differentiable, 
with support ${\rm Supp}(\Omega)\subseteq\bigl[ B^{-1} , B\bigr]$; and let 
the function $\varphi : (0,\infty)\rightarrow{\Bbb C}$ be given by 
$$\varphi(r)=\Omega\!\left( X r^2\right) r^{2it}\qquad\qquad\hbox{($r>0$).}\eqno(9.43)$$
Then $\varphi$ is infinitely differentiable on $(0,\infty)$, satisfies 
$$\varphi^{(j)}(r)\ll_{\Omega,j} (1+|t|)^j r^{-j}
\qquad\qquad\hbox{($j\in{\Bbb N}\cup\{ 0\}$ and $r>0$),}\eqno(9.44)$$ 
and has 
$${\rm Supp}(\varphi)\subseteq\bigl[\,B^{-1/2} X^{-1/2}\,,
\,B^{1/2} X^{-1/2}\,\bigr]\;.\eqno(9.45)$$ 
Moreover, the function $f : {\Bbb C}^{*}\rightarrow{\Bbb C}$ given by 
$$f(z)=\varphi\left( |z|\right)\qquad\qquad\hbox{($z\in{\Bbb C}^{*}$),}\eqno(9.46)$$ 
is even, smooth and compactly supported in $\,{\Bbb C}^{*}$. 

\goodbreak 
\noindent{\bf Proof.}\quad 
Note firstly that the function $R\mapsto\Omega(R^2)$ is infinitely 
differentiable on $(0,\infty)$. Indeed, this function 
is the composition of functions $\Omega\circ q$, where $q : (0,\infty)\rightarrow(0,\infty)$ is 
given by 
$$q(R)=R^2\qquad\qquad\hbox{($R>0$),}\eqno(9.47)$$
and so (given that the functions $x\mapsto{\rm Re}(\Omega(x))$ and 
$x\mapsto{\rm Im}(\Omega(x))$ are infinitely differentiable)  
it follows by the case $m=1$, $U=V=(0,\infty)$ and $f=q$ of Lemma~8.1 
that $\Omega\circ q$ is infinitely differentiable on $(0,\infty)$. 
Since the function $r\mapsto X^{1/2}r$ is infinitely differentiable 
(and positive valued) on $(0,\infty)$, it similarly follows that 
the function 
$$r\mapsto(\Omega\circ q)\bigl( X^{1/2}r\bigr)=\Omega\!\left( X r^2\right)$$ 
is infinitely differentiable on $(0,\infty)$. 
By the chain-rule of differential calculus (and 
the principal of induction), it may moreover be deduced that, 
for $r>0$ and $j=0,1,2,\ldots\ $, one has: 
$$\eqalignno{ 
\biggl|{{\rm d}^j\over{\rm d}r^j}\,\Omega\bigl( X r^2\bigr)\biggr| 
 &=\biggl|{{\rm d}^j\over{\rm d}r^j}\,(\Omega\circ q)\bigl( X^{1/2} r\bigr)\biggr| = {}\cr 
 &=\left|\,X^{j/2}\left(\Omega\circ q\right)^{(j)}\!\bigl( X^{1/2} r\bigr)\right| \leq {}\cr 
 &\leq\left(\max_{R>0} R^j\bigl|\left(\Omega\circ q\right)^{(j)}(R)\bigr|\right) r^{-j} 
=\biggl(\,\max_{1/\sqrt{B}\leq R\leq\sqrt{B}} 
\,R^j\bigl|\left(\Omega\circ q\right)^{(j)}(R)\bigr|\,\biggr) r^{-j}\;.  &(9.48)}$$ 

Since the function $r\mapsto\Omega(Xr^2)$ is infinitely differentiable, and 
since  
$${{\rm d}^j\over{\rm d}r^j}\,\left( r^{2it}\right) 
=(2it)(2it-1)\cdots(2it-j+1) r^{2it-j}\qquad\qquad  
\hbox{($j\in{\Bbb N}\cup\{ 0\}$, $r>0$),}$$
we may 
deduce from (9.43) and (9.48), 
by Leibniz's rule for higher order derivatives of a product,   
that the function $\varphi$ is 
infinitely differentiable on $(0,\infty)$, and is such that, 
for all $j\in{\Bbb N}\cup\{ 0\}$, and all $r>0$, one has:   
$$\left|\varphi^{(j)}(r)\right| 
=\left|\sum_{k=0}^j \pmatrix{j\cr k} 
\Biggl(\,\prod_{0\leq\ell<k}(2it-\ell)\Biggr) r^{2it-k} 
\,{{\rm d}^{j-k}\over{\rm d}r^{j-k}}\,\Omega\left( X r^2\right)\right| 
\ll_j \sum_{k=0}^j\,(1+|t|)^k r^{-k}\,O_{\Omega,j-k}\!\left( r^{-(j-k)}\right)\;.$$ 
We consequently obtain the bound stated in (9.44). 

By (9.43), ${\rm Supp}(\varphi)=\{\sqrt{x/X} : x\in{\rm Supp}(\Omega)\}$. 
The result (9.45) is therefore an immediate corollary of the hypothesis 
that ${\rm Supp}(\Omega)\subseteq [B^{-1} , B]$. 

I order to complete this proof we have now only to 
verify the assertions of the lemma concerning the 
function $f : {\Bbb C}^{*}\rightarrow{\Bbb C}$ given by (9.46). 

Firstly, we may note that, since $|{-z}|=|z|$ for $z\in{\Bbb C}$, 
it is ensured by the definition (9.46) that the function $f$ is even. 

Secondly, we observe that, by the relations (9.45) and (9.46), one has  
${\rm Supp}(f)\subseteq{\cal A}$, where 
$${\cal A}
=\left\{ z\in{\Bbb C} : B^{-1/2}X^{-1/2}\leq |z|\leq B^{1/2}X^{-1/2}\right\}
\subset {\Bbb C}^{*}\;.$$
We claim that the set ${\rm Supp}(f)$ is, therefore, a closed and bounded subset 
of ${\Bbb C}$, and so is compact (with respect to the usual topology on ${\Bbb C}$). 
The boundedness of ${\rm Supp}(f)$ follows immediately from our 
observation that ${\rm Supp}(f)$ is contained within the annular region ${\cal A}\,$  
(itself clearly bounded). To see that ${\rm Supp}(f)$ is also a closed 
subset of ${\Bbb C}$, we begin with the observation that 
${\rm Supp}(f)$ is, by definition, a closed subset with respect to the relative 
topology on ${\Bbb C}^{*}$. There is, consequently, 
some set ${\cal Z}\subseteq{\Bbb C}$ which is closed in ${\Bbb C}$ and  
satisfies ${\cal Z}\cap {\Bbb C}^{*}={\rm Supp}(f)$. 
Therefore, and since ${\rm Supp}(f)\subseteq{\cal A}\subset{\Bbb C}^{*}$, it  
follows that ${\rm Supp}(f)$ is the intersection of two closed subsets 
of ${\Bbb C}$ (namely ${\cal A}$ and ${\cal Z}$), and so is itself a 
closed subset of ${\Bbb C}$. As the above has verified 
our claims concerning the set ${\rm Supp}(f)$, we may conclude that 
$f$ is indeed compactly supported in ${\Bbb C}^{*}$. 

Thirdly (and finally), we note that one has, by (9.46), 
$$f(x+iy)=\varphi\left(\bigl( x^2+y^2\bigr)^{\!1/2}\right) 
= \bigl(\varphi\circ q^{-1}\bigr)\bigl( x^2+y^2\bigr)
\qquad\qquad\hbox{($(x,y)\in{\Bbb R}^2 -\{ (0,0)\}$),}\eqno(9.49)$$ 
where $q^{-1} : (0,\infty)\rightarrow(0,\infty)$ is the inverse of 
the function $q$ given by (9.47). 
Hence, and since all three of the functions $u\mapsto q^{-1}(u)=\sqrt{u}$, 
$\,v\mapsto{\rm Re}(\varphi(v))$ and $v\mapsto{\rm Im}(\varphi(v))$ 
are infinitely differentiable on $(0,\infty)$, it follows by the 
case $m=1$, $U=V=(0,\infty)$, $f=q^{-1}$ of Lemma~8.1 that 
the functions 
$u\mapsto{\rm Re}((\varphi\circ q^{-1})(u))$ and 
$u\mapsto{\rm Im}((\varphi\circ q^{-1})(u))$ 
are infinitely differentiable on $(0,\infty)$. 
Therefore, by two further applications of Lemma~8.1 
(both with $m=2$, $U={\Bbb R}^2-\{ (0,0)\}$, 
$V=(0,\infty)$ and $f : U\rightarrow V$ given 
by $f({\bf u})=u_1^2+u_2^2$, for ${\bf u}\in U$), 
it may be deduced from (9.49) that the function $f : {\Bbb C}^{*}\rightarrow{\Bbb C}$ 
given by (9.46) is smooth (in the sense defined at the start of Subsection~1.2)  
\ $\blacksquare$ 

\bigskip 

\goodbreak\proclaim Lemma~9.5. Let $\eta >0$. Then there exists an infinitely 
differentiable function 
$\Omega : (0,\infty)\rightarrow[0,1]$ which has 
$${\rm Supp}(\Omega)\subseteq\left[ 2^{-\eta -2}\,,\,2^{\eta +2}\right]\eqno(9.50)$$
and satisfies 
$$\Omega(u)=1\quad\hbox{for $2^{-\eta}\leq u\leq 2^{\eta}$.}\eqno(9.51)$$

\goodbreak 
\noindent{\bf Proof.}\quad 
It will suffice to construct an infinitely differentiable function 
$\Psi : {\Bbb R}\rightarrow [0,1]$ which has 
$${\rm Supp}(\Psi)\subseteq\left[ {-\eta -2}\,,\,{\eta +2}\right]\eqno(9.52)$$
and satisfies 
$$\Psi(y)=1\quad\hbox{for ${-\eta}\leq y\leq {\eta}$.}\eqno(9.53)$$
For then the function $\Omega : (0,\infty)\rightarrow[0,1]$ given by 
$$\Omega(u)=\Psi\!\left({\log u\over\log 2}\right)\qquad\qquad\quad\hbox{($u>0$)}$$
will be such that the conditions~(9.50) and~(9.51) are 
satisfied, and (by the case $m=1$, $U=(0,\infty)$, $V={\Bbb R}$, 
$f(u)=(\log u)/(\log 2)$ of Lemma~8.1) will, moreover, be an 
infinitely differentiable function on~$(0,\infty)$. 

We claim that a suitable function $\Psi : {\Bbb R}\rightarrow[0,1]\,$   
(infinitely differentiable, and such that (9.52) and (9.53) hold) 
is given by: 
$$\Psi(y)=\left( X(1)\right)^{-1}\left( 
X(y+\eta +1)-X(y-\eta -1)\right)\qquad\qquad\hbox{($y\in{\Bbb R}$),}\eqno(9.54)$$
with
$$X(x)=\int_{-\infty}^x \Phi(t)\,{\rm d}t\qquad\qquad\ \hbox{($x\in{\Bbb R}$),}\eqno(9.55)$$
where (as in the proof 
of Theorem~4) the function $\Phi$ is the infinitely differentiable real function 
defined  by the 
second of the equations below~(3.5). 

In order to verify this claim, it 
suffices to show that (9.55) defines an infinitely differentiable function 
$X : {\Bbb R}\rightarrow{\Bbb R}$ which is zero on $(-\infty,-1]$, 
strictly increasing on $(-1,1)$, and constant on $[1,\infty)$. 
For, if that is the case, then $X$ is an increasing 
and infinitely differentiable function 
on ${\Bbb R}$, with  
$$X(x)=X(-1)\quad{\rm for}\quad x\leq -1\;,\qquad\quad   
X(x)=X(1)\quad{\rm for}\quad x\geq 1\;,\qquad{\rm and}\qquad X(1)>X(-1)=0\;;$$  
and so it then follows, by (9.54), that the function  $\Psi$ is infinitely differentiable on 
${\Bbb R}$, with range contained in the interval $[0/X(1),(X(1)-0)/X(1)]=[0,1]$, 
and with  
$$\Psi(y)=\cases{
(X(-1)-X(-1))/X(1)=(0-0)/X(1)=0
&if $y\leq -\eta -2$;\cr
(X(1)-X(-1))/X(1)=(X(1)-0)/X(1)=1
&if $-\eta\leq y\leq\eta$;\cr
(X(1)-X(1))/X(1)=0/X(1)=0
&if $y\geq \eta +2$.}$$ 

Observe now that (9.55) does indeed define a real function $X$ with all 
of the properties just mentioned. Indeed, since $\Phi$ is continuous 
on ${\Bbb R}$, since 
the range of $\Phi$ is contained in $[0,\infty)$, 
and since $\Phi(t)>0$ if and only if $-1<t<1$, 
the integral on 
the right-hand side of the equation~(9.55) equals zero for $x\leq -1$, 
and (by the first fundamental theorem of integral calculus) 
is strictly increasing on the interval $(-1,1)$, and constant on $[1,\infty)$; 
one has, moreover, $X'(x)=\Phi(x)\ $  ($x\in{\Bbb R}$),  
and so, given that $\Phi$ is infinitely differentiable 
on ${\Bbb R}$, it follows that $X$ is infinitely  differentiable 
on ${\Bbb R}$.  
This completes the proof of the lemma\ $\blacksquare$ 
 
\bigskip 
\bigskip 

\goodbreak 
\centerline{\it The Proof of Theorem~10.} 

\bigskip 

By the case $\eta =4$ of Lemma~9.5, we may choose (once and for all) 
an infinitely differentiable function $\Omega : (0,\infty)\rightarrow[0,1]$ 
such that  
$${\rm Supp}(\Omega)\subseteq\left[ 2^{-6}\,,\, 2^6\right]\eqno(9.56)$$ 
and 
$$\Omega(u)=1\qquad\hbox{for $2^{-4}\leq u\leq 2^4$.}\eqno(9.57)$$
We choose this function $\Omega$ independently of all other parameters 
(as might, for example, be achieved by defining $\Omega(u)=\Psi((\log u)/(\log 2))\ $  
($u>0$), where $\Psi$ is the real function constructed in the proof of Lemma~9.5). 

Let the hypotheses of Theorem~10 be satisfied. Then, by  and the definition (1.4.11) 
and the results (1.4.3),  
(1.4.4) (established by Lemma~9.2), 
we have, for $(r,s)\in{\cal B}(R,S)$ and $n,\ell\in{\frak O}$, 
$$K_{r,s}(n,\ell) 
=\sum_{\scriptstyle 0\neq p\in{\frak O}\atop\scriptstyle (p,r)\sim 1} 
g_{r,s}\!\left( |p|^2\right) S_{1/s,\infty}\!\left( n,\ell;ps\sqrt{r}\right) 
=\sum_{c\in{}^{1/s}{\cal C}^{\infty}}^{\left(\Gamma_0(rs)\right)} 
g_{r,s}\!\left( {|c|^2\over |s|^2 |r|}\right) 
S_{1/s,\infty}\!\left( n,\ell;c\right)\;,$$
where ${}^{\frak a}{\cal C}^{\frak b}$ and 
$S_{{\frak a},{\frak b}}\!\left(\omega,\omega';c\right)$ are given by 
(1.1.13)-(1.1.15)$\,$ (with $\Gamma =\Gamma_0(rs)$ there). Since it is 
moreover the case that ${\rm Supp}(g_{r,s})\subseteq [P/2,P]$, 
we therefore have: 
$$K_{r,s}(n,\ell) 
=\sum_{\scriptstyle c\in{}^{1/s}{\cal C}^{\infty}\atop\scriptstyle 
D(r,s)/2<|c|^2<D(r,s)}^{\left(\Gamma_0(rs)\right)} 
g_{r,s}\!\left( {|c|^2\over |s|^2 |r|}\right) 
S_{1/s,\infty}\!\left( n,\ell;c\right)\qquad\qquad\hbox{($(r,s)\in{\cal B}(R,S)$, 
$\,n,\ell\in{\frak O}$),}\eqno(9.58)$$
where 
$$D(r,s)=P |s|^2 |r|\;.\eqno(9.59)$$

To prepare for an application of the Corollary to Theorems~1 and~2, 
we observe now that, if $r$, $s$, $\ell$, $n$ and $c$ satisfy the 
conditions of summation in (1.4.10) and (9.58), then, by 
(1.4.8) and (9.59), one has  
$$\left| {2\pi\sqrt{n\ell}\over c}\right|^2 
={4\pi^2 |n\ell|\over |c|^2}
\in\left( {4\pi^2 (N/4)^{1/2} (L/2)^{1/2}\over P |s|^2 |r|}\ ,\,   
{4\pi^2 N^{1/2} L^{1/2}\over\left( P |s|^2 |r|/2\right)}\right) 
\subseteq\left( 2^{-3/2} X^{-1}\,,\,2^{5/2} X^{-1}\right) ,$$
where the parameter $X\geq 2$ is that given by the equation~(1.4.7). 
The relations in (9.57)-(9.59) therefore imply that 
$$K_{r,s}(n,\ell) 
=\sum_{\scriptstyle c\in{}^{1/s}{\cal C}^{\infty}\atop\scriptstyle 
D(r,s)/2<|c|^2<D(r,s)}^{\left(\Gamma_0(rs)\right)} 
\Omega\!\left(\left| {2\pi\sqrt{n\ell}\over c}\right|^2 X\right) 
g_{r,s}\!\!\left( {|c|^2\over |s|^2 |r|}\right) 
S_{1/s,\infty}\!\left( n,\ell;c\right)$$
whenever $r$, $s$, $\ell$ and $n$ satisfy the 
conditions of summation in (1.4.10).  Now the conditions  
$D(r,s)/2<|c|^2$ and $|c|^2<D(r,s)$ are redundant in the above summation, 
for (as is implicit in our derivation of (9.58)) the factor 
$g_{r,s}(|c|^2/|s|^2 |r|)$ is equal to zero whenever those 
conditions of summation are not both satisfied. Therefore any 
weakening of those conditions has no effect on the value of the sum. 
Hence, given the definitions in (9.59) and (1.4.8), it is 
certainly the case that if  $r$, $s$, $\ell$ and $n$ satisfy the 
conditions of summation in (1.4.10) then 
$$K_{r,s}(n,\ell) 
=\sum_{\scriptstyle c\in{}^{1/s}{\cal C}^{\infty}\atop\scriptstyle 
2^{-8} D^{*}<|c|^2< 2^8 D^{*}}^{\left(\Gamma_0(rs)\right)} 
\Omega\!\left(\left| {2\pi\sqrt{n\ell}\over c}\right|^2 X\right) 
g_{r,s}\!\!\left( {|c|^2\over |s|^2 |r|}\right) 
S_{1/s,\infty}\!\left( n,\ell;c\right) ,\eqno(9.60)$$
where 
$$D^{*}=D\bigl(\sqrt{R}\,,\,\sqrt{S}\;\bigr) =PSR^{1/2}\;.\eqno(9.61)$$

For $(r,s)\in{\cal B}(R,S)$, the function 
$g_{r,s} : (0,\infty)\rightarrow{\Bbb C}$ is, by hypothesis,  
both infinitely differentiable and compactly supported, and so,  
by Mellin's inversion formula [10, Appendix, Equation~(A.2)], one has 
$$g_{r,s}(x) 
={1\over 2\pi i}\int\limits_{\sigma -i\infty}^{\sigma+i\infty} 
G_{r,s}(w) x^{-w}\,{\rm d}w\qquad\qquad\hbox{($\sigma\in{\Bbb R}$ and $x>0$),}\eqno(9.62)$$
where 
$$G_{r,s}(w) 
=\int\limits_0^{\infty} x^{w-1} g_{r,s}(x)\,{\rm d}x\qquad\qquad\quad\hbox{($w\in{\Bbb C}$).}\eqno(9.63)$$ 
Note that the integral on the right-hand side of (9.62) is absolutely convergent. 
Indeed, 
if $(r,s)\in{\cal B}(R,S)$ then, given the definition (9.63), and given our hypotheses 
concerning the function $g_{r,s}$ (which include the bounds in (1.4.9)), 
we may use  repeated integrations by parts to obtain:  
$$\eqalignno{
G_{r,s}(\sigma +it) 
 &=(-1)^j\int\limits_0^{\infty} {x^{j-1+\sigma +it}\over 
(j-1+\sigma +it)(j-2+\sigma +it)\cdots(\sigma +it)}\,g_{r,s}^{(j)}(x)\,{\rm d}x = {}
\qquad\qquad\qquad\qquad\quad\cr 
 &={(-1)^j\over (j-1+\sigma +it)(j-2+\sigma +it)\cdots(\sigma +it)}
\int\limits_{P/2}^P O_j\!\left( x^{\sigma -1}\right) {\rm d}x \ll_{\sigma ,j} {}\cr  
 &\ll_{\sigma , j}\ P^{\sigma}\prod_{0\leq m<j}^{{\over}}|m+\sigma +it|^{-1}\qquad\quad    
\hbox{for $\quad\sigma\in{\Bbb R}$, $\,0\neq t\in{\Bbb R}\,$ and $\,j\in{\Bbb N}\cup\{ 0\}$.}
 &(9.64)}$$ 

We apply the case $\sigma =1$ of the identity (9.62) 
to the factor $g_{r,s}(|c|^2/|s|^2 |r|)$  of the 
summand on the right-hand side of (9.60). By following that with a 
change in the order of summation and integration 
(justified by the finiteness of the sum concerned), we find that on the right-hand side 
of the equation (1.4.10) one has 
$$K_{r,s}(n,\ell) 
={1\over 2\pi }\int\limits_{-\infty}^{\infty} 
G_{r,s}(1+it)\left| s\sqrt{r}\right|^{2+2it} 
\!\!\!\!\!\sum_{\scriptstyle c\in{}^{1/s}{\cal C}^{\infty}\atop\scriptstyle 
2^{-8} D^{*}<|c|^2< 2^8 D^{*}}^{\left(\Gamma_0(rs)\right)} 
\!\!\Omega\!\left(\left| {2\pi\sqrt{n\ell}\over c}\right|^2 X\right) 
 \left|{c}\right|^{-2-2it}
S_{1/s,\infty}\!\left( n,\ell;c\right) {\rm d}t\;.$$ 
From this it follows (similarly) that, when $(r,s)\in{\cal B}(R,S)$, one has: 
$$\sum_{{\textstyle{N\over 4}}<|n|^2\leq N} a_n 
\sum_{{\textstyle{L\over 2}}<|\ell|^2\leq L} A(\ell) 
K_{r,s}(n,\ell) 
={1\over 2}\int\limits_{-\infty}^{\infty} 
G_{r,s}(1+it)\left| s\sqrt{r}\right|^{2+2it} (2\pi)^{-2it} 
\kappa(r,s;t)\,{\rm d}t\;,\eqno(9.65)$$
where 
$$\kappa(r,s;t) 
=\sum_{{\textstyle{N\over 4}}<|n|^2\leq N} a_n |n|^{-it}
\sum_{{\textstyle{L\over 2}}<|\ell|^2\leq L} A(\ell) |\ell|^{-it} 
\!\!\!\!\!\sum_{\scriptstyle c\in{}^{1/s}{\cal C}^{\infty}\atop\scriptstyle 
2^{-8} D^{*}<|c|^2< 2^8 D^{*}}^{\left(\Gamma_0(rs)\right)} 
\!\!\!{S_{1/s,\infty}\!\left( n,\ell;c\right)\over\pi |c|^2}
\ \varphi_t\!\!\left(\left|{2\pi\sqrt{n\ell}\over c}\right|\right) ,\eqno(9.66)$$
with 
$$\varphi_t(\rho)=\Omega\!\left( X \rho^2\right) \rho^{2it}\qquad\qquad\quad\hbox{($\rho >0$).}\eqno(9.67)$$ 

By (1.4.3), the innermost sum on the right-hand side of the equation~(9.66) is 
finite,  so it is certainly the case that, for $(r,s)\in{\cal B}(R,S)$, the 
function $t\mapsto\kappa(r,s;t)$ is both continuous and bounded on ${\Bbb R}$. 
Therefore, by (9.65) and the case $\sigma =1$ of the bounds in (9.64), 
it follows that if $(r,s)\in{\cal B}(R,S)$ and $j>1$ then 
$$\sum_{{\textstyle{N\over 4}}<|n|^2\leq N} a_n 
\sum_{{\textstyle{L\over 2}}<|\ell|^2\leq L} A(\ell) 
K_{r,s}(n,\ell) 
\ll_j\, P |s|^2 |r|\int\limits_{-\infty}^{\infty} 
\left|\kappa(r,s;t)\right|\,{{\rm d}t\over (1+|t|)^j}\;.$$ 
Consequently, given the definitions (1.4.10), (1.4.8) of 
${\Lambda}$ and ${\cal B}(R,S)$ (which ensure, amongst other things,  
that the latter is a finite set), 
one has the bounds:  
$${\Lambda}\ll_j P S R^{1/2} \int\limits_{-\infty}^{\infty} 
\Biggl(\ \sum_{(r,s)\in{\cal B}(R,S)} 
\left|b(r,s)\,\kappa(r,s;t)\right|\Biggr)\,{{\rm d}t\over (1+|t|)^j}\qquad\qquad  
\hbox{($j>1$).}\eqno(9.68)$$

The next step is (in effect) to apply the `Kloosterman to spectral' sum formula 
(Theorem~1), in order to express $\kappa(r,s;t)$ (as given by (9.66)) in 
terms of spectral data associated with the space 
$L^2(\Gamma_0(rs)\backslash G)$. 
However we save some time and space by using, in place of Theorem~1, 
the Corollary to Theorems~1 and~2. 

Note firstly that, given the constraint (9.56) on ${\rm Supp}(\Omega)$, 
it is shown by Lemma~9.4 that, for each $t\in{\Bbb R}$, the 
hypotheses of the Corollary to Theorems~1 and~2 concerning 
$\varphi : (0,\infty)\rightarrow{\Bbb C}$, $X$, $A$ and 
$f : {\Bbb C}^{*}\rightarrow{\Bbb C}$ are satisfied 
when one has: $\varphi =\varphi_t$ (the function defined in (9.67)), 
$X$ as defined in (1.4.7), $A=8$, and $f(z)=\varphi_t(|z|)\ $ ($z\in{\Bbb C}^{*}$). 
We have, in particular, 
$${\rm Supp}\left(\varphi_t\right)\subseteq\left[ 2^{-3} X^{-1/2}\,,\,2^3 X^{-1/2}\right]\qquad\qquad\quad   
\hbox{($t\in{\Bbb R}$),}\eqno(9.69)$$ 
and so, given (1.4.7) and (9.61), the conditions $2^{-8} D^{*}<|c|^2<2^8 D^{*}$ 
constraining the innermost summation on the right-hand side of the equation~(9.66) are, 
in effect, superfluous: for if $N/4<|n|^2\leq N$, and if 
$L/2<|\ell|^2\leq L$, then $\varphi_t(|2\pi\sqrt{n\ell}/c|)=0$ for all 
$c\in{}^{1/s}{\cal C}^{\infty}$ that satisfy either 
$|c|^2\leq 2^{-15/2} D^{*}$ or $|c|^2\geq 2^6 D^{*}$. It therefore   
follows by the Corollary to Theorems~1 and~2 that, when  
$(r,s)\in{\cal B}(R,S)$, $t\in{\Bbb R}$ and $\eta >0$, one has: 
$$\eqalignno{
\kappa(r,s;t) 
 &=\sum_{\scriptstyle V\atop\scriptstyle\nu_V>0}^{\left(\Gamma_0(rs)\right)} 
{\bf K} f_t\left(\nu_V,0\right) 
\overline{\sum_{{\textstyle{N\over 4}}<|n|^2\leq N} \overline{a_n}\,|n|^{it} 
c_V^{1/s}\left( n;\nu_V,0\right)} 
\sum_{{\textstyle{L\over 2}}<|\ell|^2\leq L} A(\ell)|\ell|^{-it}  
c_V^{\infty}\left( \ell;\nu_V,0\right) + {} 
 &(9.70)\cr 
 &\quad\ \, +O_{\eta}\!\!\left(\!(\log X)\,Y_t 
\left( 1+|\mu(1/s)| N^{\eta +1/2}\right) 
\!\left( 1+|\mu(\infty)| L^{\eta +1/2}\right) 
\left\| {\bf a}_N\right\|_2 
\!\Biggl(\,\sum_{{\textstyle{L\over 2}}
<|\ell|^2\leq L}\!\!\!\!|A(\ell)|^2\!\Biggr)^{\!\!\!1/2}\right)\!,\qquad}$$
where 
$$f_t(z)=\varphi_t(|z|)\qquad\qquad\quad\hbox{($z\in{\Bbb C}^{*}$)}\eqno(9.71)$$
and 
$$Y_t=X^{-3/2}\max_{\rho >0}\left|\varphi_t^{(3)}(\rho)\right|\eqno(9.72)$$ 
(while the definitions of the cuspidal subspaces $V$, the spectral parameters $\nu_V$, 
the ${\bf K}$-transform, the Fourier coefficients $c_V^{\frak a}(n;\nu,p)$, 
the factors $\mu({\frak a})$ and the norm $\|{\bf a}_N\|_2$ may, in each case, be found
either in Subsection~1.1, 
or else within the statements of Theorems~1 and~2, in Subsection~1.2). 

By Lemma~9.4 it follows that, when $t\in{\Bbb R}$, the function 
$\varphi_t : (0,\infty)\rightarrow{\Bbb C}$ given by (9.67) is 
infinitely differentiable on $(0,\infty)$, and moreover satisfies 
$$\varphi_t^{(j)}(\rho)\ll_j (1+|t|)^j \rho^{-j}\qquad\quad\hbox{for $\quad j\in{\Bbb N}\cup\{ 0\}\ $ 
and $\ \rho >0$}\eqno(9.73)$$
(although the implicit constant here does depend on $\Omega$,
as well as on $j$, we 
are nevertheless correct in omitting to indicate this dependence on $\Omega$ in (9.73), 
for our choice of $\Omega$ was not dependent on anything else, and remains fixed). 
Hence, and by (9.69), (9.71) and our hypothesis (1.4.7) that $X\geq 2$, it follows 
that, when $t\in{\Bbb R}$, we may apply the bound in (1.2.17) for $f=f_t$, 
$A=8$, and $\varphi=\varphi_t$, and so deduce, by (9.69) and (9.73) (for $j=0$), that  
$${\bf K} f_t(\nu ,0)\ll\quad\int\limits_{2^{-3}X^{-1/2}}^{2^3 X^{-1/2}} 
O(1)\,{{\rm d}\rho\over\rho}\,\min\left\{\log X\,,\,\nu^{-1}\right\} X^{\nu} 
\ll (\log X) X^{\nu}\qquad\qquad\hbox{($0<\nu\leq 1/2$).}\eqno(9.74)$$
By (9.69), (9.72) and the case $j=3$ of (9.73), we have also 
$$Y_t=X^{-3/2}\max_{\rho \geq 2^{-3} X^{-1/2}} O\left( (1+|t|)^3 \rho^{-3}\right) 
\ll (1+|t|)^3\qquad\qquad\hbox{($t\in{\Bbb R}$).}$$
Moreover, it follows by the definition (1.4.8) that, 
when $(r,s)\in{\cal B}(R,S)$, $\,q=rs$ 
and $\Gamma=\Gamma_0(q)$, one has  
$$|\mu(1/s)|=|\mu(\infty)|=|q|^{-1}$$
(for, as is observed towards the 
end of the proof of Theorem~11, if $\Gamma=\Gamma_0(q)$, and if $s\mid q$ is 
such that $q$ and $q/s$ are coprime, then $1/\mu(1/s)\sim1/\mu(\infty)\sim q$). 

Given the bound in (9.74), and given what has just been noted concerning $Y_t$, 
$\mu(1/s)$ and $\mu(\infty)$, we may deduce from the equation~(9.70)  
and Theorem~3 that, for $(r,s)\in{\cal B}(R,S)$, $t\in{\Bbb R}$ and 
$\eta =\varepsilon /16>0$, 
$$\eqalign{
\kappa(r,s;t) 
 &\ll (\log X)\sum_{\scriptstyle V\atop\scriptstyle\nu_V>0}^{\left(\Gamma_0(rs)\right)} 
X^{\nu_V}\left|\Biggl(\,\sum_{{\textstyle{N\over 4}}<|n|^2\leq N} 
\!\!\overline{a_n}\,|n|^{it} 
c_V^{1/s}\left( n;\nu_V,0\right)\Biggr)  
\Biggl(\,\sum_{{\textstyle{L\over 2}}<|\ell|^2\leq L} 
\!\!A(\ell)|\ell|^{-it}  
c_V^{\infty}\left( \ell;\nu_V,0\right)\!\Biggr)\right| + {}\cr 
 &\qquad +O_{\varepsilon}\!\!\left(\!(\log X)\,(1+|t|)^3  
\left( 1+{N^{\eta +1/2}\over |rs|}\right) 
\!\left( 1+{L^{\eta +1/2}\over |rs|}\right) 
\left\| {\bf a}_N\right\|_2 
\!\Biggl(\,\sum_{{\textstyle{L\over 2}}<|\ell|^2\leq L}
\!\!\!\!|A(\ell)|^2\!\Biggr)^{\!\!\!1/2}\right) .}$$
By this and the case $j=13$ of the bound (9.68), either  
$${\Lambda}\ll_{\varepsilon} (\log X) P S R^{1/2} 
\left\| {\bf a}_N\right\|_2 
\!\Biggl(\,\sum_{{\textstyle{L\over 2}}<|\ell|^2\leq L}
\!\!\!|A(\ell)|^2\!\Biggr)^{\!\!\!1/2}
\!\!\!\!\!\!\sum_{(r,s)\in{\cal B}(R,S)}\!\!\!|b(r,s)| 
\!\left( 1+{N^{\varepsilon /16+1/2}\over |rs|}\right) 
\!\!\left( 1+{L^{\varepsilon /16 +1/2}\over |rs|}\right)\!,\eqno(9.75)$$ 
or else  
$${\Lambda}\ll (\log X) P S R^{1/2} \int\limits_{-\infty}^{\infty} 
\bigl( S_{-t/2,t/2}^{\infty , *}\bigr)\,{{\rm d}t\over (1+|t|)^{13}}\;,\eqno(9.76)$$
with $S_{u,y}^{\infty , *}=S_{u,y}^{\infty , *}(R,S;X;L,N)$ being the 
bounded and continuous non-negative real valued 
function of $(u,y)\in{\Bbb R}\times{\Bbb R}$ that is given 
by the equation~(9.25), in Lemma~9.3.

\goodbreak 
In the latter of the two cases just described, 
the results (1.4.12) and (1.4.14) of Theorem~10 
follow immediately from (9.76) 
and the bounds (9.26), (9.27) of Lemma~9.3 
(the relevant calculation is straightforward if one notes that $RS=Q$, that 
$1\leq\delta^{-1}+|t/2|\leq (1+|t|)\delta^{-1}$ 
for $0<\delta\leq 1$ and $t\in{\Bbb R}$, and that  
$\int_{{\Bbb R}} (1+|t|)^{\sigma}{\rm d}t<\infty$ for $\sigma <-1$). 

In the former case (in which the bound (9.75) holds) we may note that, since 
$$\sum_{{\textstyle{L\over 2}}<|\ell|^2\leq L} |A(\ell)|^2 
=\sum_{{\textstyle{L\over 2}}<|\ell|^2\leq L} O(1)\ll L$$
(by the case $j=k=0$ of the hypotheses in (1.4.5)), 
and since one has also 
$$\sum_{(r,s)\in{\cal B}(R,S)}  
\left( 1+{N^{\eta +1/2}\over |rs|}\right)^{\!\!2} 
\!\left( 1+{L^{\eta +1/2}\over |rs|}\right)^{\!\!2}
\ll |{\cal B}(R,S)| \left( 1+{N\over RS}\right) 
\left( 1+{L\over RS}\right) Q^{8\eta}\qquad\qquad\hbox{($\eta >0$)}$$
and 
$$|{\cal B}(R,S)|\leq\left( 2 R^{1/2}\right)^2 \left( 2 S^{1/2}\right)^2 =16 RS=16Q\eqno(9.77)$$
(by the definition of ${\cal B}(R,S)$ in (1.4.8), and the conditions in~(1.4.6)),
it therefore follows by a combination of the bound (9.75) and 
the Cauchy-Schwarz inequality 
that    
$$\eqalignno{{\Lambda}^2
 &=O_{\varepsilon}\left(  
(\log X)^2 P^2 S^2 R \left\|{\bf a}_N\right\|_2^2 L \|\,b\|_2^2 
\,|{\cal B}(R,S)| \,Q^{\varepsilon /2}\left( 1+{N\over Q}\right)
\left( 1+{L\over Q}\right)\right) \ll_{\varepsilon} {}\cr 
 &\ll_{\varepsilon} (\log X)^2 P^2 S \left\|{\bf a}_N\right\|_2^2 L \,\| b\|_2^2 
\,Q^{\varepsilon /2}\left( Q+N\right)\left( Q+L\right)  
 &(9.78)}$$ 
where $\| b\|_2$ is as stated in (1.4.13); and it moreover follows 
by (9.78), (1.4.13) and (9.77) that one has also:  
$$\eqalignno{{\Lambda}^2 
 &=O_{\varepsilon}\left( (\log X)^2 P^2 S \left\|{\bf a}_N\right\|_2^2 L 
\,\| b\|_{\infty}^2 
\,|{\cal B}(R,S)| 
\,Q^{\varepsilon /2}\left( Q+N\right)\left( Q+L\right)\right) \ll_{\varepsilon} {}\cr 
 &\ll_{\varepsilon} (\log X)^2 P^2 S \left\|{\bf a}_N\right\|_2^2 L \,\| b\|_{\infty}^2 
\,Q^{1+\varepsilon /2}\left( Q+N\right)\left( Q+L\right) \;,
 &(9.79)}$$
where $\| b\|_{\infty}$ is given by (1.4.15). 

The bound in (9.78) implies the result stated in (1.4.12)  
(where, by definition, one has $\vartheta\geq 0$, and where, by hypothesis, 
one also has $\delta^{-1}\geq 1$). Moreover, when 
the hypotheses of Theorem~9 
concerning $a_n$ ($n\in{\frak O}-\{ 0\}$), $N,H,K\in[1,\infty)$ and 
$\alpha,\beta : {\Bbb C}\rightarrow{\Bbb C}$ are satisfied, 
one has (see (8.27), in the proof of Lemma~8.4) the bound 
$\left\|{\bf a}_N\right\|_2^2 =O_{\varepsilon}(N^{1+\varepsilon /4})$, 
so that the bound (9.79) implies the result stated in (1.4.14) 
(given that, by the condition (1.4.6), 
one has $N^{\varepsilon /4}\leq Q^{\varepsilon /2}$).

By the conclusions reached in the last three paragraphs,  
we have obtained, in all relevant cases, what is stated in Theorem~10 
\ $\blacksquare$ 

\bigskip 
\bigskip 
 
\goodbreak\centerline{\bf References} 
\bigskip 
\parskip = 2 pt 
\item{[1]} T. M. Apostol, {\it Mathematical Analysis\/}, 2nd edition, 
World Student Series, Addison-Wesley, Reading MA, 1974.

\item{[2]} R. W. Bruggeman, Fourier coefficients of cusp forms, {\it Invent. Math.}, 
{\bf 45} (1978), 1-18.

\item{[3]} R. W. Bruggeman and R. J. Miatello,
Estimates of Kloosterman sums for groups of real rank one,
{\it Duke Math. J.}, {\bf 80} (1995), 105-137.

\item{[4]} R. W. Bruggeman and Y. Motohashi, 
Sum formula for Kloosterman sums and fourth moment of the Dedekind
zeta-function over the Gaussian number field, {\it Functiones et Approximatio},
{\bf 31} (2003), 23-92.

\item{[5]} J.-M. Deshouillers and H. Iwaniec,
Kloosterman sums and Fourier coefficients of cusp forms,
{\it Invent. Math.}, {\bf 70} (1982), 219-288.

\item{[6]} J.-M. Deshouillers and H. Iwaniec,
Power mean values of the Riemann zeta-function,
{\it Mathematika}. {\bf 29} (1982), 202-212.

\item{[7]} J. Elstrodt, F. Grunewald and J. Mennicke, 
{\it Groups Acting on Hyperbolic Space\/}, Springer Monographs in Mathematics, 1997.

\item{[8]} G. H. Hardy and E. M. Wright, An Introduction to the Theory of Numbers
(5th edition) Oxford University Press 1979 (reprinted with corrections in 1983 and 1984).

\item{[9]} M. N. Huxley, the large sieve inequality for algebraic number fields, 
{\it Mathematika}, {\bf 15} (1968), 178-187.

\item{[10]} A. Ivi\'{c}, The Riemann Zeta-Function: Theory and Applications,
Dover Publications, Inc., 2003.

\item{[11]} H. Iwaniec, Spectral Theory of Automorphic Functions, 
lecture notes (manuscript), Rutgers University, 1987.

\item{[12]} H. Iwaniec, Topics in Classical Automorphic Forms, 
Graduate Studies in Mathematics {\bf 17}, AMS, Providence RI, 1997.

\item{[13]} H. H. Kim, On local $L$-functions and normalized intertwining operators, 
{\it Canad. J. Math.\/}, {\bf 57} no. 3 (2005), 535-597. 

\item{[14]} H. H. Kim and F. Shahidi, Cuspidality of symmetric powers with applications, 
{\it Duke Math. J.\/} {\bf 112} (2002), 177-197. 

\item{[15]} N. V. Kuznetsov, Petersson hypothesis for forms of weight zero and Linnik hypothesis, 
Preprint No.~02, Khabarovsk Complex Res. Inst., East Siberian Branch Acad. Sci. USSR, 
Khabarovsk, 1977 (in Russian).

\item{[16]} N. V. Kuznetsov, Petersson hypothesis for parabolic forms of weight zero 
and Linnik hypothesis. Sums of Kloosterman sums, {\it Mat. Sbornik.,} {\bf 111 (153)} No.~3 (1980), 
334-383.

\item{[17]} S. Lang, {\it Real Analysis\/}, (2nd edn.),
Addison-Wesley, Reading MA, 1983.

\item{[18]} R. P. Langlands, {\it On the functional equations satisfied by Eisenstein series\/}, 
Lecture Notes in Math., {\bf 544}, Springer-Verlag, Berlin, 1976.

\item{[19]} H. Lokvenec-Guleska, {\it Sum Formula for SL${}_2$ over
Imaginary Quadratic Number Fields\/}, thesis (in English, with
summaries in Dutch and Macedonian), University of Utrecht, 2004.

\item{[20]} O. Ramar\'{e} and R. Rumely, 
Primes in arithmetic progressions, {\it Math. Comp.}, {\bf 65}, No. 213, 
Jan. 1996, 397-425.

\item{[21]} N. Watt, Kloosterman sums and a mean value for Dirichlet
polynomials, {\it J. Number Theory}, {\bf 53} (1995), 179-210.

\item{[22]} N. Watt,
Spectral large sieve inequalities for Hecke congruence subgroups of
$SL(2,{\Bbb Z[i]})$, preprint submitted for publication, 
arXiv:1302.3112v1 [math.NT].

\item{[23]} N. Watt,
Weighted fourth moments of Hecke zeta functions
with groessencharakters, preprint in preparation. 

\item{[24]} E. T. Whittaker and G. N. Watson, 
{\it A Course of Modern Analysis\/}, Fourth edn. Cambridge University Press, 1927.
\parskip=0 pt

\bigskip\medskip    

\bye